\documentclass[twoside,12pt]{amsart}
\usepackage{amsmath,amssymb,amsrefs}
\usepackage{graphicx}
\usepackage{bm}
\usepackage{hyperref} 
\usepackage[mathscr]{euscript}
\usepackage{arydshln}
\usepackage{enumerate}
\usepackage{xcolor}
\usepackage{enumitem}

\numberwithin{equation}{section}

\newtheorem{theorem}{Theorem}[section]
\newtheorem{corollary}[theorem]{Corollary}
\newtheorem{proposition}[theorem]{Proposition}
\newtheorem{lemma}[theorem]{Lemma}
\newtheorem{example}[theorem]{Example}
\theoremstyle{definition}
\newtheorem{definition}[theorem]{Definition}
\theoremstyle{remark}
\newtheorem{remark}[theorem]{Remark}

{\unskip\nobreak\hskip 1em plus 1fil\nobreak$\Box$
\parfillskip=0pt%
\endtrivlist}

\newcommand{\R}{\mathbb{R}}
\newcommand{\C}{\mathbb{C}}
\newcommand{\Z}{\mathbb{Z}}
\newcommand{\N}{\mathbb{N}}

\newcommand{\Res}{\mathcal{K}}
\newcommand{\Ressp}{\mathcal{K}_\sharp}
\newcommand{\Ione}{I_1}
\newcommand{\Jone}{J_1}
\newcommand{\bzero}{{\bf 0}}
\newcommand{\tr}{\tilde{r}}
\newcommand{\btr}{ {\bf\tilde{r}} }
\newcommand{\tIone}{\tilde{I}_1}
\newcommand{\cH}{H^\lambda_{\rm bulk}} 
\newcommand{\cHeg}{H^\lambda_\sharp} 
\newcommand{\tHeg}{\widetilde{H}_\sharp^\lambda} 
 
\newcommand{\HTBb}{H^{^{\rm TB}}_{_{\rm bulk}}} 
\newcommand{\HTBs}{H^{^{\rm TB}}_\sharp} 
\newcommand{\be}{{\bf e}}

\newcommand{\bk}{{\bf k}}
\newcommand{\bfm}{{\bf m}}
\newcommand{\bn}{{\bf n}}

\newcommand{\br}{{\bf r}}

\newcommand{\tc}{{\tilde{c}}}
\newcommand{\cp}{{c^\prime}}

\newcommand{\bK}{{\bf K}}

\newcommand{\bv}{{\bf v}}
\newcommand{\bw}{{\bf w}}
\newcommand{\bx}{{\bf x}}
\newcommand{\bz}{{\bf z}}

\newcommand{\by}{{\bf y}}
\newcommand{\bA}{{A}}
\newcommand{\bB}{{B}}

\newcommand{\bI}{{I}}
\newcommand{\bJ}{{J}}

\newcommand{\none}{{n_1}}
\newcommand{\ntwo}{{n_2}}
\newcommand{\tm}{{\tilde{m}}}

\newcommand{\tx}{{\tilde{\bx}}}
\newcommand{\ty}{{\tilde{\by}}}

\newcommand{\tOmega}{{\widetilde{\Omega}}}

\newcommand{\vtilde}{{\bm{\mathfrak{v}}}}
\newcommand{\ktilde}{{\bm{\mathfrak{K}}}}

\newcommand{\kpar}{{k_{\parallel}}}
\newcommand{\kperp}{{k_{\perp}}}

\newcommand{\D}{\partial}

\newcommand{\eps}{\varepsilon}
\newcommand{\nit}{\noindent}
\newcommand{\nn}{\nonumber}

\newcommand{\Gfree}{G_\lambda^{\rm free}}
\newcommand{\Gatom}{G_\lambda^{\rm atom}}
\newcommand{\Hatom}{H^\lambda_{\rm atom}}
\newcommand{\Thout}{\Theta_{\rm out}}
\newcommand{\Tfree}{\Theta_{\rm free}}
\newcommand{\VGamma}{V_\Gamma^{^\lambda}}
\newcommand{\HGamma}{H_\Gamma^{^\lambda}}

\newcommand{\tsigma}{{\tilde{\sigma}}}

\newcommand{\supp}{\text{supp}}

\begin{document}

\title{
Continuum Schroedinger operators for sharply terminated graphene-like structures
}

\author{C.L. Fefferman}
\address{Department of Mathematics, Princeton University, Princeton, NJ, USA}
\email{cf@math.princeton.edu}
%
%
\author{M. I. Weinstein}
\address{Department of Applied Physics and Applied Mathematics and Department of Mathematics, Columbia University, New York, NY, USA}
\email{miw2103@columbia.edu}

\date{\today}


\begin{abstract} 

We study the single electron model of a semi-infinite graphene sheet interfaced with the vacuum and terminated along a zigzag edge. The model is a Schroedinger operator acting on $L^2(\R^2)$: $H^\lambda_{\rm edge}=-\Delta+\lambda^2 V_\sharp$, with a potential $V_\sharp$ given by a sum of translates an atomic potential well, $V_0$, of depth $\lambda^2$, 
 centered on a subset of the vertices of a  discrete honeycomb structure with a zigzag edge. 
 We give a complete analysis of the low-lying energy spectrum of $H^\lambda_{\rm edge}$ in the strong binding regime ($\lambda$ large). In particular, 
 we prove scaled resolvent convergence of $H^\lambda_{\rm edge}$ acting on $L^2(\R^2)$,
 to  the (appropriately conjugated) resolvent  of a limiting discrete tight-binding Hamiltonian acting in $l^2(\N_0;\C^2)$.
We also prove 
the existence of {\it edge states}: solutions of the eigenvalue problem  for $H^\lambda_{\rm edge}$
which are localized transverse to the edge and 
 pseudo-periodic  (propagating  or plane-wave like) parallel to the edge.   These edge states arise 
 from a ``flat-band'' of eigenstates the tight-binding  Hamiltonian. 
\end{abstract}

\maketitle   

\pagestyle{myheadings}
\thispagestyle{plain}
\markboth{Strong binding for sharply terminated periodic structure}{C.L. Fefferman, M.I. Weinstein}


\section{Introduction}\label{intro}

Tight binding models are discrete models which are central to the modeling of spatially periodic  and more general crystalline structures in condensed matter physics. 
These models apply when the the quantum state of the system is well-approximated by 
superpositions of translates of highly-localized quantum states (orbitals) within deep atomic potential wells centered at lattice sites \cite{Ashcroft-Mermin:76}.
An important example is the tight-binding model of graphene, a planar honeycomb arrangement of carbon atoms
 with two atoms per unit cell.  The two-band tight-binding model yields an explicit approximation for its lowest two dispersion surfaces, which touch conically at {\it Dirac points} over the vertices of the Brillouin zone \cite{Wallace:47}. Such Dirac points
are central to the remarkable electronic properties of graphene \cite{Kim-etal:05,geim2007rise,RMP-Graphene:09,novoselov2011nobel}
 and its artificial (electronic, photonic, acoustic, mechanical,\dots) analogues; see, for example,  \cite{artificial-graphene:11,Rechtsman-etal:13a,2014LuJoannopoulosSoljacic,irvine:15,2017KhanikaevShvets,Marquardt:17} and the  survey \cite{ozawa_etal:18}.
 The existence of Dirac points for generic honeycomb Schroedinger operators was proved in \cite{FW:12,FLW-MAMS:17}; see also \cite{berkolaiko-comech:18}. That the two-band tight-binding model gives an accurate approximation of the low-lying dispersion surfaces in the regime of strong binding   was  proved  in \cite{FLW-CPAM:17}; see also Section \ref{semi-classical}.  Other results on Dirac points for Schroedinger operators on $\R^2$ may be found in 
 \cite{Colin-de-Verdiere:91,Grushin:09,ablowitz2009conical,ACZ:12,Lee:16}, coupled oscillator models \cite{Makwana-Craster:14}  and on quantum graphs in \cite{Kuchment-Post:07,Do-Kuchment:13}.

{\it Edge states} are modes which are propagating (plane-wave like) parallel to an interface and which are localized transverse to the interface. In condensed matter physics edge states  describe the phenomenon of electrical conduction along an interface. Two types of interfaces of great physical interest are 
a sharp terminations of a bulk structure structure studied in this article (see \cite{Dresselhaus-etal:96,delplace2011zak,mong2011edge,Graf-Porta:13}) and domain wall / line-defects within the bulk 
 (see \cite{2017KhanikaevShvets,Marquardt:17,Rechtsman-etal:18,Sunku:2018qx}  and studied, for example, in  \cite{FLW-2d_edge:16,FLW-2d_materials:15,LWZ:18,D:19,DW:19} ). 
The role of edge or surface modes in the spectral theory of Schroedinger operators with potentials which model, for example, the interface between a general periodic medium and a vacuum is studied in {\it e.g.} \cite{Davies-Simon:78,Karpeshina:97}. 

{\it In this paper we study the low-lying energy spectrum (discrete and continuous spectrum) of 
 a sharply terminated honeycomb structure, corresponding to a semi-infinite sheet of graphene 
joined to the vacuum along a sharp interface. We prove convergence of the operator resolvent to that of a discrete tight-binding model and construct the continuous spectrum of edge states.}

Edge states in honeycomb structures such as graphene are of particular interest as foundational building blocks in the field topological insulators (TI). TI's are materials which are insulating in their bulk and conduction along boundaries, which is robust against large localized perturbations. When graphene is subjected to a magnetic field, its edge currents become unidirectional and acquire such robustness. This phenomenon has an explanation in terms of topological invariants
 associated with a bulk Floquet-Bloch vector bundle, which takes on non-trivial values when time-reversal symmetry is broken \cite{HK:10,Kane-Mele:05,kane2005quantum}.

    A key difference between the types of interfaces is that the sharply terminated structure has no spectral gap, resulting in certain edge orientations supporting edge states and others not. In contrast, the domain wall structures perturbations studied in \cite{FLW-2d_edge:16,FLW-2d_materials:15,LWZ:18,D:19,DW:19} have edge states which localize along arbitrary rational edges.
   For a discussion of the roles played by edge orientation and the type of symmetry breaking in the existence and robustness of edge states for domain wall / line-defects, see \cite{DW:19}. 

 Specifically, for the tight-binding model,  edge states exist at sharp terminations along a zigzag edge for a subinterval of parallel quasi-momenta, $\kpar\in[0,2\pi)$ associated with the direction of translation invariance parallel to the edge.
 They do not exist 
at the sharp termination along an armchair edge; see, for example,
 \cite{Dresselhaus-etal:96,delplace2011zak,mong2011edge,Graf-Porta:13} and Section \ref{TB}.
Such results may be interpreted as consequences of the 
non-vanishing of the  Berry-Zak phase, $\mathscr{Z}(\kpar)$,  defined as the integral of the Berry connection over the one-dimensional Brillouin zone associated with the type of edge
\cite{delplace2011zak,mong2011edge}. 

%

\subsection{Mathematical setup}\label{setup}

In this paper we initiate a study of these phenomena in the context of the underlying 
continuum equations of quantum physics, in particular the single-electron model of bulk (infinite) graphene 
and its terminations.  In particular, we study Schroedinger operators on $\R^2$ for a sharp termination 
of a honeycomb structure along a zigzag edge.
\medskip

We denote the equilateral triangular lattice in $\R^2$ by . 
\begin{equation} \Lambda\ =\ \Z\vtilde_1\oplus\Z\vtilde_2\ , \label{e-lat}
\end{equation}
where $\vtilde_1$ and $\vtilde_2$ are given by
\begin{align}
\vtilde_1 &=\ \left( \begin{array}{c} \frac{\sqrt{3}}{2} \\ {}\\  \frac{1}{2}\end{array} \right),\ \ 
\vtilde_2 =\ \left(\begin{array}{c} 0 \\ {}\\ 1 \end{array}\right)\ .\label{c-v12-def}
\end{align}
The dual lattice, $\Lambda^*$, is given by 
\begin{equation} \Lambda^*\ =\ \Z\ktilde_1\oplus\Z\ktilde_2\ , \label{e-dulat}
\end{equation}
where $\ktilde_1$ and $\ktilde_2$ are given by 
\begin{align}
 &  \ktilde_1=\ 2\pi\left(\begin{array}{c} \frac{ 2\sqrt{3} }{3}\\ {}\\ 0\end{array}\right),\ \ \
 \ktilde_2 = 2\pi\left(\begin{array}{c} -\frac{\sqrt{3} }{3}\\ {}\\ 1\end{array}\right).\ 
 \label{c-k12-def}
 \end{align}
 Note that \begin{equation}
\ktilde_l\cdot\vtilde_m=2\pi\delta_{lm}.
\label{kv-orth}\end{equation}

To generate the honeycomb structure, we first fix base points in $\R^2$:\ 
 \begin{equation}
 \bv_\bA=(0,0),\qquad \bv_\bB=\left(1/2,1/(2\sqrt3)\right).
 \label{vAB}
 \end{equation}
  The honeycomb structure, $ \mathbb{H}$, is the union of the two interpenetrating sublattices 
  \begin{equation}
  \Lambda_\bA=\bv_\bA+\Lambda,\qquad \Lambda_\bB=\bv_\bB+\Lambda\ :
  \label{Lam-AB}
  \end{equation}
\begin{equation}
\mathbb{H}\ =\ \Lambda_\bA\ \cup\ \Lambda_\bB\ .
\label{honey}\end{equation}

Let  $V_0(\bx)$ be an {\it atomic potential well} which may be considered, for the present discussion,
 to be radially symmetric, 
compactly supported  with ${\rm supp}\ V_0\subset B_{r_0}(0)$, the open disc of radius $r_0$ about $0$.  We discuss more general and physically reasonable conditions on $V_0$  below in Section \ref{setup}. \medskip 

Our bulk Hamiltonian is the honeycomb Schroedinger operator:
\begin{equation}
H_{_{\rm bulk}}^\lambda\ =\ -\Delta + \lambda^2 V(\bx)\ \textrm{acting on $L^2(\R^2)$}, 
 \label{Hbk}
 \end{equation}
where $V(\bx)$ is a superposition identical atomic potential wells, centered at the vertices of $\mathbb{H}$:
 \begin{equation}
 V(\bx)=\sum_{\bv\in \mathbb{H}}\ V_0(\bx-\bv)\ ,\ \bx\in\R^2 \ ,
 \label{Vbk}\end{equation}

 The potential $V(\bx)$  satisfies the conditions of a {\it honeycomb lattice potential} 
 in the sense of Definition 2.1 of \cite{FW:12}. For all but a discrete subset of values of $\lambda$ (including $\lambda=0$), the operator $H_{^{\rm bulk}}^\lambda$ has Dirac points at energy / quasi-momentum pairs, $(E_D^\lambda,\bK_\star)$, where  $\bK_\star$, varies over the  vertices of the Brillouin zone
  \cite{FW:12,FLW-MAMS:17}; see also \cite{berkolaiko-comech:18}. Moreover, for  $\lambda$ large (strong binding), the  low-lying Floquet-Bloch dispersion surfaces of $H_{_{\rm bulk}}^\lambda$, when rescaled,  
 are uniformly approximated by the dispersion surfaces  of the two-band tight-binding model
 \cite{FLW-CPAM:17}. 

\begin{figure}
\centering 
\includegraphics[width=.65\textwidth]{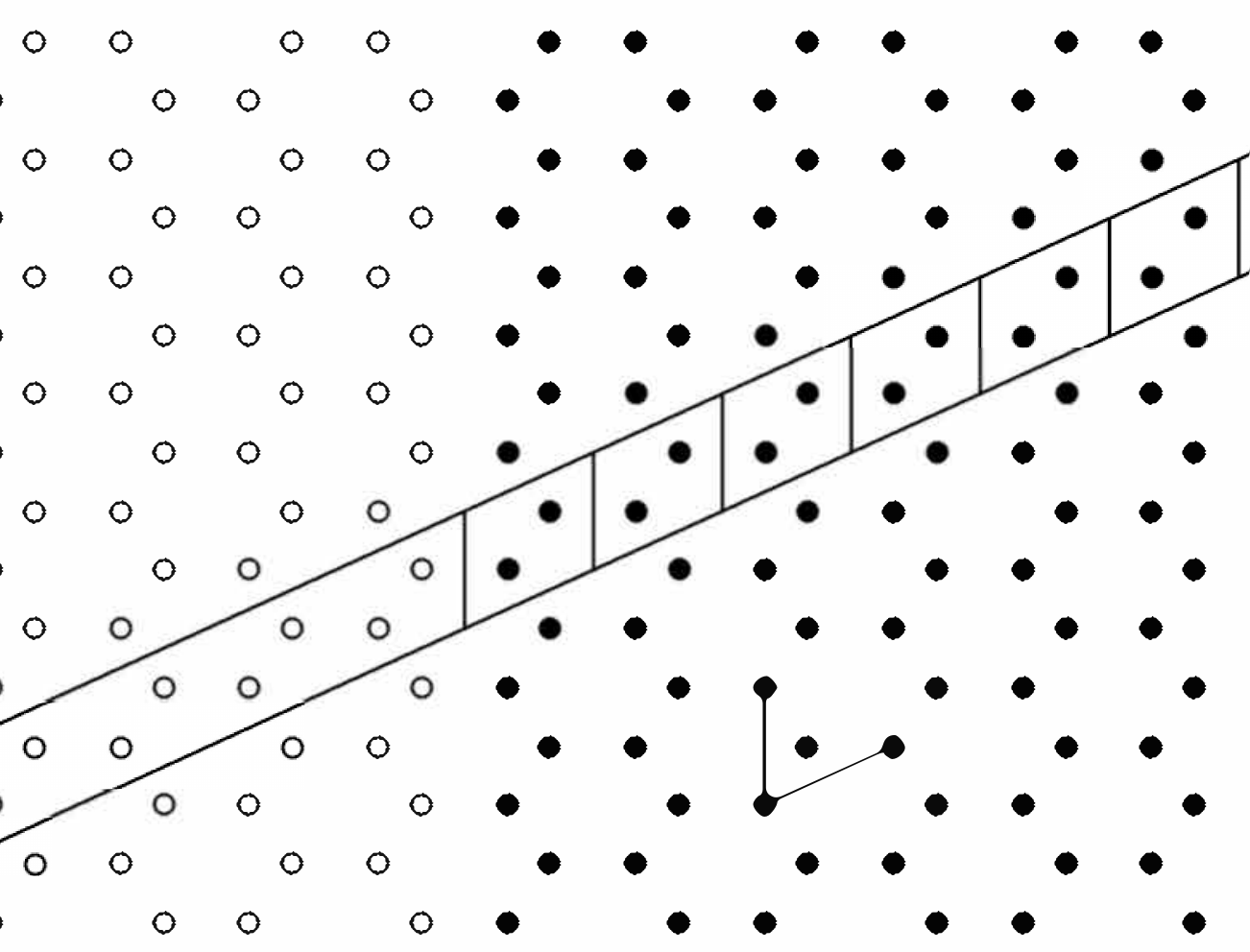}
\caption{\small  (a) $\mathbb{H}$: Bulk honeycomb structure consists of all vertices (circles, light and dark). (b) $\mathbb{H}_\sharp$: Honeycomb structure terminated along a zigzag edge consists of vertices indicated by dark circles; see \eqref{bbH+}. (c) 
$\Omega_\Sigma$: Indicated strip is a choice of fundamental cell for the cylinder $\Sigma=\R^2/\Z\vtilde_2$. $\Omega_\Sigma=\Omega_{-1}\cup\Omega_0\cup\Omega_1\cup\cdots\cup\Omega_n\cup\cdots$.
Sites: $\bv_A^n, \bv_B^n$ in finite parallelograms $ \Omega_n$, $n\ge0$, are sites in $\mathbb{H}_\sharp$. 
$\Omega_{-1}$ denotes the infinite parallelogram containing no vertices of the terminated structure, $\mathbb{H}_\sharp$.\ }
\label{zz-dimers}
\end{figure}

Consider now a ``half-plane'' of  vertices $\mathbb{H}_\sharp\subset\mathbb{H}$, whose extreme points 
trace out  a zigzag pattern:
\begin{equation}
  \mathbb{H}_\sharp\ \equiv \{\bv_\bA\ +\ \N_0\vtilde_1\oplus\Z\vtilde_2\}\ \ \cup\ \  \{\bv_\bB\ +\ \N_0\vtilde_1\oplus\Z\vtilde_2\},\qquad \N_0=\{0,1,2,\dots\},
\label{bbH+}
  \end{equation}
  The set $\mathbb{H}_\sharp$ is  invariant with respect to translations by $\vtilde_2$
   and is the subset of  sites in $\mathbb{H}$ to the right of an infinite zigzag edge; see  Figure \ref{zz-dimers}. 
The set of zigzag edge (boundary) sites, also translation invariant by $\vtilde_2$, is given by:
$\{\bv_\bA\ +\ \Z\vtilde_2\}\ \ \cup\ \  \{\bv_\bB\ +\ \Z\vtilde_2\}$ .

We define the potential 
 \begin{equation}
 V_\sharp(\bx)=\sum_{\bv\in \mathbb{H}_\sharp}\ V_0(\bx-\bv)\ ,\qquad\ \bx\in\R^2 \ .
 \label{Vsharp1}\end{equation} 
 The operator 
 \[ H^\lambda_{_{\rm edge}}=-\Delta + \lambda^2 V_\sharp(\bx)\]
models a half-plane of graphene interfaced with the vacuum along  a zigzag edge. 
Note the translation invariance: $V_\sharp(\bx+\vtilde_2)=V_\sharp(\bx)$ for all $\bx\in\R^2$.

 Let 
$(E_0^\lambda,p_0^\lambda(\bx))$, with $p_0^\lambda>0$ and 
 $L^2-$ normalized, denote the ground state eigenpair of 
the {\it atomic Hamiltonian} 
\[ \Hatom=-\Delta + \lambda^2V_0(\bx).\] 
Let $\rho_\lambda$ denote the  {\it hopping coefficient}, given by:
\begin{equation}
  \rho_\lambda\ =\ \int_{|\by|<r_0}\ p_{_0}^\lambda(\by) \lambda^2\ |V_0(\by)|\ p_0^\lambda(\by-\be)\ d\by \ ,
  \label{rho-def}
  \end{equation}
  where $\be$ is any vector from one lattice site in $\mathbb{H}$ to a nearest neighbor in $\mathbb{H}$, {\it e.g.} $\bv_B-\bv_A$.
 The potential $V_0(\by)$ and ground state $p_{_0}^\lambda(\by)$ are localized around $\by=0$, while $p_0^\lambda(\by-\be)$, is localized at any nearest neighbor site $\be\in\mathbb{H}$. Recall that $\supp\ V_0$ is contained in the set where $|\bx|<r_0$. For $\lambda$ large $\rho_\lambda$ is exponentially small (see \eqref{rho-boundsA}) \cite{FLW-CPAM:17}. \medskip

The key accomplishments of this paper are the following:  
 \begin{enumerate}
 \item {\it Theorem \ref{res-conv} (Scaled resolvent convergence)}:\ We prove for $\lambda\ge\lambda_\star$ sufficiently large (the strong binding regime), 
 that the re-centered and scaled resolvent $\left(\ (H^\lambda_{_{\rm edge}} - E_0^\lambda)/\rho_\lambda - zI\ \right)^{-1}$ 
 has a universal limit (in the uniform operator norm) described by a discrete (tight-binding) Hamiltonian, defined on a truncated honeycomb structure.
 The band structure of this limiting operator is displayed in Figure \ref{tb-spec}.\\
 \item {\it Theorem \ref{main-thm1} (Zigzag edge states):}\ We construct a continuum of edge state modes. These are eigenstates of $H^\lambda_{_{\rm edge}}$, which  are propagating (plane-wave like) parallel to and localized transverse to the zigzag edge. Upon appropriate $\lambda-$ dependent rescaling, these edge-states are close to (and converge as $\lambda$ tends to infinity to) the flat band of zero energy edge states of the tight-binding model; see Figure \ref{tb-spec}.\\
 \item {\it Resolvent kernel bounds on arbitrary discrete sets:}\ The methods of this article go considerably beyond those our previous article on the strong binding regime \cite{FLW-CPAM:17}, which established
 convergence to the (universal) two-band tight binding spectrum for the bulk graphene-like structures.  Since Theorems \ref{res-conv} and \ref{main-thm1} involve convergence of operators and eigenstates on an \underline{infinite cylinder} (Figure \ref{zz-dimers}), we required pointwise decay properties of the resolvent kernel $H^\lambda_{_{\rm edge}}$ for energies near $E_0^\lambda$. These bounds are stated in Theorem \ref{Ksharp}. In Proposition  \ref{summary-K} we  establish these kernel estimates for potentials which are a sum of atomic potentials centered on an  \underline{arbitrary discrete set}  of lattice sites $\Gamma\subset\R^2$  (not necessarily translation invariant) whose minimal pairwise distance
is $M r_0$, where $r_0$ is the radius of the support of $V_0$ and $M>2$ is some positive constant. 
We then specialize to a translation invariant set to obtain  Theorem \ref{Ksharp}.
  We believe the technique we have developed will be quite broadly applicable.
  \end{enumerate}
 
\medskip 
 
We next introduce the {\it edge state eigenvalue problem}. 
 Associated with the translation invariance of $-\Delta + \lambda^2 V_\sharp(\bx)$ by $\vtilde_2$ is a {\it parallel quasi-momentum},
  denoted $\kpar\in[0,2\pi)$. The condition that an edge state, $\Phi$, is propagating parallel to the zigzag edge is:
\begin{equation}  \Phi(\bx+\vtilde_2)\ =\ e^{i\kpar}\ \Phi(\bx),\ \ \bx\in\R^2\ \label{prop-bc}\ .
\end{equation}
We introduce the cylinder 
\begin{equation}
\Sigma=\R^2/\Z\vtilde_2.
\label{zz-cylinder}\end{equation}
 The space $L^2(\Sigma)$ consists of functions 
which are square integrable over a fundamental cell of $\Sigma$, {\it e.g.} the strip $\Omega_\Sigma$ shown in  Figure \ref{zz-dimers},  and which satisfy the periodic boundary condition with respect to $\vtilde_2$:
$\phi(\bx+\vtilde_2)\ =\  \phi(\bx)$ for almost all $\bx\in\Omega_\Sigma$ and all $\bv\in\Lambda$. 

We enforce the condition that  (i) $\Phi$ is $\kpar-$ pseudo-periodic parallel  to the zigzag edge, \eqref{prop-bc}, and (ii) decaying  to zero transverse 
to the zigzag edge as $\bx$ tends to infinity by requiring 
\[  e^{-i\frac{\kpar}{2\pi}\ktilde_2\cdot\bx}\Phi(\bx)\in L^2(\Sigma).\]
 For such functions we write $\Phi\in L^2_\kpar(\Sigma)$ or just $\Phi\in  L^2_\kpar$. 
 We  can now formulate the\medskip

\nit {\bf $\kpar-$Zigzag Edge State Eigenvalue Problem for $H^\lambda_{_{\rm edge}}=-\Delta+V_\sharp(\bx)$}:\medskip

\begin{align}
H^\lambda_{_{\rm edge}}\Psi(\bx)\ \equiv\ \left(\ -\Delta + \lambda^2 V_\sharp(\bx)\ \right)\ \Psi(\bx)\ &=\ E\ \Psi(\bx),\ \ \bx\in\R^2,\qquad \Psi\in L^2_\kpar(\Sigma).
\label{zz-evp}\end{align}

Defining $\Psi(\bx)=e^{i\frac{\kpar}{2\pi}\ktilde_2\cdot\bx}\psi(\bx)$, we may formulate \eqref{zz-evp} equivalently as:
\begin{align}
H^\lambda_{_{\rm edge}}(\kpar)\psi\ \equiv\ \left(\ -\left(\nabla+i\frac{\kpar}{2\pi}\ktilde_2\right)^2 + \lambda^2 V_\sharp(\bx)\ \right)\ \psi(\bx)\ &=\ E\ \psi(\bx),\ \ \bx\in\R^2,\qquad \psi\in L^2(\Sigma)\ .
\label{zz-kp-evp}\end{align}
We refer to non-trivial solutions of the eigenvalue problem \eqref{zz-evp} (equivalently \eqref{zz-kp-evp}) as zigzag edge states.

Before stating our main results Theorems \ref{res-conv} and  \ref{main-thm1}, 
 we recall a key observation used in \cite{FLW-CPAM:17} to obtain the low-lying dispersion surfaces 
 (energies near the atomic ground state energy, $E_0^\lambda$) of the bulk honeycomb Schroedinger operator, $H^\lambda_{_{\rm bulk}}$. That is, for $\lambda$ large, the $\bk-$ pseudo-periodic Floquet-Bloch eigenmodes which are associated with the two lowest spectral bands of $H^\lambda_{_{\rm bulk}}$, acting in $L^2(\R^2)$,  can be uniformly approximated by appropriate linear combinations of  the two $\bk-$ pseudo-periodic functions: $P_{\bk,I}^\lambda(\bx),\ I=A,B$.
These functions  are constructed as  $\bk-$ pseudo-periodic weighted sums of translates, $p_0^\lambda(\bx+\bv)$, of the atomic ground state,  where $\bv$ varies over the sublattices: $\Lambda_I=\bv_I+\Lambda,\ I=A,B$ .

 In the present work, to study the low-lying  spectral bands associated with  eigenvalue problem 
 $H^\lambda_{_{\rm edge}}(\kpar)\psi_\kpar^\lambda= E^\lambda(\kpar) \psi_\kpar$, $\psi_\kpar^\lambda\in L^2(\Sigma)$ with $E^\lambda(\kpar) $  near $E_0^\lambda$  and we find it very natural to approximate eigenstates by superpositions of the infinite family of functions 
 \begin{equation}
 p_{\bI,\kpar}^\lambda[n](\bx)\ \in\ L^2(\Sigma),\qquad I=A,B,\qquad n\ge0\ ,
 \label{Pk-intro}
 \end{equation}  
 which are  constructed as $\kpar-$ dependent and periodized (infinite) sums of translates of the ground state $p_0^\lambda(\bx)$ over the one-dimensional sublattices: $\bv_I+ n\vtilde_1 +\Z\vtilde_2$ of $\Lambda_I,\ I=A,B$  and $n\ge0$; see \eqref{Lam-AB}. The states $p_{\bI,\kpar}^\lambda[n](\bx)$
 are introduced in Definition \ref{pzn-p0} in Section \ref{nat-basis}.
 For $\lambda$ sufficiently large, any $F\in L^2(\Sigma)$ has the expansion
 \begin{equation}
  F=\sum_{I=A,B}\sum_{n\ge0} \alpha_n^I\ p_{\bI,\kpar}^\lambda[n](\bx)+ F_\perp,
  \label{Fsplit}\end{equation} where
  $\{\alpha_n^I\}\in l^2(\N_0;\C^2)$ and $F_\perp$ is  $L^2(\Sigma)-$ orthogonal to the span of the functions
   $p_{\bI,\kpar}^\lambda[n]$; see Proposition \ref{Xdecomp}.
 The tight-binding (discrete) edge Hamiltonian, $H^{\rm TB}_\sharp(\kpar)$ acting in $l^2(\N_0;\C^2)$, arises via translation and rescaling, 
of the operator whose matrix elements are $\left\langle p_{\bJ,\kpar}^\lambda[m], H_{_{\rm edge}}^\lambda(\kpar) p_{\bI,\kpar}^\lambda[n]\right\rangle_{L^2_\kpar}$, for $J,I = A,B$ and $m,n\ge0$. The tight-binding model is studied in Section \ref{TB} and its band spectrum is displayed in Figure \ref{tb-spec}.

\subsection{Main results}\label{main-results}

The relation of $H_{_{\rm edge}}^\lambda(\kpar)$ to the tight-binding Hamiltonian 
$ H^{^{\rm TB}}_\sharp(\kpar)$ is given by the following result on scaled resolvent convergence.
\begin{theorem}[Scaled resolvent convergence]\label{res-conv} 
As in Theorem \ref{main-thm1}, assume that $E_0^\lambda$, the ground state energy of the atomic Hamiltonian, $\Hatom=-\Delta+\lambda^2V_0$, satisfies the conditions (GS) \eqref{GS} and (EG) \eqref{EG} on the ground state energy and energy-gap, respectively.

 Let $\mathscr{C}$ denote a compact subset of $\C\setminus\sigma(H^{^{\rm TB}}_\sharp(\kpar))$, the resolvent set of $H^{^{\rm TB}}_\sharp(\kpar)$. 
There exist constants $\lambda_\star$, $C_\star$ and $c$,  which are independent of $\lambda$ but which depend on  $\mathscr{C}$ and conditions (GS) and (EG), such that for all $\lambda>\lambda_\star$ the following holds:\\
 
 Let 
$J_\kpar: L^2(\Sigma) \mapsto l^2(N_0;\C^2)\oplus \textrm{span}\{p_{I,\kpar}^\lambda[n]\}^\perp$ be defined,
 via \eqref{Fsplit}, $F\mapsto \left(\{\alpha_n^I[F]\},F_\perp\right)^\top$.
 
\nit Then, uniformly in  $\kpar\in[0,2\pi]$, we have  
\begin{equation}
 \Big\|\ \Big(\ \rho_\lambda^{-1}\left(H_{_{\rm{edge}}}^\lambda(\kpar)-E_0^\lambda\right)\ -\ z{\rm Id}\ \Big)^{-1}\ 
- \ J_\kpar^*\Big(\ H^{^{\rm TB}}_\sharp(\kpar)\ -\ z{\rm Id}\ \Big)^{-1}\ J_\kpar\ \Big\|_{_{L^2(\Sigma)\to L^2(\Sigma)}}\ \le\ C_\star\ e^{-c\lambda}.\label{s-conv}
\end{equation}
\end{theorem}

\medskip

In preparation for our theorem on edge states, we introduce the functions:
 \begin{equation}
  \zeta(\kpar)=1+e^{i\kpar},\quad  \delta_{\rm gap}(\kpar)=\Big|1-|\zeta(\kpar)|\Big|\ge0,\quad \delta_{\rm max}(\kpar)=1+|\zeta(\kpar)|\ .\label{zdgdm}\end{equation} 
We note that for $\kpar\in[0,2\pi]$ that $\delta_{\rm gap}(\kpar)=0$ if and only if $\kpar\in\{2\pi/3,4\pi/3\}$.

\begin{theorem}[Zigzag Edge States]\label{main-thm1}
Assume that $E_0^\lambda$, the ground state energy of the atomic Hamiltonian, $\Hatom=-\Delta+\lambda^2V_0$, satisfies the conditions (GS) \eqref{GS} and (EG) \eqref{EG} on the ground state energy and energy-gap, respectively.
Let $\mathscr{I}$ denote an arbitrary compact subinterval of quasi-momenta:
\begin{equation}
\mathscr{I}\subset\subset (2\pi/3,4\pi/3)\setminus\{\pi\} .
\label{kpar-int}
\end{equation}
Thus, $\min_{\kpar\in\mathscr{I}}\delta_{\rm gap}(\kpar)>0$.\medskip

\nit There exists $\lambda_\star=\lambda_\star(\mathscr{I})>0$ sufficiently large, such that for all $\lambda>\lambda_\star$  the following holds:
\begin{enumerate}
\item 
There is a mapping $\kpar\in\mathscr{I}\ \mapsto\ (E^\lambda(\kpar),\psi^\lambda_\kpar)$, from parallel quasimomenta $\kpar$ to simple eigenpairs of the family of the $\kpar-$ edge state eigenvalue problem
   \eqref{zz-evp}:
 \begin{align}
 H_{_{\rm edge}}^\lambda(\kpar) \psi_\kpar\ &=\ E^\lambda(\kpar)\  \psi_\kpar^\lambda,\quad \psi_\kpar\in L^2(\Sigma)\label{es-evp}\\
 E^\lambda(\kpar)\ &=\ E_0^\lambda+\rho_\lambda\ \Omega^\lambda(\kpar),
 \nn\end{align} where
$ \left|\ \Omega^\lambda(\kpar)\ \right|\ \lesssim\ e^{-c\lambda}$  with $c>0$ independent of $\lambda$. 
Correspondingly, the eigenvalue problem \eqref{zz-evp} is solved by the states $\Psi^\lambda_\kpar(\bx)=e^{i\frac{\kpar}{2\pi}\ktilde_2\cdot\bx}\psi_\kpar^\lambda(\bx)$. 
 \medskip
 %
 %
 \item The edge states $\psi^\lambda_\kpar\in L^2_\kpar(\Sigma)$ are approximated to within $\mathcal{O}(e^{-c\lambda})$  error in $L^2(\Sigma)$ as:
 \begin{equation}
  \psi^\lambda_\kpar(\bx)\ =\ \sum_{n\ge0}\ \alpha_A^n\ p_{A,\kpar}^\lambda[n](\bx)\ +\
 \sum_{n\ge0}\ \alpha_B^n\ p_{B,\kpar}^\lambda[n](\bx)\ +\ \mathcal{O}_{L^2(\Sigma)}(e^{-c\lambda}),
 \label{edge-exp}\end{equation}
 where $c>0$ is independent of $\lambda$. 
Here, 
$
\psi_\kpar^{\rm TB,bd}\ \equiv\ \{\ 
 \left(\alpha^n_A , \alpha^n_B\right)^\top\ \}_{n\ge0}\ \in\ l^2(\N_0;\C^2)$,\ $ \|\psi_\kpar^{\rm TB,bd}\|_{_{l^2(\N_0;\C^2)}}\ =\ 1$
 is a zero energy normalized eigenstate of the limiting tight-binding  edge Hamiltonian;\ $ H^{^{\rm TB}}_\sharp(\kpar)\ \psi_\kpar^{\rm TB,bd}=0$.
 See Theorem \ref{zz-spec} in Section \ref{TB}.
 \end{enumerate}
\end{theorem}

\begin{remark}[Symmetry of edge state curves]\label{es-curves}
 \nit Let $\kpar\in[0,\pi]$. If $\left(E^\lambda(\kpar),\Psi^\lambda_\kpar(\bx)\right)$ is an eigenpair of the $\kpar-$ edge state eigenvalue problem, then $\left(E^\lambda(\kpar),\overline{\Psi_\kpar^\lambda(\bx)}\right)$ is an eigenpair of the $2\pi-\kpar$ edge state eigenvalue problem. 
  \end{remark}
  \begin{remark}[Non-flatness of band]\label{noflatband}
  The large $\lambda$ edge states of eigenfrequencies, $E^\lambda(\kpar)$, in Theorem \ref{main-thm1} arise from the {\it flat band} of edge states, $\Omega(\kpar)=0$ for $2\pi/3<\kpar<4\pi/3$, of the tight-binding Hamiltonian, $H^{^{\rm TB}}_\sharp(\kpar)$. Although $E^\lambda(\kpar)$ has only exponentially small variation, we do not expect $E^\lambda(\kpar)$ to be identically constant. Indeed, numerical simulations illustrate the weak variation in $\kpar$ \cite{TWL:19}.
  \end{remark}
 \begin{remark}[Regularity]\label{regularity}
 We do not address the question of smoothness of $\kpar\in\mathscr{I}\mapsto \left(E^\lambda(\kpar),\psi^\lambda_\kpar\right)\in\R\times L^2(\Sigma)$ in the present article. We believe however that the methods of \cite{FLW-CPAM:17} may be adapted to show that this mapping extends as an analytic mapping in a complex neighborhood of $\mathscr{I}$ from which derivative bounds, {\it e.g.} on $E^\lambda(\kpar)$  ($\kpar\in\mathscr{I}$) can be derived via Cauchy estimates.
 \end{remark}
%
%
\begin{remark}\label{Psi-support}
In Theorem \ref{zz-spec} we find: $\psi_\kpar^{\rm TB,bd}=\sqrt{1-|\zeta(\kpar)|^2}\Big(\ [-\zeta(\kpar)]^n,0\ \Big)^\top$. Therefore, at leading order, $\Psi^\lambda_\kpar(\bx)$ is concentrated about the $A-$ sublattice, $\Lambda_A$:
\begin{equation}
\Psi^\lambda_\kpar(\bx)\ =\ \sqrt{1-|\zeta(\kpar)|^2}\ \sum_{n\ge0}\ [-\zeta(\kpar)]^n\ P_{A,\kpar}^\lambda[n](\bx)\ +\  \mathcal{O}_{L^2_\kpar}(e^{-c\lambda}). 
\label{edge-exp1}\end{equation}
\end{remark}
\begin{remark}\label{zak}
As noted in our discussion of the tight-binding model in Section \ref{TB} (Remark \ref{zak1}) the constraint of Theorem \ref{main-thm1} on parallel quasimomenta:
$\kpar\in(2\pi/3,4\pi/3)$ ($|\zeta(\kpar)|<1$) corresponds to the non-vanishing of the {\it Zak phase}. This is discussed further in Remark \ref{zak1}.
\end{remark}
\begin{remark}\label{armchair}
In work in progress we show, for a sharp termination of the bulk honeycomb structure along an {\it armchair edge},
that there are no edge states in an energy range about $E_0^\lambda$.
 In this case, the relevant Zak phase  vanishes 
for all $\kpar\in [0,2\pi]$.
\end{remark}

\subsection{Relation to previous work}\label{semi-classical}
Tight-binding limits arising from general distributions 
of potential wells has been discussed in the book of  \cite{DiMassi-Sjoestrand:99}
  as well as \cite{Outassourt:87,Carlsson:90}. 
There is extensive related earlier work on the semiclassical limits and methods 
{\it e.g.} \cite{Simon:83,Simon:84a,Simon:84b,Helffer-Sjoestrand:84,Helffer:88,Chantelau:90,Mohamed:91,Daumer:94,Daumer:96}.
The above works are based on detailed semiclassical (WKB) approximations
for potential wells which are assumed to have non-degenerate local minima.
In contrast, in the present article our essential assumptions are only on the 
ground state energy (GS) and spectral gap (EG) of the atomic Hamiltonian, 
$H^\lambda_{atom}$ for large $\lambda$. The relation of the continuum periodic Schroedinger operator 
with a magnetic field to 
 tight-binding models, such as the Harper model, is studied for example in \cite{HS:88}.
%

For convenience we have restricted attention here to $C^\infty$ potentials. However, we believe it will be easy to make small changes in our proofs, to apply our results to nonsmooth potentials of interest. In particular, the atomic potential may be taken to have a Coulomb singularity at the origin, or to have the form $V_0(x)=-1$ for $\bx$ in a ball $B$, $V_0(x)=0$ otherwise. Here the radius of $B$ is taken small enough to satisfy the hypothesis $(PW_2)$  of Section \ref{atomic-gs} below. Examples of artificial graphene, in which  experiments are performed, are  periodic honeycomb arrays of identical microfeatures, say small discs, with one dielectric constant inside the discs and a second dielectric constant outside the discs. Hence, compactly supported atomic potentials are a natural model; see, for example, \cite{artificial-graphene:11,Rechtsman-etal:13a,2014LuJoannopoulosSoljacic,irvine:15,2017KhanikaevShvets,Marquardt:17,ozawa_etal:18}.

For smooth atomic potentials $V_0$ with nondegenerate minima, the general semiclassical works  in  \cite{DiMassi-Sjoestrand:99,Outassourt:87,Carlsson:90}  lead to an ``interaction matrix'', which defines an  operator. In the case of periodic potentials, this can be used to compute relevant dispersion surfaces modulo exponentially small errors. These works do not assert that Dirac points form; indeed, much of the work is in the setting of a square lattice, which does not give rise to Dirac points.  However, we believe that these methods  are powerful enough to deal with Dirac points of honeycomb lattice potentials, when they are combined with the consequences of special symmetry properties of the honeycomb. The essential requirement for the semiclassical analysis approach is that the atomic potential is smooth and has a nondegenerate minimum. Another aspect of the general work of this semiclassical work is that atomic potentials are not assumed to be of compact support and the interaction matrix (hopping coefficients) are obtained in terms
 of the Agmon metric. Finally, the consideration of edge states and the spectrum for honeycombs with line defects  is not within the scope of \cite{DiMassi-Sjoestrand:99,Outassourt:87,Carlsson:90}.

%
\begin{remark}\label{domain-wall}
A different class of line-defects of  great interest in the study of topologically protected edge states is the class of  {\it domain walls}.
In our previous work,  motivated by \cite{HR:07,RH:08,Soljacic-etal:08}, domain walls are realized by starting with two periodic structures at `` $ +\infty$ '' and `` $-\infty$ '',  with a common spectral gap and phase-shifted from one another, and connecting them across a line-defect at which there is no phase-distortion. See the analytical work in 1D \cite{FLW-PNAS:14,FLW-MAMS:17,DFW:18} and 2D  \cite{FLW-2d_edge:16,FLW-2d_materials:15,LWZ:18} as well as theoretical and experimental work on photonic realizations \cite{Thorp-etal:15,CLEO-Thorp_etal:16,CLEO-Poo_etal:16}.
\end{remark}
\begin{remark}\label{q-graph}
Quantum graphs \cite{BK:13} are another class of discrete models in condensed matter, electromagnetic and other systems; see also, for example, \cite{Shipman:17,Becker-Zworski:18,BHJ:18}. An extensive discussion of edge states for nanotube structures in the setting of quantum graphs is given in \cite{Kuchment-Post:07,Do-Kuchment:13}. It would be of interest to investigate a relation between the edge modes of these models and continuum models.
\end{remark}

 \begin{figure}
\centering 
\includegraphics[width=.75\textwidth]{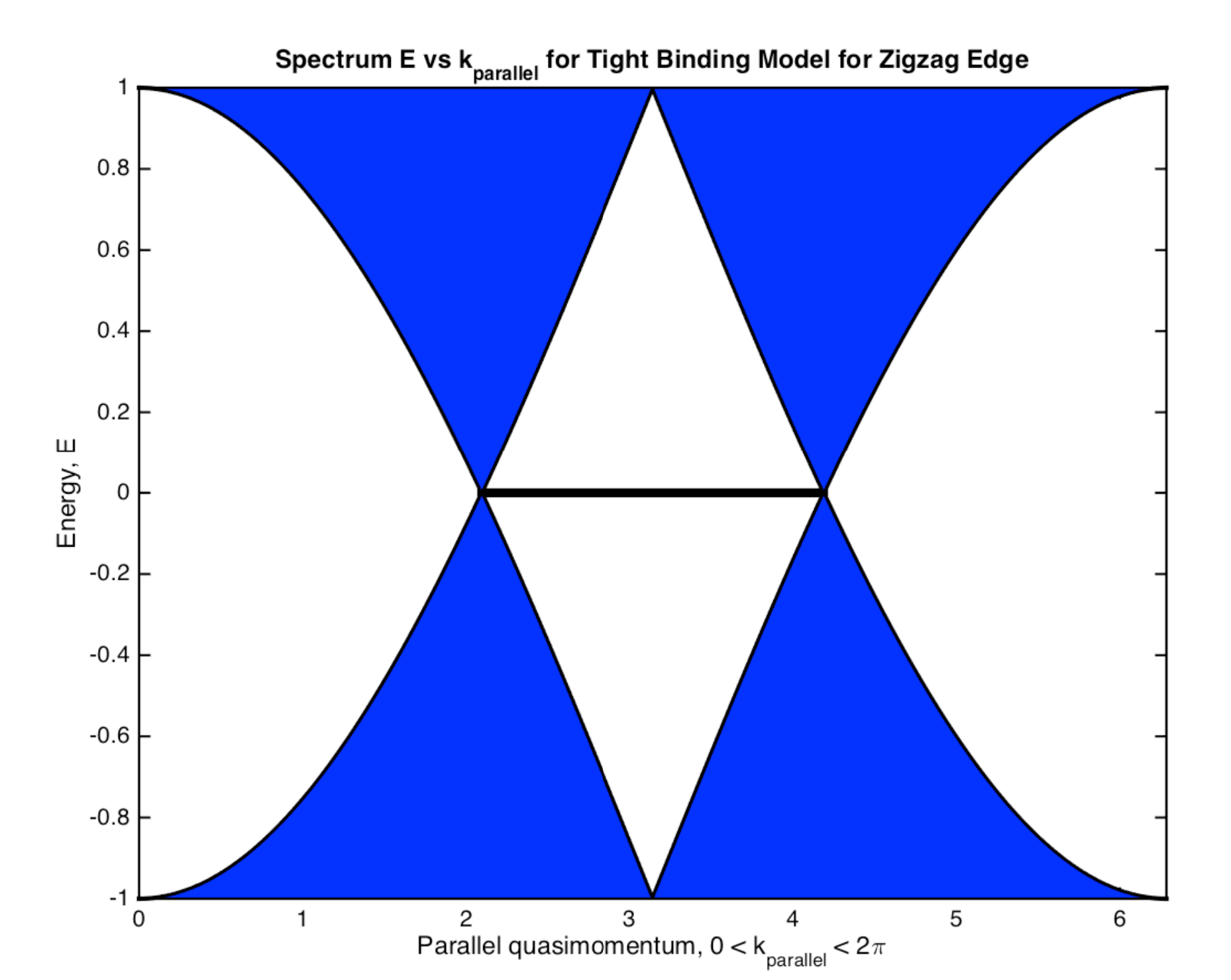}
\caption{\small Spectrum of tight-binding Hamiltonian  $\HTBs(\kpar)$, for $0\le \kpar\le2\pi$, described in Theorem \ref{zz-spec}.  This spectrum contains a {\it flat band} of zero energy states; $\cHeg(\kpar)$ has an isolated simple $0-$ energy eigenstate for $2\pi/3\le \kpar\le 4\pi/3$.  Shaded regions consist of continuous spectrum. For sufficiently large $\lambda$, the low-lying part of the spectrum of $-\Delta+\lambda^2V_\sharp-E_0^\lambda$,  after rescaling by $\rho_\lambda$, is approximated by 
spectrum of the 2-band model $\HTBs$; see Theorem \ref{main-thm1}.}
\label{tb-spec}
\end{figure}

\subsection{Outline of the paper}\label{quick} We present a brief outline.
\begin{enumerate}
\item[] Section \ref{TB} discusses tight binding models; first, the tight binding model for bulk,  and then the tight binding model for a honeycomb structure terminated along a zigzag edge. \\ 
\item[] Section \ref{setup} first introduces the atomic Hamiltonian $\Hatom=-\Delta+V_0$, where $V_0$ is a potential well whose support in a sufficiently small disc about the origin, and such that $V_0$ satisfies
some basic general assumptions $(PW_1)-(PW_4)$. The bulk honeycomb structure is defined by 
 $\cH=-\Delta+\lambda^2V$, where $V$ is the periodic potential by summing $V_0$ over translates of the atomic potential, $V_0$, over the  honeycomb structure. Thus $V$ consists of a potential well $V_0$ centered at each site of the honeycomb.  Finally the edge Hamiltonian, $H_{\rm edge}=-\Delta+\lambda^2V_\sharp$, which acts on $L^2(\R^2)$, has potential $V_\sharp$ which is identically equal to $V$ on a half-space with a zigzag edge 
  and zero on the other side of this zigzag edge. (We shall also work with the translated edge Hamiltonian $\cHeg=H_{\rm edge}-E_0^\lambda$.) The edge state eigenvalue problem for parallel quasi-momentum $\kpar$ is then stated on $L^2(\Sigma)$, where $\Sigma$ is the infinite cylinder \eqref{zz-cylinder}. \\
  \item[]  Section \ref{nat-basis} introduces a natural basis for approximating the 2 lowest lying bands of $\cHeg$ for $\lambda$ sufficiently large. 
  This basis consists of functions, $\{p_{\bI,\kpar}^\lambda[n](\bx):\bI = A,B, n\ge0\}$ on $\Sigma$, which are pseudo-periodic (with respect to the direction parallel to the edge) infinite sums of atomic orbitals. \\
  \item[]  Section \ref{en-est} establishes energy estimates on $\cHeg$ which imply invertibility of 
 $\cHeg$  on $\mathscr{X}_{AB}(\kpar)$, the orthogonal complement of the  
   orbital subspace: ${\rm span}\{ p^\lambda_{I,\kpar}[n]:n\ge0, I=A,B\}$. This implies that the resolvent of
    $\cHeg$ is well-defined and bounded on $\mathscr{X}_{AB}(\kpar)$ .\\
   \item[]  Section \ref{ls-red} implements a Lyapunov-Schmidt / Schur complement reduction strategy: The spectral problem on $L^2(\Sigma)={\rm span}\{ p^\lambda_{I,\kpar}[n]:n\ge0, I=A,B\}\oplus\mathscr{X}_{AB}(\kpar)$ is reduced, using the resolvent bounds on $\mathscr{X}_{AB}(\kpar)$, to an equivalent problem on the  space ${\rm span}\{ p^\lambda_{I,\kpar}[n]:n\ge0, I=A,B\}$. This problem depends nonlinearly on the eigenvalue parameter $E=E_0^\lambda+\rho_\lambda\Omega$ and is of the form of an infinite algebraic system:
   \[\sum_{I=A, B}\sum_{n\ge0}\mathcal{M}^{\lambda,\kpar}_{JI}[m,n](\Omega,\kpar)\ \alpha_n^I\ =\ 0;\ \ J=A, B,\ \ m\ge0\]
 for $(\Omega,\alpha)$,  where
   $\alpha=\{\alpha_n^I\}_{n\ge0, I=A,B}\in l^2(\N_0;\C^2)$ are coordinates relative to the basis
   $\{ p^\lambda_{I,\kpar}[n]:n\ge0, I=A,B\}$.
    \\
    \item[] Section \ref{lin-els}  summarizes the required properties of $\mathcal{M}^{\lambda}(\Omega,\kpar)$ acting in $ l^2(\N_0;\C^2)$. We write $\mathcal{M}^{\lambda,\kpar}(\Omega,\kpar)=\mathcal{M}_{\rm lin}^{\lambda}(\Omega,\kpar)-\mathcal{M}_{\rm nlin}^{\lambda}(\Omega,\kpar)$, separating matrix elements contibutions which are linear in $\cHeg(\kpar)$ and those which are nonlinear in $\cHeg(\kpar)$.  We have $\mathcal{M}_{\rm lin}^{\lambda}(\Omega,\kpar)=\rho_\lambda \HTBs(\kpar)\ +\ \mathcal{O}_{l^2\to l^2}(\rho_\lambda e^{-c\lambda})$ (Proposition \ref{zz-els}) and $\mathcal{M}_{\rm nlin}^{\lambda}(\Omega,\kpar) =\mathcal{O}_{l^2\to l^2}(\rho_\lambda e^{-c\lambda})$ (Proposition \ref{zz-nl-els}). These propositions
     are proved in later sections.\\
        \item[] Section \ref{zz-sb} proves Theorem \ref{main-thm1}, the existence of edge states, bifurcating from the flat band of eigenstates of $\HTBs$, via our formulation of the eigenvalue 
         in $ l^2(\N_0;\C^2)$. \\
    \item[] Section \ref{sc-res} proves Theorem \ref{res-conv}, the convergence of a translation and scaling of the resolvent of $\cHeg(\kpar)$ to that of $\HTBs(\kpar)$.\\
         \item[]  Section \ref{res-kernel} is the most technically involved and introduces  techniques not present in our earlier work. Theorem \ref{Ksharp} is  a \underline{pointwise} estimate on the resolvent kernel of $\cHeg-z=H_{\rm bulk}^\lambda-(E_0^\lambda+z)$, $z$ small,  when restricted to the orthogonal complement of ${\rm span}\{p^\lambda_{I,\kpar}[n]: n\ge0, I=A, B\}$.
These bounds are stated in Theorem \ref{Ksharp}. We first, in Proposition  \ref{summary-K}, establish these kernel estimates for potentials which are a sum of atomic potentials centered on an  arbitrary discrete set of lattice sites $\Gamma\subset\R^2$  (not necessarily translation invariant), whose minimal pairwise distance
is $M r_0$, where $r_0$ is the radius of the support of $V_0$ and $M>2$ is some positive constant. 
We then specialize to a translation invariant set to obtain  Theorem \ref{Ksharp}.
   \\
              \item[] Section \ref{m-els-est} expands the linear matrix elements, $\mathcal{M}_{\rm lin}^{\lambda}(\Omega,\kpar)$,  in terms
               of $\HTBs(\kpar)$ and estimates the corrections, proving Proposition \ref{zz-els}.\\
     \item[] Section \ref{nl-els} estimates the nonlinear matrix elements, $\mathcal{M}_{\rm nlin}^{\lambda}(\Omega,\kpar)$, proving Proposition \ref{zz-nl-els}.\\
     \item[] Finally, there are two appendices. Appendix \ref{app:error-main} introduces a technical tool used to construct  the resolvent of $\cHeg$ on $\mathscr{X}_{AB}^\lambda(\kpar)$.  Appendix \ref{overlap} contains general results on overlap integrals enabling expansion of $\mathcal{M}_{\rm lin}^{\lambda}(\Omega,\kpar)$, for $\lambda$ large, estimate corrections. 
\end{enumerate}

\subsection{Notation}\label{notation}

\begin{enumerate}
\item $\mathbb{N}=\{1,2,3,\dots\}$,\ \ $\mathbb{N}_0=\{0,1,2,3,\dots\}$\ .
\item When we write the expression $g_\eps= \mathcal{O}_{X}(\gamma_\eps)$ as $\eps\to\eps_0\in\R\cup\{\infty\}$, we  mean that there exists $C>0$, independent of $\eps$,  such that $\|g_\eps\|_X\le C\gamma_\eps$ as $\eps\to\eps_0$.
\item We shall be concerned with the asymptotic behavior of many expressions, $a(\lambda), b(\lambda),\dots$, in the regime where the parameter $\lambda$ sufficiently large.
The relation $a(\lambda)\lesssim b(\lambda)$ means that there is a constant $C$, which can be taken to be independent of $\lambda$, such that for all $\lambda$ sufficiently large: $a(\lambda)\le C b(\lambda)$.
\item 
$ \Lambda=\Z\vtilde_1 \oplus \Z\vtilde_2,$
 the equilateral triangular lattice, is  generated by the basis vectors $\vtilde_1$ and $\vtilde_2$, displayed in \eqref{c-v12-def}.
 \item $\bfm\vec\vtilde=m_1\vtilde_1+m_2\vtilde_2$, where $\bfm=(m_1,m_2)\in\Z^2$.
\item   $\Lambda^* =\ \Z \ktilde_1\oplus \Z\ktilde_2$, the dual lattice, spanned by the dual basis vectors
$\ktilde_1$ and $\ktilde_2$, displayed in \eqref{c-k12-def}.
 Note that $\ktilde_{\ell}\cdot \vtilde_{{\ell'}}=2\pi\delta_{\ell{\ell'}}$.
 \item We remark that alternative bases for $\Lambda$ and  $\Lambda^*$
  (used for example in \cite{FW:12,FLW-CPAM:17}) are:
 \begin{align}
 \bv_1 &=\vtilde_1,\quad \bv_2=\vtilde_1-\vtilde_2\nn\\
 \bk_1 &=\ktilde_1+\ktilde_2,\quad \bk_2=-\ktilde_2 .\nn
 \end{align}
We have $\Lambda=\Z\bv_1 \oplus \Z\bv_2$, $\Lambda^* =\ \Z\bk_1\oplus \Z\bk_2$
 and $\bk_{\ell}\cdot \bv_{{\ell'}}=2\pi\delta_{\ell{\ell'}}$.
 \item $\mathbb{H}$, Honeycomb structure; see \eqref{honey}.
\item  $\mathbb{H}_\sharp$, Zigzag-truncated honeycomb structure; see \eqref{bbH+}.
 \item $\Sigma=\R^2/\Z\vtilde_2$, the cylinder with $\Omega_\Sigma$, a choice of  fundamental cell for $\Sigma$; see Figure \ref{zz-dimers}.
 \item $L^2_\kpar=L^2_\kpar(\Sigma)$, functions $f$ such that $f(\bx+\vtilde_2)=e^{i\kpar}f(\bx)$ for almost all $\bx$, and
  \[ \|f\|_{_{L^2_\kpar}}^2\ =\ \int_{_{\Omega_\Sigma}}|f|^2<\infty.\]
  \item $\mathscr{H}^{(\omega)}\equiv L^2(\R^2;e^{\gamma |\bx-\omega|}\ d\bx)$, exponentially weighted $L^2$ space. 
   \item $\mathcal{B}(X)$ denotes the space of bounded  linear operators on $X$. 
 \item $\Gfree(\bx,\by)$ denotes the free Green's function defined in \eqref{Gfree-def}. 
\item  $\Gatom(\bx,\by)$ denotes the atomic Green's function defined in \eqref{Gat-def}.
 \item Hamiltonians:\\ 
$\Hatom=-\Delta+\lambda^2 V_0(\bx)$, the atomic Hamiltonian with ground state energy $E_0^\lambda$ \\
$H^\lambda_{_{\rm bulk}}=-\Delta +\lambda^2 V(\bx)$ and 
  $H^\lambda_{_{\rm edge}}=-\Delta +\lambda^2 V_\sharp(\bx)$,
denote  bulk and edge Hamiltonians acting in $L^2(\R^2)$ \\
   $\cHeg=H^\lambda_{_{\rm edge}}\ -\ E_0^\lambda$, the centered edge Hamiltonian, acting in $L^2_\kpar$\\
      $\tHeg\ =\ (\rho_\lambda)^{-1}\ \cHeg$, the scaled and centered edge Hamiltonian acting in $L^2_\kpar$\\\
 $\HTBs(\kpar)$, the tight-binding edge Hamiltonian, acting in $l^2(\N_0;\C^2)$; see Definition \ref{HTBeg}.
\end{enumerate}

\bigskip

\nit{\bf Acknowledgements:}\ The authors wish to thank Gian Michele Graf and Alexis Drouot for very stimulating discussions. We would also like to thank Bernard Helffer for correspondence concerning previous general results on tight-binding limits.
Part of this research was done while MIW was Bergman Visiting Professor at Stanford University. CLF and MIW wish to thank the Department of Mathematics at Stanford University for its  hospitality. 
This research was supported in part by 
National Science  Foundation grants DMS-1265524 (CLF) and  DMS-1412560, DMS-1620418 and Simons Foundation Math + X Investigator Award \#376319 (MIW).

\section{Tight-binding}\label{TB}

Consider a tiling of the entire plane, $\R^2$, by parallelograms of the sort shown in Figure \ref{zz-dimers}
Each parallelogram has exactly two points of $\mathbb{H}$. This is a particular {\it dimerization} of $\mathbb{H}$. 
 We assign the label $(n_1,n_2)$
 to the parallelogram which contains $\bv_A^{(n_1,n_2)}=\bv_A+n_1\vtilde_1+n_2\vtilde_2$ and $\bv_B^{(n_1,n_2)}=\bv_B+n_1\vtilde_1+n_2\vtilde_2$. To the sites $\bv_A^{(n_1,n_2)}$ and $\bv_B^{(n_1,n_2)}$ we assign  complex amplitudes
$ \psi^A_{_{n_1,n_2}}$ and 
$ \psi^B_{_{n_1,n_2}}$ and form the tight binding wave function:

\[
\psi_{_{n_1,n_2}}=\begin{pmatrix} \psi^A_{_{n_1,n_2}} \\  \psi^B_{_{n_1,n_2}}\end{pmatrix}
\]

\subsection{$\HTBb$, the tight-binding bulk Hamiltonian}\label{TB-bulk}
%
%

 The bulk tight binding Hamiltonian can be represented with respect to the above dimerization.  Starting with any dimerization  would give a unitarily equivalent operator on $l^2(\Z^2;\C^2)$.  
The nearest neighbor tight binding bulk Hamiltonian, relative to the dimerization of $\mathbb H$ in Figure \ref{zz-dimers} is:
\begin{align}\label{zz-tbevp}
\left[\ \HTBb\psi\ \right]_{_{\none,\ntwo}}\ =\ 
\begin{pmatrix} 
\left(\ \HTBb\psi\ \right)^A_{_{\none,\ntwo}}\\
 \left(\ \HTBb\psi\ \right)^B_{_{\none,\ntwo}}
 \end{pmatrix}
= \begin{pmatrix}
\psi_{_{n_1-1,n_2}}^B\ +\ \psi_{_{n_1,n_2-1}}^B\ +\ \psi_{_{n_1,n_2}}^B\\
\psi_{_{n_1+1,n_2}}^A\ +\ \psi_{_{n_1,n_2+1}}^A\ +\ \psi_{_{n_1,n_2}}^A
\end{pmatrix}
\end{align}
where $n_1,n_2\in\mathbb{Z}$. The operator $H_{_{\rm bulk}}^{\rm TB}$ is a bounded self-adjoint 
 linear operator on $l^2(\Z^2;\C^2)$ and was 
introduced in \cite{Wallace:47}. The spectrum of $\HTBb$  consists of two spectral bands which touch conically at Dirac points over the vertices of $\mathcal{B}$. The approximation and convergence as $\lambda$ increases of the low-lying dispersion surfaces and the resolvent $H_{_{\rm bulk}}^\lambda$ acting on $L^2(\R^2)$ to those of $H_{_{\rm bulk}}^{\rm TB}$ acting on $l^2(\Z^2;\C^2)$
  was studied in  \cite{FLW-CPAM:17} .

\subsection{Tight-binding Hamiltonian for the zigzag edge }\label{tb-zz}

{\ }\medskip

Our goal in this section is to introduce a tight-binding edge Hamiltonian which will act on functions 
$\psi\in l^2\left( (\N_0\times\Z) ; \C^2 \right)$ defined on the vertices of $\mathbb{H}_\sharp$. We shall do this by first expressing $ H_{\rm bulk}^{^{\rm TB}}$, as a direct integral over $\kpar$ of fiber operators $ H_{\rm bulk}^{^{\rm TB}}(\kpar)$
  acting on states which are ``$\kpar$- pseudo-periodic'' with respect to one lattice direction  and square-summable with respect to the other lattice direction. The edge Hamiltonian $H^{^{\rm TB}}_\sharp$ is then obtained from $ H_{\rm bulk}^{^{\rm TB}}(\kpar)$ by appropriate restriction to functions defined on $\mathbb{H}_\sharp$. 

Since the truncated structure $\mathbb{H}_\sharp$ and its subset  edge vertices are invariant with respect to translation  by $\vtilde_2$, we introduce $\kpar\in S^1=\R/2\pi\Z$, the parallel quasi-momentum associated with this translation invariance. 
For each $\kpar\in[0,2\pi]$, we  refer to a state as being $\kpar-$ pseudo-periodic if:
 \begin{equation}
  \psi_{n_1,n_2+1}=e^{i\kpar}\psi_{n_1,n_2},\ \ n_1\ge0,\ \ n_2\in\Z.
\label{kpseu}  \end{equation}

Functions  $\psi=\{\psi_{_{n_1,n_2}}\}\in l^2(\Z;\C^2)$  may be expressed via the discrete Fourier transform as 
\begin{equation}
\psi_{_{n_1,n_2}}\ =\ (2\pi)^{-1}\int_0^{2\pi}\ e^{i\ntwo\kpar}\psi_{_{\none}}(\kpar)\ d\kpar,
\label{fib-int}\end{equation}
as a superposition over states  $\{e^{i\ntwo\kpar}\psi_{_{\none}}(\kpar)\}$ which are square-summable over $\Z$ with respect to $n_1$ and which satisfy \eqref{kpseu}.

Therefore, the tight binding bulk Hamiltonian $H^{^{\rm TB}}_{\rm bulk}$ may be reduced to the $\kpar-$ dependent fiber (Bloch) Hamiltonians, $\HTBb(\kpar):l^2(\Z;\C^2)\to l^2(\Z;\C^2)$, defined by
\begin{align}\label{zz-ham}
\left[ \HTBb(\kpar)\psi \right]_{_{\none}}\ &\equiv 
\begin{pmatrix}
 \psi_{_{n_1-1}}^B\ +\  \left(1+e^{-i\kpar}\right)\ \psi_{_{n_1}}^B\\
\psi_{_{n_1+1}}^A\ +\  \left(1+e^{+i\kpar}\right)\ \psi_{_{n_1}}^A 
\end{pmatrix},\nn\\
&= \begin{pmatrix}0&1\\0&0\end{pmatrix}\begin{pmatrix}\psi_{_{n_1-1}}^A\\ \psi_{_{n_1-1}}^B\end{pmatrix}
\ +\ \begin{pmatrix}0&1+e^{-i\kpar}\\1+e^{+i\kpar}&0\end{pmatrix}\begin{pmatrix}\psi_{_{n_1}}^A\\ \psi_{_{n_1}}^B\end{pmatrix}
\ +\ \begin{pmatrix}0&0\\1&0\end{pmatrix}\begin{pmatrix}\psi_{_{n_1+1}}^A\\ \psi_{_{n_1+1}}^B\end{pmatrix}\ .
\end{align}

\bigskip

Finally, we define the tight-binding edge Hamiltonian, $H_\sharp^{\rm TB}$.
 For $\psi=(\psi_0,\psi_1,\psi_2,\dots)\in l^2(\N_0;\C^2)$, introduce the extension operator: 
\begin{align}
&\iota\ :\ l^2(\N_0;\C^2)\to l^2(\Z;\C^2)\nn\\
& \iota\psi\ =\ (\dots,0,0,0,\psi_0,\psi_1,\psi_2,\dots)\in l^2(\Z;\C^2)\ .
\nn\end{align}
The adjoint of $\iota$ is the restriction operator defined on 
$\phi=(\dots,\phi_{-2},\phi_{-1},\phi_0,\phi_1,\phi_2,\dots)\in l^2(\Z;\C^2)$ by:
\begin{align}
&\iota^*\ :\ l^2(\Z;\C^2)\to l^2(\N_0;\C^2)\ ,\nn\\
& \iota^*\phi\ =\ (\phi_0,\phi_1,\phi_2,\dots)\in l^2(\N_0;\C^2)\ .
\nn\end{align}

\begin{definition}\label{HTBeg}
The tight-binding edge fiber operators, $H_\sharp^{^{\rm TB}}(\kpar)$, and edge Hamiltonian $H_\sharp^{\rm TB}$ are given by
\begin{equation}
\HTBs(\kpar)\ =\ \iota^*\ \HTBb(\kpar)\ \iota\ :\ l^2(\N_0;\C^2)\to l^2(\N_0;\C^2)\ 
\label{HTBegk}
\end{equation}
and
\begin{equation}
\HTBs =\  \int_{[0,2\pi]}^\oplus\ \HTBs(\kpar)\ d\kpar\ :\ l^2(\N_0\times\Z)\to l^2(\N_0\times\Z)\ .
\label{HTBeg}
\end{equation}
\end{definition}

\subsection{Spectrum of $\HTBs(\kpar)$}\label{TB-spec}

Define, for $\kpar\in[0,2\pi]$,  the functions 
\begin{align}
\zeta(\kpar)&\equiv 1+e^{i\kpar}\quad ,\label{zeta-def	1}\\
\delta_{\rm gap}(\kpar)&\equiv \min_{k_\perp\in[0,2\pi]}\Big| 1+e^{i\kpar}+e^{i\kperp} \Big|\ =\ \Big|\ 1-|\zeta(\kpar)|\ \Big|,\label{d-gap}\\
\delta_{\rm max}(\kpar)&\equiv 1+|\zeta(\kpar)| .
\label{d-max}\end{align}
Note $\delta_{\rm gap}(2\pi/3)=\delta_{\rm gap}(4\pi/3)=0$, $\delta_{\rm gap}(\kpar)>0$ otherwise in $[0,2\pi]$, and that
 $|\zeta(\kpar)|<1$ for $\kpar\in(2\pi/3,4\pi/3)$. 
We next prove that the spectrum of $\HTBs(\kpar)$ is as displayed in Figure \ref{tb-spec}.
Let us enumerate the coordinates of the vector in $l^2(\N_0;\C^2)$, $\psi=\Big\{\begin{pmatrix}\psi_n^A\\ \psi_n^B\end{pmatrix}\Big\}_{n\ge0}$, by $\psi=(\psi_0 ^A,\psi_0^B,\psi_1 ^A,\psi_1^B,\dots)^\top$. We denote the corresponding unit vectors by  $\hat\be_1=(1,0,0,\dots)$,\ 
$\hat\be_2=(0,1,0,\dots)$,\ {\it etc}. 
\begin{theorem}[\ $\sigma(\HTBs(\kpar))$, the spectrum of $\HTBs(\kpar)$ in $l^2(\N_0;\C^2)$] \label{zz-spec}
{\ }\medskip

For each $\kpar\in[0,2\pi]$, 
 $\sigma(\HTBs(\kpar))=\sigma_{\rm pt}(\sigma(\HTBs(\kpar)))\ \cup\ \sigma_{\rm ess}(\sigma(\HTBs(\kpar)))$.

\begin{enumerate}
\item Point spectrum of $\HTBs(\kpar)$:\ 
 \begin{align*}
 \sigma_{\rm pt}(\HTBs(\kpar))\ =\
 \begin{cases} 
 \{0\}\ & \textrm{if}\ \kpar\in(2\pi/3,4\pi/3)\\
 \{-1,0,1\}\ &\ \textrm{if}\ \kpar=\pi\\
\  \ \emptyset\ & \textrm{if}\ \kpar\in[0,2\pi]\setminus(2\pi/3,4\pi/3):
 \end{cases}
 \end{align*}
 In particular, 
 \[\textrm{$\HTBs$ has a zero energy ``flat-band'' of eigenstates over the range $2\pi/3<\kpar<4\pi/3$.}\]

  For $\kpar\in(2\pi/3,4\pi/3)\setminus\{\pi\}$ the point spectrum, which consists eigenvalue $E=0$ is simple. The corresponding  normalized $0-$ energy eigenstate, 
  $\psi^{\rm TB, bd}=\left\{\psi^{\rm TB, bd}_n\right\}_{n\ge0}$,  is given by
\begin{align}
\psi^{\rm TB, bd}_n(\kpar)\ &=\ \sqrt{1-|\zeta(\kpar)|^2}\ \begin{pmatrix} \left(-\zeta(\kpar)\right)^n \\ 0\end{pmatrix},\quad n\ge0\ .
\label{zz-estate1}
\end{align}

For $\kpar=\pi$, $E=0$ is a simple eigenvalue with corresponding normalized $0-$ energy eigenstate is given by:
\begin{equation}
\psi^{\rm TB, bd}_{_{0}}(\pi)\ =\  \begin{pmatrix} 1 \\ 0\end{pmatrix},\quad \psi^{\rm TB, bd}_{_{n}}(\pi)\ =\  \begin{pmatrix} 0 \\ 0\end{pmatrix},\qquad n\ge1\ .\label{pi-state1}
\end{equation}
The eigenvalues $E=+1$ and $E=-1$ have infinite multiplicity.The corresponding eigenspaces are:
\begin{align*}
 {\rm kernel}(\HTBs(\pi)-Id)\ &=\ \Big\{\frac{1}{\sqrt2}\left(\hat\be_{2j+1}+\hat\be_{2j+2}\right):j=0,1,2,\dots\Big\},\\
{\rm kernel}(\HTBs(\pi)+Id)\ &=\ \Big\{\frac{1}{\sqrt2}\left(\hat\be_{2j+1}-\hat\be_{2j+2}\right):j=0,1,2,\dots\Big\}
\end{align*}

\item Essential spectrum of $\HTBs(\kpar)$: \begin{equation}\label{s-gap}
\sigma_{\rm ess}(H^{^{{\rm TB}}}_\sharp(\kpar))\ =\ 
\begin{cases}
\Big\{\  z\in\R : \delta_{\rm gap}(\kpar)\le |z|\le \delta_{\rm max}(\kpar)\ \Big\}, & \kpar\in[0,2\pi]\setminus\{\pi\}\\
\qquad\qquad \emptyset, & \kpar=\pi \ .
\end{cases}
\end{equation} 
\item Resolvent expansion: 
\subitem(a) Let $\kpar\in (2\pi/3,4\pi/3)\setminus\{\pi\}$. Then, for $z\in \C\setminus\sigma_{\rm ess}(H^{^{{\rm TB}}}_\sharp(\kpar))$ 
and $z\ne0$ we have
\begin{equation}
\left(\HTBs(\kpar)-z I\right)^{-1}f\ =\ \frac{1}{z} \left\langle \psi^{\rm TB, bd}(\kpar),f\right\rangle_{_{l^2(\N_0;\C^2)}}\  \psi^{\rm TB, bd}(\kpar)\ +\ \mathscr{G}_{\rm reg}(z;\kpar)f\ .
\label{laurent1}
\end{equation}
Here,  $z\mapsto \mathscr{G}_{\rm reg}(z;\kpar)$ is an analytic mapping from 
$\C\setminus\sigma_{\rm ess}(H^{^{{\rm TB}}}_\sharp(\kpar))$
to the space of bounded linear operators on $l^2(\N_0;\C^2)$. If $(z,\kpar)$ varies over a compact set 
$\Upsilon\subset\subset\R\times[0,2\pi]$ for which ${\rm distance}\left(z,\sigma_{\rm ess}\left(\HTBs(\kpar)\right)\right)\ge b>0$, 
where $b$ is a positive constant depending on $\Upsilon$, then $\|\mathscr{G}_{\rm reg}(z;\kpar)\|_{_{\mathcal{B}(l^2(\N_0;\C^2))}}<B(b)<\infty$. 
\subitem(b) Let $\kpar=\pi$. Then, $\left(\HTBs(\kpar)-z I\right)^{-1}f$ has an expression analogous to \eqref{laurent1} with poles at $z=0$, $z=+1$ and $z=-1$.
\subitem(c) Let $\kpar\in [0,2\pi]\setminus(2\pi/3,4\pi/3)$. Then, for $z\in \C\setminus\sigma_{\rm ess}(H^{^{{\rm TB}}}_\sharp(\kpar))$ 
 we have \begin{equation}\label{reg-res}
\left(\HTBs(\kpar)-z I\right)^{-1}f\ =\  \mathscr{G}_{\rm reg}(z;\kpar)f,
\end{equation}
where $z\mapsto \mathscr{G}_{\rm reg}(z;\kpar)$ is as in part (a).
\item For $\kpar\in(2\pi/3,4\pi/3)$, the equation 
  $\HTBs(\kpar)\psi = f$, where  $f\in l^2(\N_0;\C^2)$,  is solvable for $\psi\in l^2(\N_0;\C^2)$ if and only if $ \left\langle \psi^{\rm TB, bd}(\kpar),f\right\rangle_{_{l^2(\N_0;\C^2)}}=0$.
\end{enumerate}
\end{theorem}
\begin{remark}\label{zak1}   
We remark on the connection between the condition $\kpar\in(2\pi/3,4\pi/3)$ (equivalently $|\zeta(\kpar)|<1$)
and the non-vanishing of a winding number, known as the {\it Zak phase}.
 For fixed $\kpar$, consider the normalized bulk Floquet-Bloch modes of $\HTBb(\kpar)$; see \eqref{zz-ham}. There are two families
 of eigenpairs: $\left(\mu^\pm(\kpar),U^\pm_\none(k_\perp;\kpar)\right)$, where
\begin{align}
 \mu^\pm(\kpar)\ &=\pm|\zeta(\kpar)+e^{ik_\perp}|,\qquad \textrm{ ($\zeta(\kpar)=1+e^{i\kpar}$)}, \nn\\ 
U^\pm_\none(k_\perp;\kpar)\ &=e^{i k_\perp\none}\xi^\pm(k_\perp;\kpar),\ \ \xi^\pm(k_\perp;\kpar)=\frac{1}{\sqrt2}\begin{pmatrix}1\\ \pm j(k_\perp)\ \end{pmatrix},\nn\\
j(e^{ik_\perp})&=\frac{\zeta(\kpar)+e^{ik_\perp}}{|\zeta(\kpar)+e^{ik_\perp}|},\quad j(z)\overline{j(z)}=1 .\nn
\end{align}

For either family of modes (say $+$), we consider the {\it Berry connection} defined by $A(k_\perp;\kpar)\ \equiv\ \left\langle\ \xi(k_\perp;\kpar),\frac{1}{i}\D_{k_\perp}\xi(k_\perp;\kpar)\right\rangle$ and the {\it Zak phase} by $\mathscr{Z}(\kpar)\ \equiv\ \int_0^{2\pi} A(k_\perp;\kpar)\ dk_\perp$. We have
\begin{align*}
\mathscr{Z}(\kpar)\ &=\ -i\int_0^{2\pi}\overline{j(e^{ik_\perp};\kpar)}\ \frac{\D}{\D{k_\perp}}j(e^{ik_\perp};\kpar)\ dk_\perp\\
&\qquad =\ -i\int_{|w|=1}\overline{j(w;\kpar)}\ \D_zj(w;\kpar)\ dw\\
&\qquad =\ -i\int_{|w|=1}\frac{\D_w j(w;\kpar)}{j(w;\kpar)}\ dw\\
&\qquad =\  2\pi \times \textrm{Winding number of $w\in S^1\mapsto j(w;\kpar)\in\C$}\ .
\end{align*}
If $|\zeta(\kpar)|<1$, then  $\mathscr{Z}(\kpar)=2\pi$ and if $|\zeta(\kpar)|>1$, then $\mathscr{Z}(\kpar)=0$.
This is an example 
of the {\it bulk-edge correspondence} (see, for example, \cite{mong2011edge,delplace2011zak,Graf-Porta:13}) and Theorem \ref{main-thm1} establishes its validity in the strong-binding regime.
\end{remark}

\nit {\it Proof of Theorem \ref{zz-spec}:} Fix $\kpar\in[0,2\pi)$ and set $\zeta=\zeta(\kpar)=1+e^{i\kpar}$. We study the operator $\HTBs(\kpar)$ in the Hilbert space $l^2(\N_0;\C^2)$.  
An energy $z$ is in the point spectrum of $\HTBs(\kpar)$ if there exists $\psi\ne0$, $\psi\in l^2(\N_0;\C^2)$ such that
$\HTBs(\kpar)\psi=z\psi$. Written out componentwise, the eigenvalue problem is:
\begin{align}
\psi_{n-1}^B+\zeta^*\psi_n^B\ &=\ z\psi_n^A,\ n\ge0 , \label{tbevp-a}\\
 \psi_{n+1}^A+\zeta\psi_n^A\ &=\ z\psi_n^B,\ n\ge0 ,\label{tbevp-b}
 \end{align}
and $\psi_n= \begin{pmatrix}\psi_n^A\\ \psi_n^B\end{pmatrix}=\begin{pmatrix}0\\ 0\end{pmatrix}$ for all $n\le-1$. 
\medskip

We begin by showing that  for $\kpar\in(2\pi/3,4\pi/3)$, we have that $0\in\sigma_{\rm pt}(\HTBs(\kpar))$ and that for   $\kpar\in[0,2\pi]\setminus(2\pi/3,4\pi/3)$, $z=0$ is not in the point spectrum. Set $E=0$ and observe that equations \eqref{tbevp-a} and \eqref{tbevp-b} become decoupled first order difference equations:  $\psi_{n+1}^A=(-\zeta)\psi_n^A,\ n\ge0$ and $\psi_{n-1}^B= (-\zeta^*)\psi_n^B,\ n\ge0$. 

The equation for $\psi^A$ has the solution:   $\psi_n^A=(-\zeta)^n\psi^A_0$, $n\ge0$, where $\psi^A_0$ can be set arbitrarily. If $\kpar\in(2\pi/3,4\pi/3)$, then $|\zeta(\kpar)|<1$ and hence $\psi_n^A\to0$ exponentially as $n\to\infty$. 
Turning to $\psi^B$, let us first assume that $\kpar\ne\pi$ so that $\zeta(\kpar)\ne0$. In this case,
$ \psi_n^B=(-\zeta^*)^{-1}\psi_{n-1}^B\ n\ge0$. Since $\psi_{-1}^B=0$, we have $ \psi_n^B=0$ for all $n\ge0$.
 If $\kpar=\pi$ then we have from \eqref{tbevp-a} that $\psi_{n-1}^B=0$ for all $n\ge0$.
 
  Now suppose $\kpar\in[0,2\pi]\setminus(2\pi/3,4\pi/3)$. Then, the above discussion also implies that  if $\psi\in l^2(\N_0;\C^2)$ solves the eigenvalue equation with $z=0$,  then $\psi\equiv0$. 
  
 We conclude:
 {\it  $E=0$ is a point eigenvalue of  $\HTBs(\kpar)$ acting in $l^2(\N_0;\C^2)$ if and only if  $\kpar\in(2\pi/3,4\pi/3)$. 
 For $\kpar\in (2\pi/3,4\pi/3)\setminus\{\pi\}$, the $l^2(\N_0;\C^2)$- normalized eigenstate is given by:
  \begin{align}
\psi^{\rm TB, bd}_n(\kpar)\ &=\ \sqrt{1-|\zeta(\kpar)|^2}\ \begin{pmatrix} \left(-\zeta(\kpar)\right)^n \\ 0\end{pmatrix},\quad n\ge0
\label{zz-estate}\\
\zeta(\kpar)\ & \equiv\ 1+e^{i\kpar}\ .
\label{zeta-def}
\end{align}
 For  $\kpar=\pi$ ($\zeta(\kpar)=0$), the eigenstate is given by the expression:
\begin{equation}
\psi^{\rm TB, bd}_{_{0}}(\pi)\ =\  \begin{pmatrix} 1 \\ 0\end{pmatrix},\quad \psi^{\rm TB, bd}_{_{n}}(\pi)\ =\  \begin{pmatrix} 0 \\ 0\end{pmatrix},\qquad n\ge1\ ,\label{pi-state}\end{equation}
and is supported strictly at the edge.
}

We now assume that $z$ is complex and $z\ne0$, 
 and explore the invertibility of $\HTBs(\kpar)-z\ I$ on $l^2(\N_0;\C^2)$.
 Written out componentwise, the system $(\HTBs(\kpar)-z\ I)\psi=f$, where $f\in l^2(\N_0;\C^2)$ is:
\begin{align}
&\psi_{n-1}^B+\zeta^*\psi_n^B\ =\ z\psi_n^A\ +\ f_n^A,\ n\ge0 \label{tb-inhom-a}\\
&\psi_{n+1}^A+\zeta\psi_n^A\ =\ z\psi_n^B\ +\ f_n^B,\ n\ge0 ,\label{tb-inhom-b}\\
&\psi_n= \begin{pmatrix}\psi_n^A\\ \psi_n^B\end{pmatrix}
=\begin{pmatrix}0\\ 0\end{pmatrix},\quad f_n= \begin{pmatrix}f_n^A\\ f_n^B\end{pmatrix}=\begin{pmatrix}0\\ 0\end{pmatrix},\quad
 \textrm{for all $n\le-1$}\label{psi-f-bc}\\
 &\textrm{and}\ |\psi_n|\to0 \quad \textrm{as}\ n\to\infty.
\label{n2inf} \end{align}
We focus on the case $\kpar\in[0,2\pi]\setminus\{\pi\}$, so that $\zeta(\kpar)=1+e^{i\kpar}\ne0$.
\begin{remark}\label{pi-rmk}
For $\kpar=\pi$, the  system \eqref{psi-f-bc} is of the form $(\HTBs(\pi)-z)\psi=f$,
where  $\psi=(\psi_0 ^A,\psi_0^B,\psi_1 ^A,\psi_1^B,\dots)^\top$,
 $f=(f_0 ^A,f_0^B,f_1 ^A,f_1^B,\dots)^\top$ and $\HTBs(\pi)$ is a block-diagonal matrix 
 consisting of a $1\times1$ block, $0$ in the $(1,1)$ entry, followed by an infinite sequence  of identical $2\times2$ blocks, each  equal to $\sigma_1=\begin{pmatrix} 0&1\\ 1&0\end{pmatrix}$, filling out the diagonal. The statements in Theorem \ref{zz-spec} 
on the spectrum of $\HTBs(\pi)$ and the mapping $z\mapsto (\HTBs(\pi)-z)^{-1}$ 
 are easily verified.
\end{remark}

For $\kpar\ne\pi$, we next rewrite \eqref{tb-inhom-a}-\eqref{tb-inhom-b} as a first order recursion.
 Consider \eqref{tb-inhom-a} with $n$ replaced by $n+1$:
\begin{equation} \psi_{n}^B+\zeta^*\psi_{n+1}^B\ =\ z\psi_{n+1}^A\ +\ f_{n+1}^A,\ n\ge-1\ .
\label{B-shift}\end{equation}
For $n=-1$,  equation \eqref{B-shift} implies the boundary condition at site $n=0$:
\begin{equation} \zeta^*\psi_0^B\ -\ z\psi_0^A\ =\ f_0^A .\label{BC}\end{equation}
For $n\ge0$, we use $\zeta\ne0$ and \eqref{tb-inhom-b} in \eqref{B-shift} and obtain:
\begin{align}
\psi_{n+1}^B\ &=\ \left(-\frac{\zeta}{\zeta^*}\right)\ z\ \psi_n^A\ +\ \frac{z^2-1}{\zeta^*}\ \psi_n^B
  +\ \frac{z}{\zeta^*} f_n^B\ +\ \frac{1}{\zeta^*} f_{n+1}^A,\ \ n\ge0 \label{psiBn+1}
  \end{align}

  Summarizing, we have that the system: \eqref{tb-inhom-a}, \eqref{tb-inhom-b} and \eqref{psi-f-bc} is equivalent to the first order system \eqref{tb-inhom-b}, \eqref{psiBn+1} for $\psi_n=\begin{pmatrix} \psi_n^A\\ \psi_n^B\end{pmatrix}$, $n\ge0$,  with the boundary condition \eqref{BC} at $n=0$. We write this more compactly as:
  \begin{align}
 \psi_{n+1}\ &=\ M(z,\zeta)\ \psi_n\ +\ F_n(z,\zeta),\ n\ge0 ,\label{diff-eqn}\\
 &\begin{pmatrix} -z \\ \ \ \zeta^*\ \end{pmatrix}^\top\ \psi_0\ \equiv\ \begin{pmatrix} -z \\ \ \ \zeta^*\ \end{pmatrix}^\top\ \begin{pmatrix} \psi_0^A\\ \psi_0^B\end{pmatrix}\ =\ f_0^A \ ,
 \label{BC0}\\
 & |\psi_m|\ \to\ 0,\ \ m\to\infty. \label{psi-decay}\end{align}
where
 \begin{align}
 M(z,\zeta)\ &=\ \begin{pmatrix} -\zeta & z\\ -\frac{\zeta}{\zeta^*}z & \frac{z^2-1}{\zeta^*}\end{pmatrix}\ ,
 \label{MEz}\\
 F_n(z,\zeta;f)\ &=\ \begin{pmatrix} f_n^B \\ \frac{z}{\zeta^*}f_n^B\ +\ \frac{1}{\zeta^*}f_{n+1}^A\end{pmatrix}\ ,\ n\ge0.
\label{FEz} \end{align}
 We next solve \eqref{diff-eqn}-\eqref{BC0} by diagonalizing the matrix $M(z,\zeta)$.
 
The eigenvalues $\lambda$  of $M(z,\zeta)$ are solutions of the quadratic equation
\begin{equation}
\zeta^* \lambda^2\ +\ \left(1+|\zeta|^2-z^2\right) \lambda\ +\ \zeta\ =\ 0,
\label{lam-eqn}\end{equation}
whose solutions are:
\begin{align}
\lambda_1(z,\zeta)\ &=\ \frac{-\left(1+|\zeta|^2-z^2\right)+\sqrt{\left(1+|\zeta|^2-z^2\right)^2-4|\zeta|^2}}{2\zeta^*}\ \label{lam1}\\
\lambda_2(z,\zeta)\ &=\ \frac{-\left(1+|\zeta|^2-z^2\right)-\sqrt{\left(1+|\zeta|^2-z^2\right)^2-4|\zeta|^2}}{2\zeta^*}\ .\label{lam2}
\end{align}
When convenient, we suppress the dependence of $\lambda_1$ and $\lambda_2$ on $\zeta$ and $E$
and occasionally write $\lambda_j$ or $\lambda_j(z)$.
These expressions depend on $\kpar$ through $\zeta(\kpar)=1+e^{i\kpar}$.

\nit Note that  $|\lambda_1\ \lambda_2|\ =\ |\det M(z,\zeta)|\ =\ |\zeta/\zeta^*|=1$ and hence 
$M(z,\zeta)$ may have at most one eigenvalue strictly inside the unit circle in $\C$. 

Recall the definitions: $\delta_{\rm gap}(\kpar)\equiv \Big|1-|\zeta(\kpar)|\Big|$ and $\delta_{\rm max}(\kpar)\equiv1+|\zeta(\kpar)|$. 
\begin{remark}\label{prelude-roots}
We shall see just below that for fixed $\kpar\ne 2\pi/3,\pi$ or $4\pi/3$: if  (a) $|z|<\delta_{\rm gap}(\kpar)$ or (b) $|z|>\delta_{\rm max}(\kpar)$
then the discriminant in \eqref{lam1}-\eqref{lam2},  $(1+|\zeta(\kpar)|^2-z^2)^2-4|\zeta(\kpar)|^2$, is strictly positive and uniformly 
bounded away from zero. Therefore, in each of these cases  the expressions in \eqref{lam1}-\eqref{lam2} define single-valued functions $\lambda_1(z,\zeta)$ and $\lambda_2(z,\zeta)$. This property continues to hold for $\kpar\in\mathscr{I}_1\subset\subset[0,2\pi]\setminus\{2\pi/3,\pi,4\pi/3\}$ and either
(a$^\prime$)\ $ |\Re z|<\delta_{\rm gap}(\kpar)\ \ {\rm and}\ \ |\Im z|<\eta(\mathscr{I}_1)$ or 
(b$^\prime$)\ $|\Re z|>\delta_{\rm max}(\kpar)\ \ {\rm and}\ \ |\Im z|<\eta(\mathscr{I}_1)$, for some $\eta(\mathscr{I}_1)>0$ chosen sufficiently small.
 In the case where $z$ is real and $\delta_{\rm gap}(\kpar)\le |z|\le \delta_{\rm max}(\kpar)$ the discriminant is nonpositive and we do not distinguish between the roots of \eqref{lam-eqn}; they comprise a two element set on the unit circle in $\C$. 
\end{remark}
\begin{lemma}\label{whereroots}\ 
Assume $0<|\zeta(\kpar)|\ne1$, {\it i.e.} $\kpar\ne 2\pi/3,\pi$ or $4\pi/3$.    Then, the following hold:
\begin{enumerate}
\item[(1)]    Let  $z\in\R$ and assume that either 
\begin{equation}
\textrm{$|z|<\delta_{\rm gap}(\kpar)$ or $|z|>\delta_{\rm max}(\kpar)$.}
\label{en-range}\end{equation}
 Then, $M\left(z,\zeta(\kpar)\right)$ 
 has one eigenvalue inside the unit circle and one eigenvalue outside the unit circle. 
 
\item[(2)]  Let $\lambda_1(z)$ and $\lambda_2(z)$ denote be the expressions for the eigenvalues of $M\left(z,\zeta(\kpar)\right)$ displayed in \eqref{lam1}-\eqref{lam2}. 
\subitem(i) If $z\in\R$ and $|z|<\delta_{\rm gap}(\kpar)$,  then $|\lambda_1(z;\kpar)|<1<|\lambda_2(z;\kpar)|$.
 
\subitem(ii) If $z\in\R$ and $|z|>\delta_{\rm max}(\kpar)$, then $|\lambda_2(z;\kpar)|<1<|\lambda_1(z;\kpar)|$.
\subitem(iii) If $z\in\R$ and $\delta_{\rm gap}(\kpar)\le |z|\le \delta_{\rm max}(\kpar)$, then equation \eqref{lam-eqn} has two roots, $\lambda$, satisfying $|\lambda|=1$.

\item[(3)]  
Let $\mathscr{I}_1$ denote a compact subset of $[0,2\pi]\setminus\{2\pi/3,\pi,4\pi/3\}$. There exists a constant  $\eta>0$, which depends
on $\mathscr{I}_1$, such that for all $\kpar\in\mathscr{I}_1$ the following hold:\\
 (a) If $z$ is in the complex open neighborhood 
 \begin{equation}
  \mathscr{O}_0(\kpar):\quad |\Re z|<\delta_{\rm gap}(\kpar)\ \ {\rm and}\ \ |\Im z|<\eta(\mathscr{I}_1),
  \label{scrO_0}
  \end{equation}
  then \eqref{en-range} holds. Moreover,  $\lambda_1(z,\zeta)$ and $\lambda_2(z,\zeta)$ satisfy the strict inequalities of (2.i),
  and their magnitudes are uniformly bounded away from $1$, provided $z$ remains in a compact subset of  $\mathscr{O}_0(\kpar)$. \\
 (b) If $z$ is in the complex open neighborhood 
  \begin{equation}
   \mathscr{O}_+(\kpar):\quad |\Re z|>\delta_{\rm max}(\kpar)\ \ {\rm and}\ \ |\Im z|<\eta(\mathscr{I}_1),
   \label{scrO_+}
   \end{equation}
then \eqref{en-range} holds and moreover $\lambda_1(z,\zeta)$ and $\lambda_2(z,\zeta)$ 
 satisfy the inequalities of  (2.ii)  and their magnitudes are uniformly bounded away from $1$, provided $z$ remains in a compact subset of  $\mathscr{O}_+(\kpar)$.
 \end{enumerate}
 \end{lemma}
 
\nit{\it Proof of Lemma \ref{whereroots}:} Part 3 of the Lemma follows from parts (1) and (2) 
 and the expressions \eqref{lam1}, \eqref{lam2} for $\lambda_1(z;\kpar)$, and $\lambda_2(z;\kpar)$. We now proceed with the proof
  of assertions (1) and (2), which assume $z\in\R$.
 \medskip
 
  We consider the two cases delineated by the sign of the discriminant:\medskip

\nit {\bf Case 1:} $\left(1+|\zeta|^2-z^2\right)^2-4|\zeta|^2>0$\ and\ 
{\bf Case 2:} $\left(1+|\zeta|^2-z^2\right)^2-4|\zeta|^2\le0$.\medskip

\nit{\bf Case 1:}\ In this case, $\Big| 1 + |\zeta|^2-z^2\Big|>2|\zeta|$. There are two subcases:\\
  (1a) $1+|\zeta|^2-z^2>2|\zeta|$
and  (1b) $z^2-1-|\zeta|^2>2|\zeta|$. \medskip

\nit In  subcase (1a), we have  $z^2<(1-|\zeta|)^2$ and therefore $|z|<\delta_{\rm gap}(\kpar)=|1-|\zeta||$, where $\delta_{\rm gap}(\kpar)>0$ since $\kpar\ne2\pi/3, 4\pi/3$. In this subcase we also have:
 $-(1+|\zeta|^2-z^2)<-2|\zeta|<0$. Therefore, 
\begin{align*}
0> (2\zeta^*)\lambda_1 &=-\left(1+|\zeta|^2-z^2\right)+\sqrt{\left(1+|\zeta|^2-z^2\right)^2-4|\zeta|^2}\\
&>-\left(1+|\zeta|^2-z^2\right)-\sqrt{\left(1+|\zeta|^2-z^2\right)^2-4|\zeta|^2}=(2\zeta^*)\lambda_2.
\end{align*}
Let  $\lambda_1=r_1/(2\zeta^*)$ 
and $\lambda_2=r_2/(2\zeta^*)$. 
Therefore,   $|r_1|=|(2\zeta^*)\lambda_1|<|(2\zeta^*)\lambda_2|=|r_2|$. Therefore, 
$|\lambda_1| / |\lambda_2|=|r_1|/|r_2|<1$. Since $|\lambda_1|\ |\lambda_2|=1$, 
\begin{equation}
\textrm{ in subcase (1a), we have $|z|<\delta_{\rm gap}(\kpar)$, and $|\lambda_1(z)|<1<|\lambda_2(z)|$.}
\label{case1a}\end{equation}

\nit In subcase (1b) we have $|z|>1+|\zeta(\kpar)|=\delta_{\rm max}(\kpar)$. Hence, $1+|\zeta|^2-z^2<1+|\zeta|^2-(1+|\zeta|)^2=-2|\zeta|<0$ since $\kpar\ne\pi$. Therefore,
 \begin{equation}
 \textrm{ in subcase (1b), we have  $|z|>\delta_{\rm max}(\kpar)$ and  $|\lambda_2(z)|<1<|\lambda_1(z)|$.}
 \label{case1b}
 \end{equation} 

 \nit {\bf Case 2:} Here we have $\delta_{\rm gap}(\kpar)\le |z| \le \delta_{\rm max}(\kpar)$. In this case, $\lambda_1= (a+ib)/(2\zeta^*)$ and $\lambda_2= (a-ib)/(2\zeta^*)$ , where $a$ and $b$ are real. Therefore, $|\lambda_1| / |\lambda_2|= 1$ and hence $|\lambda_1|=|\lambda_2|$ implying that 
\begin{equation}
 \textrm{ in case (2), we have $\delta_{\rm gap}(\kpar)\le |z| \le \delta_{\rm max}(\kpar)$ and  $|\lambda_1(z)|=|\lambda_2(z)|=1$.}
\label{case2}
\end{equation}
  
\nit We note the assertions \eqref{case1a}, \eqref{case1b} and \eqref{case2},  hold for any $\kpar\notin\{ 2\pi/3,\pi,4\pi/3\}$. \medskip

\nit The proof of Lemma \ref{whereroots} is now complete.
 \medskip

We continue now with the proof of Theorem \ref{zz-spec}. Assume that $\kpar\in[0,2\pi]\setminus\{ 2\pi/3,\pi,4\pi/3\}$, and hence 
 $0< |\zeta(\kpar)|\ \ne\ 1$,  so that 
Lemma \ref{whereroots} applies. 
 Corresponding to the eigenvalues, $\lambda_1(z)$ and $\lambda_2(z)$ of  $M(z,\zeta)$ we  can take the corresponding eigenvectors to be of the form:
\begin{align}
\xi_1(z)\ =\ \begin{pmatrix} z\\ \zeta+\lambda_1\end{pmatrix},\qquad 
\xi_2(z)\ =\ \begin{pmatrix} z\\ \zeta+\lambda_2\end{pmatrix}\ .
\label{e-vecs}\end{align}

Due to the hypothesized constraints on $\kpar$, in particular that $\kpar\ne2\pi/3, 4\pi/3$, we have $\zeta\ne0$. For small $z$ we find
the following asymptotic expansions for $\lambda_j(z,\zeta)$, which are valid uniformly in $\kpar$ varying over any prescribed compact subset, $\mathscr{I}_1$, of $[0,2\pi]\setminus\{2\pi/3,\pi,4\pi/3\}$:
\medskip

\begin{align}
&\kpar\in \mathscr{I}_1\subset\subset (2\pi/3,4\pi/3)\setminus\{\pi\}\qquad \textrm{(hence, $0<|\zeta(\kpar)|<1$)}\nn\\
& \implies
\begin{cases}
\lambda_1=\lambda_1(z,\zeta)\ =\ -\zeta + \mathcal{O}(|z|^2) \\
\lambda_2=\lambda_2(z,\zeta)\ =\ -(\zeta^*)^{-1}\ +\ \mathcal{O}(|z|^2) .
\end{cases}
\label{lam2-est}\end{align}
and  \\
\begin{align}
&\kpar\in\mathscr{I}_1\subset\subset [0,2\pi]\setminus[2\pi/3,4\pi/3]\qquad \textrm{(equivalently, $|\zeta(\kpar)|>1$)}\nn\\
&\implies
\begin{cases}
\lambda_1=\lambda_1(z,\zeta)\ =\ -(\zeta^*)^{-1} + \mathcal{O}(|z|^2) \\
\lambda_2=\lambda_2(z,\zeta)\ =\ -\zeta\ +\ \mathcal{O}(|z|^2) .
\end{cases}
\label{lam1-estB}
\end{align}


%
%

\nit{\bf The resolvent $(\HTBs(\kpar)-z\ I)^{-1}$ on $l^2(\N_0;\C^2)$:}\ 
Let us now restrict $\kpar$ to vary over the set $(2\pi/3,4\pi/3)\setminus\{\pi\}$, and assume $0<|z|<\delta_{\rm gap}(\kpar)$; 
and construct the resolvent of $\HTBs(\kpar)$ by solving \eqref{diff-eqn}, \eqref{BC0}.  The construction of the 
resolvent for $|z|>\delta_{\rm max}(\kpar)$ for all $\kpar\in[0,2\pi]$ and all $z$ such that $|z|<\delta_{\rm gap}(\kpar)$, where 
 $\kpar\in[0,2\pi]\setminus(2\pi/3,4\pi/3)$ can be carried out similarly (see remarks below).

For $\kpar\in(2\pi/3,4\pi/3)\setminus\{\pi\}$,   the expansions \eqref{lam2-est} are valid and we have
\begin{align}
\zeta+\lambda_1\ &=\ \mathcal{O}(|z|^2),\qquad \zeta+\lambda_2\ =\ \zeta-\frac{1}{\zeta^*}\ +\ \mathcal{O}(|z|^2) ,
\nn\end{align}
and we have by \eqref{e-vecs} that the  eigenvectors satisfy
\begin{align}
\frac{1}{z}\xi_1(z)\ &=\ \begin{pmatrix}1 \\ 0\end{pmatrix}\ +\ \mathcal{O}_{_{\C^2}}(|z|),\qquad  
\xi_2(z)\ =\  \left(\zeta-\frac{1}{\zeta^*}\right)\begin{pmatrix}0 \\ 1\end{pmatrix}\ +\ \mathcal{O}_{_{\C^2}}(|z|)
\label{xi2}\end{align}
for all $z$ small. Hence,  
\[
\textrm{  $\Big\{\frac{1}{z}\xi_1(z),\xi_2(z)\Big\}$ is a basis of $\C^2$ for $0<|z|<\delta_{\rm gap}(\kpar)$ and $\kpar\in(2\pi/3,4\pi/3)\setminus\{\pi\}$}
  \]
  which does not degenerate in the limit $z\to0$. 
   Indeed, by \eqref{e-vecs} for $z\ne0$  this set is linearly independent if and only if  $\lambda_1\ne\lambda_2$. 
   However, for $0<|z|<\delta_{\rm gap}(\kpar)$ we have $|\lambda_1|<1<|\lambda_2|$. 

To solve \eqref{diff-eqn}, \eqref{BC0}  we next express $F_n=F_n(z,\zeta;f)$ in the non-degenerate basis \eqref{xi2}.
 We shall, when convenient, suppress the dependence of $F_n$ on $\zeta$ and $f$:
\begin{align}
F_n(f;z,\zeta)\ &=\ \begin{pmatrix} f_n^B \\ \frac{z}{\zeta^*}f_n^B\ +\ \frac{1}{\zeta^*}f_{n+1}^A\end{pmatrix}\nn\\
\quad &=F^{(1)}_n(f;z,\zeta)\ \frac{1}{z}\ \xi_1(z)\ +\ F^{(2)}_n(f;E,\zeta)\ \xi_2(z)\ .
\label{F-expand}\end{align}
We also seek a solution as an expansion in the basis  \eqref{xi2}: 
\begin{equation}
\psi_n\ =\ \psi^{(1)}_n\ \frac{1}{z}\ \xi_1(z)\ +\ \psi^{(2)}_n\ \xi_2(z),
\label{psi-expand}
\end{equation}
where $\psi^{(1)}_n=\psi^{(1)}_n(z)$ and $\psi^{(2)}_n=\psi^{(2)}_n(z)$ are to be determined.
Then, we obtain the two decoupled first order difference equations:
\begin{align}
\psi^{(1)}_{n+1}\ &=\ \lambda_1(z)\psi^{(1)}_n+F^{(1)}_n(z),\quad\ n\ge0, \label{diag1}\\
\psi^{(2)}_{n+1}\ &=\ \lambda_2(z)\psi^{(2)}_n+F^{(2)}_n(z),\quad\ n\ge0,  \label{diag2}
\end{align}
with boundary condition \eqref{BC0} to be expressed in terms of
$\psi^{(j)}_0$, and $F^{(j)}_0$, $j=1,2$:
\begin{equation}
\frac{1}{z}\ \begin{pmatrix} -z \\ \ \ \zeta^*\ \end{pmatrix}^\top\ \xi_1(z)\ \psi_0^{(1)}\ +\ 
 \begin{pmatrix} -z \\ \ \ \zeta^*\ \end{pmatrix}^\top\ \xi_2(z)\ \psi_0^{(2)}\ =\ f_0^A
\label{BC0-a}
\end{equation}

We now proceed to solve the decoupled system \eqref{diag1}-\eqref{diag2} and then impose the boundary condition \eqref{BC0-a}.
Recall our assumption that $ 0\ <\ |\zeta|\ <\ 1$, {\it i.e.} $\kpar\in(2\pi/3,4\pi/3)\setminus\{\pi\}$ and therefore for $z$ real and $|z|<\delta_{\rm gap}(\kpar)$, we have 
that $|\lambda_1(z)|<1<|\lambda_2(z)|$. In this case, 
 the most general solution of \eqref{diag1}, which decays as $n\to+\infty$ is:
 \begin{equation}
 \psi_n^{(1)}(z)\ =\ \sum_{j=0}^{n-1}\ \left(\lambda_1(z)\right)^{n-1-j}\ F_j^{(1)}(z)\ +\ \mu\ \left(\lambda_1(z)\right)^n \ .
 \label{psi1-gen}
 \end{equation}
 where $\mu$ is an arbitrary constant to be determined and $F_j^{(1)}(f;z,\zeta),\ F_j^{(2)}(f;z,\zeta)$ are defined by \eqref{F-expand}.
 
 Furthermore, the most general solution of \eqref{diag2} which decays as $n\to+\infty$ is:
 \begin{equation}
 \psi_n^{(2)}(z)\ =\ -\sum_{j=n}^\infty\ \left(\lambda_2(z)\right)^{n-j-1}\ F_j^{(2)}(z)\ .
 \label{psi2-gen}
 \end{equation}
 
 Finally, we now turn to the boundary condition \eqref{BC0-a}. Using \eqref{psi1-gen} and \eqref{psi2-gen}
 for $n=0$ in \eqref{BC0-a} we find:
\begin{equation} 
\mu\frac{1}{z}  \begin{pmatrix} -z \\ \ \ \zeta^*\ \end{pmatrix}^\top\ \xi_1(z)\ -\ 
 \begin{pmatrix} -z \\ \ \ \zeta^*\ \end{pmatrix}^\top\ \xi_2(z)\ \sum_{j=0}^\infty\ \left(\lambda_2(z)\right)^{-j-1}\ F_j^{(2)}(z,\zeta;f)\\ =\ f_0^A .
\label{BC0-b}
\end{equation}
By \eqref{lam-eqn}, the quadratic equation for the roots $\lambda_j$, we find:
\begin{equation}
\begin{pmatrix} -z \\ \ \ \zeta^*\ \end{pmatrix}^\top\ \xi_j(z)\ =\ 
\begin{pmatrix} -z \\ \ \ \zeta^*\ \end{pmatrix}^\top\ \begin{pmatrix} z\\ \zeta+\lambda_j(z)\end{pmatrix}
\ =\ -\frac{\zeta+\lambda_j(z)}{\lambda_j(z)},\ \ j=1,2.
\label{coeff-j}
\end{equation}

\nit {\it Claim:}\ Assume $z\ne0$ and $z\in\R$. If  $\lambda(z)$ is any root of  \eqref{lam-eqn}, then
$\frac{\zeta+\lambda(z)}{\lambda(z)}\ne0$.\medskip

\nit It follows from this claim and \eqref{coeff-j} that the coefficient of $\mu$ in   \eqref{BC0-b} is non-zero and hence
\[\textrm{if $z\ne0$ we can solve \eqref{BC0-a} for $\mu=\mu(z,\zeta;f)$\ .}
\]
  To prove the above Claim we first note that $\lambda\ne0$.
Indeed, if $\lambda=0$ then \eqref{lam-eqn} would then imply $\zeta=1+e^{i\kpar}=0$; this contradicts our assumption that $\kpar\ne\pi$. Thus, $\lambda(z)\ne0$. Furthermore, we claim that $\zeta+\lambda(z)\ne0$. Again, using \eqref{lam-eqn} we have that if $\zeta+\lambda=0$ then $\zeta\ z^2=0$. This contradicts the assumptions that $z\ne0$ and $\zeta\ne0$.
 
 It follows from this discussion that for $z\ne0$ and $\kpar\ne\pi$:
 \begin{align}
 \mu(f;z,\zeta)\ = -\frac{z\ \lambda_1(z)}{\zeta+\lambda_1(z)}\ \Big[\ f^A_0\ -\ 
  \frac{\zeta+\lambda_2(z)}{\lambda_2(z)}\ \sum_{j=0}^\infty\ \left(\lambda_2(z)\right)^{-j-1}\ F_j^{(2)}(f;z,\zeta)\ \Big]
\label{mufEz} \end{align}
Therefore  if $0<|z|<\delta_{\rm gap}(\kpar)$ and  $\kpar\in(2\pi/3,4\pi/3)\setminus\{\pi\}$, we can solve for $\mu=\mu(z,\zeta;f)$.  
We obtain for any $f\in l^2(\N_0;\C^2)$,  the unique solution of \eqref{diff-eqn}, \eqref{BC0} and \eqref{psi-decay}
\[\psi=\{\psi_n\}_{n\ge0}, \quad \textrm{ with $\psi_n$ tending to zero as $n\to\infty$, is given by}\]
\begin{align}
\psi_n\ &=\ \left[\ \sum_{j=0}^{n-1}\ \left(\lambda_1(z,\zeta)\right)^{n-1-j}\ F_j^{(1)}(f;z,\zeta)\ +\ \mu(z,\zeta;f)\ \left(\lambda_1(z,\zeta)\right)^n  \right]\ \frac{1}{z}\xi_1(z,\zeta)\nn\\
&\quad  -\ \left[\sum_{j=n}^\infty\ \left(\lambda_2(z,\zeta)\right)^{n-j-1}\ F_j^{(2)}(f;z,\zeta)\right]\ \xi_2(z,\zeta),\qquad n\ge0\ ,
\label{psi-resf}\end{align} 
where $\mu=\mu(z,\zeta;f)$ is obtained from \eqref{BC0-b}.  By \eqref{F-expand}, we may express $F_j^{(1)}$ and $F_j^{(2)}$
as 
\begin{equation}
F_j^{(1)}\ =\ \alpha_1(z,\zeta)\ f_j^B\ +\ \alpha_2(z,\zeta)\ f_{j+1}^A,\quad 
F_j^{(2)}\ =\ \beta_1(z,\zeta)\ f_j^B\ +\ \beta_2(z,\zeta)\ f_{j+1}^A,
\label{albe-def}\end{equation}
 where the coefficients are bounded and smooth over the ranges of $z$ and $\kpar$ under consideration. 

Next, introduce the discrete vector-valued kernel, depending on parameters $\alpha$ and $\beta$:
\begin{equation}
\mathscr{K}(n,j;\alpha,\beta) =\ 
\begin{cases}
 \alpha\ \lambda_1(z,\zeta)^{n-1-j}\ \frac{1}{z}\ \xi_1(z,\zeta)\ ,\ & 0\le j\le n-1\\ 
 &\\ 
-\beta\ \lambda_2(z,\zeta)^{n-1-j}\ \xi_2(z,\zeta)\ ,\ & n\le j<\infty\ .
\end{cases}
\label{k-kernel}\end{equation}
Then, we have
\begin{align}
\psi_n\ &=\ \sum_{j=0}^\infty\ \mathscr{K}(n,j;\alpha_1,\beta_1) f_j^B\ +\ \sum_{j=0}^\infty\ \mathscr{K}(n,j;\alpha_2,\beta_2) f_{j+1}^A\nn\\
&\qquad +\ \mu(f;z,\zeta)\ (\lambda_1(z,\zeta))^n\ \frac{1}{E}\xi_1(z,\zeta),
\label{psi_n}\end{align}
where $\mu(f;z,\zeta)$ is given by the linear functional of $f$, displayed in \eqref{mufEz}.

  \begin{proposition}\label{TB-res}
Let $\mathscr{I}_1$ denote a compact subset of $(2\pi/3,4\pi/3)\setminus\{\pi\}$ and let  $\eta(\mathscr{I}_1)>0$,
denote the constant appearing in part (3) of Lemma \ref{whereroots}. 
\begin{enumerate}
\item There is a constant, $C$, depending on
$\mathscr{I}_1$ such that for all complex energies, $z\in\mathscr{O}_0(\kpar)\setminus\{0\}$ (\ see \eqref{scrO_0}\ ),
  the resolvent operator:
 \begin{equation}
 f\in l^2(\N_0;\C^2)\mapsto \psi=\{\psi_n\}_{n\ge0}\ \equiv\ \left(\HTBs(\kpar)-z\right)^{-1}f,\label{tb-res-def}
 \end{equation}
given by the expression in \eqref{psi_n}, 
  defines a bounded linear operator on $l^2(\N_0;\C^2)$ with 
\begin{equation}
 \Big\|\left(\HTBs(\kpar)-z\right)^{-1}f\Big\|_{l^2(\N_0;\C^2)}\le C \frac{1}{|z|}\ \|f\|_{l^2(\N_0;\C^2)}\ ,
 \label{TB-res-bd}
 \end{equation}
where the constant, $C$,  is independent of depends on the compact set $\mathscr{I}_1$.
\item The mapping $z\mapsto\left(\HTBs(\kpar)-z\right)^{-1}$ is meromorphic for $z$ varying in the open set $\mathscr{O}_0(\kpar)$
 into $\mathcal{B}( l^2(\N_0;\C^2))$, the space of bounded linear operators on  $l^2(\N_0;\C^2)$,  with only pole at $z=0$. For $z\in \mathscr{O}_0(\kpar)\setminus\{0\}$ we have
\begin{equation}
\left(\HTBs(\kpar)-z I\right)^{-1}f\ =\ \frac{1}{z} \left\langle \psi^{\rm TB, bd}(\kpar),f\right\rangle_{_{l^2(\N_0;\C^2)}}\  \psi^{\rm TB, bd}(\kpar)\ +\ \mathscr{G}_{\rm reg}(z;\kpar)f,
\label{laurent}\end{equation}
where $z\mapsto \mathscr{G}_{\rm reg}(z;\kpar)$ is an analytic map from $\mathscr{O}_0(\kpar)$
to $\mathcal{B}( l^2(\N_0;\C^2))$.
\item $\HTBs(\kpar)\psi = f\in l^2(\N_0;\C^2)$  has a solution in the space  $  l^2(\N_0;\C^2)$  if and only if $ \left\langle \psi^{\rm TB, bd}(\kpar),f\right\rangle_{_{l^2(\N_0;\C^2)}}=0$.
\end{enumerate}
\end{proposition}

\nit {\it Proof of Proposition \ref{TB-res}:}  We fix $\mathscr{I}_1\subset\subset(2\pi/3,4\pi/3)\setminus\{\pi\}$
 and take $E\in\mathscr{O}_0(\kpar)\setminus\{0\}$. To bound the resolvent we estimate the expression in $\{\psi_n\}_{n\ge0}$ displayed in \eqref{psi_n} in $l^2(\N_0;\C^2)$. 
We begin with an estimate of the latter term in \eqref{psi_n}: $\mu(f;z,\zeta)\ (\lambda_1(z,\zeta))^n\ \frac{1}{z}\xi_1(z,\zeta)$. From the expression for $\mu$ in \eqref{mufEz} and 
 the definition of  $F_j^{(2)}$ in \eqref{F-expand} (recall $F_j^{(1)}$ and $F_j^{(2)}$
  are coordinates of $F_j\in\C^2$, also given in \eqref{F-expand}) with respect to the basis $\{\frac{1}{z}\xi_1(z),\xi_2(z)\}$),
   we have that $|\mu(f;z,\zeta)|\lesssim |f_0^A|+\sum_{j=0}^\infty|\lambda_2|^{-j-1}\left(\ |f_j^B|\ +\ |f_{j+1}^A|\ \right)\le C_1(z,\zeta)\ \|f\|_{l^2(\N_0;\C^2)}$, where $C_1(z,\zeta)$ is a finite constant which depends on $z$ and $\zeta$ in the ranges specified above.  The constant $C_1(z,\zeta)$ is bounded for $z$  bounded away from $z=0$ and $\kpar\in\mathscr{J}_1$. As we shall see below, for $\kpar\in\mathscr{J}_1$, there is pole of order one as $E\to0$.

Therefore, applying Young's inequality to the first two terms in \eqref{psi_n} we obtain:
\[\|\psi\|_{_{l^2(\N_0;\C^2)}}\ \le\ \Big(\ C(\mathscr{K},z,\zeta)\  +\ C_1(z,\zeta)\ \Big)\ \|f\|_{_{l^2(\N_0;\C^2)}},\]
 where
\begin{equation}
C(\mathscr{K},z,\zeta)= \max_{r=1,2}\Big(\ \sup_{n\ge0}\sum_{j=0}^\infty\ |\mathscr{K}(n,j,\alpha_r,\beta_r)|\ +\ \sup_{j\ge0}\sum_{n=0}^\infty\ |\mathscr{K}(n,j,\alpha_r,\beta_r)| \Big)\ ,
\label{y-norm}
\end{equation}
and  we recall from \eqref{albe-def} that $\alpha_r$ and $\beta_r$ are smooth and bounded functions of $z$ and $\zeta$.
Estimating the first sum in \eqref{y-norm}, we have for $r=1, 2$: 
\begin{align} \sum_{j=0}^\infty |\mathscr{K}(n,j,\alpha_r,\beta_r)|\ &\lesssim\ |\alpha_r(z,\zeta)| \sum_{j=0}^{n-1}|\lambda_1(z,\zeta)|^{n-1-j}
\ +\ |\beta_r(z,\zeta)|\ \sum_{j=n}^{\infty}|\lambda_2(z,\zeta)|^{n-1-j}\nn\\
& \lesssim\ |\alpha_r(z,\zeta)|\ (1-|\lambda_1(z,\zeta)|)^{-1}+|\beta_r(z,\zeta)|\ (|\lambda_2(z,\zeta)|-1)^{-1}
\label{Kbd}\end{align}
The bound \eqref{Kbd} holds, for $r=1, 2$ and any fixed $z\in\mathscr{O}_0(\kpar)\setminus\{0\}$, uniform in $\kpar\in\mathscr{I}_1$. 
  The second sum in \eqref{y-norm} is bounded similarly. Therefore,  we have
  for all $\kpar\in\mathscr{I}_1$ and any $z\in\mathscr{O}_0(\kpar)$, the resolvent operator:  $f\mapsto \left(\HTBs(\kpar)-z\right)^{-1}f$ (see \eqref{tb-res-def}) is a bounded linear operator on   $l^2(\N_0;\C^2)$.
The next step in the proof of Proposition \ref{TB-res} requires us to consider the resolvent for
  small complex $z$ in  $\mathscr{O}_0(\kpar)\setminus\{0\}$.

\subsection{The resolvent $\left(\HTBs(\kpar)-z\ I\right)^{-1}$ for  $z$ near zero energy}

Since there is a simple zero energy eigenstate for each $\kpar\in (2\pi/3,4\pi/3)$, we expect a simple pole of the resolvent at $z=0$. We now make this explicit by expanding the resolvent in a neighborhood of $z=0$ for $\kpar\in (2\pi/3,4\pi/3)$.
 In order to work with the above detailed calculations, we restrict our discussion to the case where $\kpar\ne\pi$ ($\zeta\ne0$).
Consider first the relation \eqref{BC0-b}, which determined the free parameter $\mu=\mu(f;z,\zeta)$. We shall simplify  \eqref{BC0-b} using the following 
expansions which hold for $|z|$ small: 
\begin{align}
\begin{pmatrix} -z \\ \ \ \zeta^*\ \end{pmatrix}^\top\ \frac{1}{z}\ \xi_1(z)\ =\ 
\begin{pmatrix} -z \\ \ \ \zeta^*\ \end{pmatrix}^\top\ \frac{1}{z}\ \begin{pmatrix} z\\ \zeta+\lambda_1(z)\end{pmatrix}
\ =\ -\frac{1}{z}\ \frac{\zeta+\lambda_1(z)}{\lambda_1(z)}\ =\   \frac{z}{|\zeta|^2-1}\ +\mathcal{O}(|z|^3)
\label{coeff-1exp}\\
\begin{pmatrix} -z \\ \ \ \zeta^*\ \end{pmatrix}^\top\ \xi_2(z)\ =\ 
\begin{pmatrix} -z \\ \ \ \zeta^*\ \end{pmatrix}^\top\ \begin{pmatrix} z\\ \zeta+\lambda_2(z)\end{pmatrix}
\ =\ -\frac{\zeta+\lambda_2(z)}{\lambda_2(z)}\ =\ |\zeta|^2\ -\ 1\ +\ \mathcal{O}(|z|^2)
\label{coeff-2exp}
\end{align}
We also have from \eqref{F-expand} that 
\begin{align}
F_n(f;z,\zeta)\ &=\
 \begin{pmatrix} f_n^B \\ \frac{z}{\zeta^*}f_n^B\ +\ \frac{1}{\zeta^*}f_{n+1}^A
 \end{pmatrix}\nn\\
&= f_n^B\ \frac{1}{z} \xi_1(z)\ +\  
f_{n+1}^A\ \frac{1}{\zeta^*}\cdot \left(\zeta-\frac{1}{\zeta^*}\right)^{-1}\ \xi_2(z)\ +\ \mathcal{O}\left(|z|\  [\ |f_n|+|f_{n+1}|\ ]\right) .\nn
\end{align}
Therefore, for $|z|$ small
\begin{align}
F_n^{(1)}(f;z,\zeta)\ &=\  f_n^B\ +\ \mathcal{O}\left(|z|\ [\ |f_n|+|f_{n+1}|\ ]\right),\nn\\
 F_n^{(2)}(f;z,\zeta)\ &=\  \frac{1}{|\zeta|^2-1}\ f_{n+1}^A\ +\
 \mathcal{O}\left(|z|\  [\ |f_n|+|f_{n+1}|\ ]\right)\ .
 \label{F-expandE}\end{align}
Substitution of the expansions \eqref{coeff-1exp}, \eqref{coeff-2exp} and \eqref{F-expandE} into \eqref{BC0-b}, 
we obtain:
\begin{align}
&\frac{z}{|\zeta|^2-1}\ \mu\ -\ (|\zeta|^2-1)\ \sum_{j=0}^\infty\ \left(-\frac{1}{\zeta^*}\right)^{-(j+1)}\ 
\frac{1}{|\zeta|^2-1}\ f_{j+1}^A\nn\\
&\qquad +\ \mathcal{O}\left(|z|\ \|f\|_{l^2(\N_0;\C^2)}\right)\ +\ \mathcal{O}\left(|z|\ |\mu|\right)\ =\ f_0^A\ .
\label{almost1}\end{align}

Hence, 

\begin{align}
\frac{z}{|\zeta|^2-1}\ \mu\ &=\  \ f_0^A\ +\ \sum_{j=0}^\infty\ \left(-\frac{1}{\zeta^*}\right)^{-(j+1)} \ f_{j+1}^A\ +\ \mathcal{O}\left(|z|\ \|f\|_{l^2(\N_0;\C^2)}\right)\nn\\
\ &=\ \sum_{j=0}^\infty\ \left(-\zeta^*\right)^{j} \ f_j^A\ +\ \mathcal{O}\left(|z|\ \|f\|_{l^2(\N_0;\C^2)}\right)\ .
\label{almost2}\end{align}

Recall that we have assumed $\kpar\in\mathscr{I}_1\subset\subset(2\pi/3,4\pi/3)\setminus\{\pi\}$ (thus $|\zeta(\kpar)|^2-1\ne0$) and $z\in\mathscr{O}_0(\kpar)\setminus\{0\}$.   Solving \eqref{almost2} for $\mu(f;z,\zeta)$ and using the expression for 
 $\{\ \psi_j^{\rm TB, bd}(\kpar)\ \}_{j\ge0}$, the zero energy eigenstate of $\HTBs$ in \eqref{zz-estate}, 
 we obtain:
 \begin{align}
\mu(z,\zeta;f)\ &=\ 
\frac{1}{z}\ \sqrt{1-|\zeta|^2}\ \sum_{j=0}^\infty\ \overline{\psi_j^{\rm TB, bd}(\kpar)} \ f_j^A\ +\ \mathcal{O}\left( \|f\|_{l^2(\N_0;\C^2)}\right)\nn\\
&\ =\ \frac{1}{z}\ \sqrt{1-|\zeta|^2}\ \left\langle \psi^{\rm TB, bd}(\kpar),f\right\rangle_{_{l^2(\N_0;\C^2)}}\ +\ \mathcal{O}\left( \|f\|_{l^2(\N_0;\C^2)}\right),
\label{mu-Esmall}\end{align}
The error bound in \eqref{mu-Esmall} is uniform in $\kpar\in \mathscr{I}_1\setminus\{\pi\}$ and bounds an expression 
which is analytic in $z\in\mathscr{O}_0(\kpar)\setminus\{0\}$.
From the previous discussion we conclude the following. Fix any $\kpar\in\mathscr{I}_1\subset\subset(2\pi/3,4\pi/3)\setminus\{\pi\}$.
Let $\mathscr{O}_0(\kpar)$ denote the open neighborhood in $\C$ defined in \eqref{scrO_0}. 
Then, for all $z$ in  $\mathscr{O}_0(\kpar)$,   the mapping
\[ z\in\mathscr{O}_0(\kpar)\mapsto \left(\HTBs(\kpar)-z I\right)^{-1}\ \textrm{is meromorphic with values in  }\ \ 
 l^2(\N_0;\C^2)\]
with only one pole, located at $z=0$. Moreover, for $z\in \mathscr{O}_0(\kpar)\setminus\{0\}$ we have
\begin{equation}
\left(\HTBs(\kpar)-z I\right)^{-1}f\ =\ \frac{1}{z} \left\langle \psi^{\rm TB, bd}(\kpar),f\right\rangle_{_{l^2(\N_0;\C^2)}}\  \psi^{\rm TB, bd}(\kpar)\ +\ \mathscr{G}_{\rm reg}(z;\kpar)f,
\label{laurent}\end{equation}
where $z\mapsto \mathscr{G}_{\rm reg}(z;\kpar)$ is an analytic map from $\mathscr{O}_0(\kpar)$
to $\mathcal{B}(l^2(\N_0;\C^2))$. Thus we have proved part (3a) of Theorem \ref{zz-spec}, except for the case $\kpar=\pi$. 
We leave this as an exercise for the reader. 

Note that for all $\kpar\in(2\pi/3,4\pi/3)$, we have that
\begin{equation}
  \HTBs(\kpar)\psi = f\in l^2(\N_0;\C^2)\ \textrm{ is solvable in}\  l^2(\N_0;\C^2)\ \ \iff\ \  \left\langle \psi^{\rm TB, bd}(\kpar),f\right\rangle_{_{l^2(\N_0;\C^2)}}=0.\label{solvability}
  \end{equation}
Thus we have proved all assertions of Theorem \ref{zz-spec}  for  $\kpar\in \mathscr{I}_1$
( $\mathscr{I}_1$ arbitrary compact subset of $(2\pi/3,4\pi/3)$, 
and all $E$ in the open complex neighborhood $\mathscr{O}_0(\kpar)$, defined in \eqref{scrO_0}.
\medskip

 It remains to address the cases:\\
  (A) $\kpar\in[0,2\pi]\setminus(2\pi/3,4\pi/3)$ and $z\in\mathscr{O}_0(\kpar)$, defined in \eqref{scrO_0} and\\
   (B) $\kpar\in[0,2\pi]$ and $z\in\mathscr{O}_+(\kpar)$, defined in \eqref{scrO_+}.
   \medskip
   
\nit In case (A), Lemma \ref{whereroots} tells us that $|\lambda_1(z)|<1<|\lambda_2(z)|$.
 Hence, the construction of the resolvent is as above, and gives the map $f\mapsto \psi$
 defined by \eqref{psi-resf}. However now, since $z=0$ is not an eigenvalue, $\mu=\mu(f;z,\zeta)$ does not have a pole, as was the case in for for $\kpar\in(2\pi/3,4\pi/3)$; 
  see \eqref{mufEz}.

\nit In case (B), Lemma \ref{whereroots} tells us that    $|\lambda_2(z)|<1<|\lambda_1(z)|$. 
The construction of the resolvent is analogous with the roles of the eigenpairs:  $(\lambda_1,\xi_1)$ and  $(\lambda_2,\xi_2)$ interchanged. Since in $\mathscr{O}_+(\kpar)$ 
 $|z|>|\Re z|>\delta_{\rm max}(\kpar)\ge1$ and the only possible eigenvalue is at $z=0$,  the analogue of the $\mu(f;z,\zeta)-$ term in \eqref{psi-resf} does not have a pole in this case as well.
Therefore, in both cases (A) and (B) the mapping $z\mapsto\left(\ \HTBs-z I\ \right)^{-1}$ is analytic with values in $\mathcal{B}(l^2(\N_0;\C^2))$.

Finally, using part (2) of  Lemma \ref{whereroots}, one can check that $\HTBs(\kpar)-z\ I$ is not invertible for $\delta_{\rm gap}(\kpar)\le|z|\le\delta_{\rm max}(\kpar)$ since the eigenvalues of  $M(z,\zeta)$ satisfy: $|\lambda_1(z,\zeta)|=|\lambda_2(z,\zeta)|=1$.
 Such energies $z$ comprise  the essential spectrum of $\HTBs(\kpar)$, $\sigma_{\rm ess}\left(\HTBs(\kpar)\right)$.
 The details are left to the reader.
 This completes the proof of Theorem \ref{zz-spec}.

\section{Setup for the continuum problem; zigzag edge Hamiltonian and the zigzag edge-state eigenvalue problem}\label{setup}

In this section we begin our detailed formulation and discussion of the continuum edge state eigenvalue problem. For this we must first discuss the atomic, bulk and edge Hamiltonians: $\Hatom$, $\cH$ and $\cHeg$. 

\subsection{The atomic Hamiltonian and its ground state}\label{atomic-gs}

We work with the class of ``atomic potential wells '' introduced in  \cite{FLW-CPAM:17}. Fix a smooth potential  $V_0(\bx)$ on $\R^2$  with the following properties.
\begin{itemize}
\item[($PW_1$)] $-1\le V_0(\bx)\le0,\ \bx\in\R^2$.
\item[($PW_2$)] ${\rm supp}\ V_0\ \subset\ \{\bx\in\R^2: |\bx|<r_0\}$, where $r_0<r_{\rm cr}$. Here, $r_{\rm cr}$ is a universal constant defined in \cite{FLW-CPAM:17} satisfying
$ 0.33 |\be|\le r_{\rm cr}< 0.5|\be|$,  and 
$|\be|=|\bv_B-\bv_A|= 1/\sqrt3$ is the distance between one vertex in $\mathbb{H}$ and any nearest neighbor.
\item[($PW_3$)]
 $V_0(\bx)$ is invariant under a $2\pi/3$ ($120^\circ$) rotation about the origin, $\bx=0$.
\item[($PW_4$)] $V_0(\bx)$ is inversion-symmetric with respect to the origin; $V_0(-\bx)=V_0(\bx)$.
\end{itemize}

Consider the ``atomic'' Hamiltonian: 
 $\Hatom=  -\Delta + \lambda^2 V_0(\bx)$
 acting in $L^2(\R^2)$. Let $p_0^\lambda(\bx), E_0^\lambda$, respectively, be the ground state eigenfunction and its strictly negative ground state eigenvalue:
 \begin{equation}
 \left(\ -\Delta\ +\ \lambda^2V_0(\bx)\ -\ E_0^\lambda\ \right)p_0^\lambda(\bx)=0,\ \ p_0^\lambda\in L^2(\R^2), \ E_0^\lambda<0.
 \label{gs-pair}\end{equation}
  This eigenpair is simple and, by the symmetries of $V_0(\bx)$,  the ground state $p_0^\lambda(\bx)$ is invariant under a $\pi/3$ ($60^\circ$) rotation about the origin. 
 We choose $p_0^\lambda(\bx)$ so that $p_0^\lambda(\bx)>0$ for all $\bx\in\R^2$ and 
$
\int_{\R^2} |p_0^\lambda(\bx)|^2\ d\bx=1.
$
Since $V_0\in L^\infty(\R^2)$ and  $-\Delta p_0^\lambda=(E-\lambda^2V_0)p_0^\lambda$, it follows that $p_0^\lambda\in H^2(\R^2)$.
\medskip

Recall  the hopping coefficient $\rho_\lambda$ given by:
\begin{equation}
  \rho_\lambda\ =\ \int_{|\by|<r_0}\ p_{_0}^\lambda(\by) \lambda^2\ |V_0(\by)|\ p_0^\lambda(\by-\be)\ d\by \ .
  \label{rho-defA}
  \end{equation}
  By Proposition 4.1 of \cite{FLW-CPAM:17} we have, under hypotheses $(PW_1),\dots,(PW_4)$ on $V_0(\bx)$ the upper and lower bounds for large $\lambda$ :
 \begin{equation}
  e^{-c{_{-}}\lambda}\ \lesssim\ \rho_\lambda\ \lesssim\  e^{-c_{_{+}}\lambda}\ 
  \label{rho-boundsA}
  \end{equation}
  for some constants: $0<c_+<c_-$ which depend on $V_0$ but not on $\lambda$ .
  \begin{remark}\label{hatc}
  The edge states we construct will have energies $E^\lambda=E_0^\lambda+\Omega^\lambda$, 
  with $\rho_\lambda^{-1} |\Omega^\lambda| \ll1$. In preparation for our later discussion, it is useful at this stage 
  to introduce a positive constant, $\hat{c}$, such that $\hat{c}>c_-$ (see \eqref{rho-boundsA}) and to observe that 
  \[ |\Omega^\lambda | <e^{-\hat{c}\lambda}\quad\implies\ \rho^{-1}_\lambda |\Omega^\lambda | <e^{-(\hat{c}-c_-)\lambda}\downarrow0
\ \textrm{as}\ \lambda\uparrow\infty.\]

  \end{remark}

In addition to hypotheses $(PW_1),\dots,(PW_4)$ on $V_0(\bx)$, we assume the following two spectral properties of 
 $\Hatom=-\Delta + \lambda^2 V_0$ acting on $L^2(\R^2)$:
 \medskip
 
\nit{\bf (GS) Ground state energy upper bound:} For $\lambda$ large, $E_0^\lambda$, the ground state energy of  $-\Delta + \lambda^2 V_0(\bx)$, 
satisfies the upper bound 
\begin{equation} E_0^\lambda\  \le\ -c_{\rm gs}\ \lambda^2 .
\label{GS}\end{equation}
Here, $c_{\rm gs}$ is a strictly positive constant depending on $V_0$. A simple consequence of the variational characterization of $E_0^\lambda$ is the lower bound $E_0^\lambda \ge -\|V_0\|_{_{L^\infty}} \lambda^2=-\lambda^2$. However, the upper bound \eqref{GS} requires further restrictions on $V_0$.
Using the condition (GS), we can show that  $p_0^\lambda$, satisfies the
following pointwise bound:
 \begin{align}
 | p_0^\lambda(\bx) | &\le C_1\left(\ \lambda\ {\bf 1}_{_{|\bx|< r_0+\delta_0}}\ +\ e^{-c_1\lambda|\bx|}\ \right)
  \label{p0-bd-cpam}
  \end{align}
  where $\supp(V_0)\subset B(0,r_0)$,  $\delta_0>0$ is arbitrary, and $C_1$ and $c_1$ are constants that depend on 
   $V_0$, $r_0$ and $\delta_0$; see Corollary 15.5 of \cite{FLW-CPAM:17}.
 \medskip

\nit{\bf (EG) Energy gap property:} For $\lambda>0$ sufficiently large, there exists $c_{\rm gap}>0$, independent of $\lambda$, such that if $\psi\in H^2(\R^2)$ and $\left\langle p_0^\lambda,\psi\right\rangle_{_{L^2(\R^2)}} =\ 0$, then
\begin{equation}
\left\langle\ \left(-\Delta + \lambda^2 V_0 - E_0^\lambda\right)\psi,\psi\ \right\rangle_{_{L^2(\R^2)}}\ \ge\ c_{\rm gap}\ \|\psi\|_{_{L^2(\R^2)}}^2.\label{EG}
\end{equation}

 In Section 4.1 of \cite{FLW-CPAM:17} we discuss examples of potentials for which $-\Delta+\lambda^2 V_0$ satisfies (GS) and (EG).

\subsection{Review of terminology and formulation}\label{review}

We conclude this section with a review of some terminology and the formulation of the edge state eigenvalue problem. 

\begin{enumerate}
\item {\it Continuum bulk Hamiltonian, $\cH$:}\ 
\begin{equation}
\cH \equiv -\Delta + \lambda^2 V(\bx)\qquad  \textrm{acting on}\ \ L^2(\R^2)\ .
\end{equation}
Here, $V(\bx)$, the {\it bulk periodic potential}, is  defined to be the sum of all translates of atomic wells, $V_0(\bx-\bv)$, where $\bv$ ranges over $\mathbb{H}$:
$
 V(\bx)\ =\ \sum_{\bv\in\mathbb{H}} V_0(\bx-\bv) ;
$
see \eqref{Vbk}.

The potential $V(\bx)$ is a honeycomb lattice potential in the sense of Definition 2.1 of \cite{FW:12};
  $V$ is real-valued, and with respect to an origin placed at the center of a regular hexagon of the tiling of $\R^2_\bx$: $V$ is inversion symmetric and rotationally invariant by $2\pi/3$. 
  \item {\it Continuum zigzag edge Hamiltonian, $H^\lambda_{_{\rm edge}}$:}\label{eg-Ham}
The potential for a honeycomb structure interfaced with the vacuum along a sharp interface with direction $\vtilde_2\in\Lambda$
 (parallel to the zigzag edge) is obtained by summing translates of $V_0$  over the truncated structure, $\mathbb{H}_\sharp$,
defined in  \eqref{bbH+}:
\begin{align} 
V_{\sharp}(\bx)\ &=\ \sum_{\bv\in\mathbb{H}_\sharp} V_0(\bx-\bv)\ .\label{Vsharp}
\end{align}
The Hamiltonian for the truncated structure is given by
\begin{equation}
H^\lambda_{_{\rm edge}}\ \equiv -\Delta + \lambda^2 V_\sharp(\bx),\qquad  \textrm{acting on}\ \ L^2(\R^2)\ ,
\label{Hedge}\end{equation}
and its centering at the ground state energy, $E_0^\lambda$, of $\Hatom$ is denoted:
\begin{equation}
\cHeg \equiv -\Delta + \lambda^2 V_\sharp(\bx) - E_0^\lambda\qquad  \textrm{acting on}\ \ L^2(\R^2)\ .
\label{Hsharp}\end{equation}
Since $H^\lambda_{_{\rm edge}}$ and  $\cHeg$ are invariant under the translation invariance: $\bx\mapsto \bx+\vtilde_2$, these operators  
 act in $L_\kpar^2(\Sigma)$, $\Sigma =\R^2/\Z\vtilde_2$.
\item The {\it $\kpar-$ dependent Edge Hamiltonian}, $\cHeg(\kpar)$,  acting in $L^2(\Sigma)$  is given by:
\begin{equation}
\cHeg(\kpar) \equiv -\left(\nabla+i\frac{\kpar}{2\pi}\ktilde_2\right)^2 + \lambda^2 V_\sharp(\bx)\ -\ E_0^\lambda\ .  
\label{Hsharp-k}\end{equation}
\end{enumerate}

Finally we recall that  the {\it Zigzag Edge state Eigenvalue Problem} is given by \eqref{zz-evp}, or equivalently, 
\eqref{zz-kp-evp}. With $E=E_0^\lambda+\Omega$, we have:
\begin{align}
 \left(\cHeg(\kpar)\ -\ \Omega\right)\psi\ =\ 0\ ,\quad \psi\in L^2_\kpar.
\label{kpar-evp}\end{align}

\section{A natural subspace of $L^2_{_\kpar}(\Sigma)$}\label{nat-basis}

Define, for all $n\ge0$
\footnote{The labeling convention of $\bA-$ points and $\bB-$ sublattice points used in the present article differs from  that used in \cite{FLW-CPAM:17}. This has no effect on the results in this article or in  \cite{FLW-CPAM:17}. }
\begin{equation}
\bv_\bA^n\ \equiv\ \bv_\bA+n\vtilde_1,\qquad \bv_\bB^n\ \equiv\ \bv_\bB+n\vtilde_1,
\label{vAB_n}
\end{equation}
where $\bv_A^0=\bv_A$ and $\bv_B^0=\bv_B$.
The cylinder $\Sigma=\R^2/ \Z\vtilde_2$ has fundamental domain $\Omega_\Sigma\subset\R^2$, which may be expressed as the union of paralleograms: 
\begin{equation} \Omega_\Sigma = \cup_{n\ge0}\ \Omega_n\ \cup\ \Omega_{-1}\quad \textrm{as in Figure \ref{zz-dimers}\ .} 
\label{OmSig}
\end{equation}
Each parallelogram $\Omega_n$ with $n\ge0$ contains two atomic sites: $\bv_\bA^n$ and $\bv_\bB^n$.  The infinite parallelogram, $\Omega_{-1}$, contains no atomic sites. A fundamental cell of the cylinder $\Sigma$, $\Omega_\Sigma$,  and its decomposition into parallelograms $\Omega_n$, for $n\ge-1$ is depicted in Figure \ref{zz-dimers}.
The zigzag sharp truncation of $\mathbb{H}$ may be expressed as a union over ``vertical translates'' (translates with respect to $\vtilde_2$) of sites within $\Omega_\Sigma$:
  \[ \mathbb{H}_\sharp\ =\ \cup_{n_2\in\Z}\ \cup_{\substack{n_1\ge0}}\ \
   \Big\{\ \bv^{n_1}_{\bA}+n_2\vtilde_2\ ,\  \bv^{n_1}_{\bB}+n_2\vtilde_2\   \Big\}\ .
  \]

\nit We next introduce  approximate $\kpar-$ pseudo-periodic solutions of $\cHeg\Psi=0$ via $\kpar-$ pseudo-periodization  of the atomic ground state, $p_0^\lambda$:

\begin{definition}\label{pzn-p0}
Fix $\kpar\in [0,2\pi]$ and $\bI=\bA, \bB$. 
For each $n\in\N_0\equiv\{0,1,2,\dots\}$, define 
\begin{align}\label{pkpar}
p_{_{\kpar,I}}^{^{\lambda}}[n](\bx)\ &\equiv\ \sum_{m_2\in\Z}\ 
p_{_0}^\lambda(\bx-\bv_\bI^n-m_2\vtilde_2)\ e^{-i\frac{\kpar}{2\pi}\ktilde_2\cdot(\bx-\bv_\bI^n-m_2\vtilde_2)}\\
&=\ e^{-i\frac{\kpar}{2\pi}\ktilde_2\cdot(\bx-\bv_\bI)} \ \sum_{m_2\in\Z}\ 
e^{i\kpar m_2}\ p_{_0}^\lambda(\bx-\bv_\bI^n-m_2\vtilde_2)\ 
\nn\end{align}
and 
\begin{align}\label{Pkpar}
P_{_{\kpar,I}}^\lambda[n](\bx)\ &\equiv\ e^{ i\frac{\kpar}{2\pi}\ktilde_2\cdot(\bx-\bv_\bI) }\ p_{_{\kpar,I}}^\lambda[n](\bx)
\ =\ \sum_{m_2\in\Z}\ e^{i\kpar m_2}\
p_{_0}^\lambda(\bx-\bv^n_\bI-m_2\vtilde_2)\ .
\end{align}
\end{definition}
\nit   The function $\bx\mapsto p_{_{\kpar,I}}^{^{\lambda}}[n](\bx)$ is defined on the cylinder $\Sigma$, {\it i.e.}
$ p_{_{\kpar,I}}^{^{\lambda}}[n](\bx+\vtilde_2)=p_{_{\kpar,I}}^{^{\lambda}}[n](\bx)$. 
  To see this,  
 replace $\bx$ by $\bx+\vtilde_2$ and redefine the summation index. Furthermore, we note that:
$P_{_{\kpar,I}}^\lambda[n](\bx+\vtilde_2)=e^{i\kpar}P_{_{\kpar,I}}^\lambda[n](\bx)$. 
 
 The functions: $p_{_{\kpar,I}}^{\lambda}[n]$, $\bI=\bA,\bB$, $n\ge0$,  form a nearly orthonormal set in $L^2(\Sigma)$ for large $\lambda$. In particular, we have:
\begin{proposition}\label{prop:pkpar} Fix $\kpar\in[0,2\pi]$ and $\lambda>0$. 
\begin{enumerate}
\item For all $n\in\N_0$, we have
$p_{_{\kpar,I}}^{\lambda}[n]\in L^2(\Sigma)$ and $P_{_{\kpar,\bI}}^{\lambda}[n]\in L^2_\kpar$.
\end{enumerate}
 Furthermore, there exist constants $\lambda_\star, c>0$ such that for all $\lambda\ge\lambda_\star$:
\begin{enumerate}
\item[(2)] For $n\in\N_0$, $\bI=\bA,\bB$ 
\begin{equation}\label{near-orth1}
\Big|\ \left\langle p_{_{\kpar,\bI}}^{^{\lambda}}[n],p_{_{\kpar,\bJ}}^{^{\lambda}}[n]\right\rangle_{_{L^2(\Sigma)}}\ -\ \delta_{_{\bI\bJ}}\ \Big| \
\lesssim\ e^{-c\lambda}\ ,
\end{equation}
where $\delta_{_{\bI\bJ}}$ denotes the Kronecker delta symbol.
\item[(3)] For $\bI=\bA,\bB$, $m,n\in\N_0$ with $m\ne n$ and all $\lambda>0$ sufficiently large: 
\begin{equation}\label{near-orth2}
\Big|\ \left\langle p_{_{\kpar,\bI}}^{\lambda}[m],p_{_{\kpar,\bJ}}^\lambda[n]\right\rangle_{_{L^2(\Sigma)}}\ \Big|
 \lesssim\ e^{-c\lambda|m-n|}\ .
\end{equation}
\end{enumerate}
Assertions \eqref {near-orth1} and \eqref{near-orth2} hold as well with $p_{_{\kpar,\bI}}^{\lambda}[m]$ replaced by $P_{_{\kpar,\bI}}^{\lambda}[m]$,
defined in \eqref{Pkpar}, and with $L^2(\Sigma)$ replaced by $L^2_\kpar(\R^2)$. Here, $\lambda_\star$ depends only on $V$.
\end{proposition}
This proposition follows from the normalization and decay properties of the atomic ground state, $p_0^\lambda$; the details are omitted.

We conclude this section by showing that the functions $p_{_{\kpar,\bI}}^{\lambda}[n]$, $\bI=\bA,\bB$, $n\ge0$,
 are nearly annihilated by $\cHeg(\kpar)$.
\begin{proposition}\label{Hp_ap0}
There exist positive constants $\lambda_\star$ (large) and $c>0$, such that  for all $\lambda>\lambda_\star$ and all $\bI=\bA, \bB$ and $n\ge0$:
\begin{align}
\left|\ \cHeg(\kpar) p_{_{\kpar,\bI}}^{\lambda}[n](\bx)\ \right|\  &\lesssim\ e^{-c|\bx-n\vtilde_1|}\ e^{-c\lambda},\ \ \bx\in\Omega_\Sigma\label{Hpsi-decay}\\
\left\| \cHeg(\kpar) p_{_{\kpar,\bI}}^{\lambda}[n] \right\|_{L^2(\Sigma)}\ &\lesssim\ e^{-c\lambda}\ .
\label{Hpsi-small}
\end{align}
\end{proposition}

\nit{\it Proof of Proposition \ref{Hp_ap0}:} We first note that \eqref{Hpsi-small} follows from \eqref{Hpsi-decay}
by integrating the square of bound  \eqref{Hpsi-decay} over a fundamental domain (strip), $\Omega_\Sigma$. Thus we focus on the pointwise bound \eqref{Hpsi-decay}.   The identity $\nabla_\bx=e^{i\frac{\kpar}{2\pi}\ktilde_2\cdot\bx}\ (\nabla+i\frac{\kpar}{2\pi}\ktilde_2) 
e^{-i\frac{\kpar}{2\pi}\ktilde_2\cdot\bx}$ and  \eqref{gs-pair} imply that for arbitrary $\hat\bv\in\R^2$:
\begin{equation}
\left(\ -\left(\nabla+i\frac{\kpar}{2\pi}\ktilde_2\right)^2\ +\ \lambda^2V_0(\bx-\hat\bv)\ -\ E_0^\lambda\ \right) e^{-i\frac{\kpar}{2\pi}\ktilde_2\cdot(\bx-\hat\bv)}p_0^\lambda(\bx-\hat\bv)=0\ ;
\label{p0-shift}
\end{equation}
we shall apply \eqref{p0-shift} for $\hat\bv\in\mathbb{H}_\sharp$.

As a first step toward obtaining the bound \eqref{Hpsi-decay} for $\cHeg(\kpar) p_{_{\kpar,\bI}}^{\lambda}[n](\bx)$,
we observe that
\[
\ \textrm{for}\ \bx\in\Omega_\Sigma,\qquad  V_\sharp(\bx)\ =\ \sum_{\bJ=\bA,\bB}\ \sum_{n_1\ge0}\ V_0(\bx-\bv_\bJ-n_1\vtilde_1).\]
Therefore, for $\bx\in\Omega_\Sigma$ we have
\begin{align}
&\cHeg(\kpar) p_{_{\kpar,\bI}}^{\lambda}[n](\bx)\ 
=\ \sum_{m_2\in\Z}\ \cHeg(\kpar)\ e^{-i\frac{\kpar}{2\pi}\ktilde_2\cdot(\bx-\bv_\bI^n-m_2\vtilde_2)}\ p_0^\lambda(\bx-\bv^n_\bI-m_2\vtilde_2)\nn\\
&=\ \cHeg(\kpar)\ e^{-i\frac{\kpar}{2\pi}\ktilde_2\cdot(\bx-\bv_\bI^n)}\ p_0^\lambda(\bx-\bv^n_\bI)\nn\\
&\qquad\ +\
 \sum_{m_2\in\Z\setminus\{0\}}\ \left(\ -\left(\nabla+i\frac{\kpar}{2\pi}\ktilde_2\right)^2\ -\ E_0^\lambda\ \right) e^{-i\frac{\kpar}{2\pi}\ktilde_2\cdot(\bx-\bv_\bI^n-m_2\vtilde_2)}\ p_0^\lambda(\bx-\bv^n_\bI-m_2\vtilde_2)\nn\\
&\qquad\ +\  \sum_{m_2\in\Z\setminus\{0\}}\ \lambda^2 V_\sharp(\bx)\ e^{-i\frac{\kpar}{2\pi}\ktilde_2\cdot(\bx-\bv_\bI^n-m_2\vtilde_2)}\ p_0^\lambda(\bx-\bv^n_\bI-m_2\vtilde_2)\ .\nn
\end{align}
In the second equality just above we have split off the $m_2=0$ and $m_2\ne0$ contributions.
The first term of the $m_2\ne0$ contribution vanishes identically for $\bx\in\Omega_\Sigma$.
 Indeed,  equation \eqref{p0-shift} for $p_0^\lambda$ implies that this term is  a sum of terms, each containing a factor 
  $\lambda^2V_0(\bx-\bv^n_\bI-m_2\vtilde_2)$ for some $m_2\in\Z\setminus\{0\}$. Each of these terms  vanishes since the constraint: $m_2\ne0$ implies they are all supported outside of $\Omega_\Sigma$. Therefore,
  \begin{align}
\cHeg(\kpar) p_{_{\kpar,\bI}}^{\lambda}[n](\bx)\ 
&=\ \cHeg(\kpar)\ e^{-i\frac{\kpar}{2\pi}\ktilde_2\cdot(\bx-\bv_\bI^n)}\ p_0^\lambda(\bx-\bv^n_\bI)\nn\\
&\qquad\ +\  \sum_{m_2\in\Z\setminus\{0\}}\ \lambda^2 V_\sharp(\bx)\ e^{-i\frac{\kpar}{2\pi}\ktilde_2\cdot(\bx-\bv_\bI^n-m_2\vtilde_2)}\ p_0^\lambda(\bx-\bv^n_\bI-m_2\vtilde_2)\ .\label{Heg-cncl}
\end{align}
We may now use \eqref{p0-shift} with $\hat\bv=\bv_\bI^n=\bv_\bI+n\vtilde_1$ to simplify the first term on the right hand side of the previous equation. For all $\bx\in\Omega_\Sigma$ with $n\ge0$ and $\bI, \bJ\in\{\bA,\bB\}$ with $\bI\ne\bJ$, we obtain:
  \begin{align}
\cHeg(\kpar) p_{_{\kpar,\bI}}^{\lambda}[n](\bx)\ 
&=\ \left(\ \lambda^2\sum_{\substack{n_1\ge0\\ n_1\ne n}} V_0(\bx-\bv_\bI^{n_1})\ \right)
 e^{-i\frac{\kpar}{2\pi}\ktilde_2\cdot(\bx-\bv_\bI^n)}\ p_0^\lambda(\bx-\bv_\bI^n) \nn\\
 &\qquad\ +\ \left(\ \lambda^2\sum_{n_1\ge0} V_0(\bx-\bv^{n_1}_\bJ)\ \right)
 e^{-i\frac{\kpar}{2\pi}\ktilde_2\cdot(\bx-\bv_\bI^n)}\ p_0^\lambda(\bx-\bv_\bI^n)\nn\\
&\qquad\ +\  \sum_{m_2\in\Z\setminus\{0\}}\ \lambda^2 V_\sharp(\bx)\ e^{-i\frac{\kpar}{2\pi}\ktilde_2\cdot(\bx-\bv_\bI^n-m_2\vtilde_2)}\ p_0^\lambda(\bx-\bv^n_\bI-m_2\vtilde_2)\ .\nn
\end{align}
Thus,
 \begin{align}
&\left|\ \cHeg(\kpar) p_{_{\kpar,\bI}}^{\lambda}[n](\bx)\ \right|\nn\\ 
&\le\ \left(\ \lambda^2\sum_{\substack{n_1\ge0\\ n_1\ne n}} |V_0(\bx-\bv_\bI^{n_1})|\ \right)
\ p_0^\lambda(\bx-\bv_\bI^n)\ +\ \left(\ \lambda^2\sum_{n_1\ge0} |V_0(\bx-\bv^{n_1}_\bJ)|\ \right)
  p_0^\lambda(\bx-\bv_\bI^n)\nn\\
  &\qquad\qquad +\  \sum_{m_2\in\Z\setminus\{0\}}\ \lambda^2 |V_\sharp(\bx)|\  p_0^\lambda(\bx-\bv^n_\bI-m_2\vtilde_2)\nn\\
  & \equiv\ T_1(\bx;n)\ +\ T_2(\bx;n)\ +\ T_3(\bx;n)\ .\label{Hp3}
\end{align}
To bound the first term of \eqref{Hp3}, we note that for $n_1\ne n$
\begin{align}
 |V_0(\bx-\bv_\bI^{n_1})|\ p_0^\lambda(\bx-\bv_\bI^n)&\le\ \|V_0\|_\infty\ {\bf 1}_{|\bx-\bv_\bI^{n_1}|<r_0}\ p_0^\lambda(\bx-\bv_\bI^n)\nn\\
 &\lesssim\ {\bf 1}_{|\bx-\bv_\bI^{n_1}|<r_0}\ e^{-c\lambda|\bx-\bv_\bI^n|}\nn\\
 &\lesssim\ {\bf 1}_{|\bx-\bv_\bI^{n_1}|<r_0}\ e^{-\frac{c}2\lambda|\bx-\bv_\bI^n|}\ e^{-\tc\lambda|n_1-n|}\ .
\nn \end{align}
Summing over $n_1\ge0$ with $n_1\ne n$ we obtain $T_1(\bx;n)\lesssim e^{-c^\prime\lambda}\ 
e^{-c^{\prime\prime}\lambda|\bx-\bv_\bI^n|}$. Very similarly we obtain: $T_2(\bx;n)\lesssim e^{-c^\prime\lambda}\ 
e^{-c^{\prime\prime}\lambda|\bx-\bv_\bI^n|}$. 
We finally consider $T_3(\bx;n)$. For $\bx\in\Omega_\Sigma$,
\begin{align}
T_3(\bx;n)\ &\lesssim\ \lambda^2\ \|V_0\|_\infty\ \sum_{m_2\in\Z\setminus\{0\}}\ e^{-c\lambda|\bx-\bv^n_\bI|}\ e^{-c\lambda |m_2|}\ \lesssim\ e^{-c^\prime\lambda}\ e^{-c^{\prime\prime}\lambda|\bx-\bv^n_\bI|}\ .
\nn\end{align}
This completes the proof of Proposition \ref{Hp_ap0}.

\subsection{The subspace  $\mathscr{X}^\lambda_{AB}(\kpar)$}\label{XAB}
We introduce the closed subspace of $L^2(\Sigma)$:
\begin{equation}
\mathscr{X}^\lambda_{AB}(\kpar)= \textrm{the orthogonal complement in $L^2(\Sigma)$ of}\ 
\textrm{span}\Big\{p_{_\kpar}^{^{\lambda,\bI}}[n]: I=A,B;\ n\ge0\Big\}\ .
\label{XAB-def}\end{equation}
We shall sometimes suppress the dependence on $\lambda$ and write $\mathscr{X}_{AB}(\kpar)$.
The space $L^2(\Sigma)$ may be decomposed as the orthogonal sum of subspaces:
\begin{equation}
L^2(\Sigma)\ =\ \textrm{span}\Big\{p_{_{\kpar,\bI}}^{^{\lambda}}[n]: I=A,B;\ n\ge0\Big\}\oplus 
\mathscr{X}_{AB}(\kpar)\ .
\end{equation}

We also introduce the orthogonal projection onto $\mathscr{X}_{AB}(\kpar)$:
   \begin{equation} \Pi_{_{AB}}=\Pi_{_{AB}}(\kpar): L^2(\Sigma)\to \mathscr{X}_{AB}(\kpar) . \label{Xkproj}\end{equation}
   
  Since the set $\Big\{p_{_{\kpar,\bI}}^{^{\lambda}}[n]: I=A,B;\ n\ge0\Big\}$ is only nearly-orthonormal for $\lambda$ large (Proposition \ref{prop:pkpar}), we make use of the following:
     
   \begin{proposition}\label{Xdecomp}
  There exists $\lambda_\star>0$ such that for all $\lambda>\lambda_\star$ the following holds.  
Fix $\kpar\in[0,2\pi]$.
  \begin{enumerate}
  \item  Then, 
  for $F\in L^2(\Sigma)$ we have that
  \[ F\equiv0\ \ \iff\ \ \Pi_{AB}(\kpar)F=0\ \ {\rm and}\ \  \left\langle p_{_{\kpar,I}}^\lambda[n],F\right\rangle_{_{L^2(\Sigma)}}=0,\ \ n\ge0,\ I=A, B\ .\]
  \item Any $\psi\in L^2(\Sigma)$ may be expressed in the form:
  \begin{equation}  \label{psi-decomp}
  \psi\ =\ \sum_{J=A,B}\ \sum_{n\ge0} \alpha_n^J\ p_{_{\kpar,J}}^\lambda[n]\ +\ \widetilde{\psi},
  \end{equation}
  where $\alpha=\{(\alpha_n^A,\alpha_n^B)^{^\top}\}_{n\ge0}\in l^2(\N_0;\C^2)$ and $\Pi_{_{AB}}(\kpar)\widetilde\psi=\widetilde\psi\in\mathscr{X}^\lambda_{AB}(\kpar)$.
  \end{enumerate}
     \end{proposition}
\nit The proof is similar to that of Lemma 8.2 on page 31 of \cite{FLW-CPAM:17} and is omitted.

\section{Energy estimates and the resolvent }\label{en-est}

The following proposition concerns the invertibility of $\Pi_{_{AB}}(\kpar)\left(\ \cHeg(\kpar)-\Omega\ \right)\Pi_{_{AB}}(\kpar)$ on $\mathscr{X}_{AB}(\kpar)$ for $\lambda$ sufficiently large. This will facilitate reduction of the edge state eigenvalue problem,
 \eqref{zz-evp} or \eqref{zz-kp-evp},  to a problem on the linear space
 ${\rm span}\Big\{p_{_{\kpar,J}}^\lambda[n]\ :\ \bI=\bA,\bB,\ n\ge0\Big\}$; see \eqref{XAB-def}.
 The proof uses
arguments analogous to those in \cite{FLW-CPAM:17}. The necessary modifications in the strategy are discussed at the end of this section. 

\begin{proposition}\label{resolvent}
There exist constants $\lambda_\star>0$ (sufficiently large) and $c^\prime>0$ (sufficiently small), such that for all $\lambda>\lambda_\star$, $\kpar\in[0,2\pi]$ and $|\Omega|\le c^\prime$
  the following hold:
\begin{enumerate}
\item  For all $\varphi\in \mathscr{X}_{_{AB}}(\kpar)$, the equation 
 \begin{equation}
\Pi_{_{AB}}(\kpar)\ \left(\ \cHeg(\kpar)-\Omega\ \right)\psi\ =\ \varphi\ ,
\label{non-hom}
\end{equation}
has a unique solution 
\[ \psi\ \equiv\ \Ressp^{\lambda}(\Omega,\kpar)\varphi\ \in \mathscr{X}_{_{AB}}\cap H^2(\Sigma).\]\\
Thus, $\Ressp^{\lambda}(\Omega,\kpar)$ is the inverse of $\Pi_{_{AB}}(\kpar)\left(\ \cHeg(\kpar)-\Omega\ \right)\Pi_{_{AB}}(\kpar)$ or equivalently $
 \Pi_{_{AB}}(\kpar)\left(\ \cHeg(\kpar)-\Omega\ \right)$ acting on $\mathscr{X}_{_{AB}}$.
 \item The mapping
$\varphi\mapsto \Ressp^{\lambda}(\Omega,\kpar)\varphi$ is a bounded linear operator :
 \begin{equation} \Ressp^{\lambda}(\Omega,\kpar):\mathscr{X}_{_{AB}}(\kpar)\to H^2(\Sigma)\cap \mathscr{X}_{_{AB}}(\kpar). \label{Res-def}\end{equation}
 \item We have the following operator norm bounds on $\Ressp^{\lambda}(\Omega,\kpar)$:
 \begin{align}
\left\|\ \Ressp^{\lambda}(\Omega,\kpar)\ \right\|_{_{\mathscr{X}_{AB}\to\mathscr{X}_{AB}}}&\lesssim1
\label{Res-bd1}\\
\lambda^{-1}\ \left\|\ \nabla_\bx\ \Ressp^{\lambda}(\Omega,\kpar)\ \right\|_{_{\mathscr{X}_{AB}\to\mathscr{X}_{AB}}}&\lesssim1
\label{Res-bd2}\\
\left\|\ \Ressp^{\lambda}(\Omega,\kpar)\ \right\|_{_{\mathscr{X}_{AB}\to H^2(\Sigma)\cap\mathscr{X}_{AB}}}
&\le C(\lambda,\kpar)\ .\ \label{Res-bd3}
\end{align}
\item  Furthermore, this mapping depends analytically on $\Omega\in\C$ for $|\Omega|< c^\prime$, 
 and for all such $\Omega$:
 \begin{equation} \ \Big\|\ \partial_\Omega\ \Ressp^{\lambda}(\Omega,\kpar)\ \Big\|_{_{ \mathscr{X}_{_{AB}}\to\mathscr{X}_{_{AB}}   }}\lesssim1.
 \label{DOmega-bd}\end{equation}
\item For real $\Omega\in(-c^\prime,c^\prime)$, $\Ressp^{\lambda}(\Omega,\kpar)$ is self-adjoint on the Hilbert space $\mathscr{X}_{AB}$, endowed with the $L^2(\Sigma)$ inner product. 
 \end{enumerate}
\end{proposition}

\nit A key step to proving Proposition \ref{resolvent}  is the following energy estimate on the space $ \mathscr{X}_{AB}(\kpar)$:

\begin{proposition}[Energy Estimate]\label{energy}
Fix $\kpar\in[0,2\pi]$. There exists $\lambda_\star>0$, independent of $\kpar$,  and a constant $C_\star>0$ such that the following holds for all $\lambda\ge\lambda_\star$. Let 
$\psi\in\mathscr{X}_{AB}(\kpar)\cap H^2(\Sigma)$. That is,
\begin{equation}
 \left\langle p^\lambda_{_{\kpar,J}}[n],\psi\right\rangle_ {L^2(\Sigma)}\ =\ 0,\ n\ge0,\ \ J=A, B\ .
 \label{p-orth}
 \end{equation}
Then, 
\begin{equation}
\|\  \cHeg(\kpar)\psi\ \|^2_{_{L^2(\Sigma)}}\ \ge c_\star\ \left(\ \left\|\psi\right\|^2_{_{L^2(\Sigma)}}\ +\
 \lambda^{-2}\left\|\nabla\psi\right\|^2_{_{L^2(\Sigma)}}\ \right)\ .
\label{main-en}
\end{equation}
The constant $c_\star$ can be taken independent of $\kpar$ but it does depend on properties of the atomic potential, $V_0$, in particular on  the constants $c_{\rm gs}$ and $c_{\rm gap}$; see \eqref{GS} and \eqref{EG}.
\end{proposition}

The proof of Proposition \ref{resolvent} follows the general structure of the proof of the energy estimates in \cite{FLW-CPAM:17} .
We now discuss the modifications in these arguments, which are required to prove Propositions \ref{energy} and \ref{resolvent}.
  We follow the discussion of Section 9 of  \cite{FLW-CPAM:17}  with $\Sigma=\R^2/\Z\vtilde_2$ playing the role of $\R^2/\Lambda$, and with the approximate eigenfunctions $p_{_{\kpar,I}}^\lambda[n]\in L^2(\Sigma)$ playing the role of $p_{_{\bk,I}}^\lambda\in L^2(\R^2/\Lambda)$ in \cite{FLW-CPAM:17} . 
   
For $n\ge0$, let  $\bx_{I}^n,\ I=A,B$   denote the two atomic sites in $\Omega_n$, where $n\ge0$. Recall $\Omega_\Sigma$ is the union, for $n\ge-1$,  over all $\Omega_n$; see Figure \ref{zz-dimers}. In place of the partitions of unity (9.11) in \cite{FLW-CPAM:17} on $\R^2/\Lambda$, we introduce here analogous partitions on $\Sigma$:
   \[ 1\ =\ \Theta_0^2\ +\ \sum_{\substack{n\ge0\\ I=A,B}}\ \Theta_{n,I}^2,\qquad 1\ =\ \widetilde\Theta_0^2\ +\ \sum_{\substack{n\ge0\\ I=A,B}}\ \widetilde\Theta_{n,I}^2\]
   where $\Theta_{n,I}$ and $\widetilde{\Theta}_{n,I}$ are supported near $\bx_I^n$. 
   All the arguments in Sections 9.1 through 9.4 of \cite{FLW-CPAM:17} go through in the above setting, with minimal changes.
    This gives Proposition \ref{energy}.\medskip
    
    We seek to show that the inverse  of 
  $\Pi_{_{AB}}(\kpar)\left(\ \cHeg(\kpar)-\Omega\ \right)\Pi_{_{AB}}(\kpar)$, is a bounded linear operator on $\mathscr{X}^\lambda_{_{AB}}(\kpar)$, satisfying the bound \eqref{Res-bd1} and \eqref{Res-bd2}
 and furthermore that  $\Ressp^{\lambda}(\Omega,\kpar)$ maps $\mathscr{X}^\lambda_{_{AB}}(\kpar)$ to $H^2(\Sigma)\cap\mathscr{X}^\lambda_{_{AB}}(\kpar)$ and satisfies the operator bound \eqref{Res-bd3}.
 
 To adapt Section 9.5 of \cite{FLW-CPAM:17} to our setting requires an additional argument which we now supply. Suppose
    we have $\Pi_{AB}(\kpar) \left[\ \cHeg(\kpar)\ -\ \Omega\ I\ \right]\psi = f$, where   $\psi\in L^2(\Sigma)\cap\mathscr{X}_{_{AB}}^\lambda(\kpar)$ and $f\in L^2(\Sigma)$. Then, for some $\{\alpha_{I,n}\}$, $(I=A,B\ n\ge0)$, in $ l^2(\N_0;\C^2)$:
    \begin{equation}
    \left[\ \cHeg(\kpar)\ -\ \Omega\ I\ \right]\psi\ = \ f\ +\ \sum_{\substack{I=A,B\\n\ge0}}\ \alpha_{I,n}\ p^\lambda_{_{\kpar,I}}[n]\ ,
 \label{inhom}   \end{equation}
where the right hand sum is convergent in $L^2(\Sigma)$  and the left hand side is interpreted as a distribution on $\Sigma$. 
 Taking the inner product in $L^2(\Sigma)$ of \eqref{inhom} with 
    $p_{_{\kpar,J}}^\lambda[m]$, we find that 
    \begin{align}\label{proj-inhom}
    &\ \sum_{\substack{I=A,B\\n\ge0}}\ \alpha_{I,n}\ \left\langle p^\lambda_{_{\kpar,J}}[m] , p^\lambda_{_{\kpar,I}}[n]\right\rangle\ 
    =\ \xi^\lambda_{_{\kpar,J}}[m],\quad \textrm{where}\nn\\
& \qquad\qquad   \xi^\lambda_{_{\kpar,J}}[m]\ \equiv\ \left\langle \cHeg(\kpar) p^\lambda_{_{\kpar,J}}[m], \psi\right\rangle\ -\ 
     \left\langle p^\lambda_{_{\kpar,J}}[m], f\right\rangle\  .
     \end{align}
     We have 
     \begin{equation}
     \Big|\xi^\lambda_{_{\kpar,J}}[m]\Big|^2\ \lesssim\ 
     \Big|\ \left\langle \cHeg(\kpar) p^\lambda_{_{\kpar,J}}[m], \psi\right\rangle\ \Big|^2\ +\ 
     \Big|\ \left\langle p^\lambda_{_{\kpar,J}}[m], f\right\rangle\ \Big|^2
     \label{pre-xi}
     \end{equation}
     and summing over $J= A, B$ and $m\ge0$ yields
      \begin{equation}
   \sum_{\substack{J=A,B\\m\ge0}}\ \Big|\xi^\lambda_{_{\kpar,J}}[m]\Big|^2\ \lesssim\ 
    \sum_{\substack{J=A,B\\m\ge0}}\ \Big|\ \left\langle \cHeg(\kpar) p^\lambda_{_{\kpar,J}}[m], \psi\right\rangle\ \Big|^2\ +\ 
    \sum_{\substack{J=A,B\\m\ge0}}\ \Big|\ \left\langle p^\lambda_{_{\kpar,J}}[m], f\right\rangle\ \Big|^2\ .
     \label{pre-xi1}
     \end{equation}
In order to bound the second term on the right in \eqref{pre-xi1}, note that the near-orthonormality of the set $\{p^\lambda_{_{\kpar,J}}[m]:J=A, B,\ m\ge0\}$ for $\lambda$ large (Proposition \ref{prop:pkpar}) implies the Bessel-type inequality:
        \[ \sum_{\substack{J=A,B\\ m\ge0}}\ \Big| \left\langle p^\lambda_{_{\kpar,J}}[m],f\right\rangle\Big|^2\ \lesssim\ 
        \|f\|_{L^2(\Sigma)}^2\ .\]
   Consider next the first term on the right in \eqref{pre-xi1}. 
   Thanks to  the pointwise bound on $ \cHeg(\kpar) p^\lambda_{_{\kpar,J}}[m](\bx)$ from 
   Proposition \ref{Hp_ap0}, a Young-type inequality yields:
   \[
    \sum_{\substack{J=A,B\\m\ge0}}\ \Big|\ \left\langle \cHeg(\kpar) p^\lambda_{_{\kpar,J}}[m], \psi\right\rangle\ \Big|^2
    \lesssim\ e^{-c\lambda}\ \|\psi\|_{_{L^2(\Sigma)}}^2
     \ .\]
    Again, by Proposition \ref{prop:pkpar}, we have
           \begin{align}
      \sum_{\substack{J=A,B\\ m\ge0}}\ |\alpha_m^I |^2 \ &\lesssim\  \sum_{\substack{J=A,B\\ m\ge0}}\ |\ \xi^\lambda_{\kpar,I}[m]\ |^2\nn\\
        &\qquad \lesssim\ e^{-c\lambda}\ \|\psi\|_{L^2(\Sigma)}^2\ +\ 
        C\ \|f\|_{L^2(\Sigma)}^2\ .
        \label{a1!}\end{align}
 And finally one more application of Proposition \ref{prop:pkpar} gives
        \begin{equation}
       \Big\|\  \sum_{\substack{I=A,B\\ n\ge0}}\ \alpha^\lambda_{\kpar,I}[n]\ p^\lambda_{\kpar,I}[n] \Big\|_{L^2(\Sigma)}\ \lesssim\ 
       e^{-c\lambda}\ \|\psi\|_{L^2(\Sigma)}\ +\ 
        C\ \|f\|_{L^2(\Sigma)}\ .\label{a2!}\end{equation}
The estimates \eqref{a1!} and \eqref{a2!} allow us to argue as in Section 9.5 of \cite{FLW-CPAM:17}, using our energy estimates,
 that the  operator 
 $\Ressp^{\lambda}(\Omega,\kpar)$, the inverse of $\Pi_{_{AB}}(\kpar)\left(\ \cHeg(\kpar)-\Omega\ \right)\Pi_{_{AB}}(\kpar)$, is a bounded linear operator on $\mathscr{X}^\lambda_{_{AB}}(\kpar)$, satisfying the bounds \eqref{Res-bd1} and \eqref{Res-bd2}.
 
 To complete the proof of Proposition \ref{resolvent} must show that $\Ressp^{\lambda}(\Omega,\kpar)$ maps $\mathscr{X}^\lambda_{_{AB}}(\kpar)$ to $H^2(\R^2)\cap \mathscr{X}^\lambda_{_{AB}}(\kpar)$.    To bound  $\|\Delta\psi\|_{_{L^2(\Sigma)}}$, we use  \eqref{inhom} to obtain an expression for $\Delta\psi$ in terms of $\psi$ and $\nabla\psi$. Then,  the energy estimate for $\|\psi\|_{_{L^2(\Sigma)}}$ and 
   $\|\nabla\psi\|_{_{L^2(\Sigma)}}$, and the bound \eqref{a2!} imply that
for $\lambda$ sufficiently large, the $L^2(\Sigma)$ norm of each term in the expression $\Delta\psi$ can be bounded by $C(\lambda)\times\|f\|_{_{L^2(\Sigma)}}$, where $C(\lambda)$ denotes a $\lambda-$ dependent constant.

\section{Lyapunov-Schmidt reduction; formulation as a problem in $\mathscr{X}^\lambda_{_{AB}}(\kpar)$}\label{ls-red}

The resolvent bounds of Proposition \ref{resolvent} ensure that on the subspace $\mathscr{X}_{_{AB}}(\kpar)$, the operator $H_\sharp^\lambda(\kpar)-\Omega$ is invertible in a neighborhood of $\Omega=0$, {\it i.e.} the spectrum of $\Pi_{_{AB}}(\kpar)H_\sharp^\lambda(\kpar)\Pi_{_{AB}}(\kpar)$
 is bounded away from zero, uniformly in $\lambda\gg1$.  
   In this section, we make use of this spectral separation to obtain a reduction of the $L^2_\kpar$ eigenvalue problem 
 to a problem on the subspace  of $L^2(\Sigma)$ given by: $\textrm{span}\Big\{p_{_\kpar}^{^{\lambda,\bI}}[n]: I=A,B;\ n\ge0\Big\}$.
 
Consider the eigenvalue problem:
\begin{equation}
 \left(\ -\left(\nabla+i\frac{\kpar}{2\pi}\ktilde_2\right)^2 + \lambda^2V_\sharp(\bx) \ \right)\psi\ =\ E\psi,\ \ 
 \psi\in H^2(\Sigma).
\label{evp-kpar} \end{equation}
Let 
\begin{equation} E= E_0^\lambda + \Omega\  \label{Edef}\end{equation}
Recall the centered edge-Hamiltonian:
\begin{equation}
H^\lambda_\sharp(\kpar)= -\left(\nabla+i\frac{\kpar}{2\pi}\ktilde_2\right)^2 + \lambda^2V_\sharp(\bx)-E^\lambda_0\ ;
\label{Hshkpar}\end{equation}
see also  \eqref{Hsharp}.
Then, the eigenvalue problem may be rewritten as:
\begin{equation}
  \left(\ H^\lambda_\sharp(\kpar)\ -\ \Omega\ \right)\psi\ =\ 0,\quad \psi\in H^2(\Sigma)\ .
  \label{evp-1}
  \end{equation}

By Proposition \ref{Xdecomp} any $\psi\in H^2(\Sigma)$ may be written in the form:
\begin{equation}  \label{psi-decomp1}
  \psi\ =\ \sum_{I=A,B}\ \sum_{n\ge0} \alpha_n^I\ p_{_{\kpar,I}}^\lambda[n]\ +\ \widetilde{\psi},
  \end{equation}
  where $\alpha=\{(\alpha_n^A,\alpha_n^B)^{^\top}\}_{n\ge0}\in l^2(\N_0;\C^2)$ and 
  $\Pi_{_{AB}}(\kpar)\widetilde\psi=\widetilde\psi$. We adopt the convention 
  \[ \alpha_n^I=0,\ \ n\le-1,\ \ I=\bA,\bB.\]
  
  Substitution of \eqref{psi-decomp1} into \eqref{evp-1} yields:
  \begin{equation}
  \sum_{I=A,B}\ \sum_{n\ge0} \alpha_n^I\ \left(\ H^\lambda_\sharp(\kpar)\ -\ \Omega\ \right)\ p_{_{\kpar,I}}^\lambda[n]
   + \left(\ H_\sharp(\kpar)\ -\ \Omega\ \right)\ \widetilde{\psi}\ =\ 0\ .
   \label{evp-pre2}
   \end{equation}
   
By part (1) of Proposition \ref{Xdecomp}, the eigenvalue problem \eqref{evp-1} is seen to be equivalent to the system obtained by:\ 
 (i) applying the orthogonal projection $\Pi_{_{AB}}(\kpar)$ to  \eqref{evp-pre2}:
\begin{align}
&\Pi_{_{AB}}(\kpar)\left(\ H^\lambda_\sharp(\kpar)\ -\ \Omega\ \right)\ \widetilde{\psi}\ +\ 
\sum_{I=A,B}\ \sum_{n\ge0} \alpha_n^I\ \Pi_{_{AB}}(\kpar)\ \left(\ H^\lambda_\sharp(\kpar)\ -\ \Omega\ \right)\ p_{_{\kpar,I}}^\lambda[n] =\ 0\label{AB-perp}
\end{align}
and
  (ii) taking the inner product of \eqref{evp-pre2} with the states: $p_{_{\kpar,J}}^\lambda[m];\ m\ge0,\ J=A,B$:
\begin{align}
&\left\langle p_{_{\kpar,J}}^\lambda[m], \sum_{I=A,B}\ \sum_{n\ge0} \alpha_n^I\  \left(\ H^\lambda_\sharp(\kpar)\ -\ \Omega\ \right)\ p_{_{\kpar,I}}^\lambda[n]
\right\rangle\  +\ \left\langle  \left(\ H^\lambda_\sharp(\kpar)\ -\ \overline\Omega\ \right)p_{_{\kpar,I}}^\lambda[m],\widetilde\psi\right\rangle\ =\ 0 \label{AB}
  \end{align}
  where $m=0,1,2,\dots$. 
  
  Using Proposition \ref{resolvent} we solve \eqref{AB-perp} for $\widetilde{\psi}$ as a function of $\alpha=(\alpha^A,\alpha^B)^\top\in l^2(\N_0;\C^2)$:
  \begin{align}
\widetilde{\psi}\ =\ 
-\sum_{I=A,B}\ \sum_{n\ge0} \alpha_n^I\ \Ressp^{\lambda}(\Omega,\kpar)\ \Pi_{_{AB}}(\kpar)\  \cHeg(\kpar)\ p_{_{\kpar,I}}^\lambda[n] \ .
\label{psi-t}  \end{align}
Here we have used that $\Pi_{_{AB}}(\kpar)\ p_{_{\kpar,I}}^\lambda[n]\ =0$.
  Substitution of \eqref{psi-t} into \eqref{AB} yields
  \begin{equation} \sum_{I=A, B}\sum_{n\ge0}\mathcal{M}^{\lambda,\kpar}_{JI}[m,n](\Omega,\kpar)\ \alpha_n^I\ =\ 0;\ \ J=A, B,\ \ m\ge0\ ,
  \label{eig-M}\end{equation}
  where 
\begin{align}
&\mathcal{M}_{JI}^{\lambda}[m,n](\Omega,\kpar)\nn\\
 &\qquad \equiv\ 
\left\langle p_{_{\kpar,J}}^\lambda[m],   \left(\ H^\lambda_\sharp(\kpar)\ -\ \Omega\ \right)\ p_{_{\kpar,I}}^\lambda[n]
\right\rangle_{_{L^2(\Sigma)}}\nn\\
&\qquad\ -\  \left\langle   \cHeg(\kpar)\  p_{_{\kpar,J}}^\lambda[m]\ ,\ \Pi_{_{AB}}(\kpar)\  \Ressp^{\lambda}(\Omega,\kpar)\ \Pi_{_{AB}}(\kpar)\  \cHeg(\kpar)\ p_{_{\kpar,I}}^\lambda[n]
\right\rangle_{_{L^2(\Sigma)}}\ .\nn\\
&\label{cM-def}
\end{align}

\begin{remark}\label{interp}
 For fixed $J=A$ or $B$ and fixed $m\ge0$, the equation \eqref{eig-M} expresses the interaction of all atomic $A-$ and $B-$ sites within the cylinder, $\Sigma$, with the  atomic site $J$  in cell $m$. In particular, 
the  $\mathcal{M}_{JA}[m,n]$ are interaction coefficients between site $J$ in $\Omega_m$ and all sites $\bv_A^n$, $n\ge0$, and $\mathcal{M}_{JB}[m,n]$ are interaction coefficients between site $J$ in cell $\Omega_m$ and all sites $\bv_B^n$, $n\ge0$.
 \end{remark}
 
Due to their dependence on the Hamilitonian, $\cHeg$, we refer to the first term on the right in \eqref{cM-def} as the linear matrix elements, $\mathcal{M}^{^{\lambda,{\rm lin}}}[m,n](\Omega,\kpar)$ and second term on the right in  \eqref{cM-def} as the non-linear matrix elements, $\mathcal{M}^{^{\lambda,{\rm nl}}}[m,n](\Omega,\kpar)$. Thus, 
 \begin{equation}
 \mathcal{M}^{\lambda}[m,n](\Omega,\kpar)\ \equiv\ 
 \mathcal{M}^{^{\lambda,{\rm lin}}}[m,n](\Omega,\kpar)\ -\ \mathcal{M}^{^{\lambda,{\rm nl}}}[m,n](\Omega,\kpar)\ .
\label{cM-def1} \end{equation}
 
   In the subsequent sections we compute highly accurate approximations to  the linear (Section \ref{lin-els}) and non-linear (Section \ref{nl-els}) matrix elements. This will enable us to recast and solve \eqref{eig-M}
as a perturbation of a tight-binding model for $\lambda$ sufficiently large (Section \ref{zz-sb}).

\section{Matrix elements $\mathcal{M}_{_{JI}}^{^{\lambda,{\rm lin}}}[m,n](\Omega,\kpar)$ and $\mathcal{M}_{_{JI}}^{^{\lambda,{\rm nl}}}[m,n](\Omega,\kpar)$}
\label{lin-els}

In this section we provide expansions of the matrix entries of 
$\mathcal{M}_{_{JI}}^{^{\lambda,{\rm lin}}}[m,n](\Omega,\kpar)$. 
 Recall that 
 \begin{align}\label{Pkpar1}
P_{_{\kpar,\bI}}^\lambda[n](\bx)\ &\equiv\ e^{ i\frac{\kpar}{2\pi}\ktilde_2\cdot(\bx-\bv_\bI) }\ p_{_{\kpar,\bI}}^\lambda[n](\bx)
\ =\ \sum_{m_2\in\Z}\ e^{i\kpar m_2}\
p_{0}^\lambda(\bx-\bv^n_\bI-m_2\vtilde_2)\ ;
\end{align} 
(see also \eqref{Pkpar})
and that
$\cHeg\ =\ -\Delta + \lambda^2 V_\sharp(\bx)\ -\ E_0^\lambda$\ .

In preparation for our expansions, introduce the {\it nearest-neighbor hopping coefficient}:
 \begin{equation}
  \rho_\lambda\ =\ \int_{B_{_{r_0}}(0)}\ p_{0}^\lambda(\by) \lambda^2\ |V_0(\by)|\ p_0^\lambda(\by+\be)\ d\by 
 \ =\ \int_{\R^2}\ p_{0}^\lambda(\by) \lambda^2\ |V_0(\by)|\ p_0^\lambda(\by+\be)\ d\by ,
  \label{rho-def1}
  \end{equation}
  where $\be=\bv_B-\bv_A$. 
  The latter equality holds since $V_0$ has compact support in $B_{_{r_0}}(0)$. 
  We further recall the bounds \eqref{rho-boundsA} :
 \begin{equation}
  e^{-c{_{-}}\lambda}\ \lesssim\ \rho_\lambda\ \lesssim\  e^{-c_{_{+}}\lambda}
  \label{rho-bounds}
  \end{equation}
 for some constants $c_-, c_+>0$  and all $\lambda>0$ sufficiently large;
 this was proved in \cite{FLW-CPAM:17}.

The main results of this section (Propositions \ref{zz-els} and \ref{zz-nl-els}) are the following two propositions which (i) isolate the dominant (nearest neighbor) behavior of the linear matrix elements and provide estimates on the corrections, and (ii)
estimate the nonlinear  matrix elements. 

 \begin{proposition}[Expansion of linear matrix elements]\label{zz-els} 
 {\ }
 
For all $\lambda>\lambda_\star$ (sufficiently large), and all $\kpar\in[0,2\pi]$, we have:
 \begin{enumerate}
 \item  For $m\ge0$, 
 {\footnotesize{
\begin{align}
\label{ip-mmBA}
\left\langle P_{_{\kpar,\bB}}^{^{\lambda}}[m],\cHeg\ P_{_{\kpar,\bA}}^{^{\lambda}}[m]\right\rangle_{_{L^2(\Sigma)}}
\ =\ \left\langle p_{_{\kpar,\bB}}^{^{\lambda}}[m],\cHeg(\kpar) p_{_{\kpar,\bA}}^{^{\lambda}}[m]\right\rangle_{_{L^2(\Sigma)}}\ =\ -\rho_\lambda\ \left(1+e^{i\kpar} \right)\ +\ \mathcal{O}(e^{-c\lambda}\ \rho_\lambda)\ ,\\
\left\langle P_{_{\kpar,\bA}}^{^{\lambda}}[m],\cHeg\  P_{_{\kpar,\bB}}^{^{\lambda}}[m]\right\rangle_{_{L^2(\Sigma)}}
\ =\ \overline{\left\langle P_{_{\kpar,\bB}}^{^{\lambda}}[m],\cHeg\ P_{_{\kpar,\bA}}^{^{\lambda}}[m]\right\rangle_{_{L^2(\Sigma)}}}\ =\ -\rho_\lambda\ \left(1+e^{-i\kpar} \right)\ +\ \mathcal{O}(e^{-c\lambda}\ \rho_\lambda)\ .
\label{ip-mmAB}
\end{align}
}}
 \item For $m\ge0$, 
 \begin{align}\label{ip-mm+1}
\left\langle P_{_{\kpar,\bB}}^{^{\lambda}}[m],\cHeg\ P_{_{\kpar,\bA}}^{^{\lambda}}[m+1]\right\rangle_{_{L^2(\Sigma)}}
\ &=\ -\rho_\lambda\ +\ \mathcal{O}(e^{-c\lambda}\ \rho_\lambda),
\end{align}
and for $ m\ge1$
\begin{align}
\left\langle P_{_{\kpar,\bA}}^{^{\lambda}}[m],\cHeg\  P_{_{\kpar,\bB}}^{^{\lambda}}[m-1]\right\rangle_{_{L^2(\Sigma)}}
\ &\ =\ -\rho_\lambda\ +\ \mathcal{O}(e^{-c\lambda}\ \rho_\lambda)\ .
\label{ip-mm-1}
\end{align}
  \item
 \begin{align}
\left\langle P_{_{\kpar,\bB}}^{^{\lambda}}[m],\cHeg\ P_{_{\kpar,\bA}}^{^{\lambda}}[n]\right\rangle_{_{L^2(\Sigma)}}
=\ \mathcal{O}\left(\  e^{-c\lambda |m-n|}\ \rho_\lambda\ \right),\ m,n\ge0,\ \ n\ne m,m+1\ ,\label{ip-mn_BA}\\
\left\langle P_{_{\kpar,\bA}}^{^{\lambda}}[m],\cHeg\ P_{_{\kpar,\bB}}^{^{\lambda}}[n]\right\rangle_{_{L^2(\Sigma)}}
=\ \mathcal{O}\left(\ e^{-c\lambda |m-n|}\ \rho_\lambda  \right),\ m,n\ge0,\ \ n\ne m,m-1\ .\label{ip-mn_AB}
\end{align}

\item 
For $m, n\ge0$ and $\bI=A$ or $B$ \ \ 
 \begin{align}\label{ip-II}
\left\langle P_{_{\kpar,\bI}}^{^{\lambda}}[m],\cHeg\ P_{_{\kpar,\bI}}^{^{\lambda}}[n]\right\rangle_{_{L^2(\Sigma)}}
=\ \mathcal{O}\left(\ \ e^{-c\lambda}\ e^{-c\lambda|m-n|}\ \rho_\lambda\ \right)\ .
\end{align}
 \end{enumerate}
 The implied constants in the $\mathcal{O}(\cdot)$ estimates and the constants $\lambda_\star$ and $c$ are independent of $\kpar$.
\end{proposition}
\bigskip

We note, by part (4) of Proposition \ref{resolvent},  that the function
\[\Omega\mapsto\
 \left\langle\ \cHeg(\kpar) p^\lambda_{\kpar,J}[n]\ ,\ \Pi^\lambda_{_{AB}}(\kpar)\ \Ressp^{\lambda}(\Omega,\kpar)\ 
 \Pi^\lambda_{_{AB}}(\kpar)\ \cHeg(\kpar) p^\lambda_{\kpar,I}[m]\ \right\rangle_{L^2(\Sigma)}
 \]
 is analytic for $|\Omega|<c^\prime$.

\begin{proposition}[Estimation of nonlinear matrix element contributions]\label{zz-nl-els} 
There exists $\lambda>\lambda_\star$ (sufficiently large), such that for all $\kpar\in[0,2\pi]$ and $|\Omega|\le e^{-c^\prime\lambda}$ ($c^\prime$, a  sufficiently small constant determined by $V_0$) the following holds for $j=0,1$:
 \begin{align}
  &\Big|\ \left\langle\  \cHeg(\kpar) p^\lambda_{\kpar,J}[n]\ ,\ \Pi^\lambda_{_{AB}}(\kpar)\ \partial_\Omega^j \Ressp^{\lambda}(\Omega,\kpar)\ \Pi^\lambda_{_{AB}}(\kpar)\ \cHeg(\kpar) p^\lambda_{\kpar,I}[m]\ \right\rangle_{L^2(\Sigma)}\ \Big|\nn\\  &\qquad\qquad\qquad\qquad \lesssim\  \rho_\lambda\ e^{-c\lambda}\ e^{-c |n-m|}\ .  \label{nl-est1} \end{align}
 The implied constants in the $\mathcal{O}(\cdot)$ estimates and the constants $\lambda_\star$ and $c$ are independent of $\kpar$.
\end{proposition}
%
%

 Proposition \ref{zz-els} is proved in Section \ref{m-els-est} and  Proposition \ref{zz-nl-els} in Section \ref{nl-els}. The proof of  Proposition \ref{zz-nl-els}  requires  detailed information on the resolvent, which we need to control in weighted spaces. We obtain this control by constructing the  resolvent kernel and obtaining pointwise bounds for it. 
 The construction is carried out in Section \ref{res-kernel}.

\section{ Existence of zigzag edge states in the strong binding regime}\label{zz-sb}

In this section we apply Propositions \ref{zz-els} and \ref{zz-nl-els} to rewrite the edge state eigenvalue problem as a perturbation of the eigenvalue problem for the  tight-binding limiting operator studied in Section \ref{TB}.
We then use this reformulation to construct zigzag edge states for arbitrary $\lambda>\lambda_\star$, where $\lambda_\star$ is fixed and  sufficiently large.

Recall our reduction, for $\kpar\in\mathscr{J}\subset\subset(2\pi/3,4\pi/3)$, of the edge state eigenvalue problem for $\cHeg(\kpar)$ to the discrete eigenvalue problem for $\{(\alpha^A_m,\alpha^B_m)\}_{m\ge0}$ in $l^2(\N_0;\C^2)$:

\begin{equation} \sum_{I=A, B}\sum_{n\ge0}\mathcal{M}^\lambda_{JI}[m,n](\Omega,\kpar)\ \alpha_n^I\ =\ 0;\ \ J=A, B,\ \ m\ge0\ ,
  \label{eig-M1}\end{equation}
  Let's cast \eqref{eig-M1} in a form in which the tight-binding operator $\HTBs(\kpar)$ is made explicit. First, \eqref{eig-M1} is equivalent to the following system for $m\ge0$:
  \begin{align}
 \sum_{n\ge0} \mathcal{M}^\lambda_{AA}[m,n](\Omega,\kpar)\ \alpha_n^A\ +\ 
  \sum_{n\ge0} \mathcal{M}^\lambda_{AB}[m,n](\Omega,\kpar)\ \alpha_n^B\ =\ 0\ ,\nn\\
   \sum_{n\ge0} \mathcal{M}^\lambda_{BA}[m,n](\Omega,\kpar)\ \alpha_n^A\ +\ 
   \sum_{n\ge0} \mathcal{M}^\lambda_{BB}[m,n](\Omega,\kpar)\ \alpha_n^B\ =\ 0\ .
   \label{AB-exact}
    \end{align}
       %
%
To isolate the dominant terms (see Propositions \ref{zz-els} and \ref{zz-nl-els}), we rearrange the expressions and obtain for $m\ge0$:
\begin{align}
&\mathcal{M}^\lambda_{AB}[m,m-1](\Omega,\kpar)\ \alpha_{m-1}^B\ +\ 
\mathcal{M}^\lambda_{AB}[m,m](\Omega,\kpar)\ \alpha_m^B\ +\ \mathcal{M}^\lambda_{AA}[m,m](\Omega,\kpar)\ \alpha_m^A\nn\\ 
&\qquad\qquad =\  -\sum_{\substack{n\ge0\\ n\ne m, m-1}} \mathcal{M}^\lambda_{AB}[m,n](\Omega,\kpar)\ \alpha_n^B\  - \sum_{\substack{n\ge0\\ n\ne m}} \mathcal{M}^\lambda_{AA}[m,n](\Omega,\kpar)\ \alpha_n^A\nn \\
&  \mathcal{M}^\lambda_{BA}[m,m](\Omega,\kpar)\ \alpha_m^A\ +\ 
\mathcal{M}^\lambda_{BA}[m,m+1](\Omega,\kpar)\ \alpha_{m+1}^A\ +\ \mathcal{M}^\lambda_{BB}[m,m](\Omega,\kpar)\ \alpha_m^B\nn\\ 
&\qquad\qquad  =\ - \sum_{\substack{n\ge0\\ n\ne m, m+1}}  \mathcal{M}^\lambda_{BA}[m,n](\Omega,\kpar)\ \alpha_n^A\   -\ \sum_{\substack{n\ge0\\ n\ne m}} \mathcal{M}^\lambda_{BB}[m,n](\Omega,\kpar)\ \alpha_n^B\ .
  \label{AB-exact1}  \end{align}
Here, $\mathcal{M}^\lambda_{JI}[m,n]$ is given by \eqref{cM-def}, where we take  $\mathcal{M}^\lambda_{BA}[m,m-1]=0$ for $m=0$.
   The system \eqref{AB-exact1} is equivalent to \eqref{eig-M1}.\medskip

Our next step will be to express the matrix elements on the left hand side of \eqref{AB-exact1}, using Proposition \ref{prop:pkpar}, Proposition \ref{zz-els} and Proposition \ref{zz-nl-els}. Since the leading order expressions
 are proportional to $\rho_\lambda$, it is natural to 
introduce  the rescaled energy: 
\begin{equation}
 \Omega\ \equiv\ \rho_\lambda\ \tOmega. \label{tOmega-def}\end{equation}
 Recall our general upper and lower bounds on $\rho_\lambda$: $e^{-c_-\lambda}\lesssim \rho_\lambda\lesssim e^{-c_+\lambda}$ 
 (see \eqref{rho-bounds} or \eqref{rho-boundsA}) and let $\hat{c}>c_->0$  denote the positive constant introduced in Remark \ref{hatc}. We now constrain $\Omega$ to satisfy
 $|\Omega|< e^{-\hat{c}\lambda}$.  Then, $|\tOmega|=| \rho_\lambda^{-1}\Omega|\le e^{-(\hat{c}-c_-)\lambda}<e^{-c^{\prime\prime}\lambda}$,
 where $c^{\prime\prime}$ is a small positive constant, for any finite $\lambda$ sufficiently large.
 
 Using Proposition \ref{prop:pkpar}, Proposition \ref{zz-els} and Proposition \ref{zz-nl-els} in \eqref{AB-exact1}  we  obtain after dividing by $-\rho_\lambda$:
     \begin{align}
  &   \left(\ -1+\mathcal{O}(e^{-c\lambda})\ \right)\ \alpha_{m-1}^B\ +\  \left(\ -(1+e^{-i\kpar})+\mathcal{O}(e^{-c\lambda})\ \right)\ \alpha_m^B\ +\ 
   \left(\ -1+\mathcal{O}(e^{-c\lambda})\ \right)\  \tOmega\   \alpha_m^A\nn\\
 &\qquad =\  \sum_{\substack{n\ge0\\ n\ne m,m-1}}\ \mathcal{O}(e^{-c\lambda}\ e^{-c|m-n|})\ \alpha_n^B\ +\ 
 \sum_{n\ge0}\ \mathcal{O}(e^{-c\lambda}\ e^{-c|m-n|})\ \alpha_n^A,\nn \\
 &\label{ab1}\\
 &\qquad\qquad \textrm{where $\alpha_{m-1}^B=0$ for $m=0$, and }\nn\\
 &\nn\\
  &   \left(\ -(1+e^{i\kpar})+\mathcal{O}(e^{-c\lambda})\ \right)\ \alpha_{m}^A\ +\  \left(\ -1+\mathcal{O}(e^{-c\lambda})\ \right)\ \alpha_{m+1}^A\ +\  \left(\ -1+\mathcal{O}(e^{-c\lambda})\ \right)\ \tOmega\  \alpha_m^B\nn\\ 
  &\qquad =\ \sum_{\substack{ n\ge0 \\ n\ne m,m+1}}\ \mathcal{O}(e^{-c\lambda}\ e^{-c|m-n|})\ \alpha_n^A\ +\ 
 \sum_{n\ge0}\ \mathcal{O}(e^{-c\lambda}\ e^{-c|m-n|})\ \alpha_n^B\ ,
 &\nn\\
 \label{ab2}    \end{align}
 where $|\tOmega|< c^{\prime\prime}$.
 \begin{remark}\label{Om-an}      
By Proposition \ref{resolvent} (part 4) and Proposition \ref{zz-nl-els}, the expressions  in \eqref{ab1}-\eqref{ab2} of the form $\mathcal{O}(g(\lambda))$ are analytic functions of $\tOmega$ for $\tOmega$ varying in the open subset of $\C$: $|\Omega|<e^{-\hat{c}\lambda}$. Moreover, these expressions are all uniformly bounded by $g(\lambda)$ for all $\tOmega$ such that  $|\tOmega|< c^{\prime\prime}$, a small positive constant.
\end{remark}
We obtain, for $m\ge0$ and $|\tOmega|<c^{\prime\prime}$:
      \begin{align}
  &    -\alpha_{m-1}^B\ -\  (1+e^{-i\kpar})\ \alpha_m^B\ -\     \tOmega\   \alpha_m^A \nn\\
 &\qquad =\  \sum_{n\ge0}\ \mathcal{O}(e^{-c\lambda}\ e^{-c|m-n|})\ \alpha_n^B\ +\ 
 \sum_{n\ge0}\ \mathcal{O}(e^{-c\lambda}\ e^{-c|m-n|})\ \alpha_n^A \ ,
 \label{TB-pertA}\\
 &\qquad\qquad \textrm{where $\alpha_{m-1}^B=0$ for $m=0$, and }\nn\\
 &\nn\\
  &   -(1+e^{i\kpar})\ \alpha_{m}^A\ -\  \alpha_{m+1}^A\  -\   \tOmega\  \alpha_m^B \nn \\ 
  &\qquad =\ 
 \sum_{n\ge0}\ \mathcal{O}(e^{-c\lambda}\ e^{-c|m-n|})\ \alpha_n^B\
 +\  \sum_{n\ge0}\ \mathcal{O}(e^{-c\lambda}\ e^{-c|m-n|})\ \alpha_n^A  \ . \label{TB-pertB}  
   \end{align}
   Again we remark, as in Remark \ref{Om-an}, that in \eqref{TB-pertA}-\eqref{TB-pertB} expressions of the form $\mathcal{O}(g(\lambda))$ are analytic in $\tOmega$ and uniformly bounded by $g(\lambda)$ for $|\tOmega|< c^{\prime\prime}$.

The system \eqref{TB-pertA}-\eqref{TB-pertB} is  of the form:
  \begin{equation}
    \left[\ \left(\ H^{^{\rm TB}}_\sharp(\kpar)\ -\ \tOmega\ \right) \begin{pmatrix}\alpha^A\\ \alpha^B\end{pmatrix}\right]_m\ =\ \left[\ \mathscr{P}(\lambda;\rho_\lambda \widetilde\Omega)\ \begin{pmatrix}\alpha^A\\ \alpha^B\end{pmatrix}\ \right]_m,\qquad \textrm{ for \ \ $m\ge0$,}
      \label{TB-pert1}
  \end{equation}
where $H^{^{\rm TB}}_\sharp(\kpar)$
is the tight binding Hamiltonian  for a zigzag termination of $\mathbb{H}$, studied in Section \ref{TB}; see, in particular, \eqref{zz-ham}, \eqref{HTBegk}
\footnote{ Actually, the operator which emerges in \eqref{TB-pertA}-\eqref{TB-pertB} is 
 $-H^{^{\rm TB}}_\sharp(\kpar)$, minus one times the operator studied in Section \ref{TB}. However, since $\sigma_2H^{^{\rm TB}}_\sharp(\kpar)\sigma_2=-H^{^{\rm TB}}_\sharp(\kpar)$, the spectrum of $H^{^{\rm TB}}_\sharp(\kpar)$ is symmetric about zero energy
and $-H^{^{\rm TB}}_\sharp(\kpar)-z{\rm Id}$ has the same invertibility properties of $H^{^{\rm TB}}_\sharp(\kpar)-z{\rm Id}$.  Hence, in this and the following section we take $H^{^{\rm TB}}_\sharp(\kpar)$ to denote the negative of the operator
studied in  Section \ref{TB}.}
Note that our definition of $H^{^{\rm TB}}_\sharp(\kpar)$ implies that the scaled spectral parameter, $\tOmega$,  appears with a plus rather than a minus sign in \eqref{TB-pert1}.

Furthermore, using that the mapping $\{\gamma_m\}_{m\ge0}\mapsto\ \Big\{ \sum_{n\ge0}e^{-c|m-n|}\ \gamma_n\Big\}_{m\ge0}$ is bounded on $l^2(\N_0)$,  we have that the mapping $\tOmega\mapsto \mathscr{P}(\lambda;\rho_\lambda\tOmega)$ is an analytic mapping for $|\tOmega|<c^{\prime\prime}$ with values in the space of bounded linear operators on  $l^2(\N_0;\C^2)$. We also have, for all 
$|\tOmega|\le c^\prime$,\ $(c^\prime<c^{\prime\prime})$:
\begin{equation}
  \|\ \mathscr{P}(\lambda;\rho_\lambda\tOmega)\ \|_{l^2\to l^2}\lesssim e^{-c\lambda}\ ,
  \label{Plam-bd}\end{equation}
where the implied constant is independent of $\tOmega$, but depends on $c^\prime$.
Recall that $\kpar$ varies in a compact subinterval of $(2\pi/3,4\pi/3)$, where $\delta_{\rm gap}(\kpar)=\Big|1-|\zeta(\kpar)|\Big|=\Big|1-|1+e^{i\kpar}|\Big|>0$.
 We will further restrict $\tOmega$ to satisfy $|\tOmega|<c^\prime<\delta_{\rm gap}(\kpar)$.

Our goal is to construct, for all $\lambda$ sufficiently large, a solution of \eqref{TB-pert1}:
  \begin{align}
   \lambda\ &\ \mapsto\vec\alpha(\lambda)\ =\ (\ \alpha^A(\lambda),\alpha^B(\lambda)\ )\in l^2(\N_0;\C^2)\nn\\
   \lambda &\mapsto\ \tOmega(\lambda),\quad \textrm{such that}\ \ |\tOmega(\lambda)|\lesssim e^{-c\lambda}\le c^{\prime}\ .
\label{sol-XAB}  \end{align}
Given the mappings \eqref{sol-XAB}, equations \eqref{psi-decomp1}, \eqref{psi-t} and the relation $E=E_0^\lambda+\rho_\lambda\tOmega$ define a solution to the $L^2_\kpar(\Sigma)$ edge state eigenvalue problem, $\Psi_\kpar^\lambda(\bx)=e^{i\frac{\kpar}{2\pi}\ktilde_2\cdot\bx} \psi^\lambda_\kpar(\bx)$, where  
\begin{align}  \label{psi-lam}
  \psi^\lambda_\kpar(\bx)\ &=\ \sum_{I=A,B}\ \sum_{n\ge0} \alpha_n^I(\lambda)\ p_{_{\kpar,I}}^\lambda[n](\bx)\ +\ \widetilde{\psi}[\vec\alpha(\lambda)](\bx)\  ,\\
  E^\lambda(\kpar)\ &=\ E_0^\lambda+\rho_\lambda\tOmega(\lambda;\kpar) ,\nn
  \end{align}
  and the map $\vec\alpha\mapsto\widetilde{\psi}[\vec\alpha](\bx)$ is given in \eqref{psi-t}.
We shall succeed in this construction for $\kpar\in\mathscr{I}\subset\subset(2\pi/3,4\pi/3)$ and 
$\lambda>\lambda_\star(\mathscr{I})$ sufficiently large. 

The first step in this construction is to  note that as  $\lambda$ tends to infinity the system \eqref{TB-pertA}-\eqref{TB-pertB} formally reduces to the edge state eigenvalue problem
for the tight-binding Hamiltonian, $H_\sharp^{^{\rm TB}}$ (see  \eqref{zz-tbevp}, \eqref{zz-ham}) given by:
 \begin{align}
  &  \alpha_{m-1}^B\ +\ (1+e^{-i\kpar})\ \alpha_m^B\ -\ \tOmega\   \alpha_m^A
      \ =\ 0,\ \ m\ge0\nn \\
  &    (1+e^{i\kpar})\ \alpha_{m}^A\ +\ \alpha_{m+1}^A\ -\ \tOmega\   \alpha_m^B\ 
= 0 \ ,\  \ m\ge0\ ,\ {\rm with}\ \   \alpha_{-1}^B=0\ .\ \label{TB-pert0}  \end{align}
 By Theorem \ref{zz-spec}, if $\kpar\in(2\pi/3,4\pi/3)$ the system
 \eqref{TB-pert0} has an isolated and simple eigenvalue at $\widetilde{\Omega}^{_{\rm TB}}=0$ with corresponding vector 
  $\vec\alpha^{^{\rm TB}}=\{\ \alpha^{^{\rm TB}}_{_{m}}\ \}_{m\ge0}\in l^2(\N_0;\C^2)$ given by: 
\begin{equation}
\alpha^{^{\rm TB}}_{_{m}}\ =\ \begin{pmatrix} \alpha^{^{{\rm TB},A}} \\ \alpha^{^{{\rm TB},B}}\end{pmatrix}_{m}\ =\ \gamma_\star\begin{pmatrix} (-1)^{m}\left(1+e^{i\kpar}\right)^{m} \\ 0\end{pmatrix},\ \ \textrm{for}\ \ m\ge0\ ,
\label{zz-estate2}\end{equation}
where we take $\gamma_\star=\sqrt{1-|\zeta(\kpar)|^2}\ne0$ so that $\vec\alpha^{^{\rm TB}}$ has $l^2(\N_0;\C^2)-$ norm equal to one. 

To prove that \eqref{TB-pert1} has a solution in $l^2(\N_0,\C^2)$ which for $\lambda$ large is approximately equal to $\vec\alpha^{^{\rm TB}}$, we seek a solution of   \eqref{TB-pert1}  of the form:
\begin{align}
\vec\alpha(\lambda) &=\ \vec\alpha^{^{\rm TB}}\ +\ \vec\beta(\lambda)\ =\  \begin{pmatrix}\alpha^{^{\rm TB,A}}\\ \alpha^{^{\rm TB,B}}\end{pmatrix}\ +\ 
 \begin{pmatrix}\beta^A(\lambda)\\ \beta^B(\lambda)\end{pmatrix},\nn\\
  \tOmega &=\ \tOmega(\lambda)\ ,\quad\textrm{where we take}\ \left\langle\vec\alpha^{^{\rm TB}},\vec\beta\right\rangle_{l^2(\N_0;\C^2)}=0.
 \label{beta-def}\end{align}
 
 Introduce the orthogonal projection $\Pi_{_0}^{^{\rm TB}}:l^2(\N_0;\C^2)\to \Big(\ {\rm span}\Big\{\vec\alpha^{^{\rm TB}}\Big\}\ \Big)^\perp$. Substituting \eqref{beta-def} into \eqref{TB-pert1}
  and projecting onto ${\rm span}\{\vec\alpha^{^{\rm TB}}\}$ and its orthogonal complement, we obtain
   the equivalent system for $\vec\beta$ and $\tOmega$:
\begin{align}
    &\left(H^{^{\rm TB}}_\sharp(\kpar)-\tOmega\right) \vec\beta\ =\  \Pi_{_0}^{^{\rm TB}}\ \mathscr{P}(\lambda;\rho_\lambda \widetilde\Omega)\ \vec\alpha^{^{\rm TB}}\ +\ \Pi_{_0}^{^{\rm TB}}\ \mathscr{P}(\lambda;\rho_\lambda \widetilde\Omega)\vec\beta\ ,
      \label{beta-eqn}\\
    &\tOmega\ +\ \left\langle \vec\alpha^{^{\rm TB}},\mathscr{P}(\lambda;\rho_\lambda \widetilde\Omega)\ \vec\alpha^{^{\rm TB}}\right\rangle\ +\ 
    \left\langle \vec\alpha^{^{\rm TB}},
         \mathscr{P}(\lambda;\rho_\lambda \widetilde\Omega)\ \vec\beta\right\rangle\ =\ 0.
    \label{tOm-eqn} 
  \end{align}
  Let $\mathscr{R}^{^{\rm TB}}(\tOmega;\kpar)$ denote the inverse of 
  $\Pi_{_0}^{^{\rm TB}}\left(H^{^{\rm TB}}_\sharp(\kpar)-\tOmega\right)\Pi_{_0}^{^{\rm TB}}$, which for $|\tOmega|<c^\prime$ is well-defined as a bounded operator on the $l^2(\N_0;\C^2)-$  orthogonal complement of ${\rm span}\{\vec\alpha^{^{\rm TB}}(\kpar)\}$. Moreover,  $\|\mathscr{R}^{^{\rm TB}}(\tOmega;\kpar)\|\lesssim1$ for $|\tOmega|<c^{\prime}<\delta(\kpar)$, by Theorem \ref{zz-spec}.  For $\lambda$ sufficiently large we may solve \eqref{beta-eqn} for $\vec\beta[\tOmega;\lambda]\in {\rm Range}\ \Pi_{_0}^{^{\rm TB}}$
  and obtain:
  \begin{align}
  \vec\beta[\tOmega;\lambda]\ &=\ \Big[\ I\ -\ 
  \mathscr{R}^{^{\rm TB}}(\tOmega;\kpar) \Pi_{_0}^{^{\rm TB}} \mathscr{P}(\lambda;\rho_\lambda \widetilde\Omega) \Big]^{-1}\ \Pi_{_0}^{^{\rm TB}} \mathscr{P}(\lambda;\rho_\lambda \widetilde\Omega)\
  \vec\alpha^{^{\rm TB}} \nn\\  
  &\equiv \mathscr{A}(\tOmega;\lambda)\ \Pi_{_0}^{^{\rm TB}} \mathscr{P}(\lambda;\rho_\lambda \widetilde\Omega)\ \vec\alpha^{^{\rm TB}}\ .
  \end{align}
  This follows by the bound $\| \mathscr{P}(\lambda;\rho_\lambda\tOmega)\ \|_{l^2\to l^2}\lesssim e^{-c\lambda}$; see \eqref{Plam-bd}. Therefore,  the construction of $\vec\beta(\lambda)$, $\tOmega(\lambda)$ (see \eqref{sol-XAB}) boils down to solving the following scalar nonlinear equation for $\tOmega$ as a function of $\lambda$:
  \begin{align}
  \tOmega\ +\ \left\langle \vec\alpha^{^{\rm TB}},\mathscr{P}(\lambda;\rho_\lambda \widetilde\Omega)\ \vec\alpha^{^{\rm TB}}\right\rangle\ +\ 
    \left\langle\ \vec\alpha^{^{\rm TB}},
           \mathscr{P}(\lambda;\rho_\lambda \widetilde\Omega)\ \mathscr{A}(\tOmega;\lambda)\ \Pi_{_0}^{^{\rm TB}} \mathscr{P}(\lambda;\rho_\lambda \widetilde\Omega)\ \vec\alpha^{^{\rm TB}}\right\rangle\ =\ 0.
    \label{tOm-eqn1} 
  \end{align}
  
  Using analyticity in $\tOmega$ and previous bounds, we may write \eqref{tOm-eqn1}  as
  \begin{equation}
  \tOmega\ +\ \left\langle \vec\alpha^{^{\rm TB}},\mathscr{P}(\lambda;0)\ \vec\alpha^{^{\rm TB}}\right\rangle\ +\
  \mathscr{G}(\tOmega;\lambda)\ =\ 0.
 \label{tOm-eqn2} \end{equation}
 Here, $  \mathscr{G}(\tOmega;\lambda)$ is analytic 
  with $|\partial_\tOmega^j \mathscr{G}(\tOmega;\lambda)|\lesssim e^{-c\lambda}$ ($j=1,2$) for all $\tOmega$ in the complex neighborhood of zero, $|\tOmega|<c^\prime$.
  Since $\Big|\left\langle \vec\alpha^{^{\rm TB}},\mathscr{P}(\lambda;0)\ \vec\alpha^{^{\rm TB}}\right\rangle\Big|\le \ e^{-c\lambda}$,
   for $\lambda$ sufficiently large, equation \eqref{tOm-eqn2} may be solved for $\tOmega(\lambda)$ by 
 using a contraction mapping argument on the disc: $|\tOmega|\le 2C e^{-c\lambda}$.
  Therefore, modulo Propositions \ref{zz-els}  and  \ref{zz-nl-els} which are proved in Sections \ref{res-kernel}, \ref{m-els-est} and \ref{nl-els}, we have proved our main result,  Theorem \ref{main-thm1}.

  \section{Resolvent convergence; proof of Theorem \ref{res-conv}}\label{sc-res}

We study the scaled resolvent: 
\[
\Big(\ \rho_\lambda^{-1}\ H^\lambda_\sharp\ -\ z {\rm Id}\ \Big)^{-1}
\ =\ \Big(\ \rho_\lambda^{-1}\left(\ -\Delta+V_\sharp-E_0^\lambda\ \right)\ -\ z {\rm Id}\ \Big)^{-1}\]
as an operator on $L^2(\R^2)$.
We consider the scaled non-homogeneous equation
\begin{equation}
\left(\ \rho_\lambda^{-1}H^\lambda_\sharp(\kpar) - z{\rm Id}\ \right)\psi\ =\ \varphi,\qquad \varphi\in L^2(\Sigma).
\label{inhom}\end{equation}
or equivalently
\begin{equation}
\left(\ H^\lambda_\sharp(\kpar) - \rho_\lambda z{\rm Id}\ \right)\psi\ =\ \rho_\lambda\varphi,\qquad \varphi\in L^2(\Sigma).
\label{inhom-a}\end{equation}
We express $\varphi$ as:
\begin{equation}
\varphi\ =\ \sum_{J=A,B}\ \sum_{n\ge0} \beta_n^J\ p_{_{\kpar,J}}^\lambda[n]\ +\ \widetilde\varphi\ ,\qquad 
\Pi_{_{AB}}(\kpar)\widetilde\varphi=\widetilde\varphi
\label{vphi}\end{equation}
and seek a solution of \eqref{inhom} in the form
  \begin{equation}  \label{psi-dcmp1}
  \psi\ =\ \sum_{I=A,B}\ \sum_{n\ge0} \alpha_n^I\ p_{_{\kpar,I}}^\lambda[n]\ +\ \widetilde{\psi},\qquad
 \Pi_{_{AB}}(\kpar)\widetilde\psi=\widetilde\psi.
  \end{equation}
  where $\alpha=\{(\alpha_n^A,\alpha_n^B)^{^\top}\}_{n\ge0}\in l^2(\N_0;\C^2)$ and
   $ \widetilde\psi=\Pi_{_{AB}}(\kpar)\widetilde\psi\in\mathscr{X}^\lambda_{AB}(\kpar)$.

Substitution of \eqref{vphi} and  \eqref{psi-dcmp1} into \eqref{inhom-a} and projecting the resulting 
equation with $\Pi_{_{AB}}(\kpar)$ and $I-\Pi_{_{AB}}(\kpar)$ (whose range is 
 ${\rm span}\{p^\lambda_{_{\kpar,I}}[n]:I=A,B,\ n\ge0\}$),  yields the coupled system for  
 $\alpha=\{\alpha_{n}^I:n\ge0, I=A,B\}\in l^2(\N_0;\C^2)$ and 
 $\widetilde\psi\in\mathscr{X}^\lambda_{AB}(\kpar)$:
 \begin{align}
& \Pi_{_{AB}}(\kpar)\left(\ H_\sharp^\lambda(\kpar)-\rho_\lambda z{\rm Id}\ \right)\widetilde\psi
 \ =\ -\ \sum_{I, n}\  \alpha_n^I\ \Pi_{_{AB}}(\kpar)\ H^\lambda_\sharp(\kpar)\ p_{_{\kpar,I}}^\lambda[n]\ +\    
  \rho_\lambda\widetilde\varphi \label{psit-eqn}\\
& \sum_{I,n} \left\langle p_{_{\kpar,J}}^\lambda[m], 
 \Big( H^\lambda_\sharp(\kpar)-\rho_\lambda z\ {\rm Id}\Big)
 p_{_{\kpar,I}}^\lambda[n] \right\rangle \alpha_n^I\ +\  
 \left\langle H^\lambda_\sharp(\kpar)p_{_{\kpar,J}}^\lambda[m], \widetilde\psi\right\rangle   \label{alph-eq}\\
 &\qquad =\ \rho_\lambda\sum_{I,n}  \left\langle  p_{_{\kpar,J}}^\lambda[m],  p_{_{\kpar,I}}^\lambda[n] \right\rangle \beta_n^I,\ \textrm{for}\ J=A, B\ {\rm and}\ m\ge0,\nn
 \end{align}
 where the sums $\sum_{I,n}$ are over $I=A, B$ and $n\ge0$. 
 
 We next use Proposition \ref{resolvent} to solve \eqref{psit-eqn} for $\widetilde\psi\in L^2(\Sigma)$ and obtain:
 \begin{equation}
 \widetilde\psi\ =\ - \sum_{I,n} \alpha_n^I\ \mathcal{K}_\sharp^\lambda(\rho_\lambda z,\kpar)H_\sharp^\lambda(\kpar)p_{_{\kpar,I}}^\lambda[n]\ +\ \rho_\lambda\ \mathcal{K}_\sharp^\lambda(\rho_\lambda z,\kpar)\widetilde\varphi\ .
 \label{tpsi-in}
 \end{equation}
 Substitution of the expression in \eqref{tpsi-in} for $\widetilde\psi$ into the left hand side of \eqref{alph-eq} yields the closed non-homogeneous system for $\alpha\in l^2(\N_0;\C^2)$:
 \begin{align}
 &\sum_{I,n}\ \mathcal{M}^\lambda_{JI}[m,n] \alpha_n^I\ =\ 
 \rho_\lambda\Big[\ \sum_{I,n}  \left\langle  p_{_{\kpar,J}}^\lambda[m],  p_{_{\kpar,I}}^\lambda[n]\right\rangle \beta_n^I\ 
  -\ \left\langle H_\sharp^\lambda(\kpar)p_{_{\kpar,J}}^\lambda[m],\mathcal{K}_\sharp^\lambda(\rho_\lambda z,\kpar)\widetilde\varphi\right\rangle\ \Big],
  \label{inhom-alph}
 \end{align}
for each $J=A, B$ and $m\ge0$. The matrix elements $\mathcal{M}^\lambda_{JI}[m,n] $ are displayed in \eqref{cM-def}. As in our study of the edge state eigenvalue problem (Section \ref{zz-sb}) we expand the 
 $\mathcal{M}^\lambda_{JI}[m,n]$ using Proposition \ref{zz-els} and obtain the following system, which is equivalent to \eqref{inhom-alph}
 \footnote{As in Section \ref{zz-sb} (see the footnote after \eqref{TB-pert1}), based
 on the observation $\sigma_2 H^{^{\rm TB}}_\sharp(\kpar)\sigma_2=-H^{^{\rm TB}}_\sharp(\kpar)$)  we let $ H^{^{\rm TB}}_\sharp(\kpar)$ denote 
 the negative of the operator studied in Section \ref{TB}.}
 :
 \begin{align}
 &   \left[\ \left(\ H^{^{\rm TB}}_\sharp(\kpar)\ -\ z\ {\rm Id}\ -\ \mathscr{P}(\lambda;\rho_\lambda z)\ \right) \begin{pmatrix}\alpha^A\\ \alpha^B\end{pmatrix}\right]_m\nn\\
    &\quad  =\ 
-\  \begin{pmatrix}
 \sum_{I,n} \left\langle  p_{_{\kpar,A}}^\lambda[m],  p_{_{\kpar,I}}^\lambda[n]\  \beta_n^I\  \right\rangle \\ \\
  \sum_{I,n} \left\langle  p_{_{\kpar,B}}^\lambda[m],  p_{_{\kpar,I}}^\lambda[n]\  \beta_n^I\ \right\rangle
  \end{pmatrix}
  +\ 
  \begin{pmatrix}
  \left\langle H_\sharp^\lambda(\kpar)p_{_{\kpar,A}}^\lambda[m],\mathcal{K}_\sharp^\lambda(\rho_\lambda z,\kpar)\widetilde\varphi\right\rangle\\ \\
    \left\langle H_\sharp^\lambda(\kpar)p_{_{\kpar,B}}^\lambda[m],\mathcal{K}_\sharp^\lambda(\rho_\lambda z,\kpar)\widetilde\varphi\right\rangle
    \end{pmatrix}
    ,\quad m\ge0.\nn\\
 &     \label{TB-inhom}
  \end{align}
  Recalling the bound $\|\mathscr{P}(\lambda;\rho_\lambda\tOmega)\ \|_{l^2\to l^2}\lesssim e^{-c\lambda}$  (see \eqref{Plam-bd}), Proposition \ref{prop:pkpar} and Proposition \ref{Hp_ap0} we solve for $\alpha^\lambda$ and find
  \begin{align}
  \alpha^\lambda\ &=\ \left(\ H^{^{\rm TB}}_\sharp(\kpar)\ -\ z\ {\rm Id}\ \right)^{-1}\beta\ +\ \alpha_1^\lambda, \quad \textrm{where}\nn\\
 \| \alpha_1^\lambda \|_{_{l^2(\N_0;\C^2)}}\ &\lesssim\ e^{-c\lambda}\Big(\ \|\beta\|_{_{l^2(\N_0;\C^2)}}\ +\ 
 \|\Pi_{_{AB}}(\kpar)\varphi\|_{_{L^2(\Sigma)}} \ \Big)
 \end{align}
We therefore have that $\psi=\left(\ \rho_\lambda^{-1}H^\lambda_\sharp(\kpar) - z{\rm Id}\ \right)^{-1}\varphi\in L^2(\Sigma)$ is given by:
\begin{align}
\left(\ \rho_\lambda^{-1}H^\lambda_\sharp(\kpar) - z{\rm Id}\ \right)^{-1}\varphi\ &=\ \sum_{I,n}\ \Big[\ \left(\ H^{^{\rm TB}}_\sharp(\kpar)\ -\ z\ {\rm Id}\ \right)^{-1}\beta\ +\ \alpha_1^\lambda\ \Big] p_{_{\kpar,I}}^\lambda[n]\nn\\
&\qquad  + \mathcal{O}_{L^2(\Sigma)}\left(\ e^{-c\lambda} \|\beta\|_{_{l^2(\N_0;\C^2)}} + e^{-c\lambda}\|\Pi_{_{AB}}(\kpar)\varphi\|_{_{L^2(\Sigma)}} \right)
\end{align}

Introduce $H^\lambda_{\sharp,\kpar}$, the restriction of $H^\lambda_\sharp$, to the space $H^2_\kpar$.
 Since $H^\lambda_\sharp$ commutes with $\bx\mapsto \bx+\vtilde_2$ it follows that $H^\lambda_{\sharp,\kpar}$ maps the space $H^2_\kpar$ into $L^2_\kpar$.  Let $\mathscr{P}_{AB,\kpar}$ denote the projection of $L^2_\kpar$ onto the orthogonal complement of the subspace of $L^2_\kpar$ spanned by the states:
 $P_{_{\kpar,I}}^\lambda[n]=e^{i\frac{\kpar}{2\pi}(\bx-\bv_I)}p_{_{\kpar,I}}^\lambda[n]\in L^2_\kpar$, where $I=A, B$ and $n\ge0$; see  \eqref{Pkpar}.
Therefore, for any $F\in L^2_\kpar$:
\begin{align}
\left(\ \rho_\lambda^{-1}H^\lambda_{\sharp,\kpar} - z{\rm Id}\ \right)^{-1}F\ &=\ \sum_{I,n}\ \Big[\ 
\left(\ H^{^{\rm TB}}_\sharp(\kpar)\ -\ z\ {\rm Id}\ \right)^{-1}\beta\ +\ \alpha_1^\lambda\ \Big] P_{_{\kpar,I}}^\lambda[n]\nn\\
&\qquad  + \mathcal{O}_{L^2(\Sigma)}\left(\ e^{-c\lambda} \|\beta\|_{_{l^2(\N_0;\C^2)}} + e^{-c\lambda}\|\mathscr{P}_{AB,\kpar}F\|_{_{L^2_\kpar}} \right).
\label{sc-res1}
\end{align}

Any $F\in L^2_\kpar$ has the  representation $F=\sum_{I,n}\alpha_n^I[F] P_{_{\kpar,I}}^\lambda[n]+F_\perp$,  where $\{\alpha_n^I[F]\}_{I,n}\in l^2(\N_0;\C^2)$ and $F_\perp\in \textrm{Range}(\mathscr{P}_{AB,\kpar})$. Define the map $J_\kpar:L^2_\kpar\to l^2(\N_0;\C^2)\oplus\textrm{Range}(\mathscr{P}_{AB,\kpar})$ by:
\begin{align}
J_\kpar:\ F\ \mapsto\ \begin{pmatrix} \{\alpha_n^I[F]\} \\ F_\perp\end{pmatrix}= 
\begin{pmatrix}\Big\{\left\langle P_{_{\kpar,I}}^\lambda[n],F\right\rangle\ +\ \mathcal{O}(e^{-c\lambda}\|F\|_{L^2_\kpar})\Big\}\ \\ F_\perp\end{pmatrix}
\end{align}

We therefore have from \eqref{sc-res1} that 
\[
\left(\ \rho_\lambda^{-1}H^\lambda_{\sharp,\kpar} - z{\rm Id}\ \right)^{-1}\ -\ 
 J_\kpar^*\ 
 \begin{pmatrix} 
 \left(\ H^{^{\rm TB}}_\sharp(\kpar)\ -\ z\ {\rm Id}\ \right)^{-1} & 0\\
 0 & 0\end{pmatrix}\ J_\kpar\ =\ \mathcal{O}_{_{L^2_\kpar\to L^2_\kpar}}(e^{-c\lambda}).
 \]

This completes the proof of Theorem \ref{res-conv}.

\section{The resolvent kernel and weighted resolvent bounds}\label{res-kernel} {\ }\medskip

It remains for us to prove Propositions \ref{zz-els} and \ref{zz-nl-els} on the expansion and estimation of matrix elements. 
 The proof of Proposition \ref{zz-els} concerning the linear matrix elements  uses the energy estimates on the resolvent obtained in Section \ref{en-est}.

 To prove Proposition  \ref{zz-nl-els} we require exponentially weighted estimates, which we obtain by constructing the resolvent kernel  and obtaining pointwise bounds on it.  We carry this out in the present section.  In Section \ref{m-els-est} we then give the proof of 
Proposition \ref{zz-els} and in Section \ref{nl-els} we prove Proposition \ref{zz-nl-els}.
\bigskip

In Section \ref{en-est} we  obtained energy estimates for $\Ressp^\lambda(\Omega,\kpar)$,  the inverse of 
{\footnotesize{
\[ \Pi_{_{AB}}(\kpar)\left(\ \cHeg(\kpar)-\Omega\ \right)\Pi_{_{AB}}(\kpar)\ =\  
 \Pi_{AB}\ \left[\ -\left(\nabla_\bx+i\frac{\kpar}{2\pi}\ktilde_2\right)^2\ +\ \lambda^2 V_\sharp(\bx)\ -\ E_0^\lambda-\Omega\ \right]\ \Pi_{AB},\]
 }}
 defined
as a bounded operator from $\mathscr{X}_{AB}(\kpar)$ to  $\mathscr{X}_{AB}(\kpar)\cap H^2(\Sigma)$; see Proposition \ref{resolvent}, 
which holds for all $|\Omega|<c^\prime$, where $c^\prime$ is a sufficiently small positive constant. We may extend $\Ressp^\lambda(\Omega,\kpar)$ to an operator acting on all of $L^2(\Sigma)$, not just $\mathscr{X}_{AB}(\kpar)$, by composing it with $\Pi_{_{AB}}(\kpar)$,
 {\it i.e.} we require $\Ressp^\lambda(\Omega,\kpar)\psi=0$ if $\Pi_{_{AB}}(\kpar)\psi=0$.
\medskip

In this section we shall prove, under the more stringent restriction on $\Omega$: $|\Omega|\le e^{-c\lambda}$ for some $c>0$ and $\lambda\gg1$, that this operator derives from a kernel 
 $\mathcal{K}_\sharp^\lambda(\bx,\by;\Omega,\kpar)$. Specifically, we have
 \bigskip
 
\begin{theorem}\label{Ksharp} There exist constants $\lambda_\star, c>0$ such that for $\lambda\ge\lambda_\star$, $|\Omega|\le e^{-c\lambda}$ and
for each $\kpar\in[0,2\pi]$ the following holds for the operator $\mathcal{K}_\sharp^\lambda(\Omega,\kpar)$, which is bounded on $L^2(\Sigma)$:\medskip

\begin{enumerate}
 \item $\mathcal{K}_\sharp^\lambda(\Omega,\kpar)$ is arises from an integral kernel
 $\mathcal{K}_\sharp^\lambda(\bx,\by;\Omega,\kpar)$:
 
\begin{equation}\mathcal{K}_\sharp^\lambda(\Omega,\kpar)[f](\bx)\ =\ \Pi_{AB}^\lambda\ \Ressp(\Omega,\kpar)\ \Pi_{AB}^\lambda[f](\bx)\ =\ \int_{\Omega_\Sigma}\ \mathcal{K}_\sharp^\lambda(\bx,\by;\Omega,\kpar)\ f(\by)\ d\by\ .
\label{Klam-rep}\end{equation}

\item The integral kernel
 $\mathcal{K}_\sharp^\lambda(\bx,\by;\Omega,\kpar)$ satisfies the following bound: there exist positive constants $R, C_1, C_2$, independent  of $\kpar$ and $\Omega$,  such that for all $\bx, \by\in\R^2$: 
 
 \begin{equation}
 \left|\ \mathcal{K}_\sharp^\lambda(\bx,\by;\Omega,\kpar)\ \right|\ \le\ C_1\left[\ \lambda^4\ +\ 
 \Big|\ \log|\bx-\by|\ \Big|\ \right]\ {\bf 1}_{_{|\bx-\by|\le R}}
  \ +\ C_2\ e^{-c\lambda}\ e^{-c\lambda|\bx-\by|}\ .
 \label{Ksharp-est} \end{equation}

\end{enumerate}
\end{theorem}

Theorem \ref{Ksharp} is at the heart of the proof of Proposition \ref{zz-nl-els}, which provides bounds on the nonlinear  matrix elements of $\mathcal{M}^\lambda(\Omega,\kpar)$. The remainder of this section is devoted to the proof of Theorem \ref{Ksharp}.
The construction and estimation $\mathcal{K}_\sharp^\lambda$ is based on a strategy, in which we piece together localized
 atomic Green's functions with appropriate corrections.

\subsection{The free Green's function and bounds on the atomic ground state}\label{atom-bnds}

Denote by $\Gfree(\bx)$ the fundamental solution of $-\Delta - E_0^\lambda$:
\begin{equation}
\left(-\Delta_\bx - E_0^\lambda\right)\ \Gfree(\bx)\ =\ \delta(\bx) ,
\label{Gfree-def}\end{equation}
where $\delta(\bx)$ is the Dirac delta function. Here, $E_0^\lambda$ denotes the ground state of $\Hatom=-\Delta+\lambda^2V_0$;
 see hypothesis (GS), \eqref{GS}.
 Note that  $\Gfree(\bx)=\ G^{\rm free}\left(\sqrt{|E_0^\lambda|}\ \bx\right)$,
where  $G^{\rm free}(\bx)$  satisfies
$\left(\ -\Delta_\bx+1\ \right)\ G^{\rm free}(\bx)=\delta(\bx),\ \bx\in\R^2$.
 $G^{\rm free}(\bx)=K_0(|\bx|)$ is the modified Bessel function of order zero, which decays to zero exponentially  as $|\bx|\to\infty$ \cites{WW:02}.
  The following lemma summarizes important standard properties of $\Gfree(\bx)$; see \cite{Simon:76,FLW-CPAM:17}
 \begin{lemma}\label{K0-prop} For $\bx\in\R^2$,
\begin{enumerate}
\item $G^{\rm free}(\bx)=G^{\rm free}(|\bx|)$ is positive and strictly decreasing for $|\bx|\ge0$.
\item There exist entire functions $f$ and $g$ and constants $C_1, c_2$,  such that
\begin{equation} G^{\rm free}(\bx)\ =\ f(|\bx|) \log |\bx|\ +\ g(|\bx|)\ ,\label{fg-Gfree}\end{equation}
where $f(0)=-1/2\pi$ and $|\partial_s^jf(s)|, |\partial_s^jg(s)|\le C_1 e^{-c_2s}$, for $j=0,1$ and all $s\in[0,\infty)$.
\item $G^{\rm free}(\bx)\ \lesssim |\bx|^{-\frac12} e^{-|\bx|}$ for $|\bx|$ large. 
  \end{enumerate}
  \end{lemma}
The bounds on $f(s)$ and $g(s)$ are proved, for the case $j=0$, in \cite{Simon:76}. This proof can be extended to a derivation of the bounds for $j=1$. Alternatively, these bounds may be deduced
directly  from the integral representation for $G^{\rm free}(\bx)$ used in the proof of Lemma 15.3 of \cite{FLW-CPAM:17}. 

 We shall apply the following consequence of  Lemma \ref{K0-prop} and  \eqref{GS}:\medskip

\nit There exist  $c, c^\prime>0$, and for each $R>0$, additional constants $C_R, C_R^\prime>0$,  such that 
\begin{align}
0<G^{\rm free}_\lambda(\bx)\ &=\ G^{\rm free}\left(\sqrt{|E_0^\lambda|}\ \bx\right)\ \le\ C_R\ e^{-c\lambda|\bx|}
\left(\ \Big|\log (\lambda |\bx|)\Big|\ {\bf 1}_{_{\{\lambda|\bx|\le R\}}}\ +\ 1\right)\ ,\ \bx\in\R^2.
\label{Gf-allx}\\
|\nabla_\bx G^{\rm free}_\lambda(\bx)|\ &\le\ C_R^\prime e^{-c^\prime\lambda|\bx|}
\left(\ \frac{1}{\lambda|\bx|}\ {\bf 1}_{_{\{\lambda|\bx|\le R\}}}\ +\ 1\ \right)\label{DGf-allx}
 \end{align}

\subsection{The atomic Green's function}\label{G-atomic}

{\ }\medskip

In this section we establish bounds (integral and then pointwise) on the Green's function 
associated with $\Hatom-E_0^\lambda=-\Delta + \lambda^2V_0(\bx) -E_0^\lambda$.
 Since $\Hatom$ has a one dimensional kernel spanned by $p_0^\lambda(\bx)$,
  and a spectral gap (see \eqref{EG}), the operator $\Hatom-E_0^\lambda$ is invertible on 
the orthogonal complement of $\textrm{span}\{p_0^\lambda\}$.

We denote by $\Gatom(\bx,\by)$ the associated Green's kernel, which solves
\begin{align}
&\left(\ -\Delta_\bx + \lambda^2V_0(\bx) -E_0^\lambda\ \right)\ \Gatom(\bx,\by)\ =\ 
\delta(\bx-\by)\ -\ p_0^\lambda(\bx)p_0^\lambda(\by)\label{Gat-def}
\end{align}
 and which satisfies
 \begin{align}
&\int_{\R^2} \Gatom(\bx,\by)\ p_0^\lambda(\by)\ d\by\ =\ 0,\ \ \textrm{for all}\ \bx\in\R^2,\label{Gat-orth}\\
& \Gatom(\bx,\by)\ =\ \Gatom(\by,\bx)\ \ \textrm{for}\ \bx,\ \by\ \in\ \R^2 \ \textrm{with}\ \ \bx\ne\by\ .\label{Gat-sym}
\end{align}
For fixed $\bx$, the function $\by\mapsto\Gatom(\bx,\by)$ belongs to $ L^2(\R^2_\by)$, and we have for any
 $f\in L^2(\R^2)$ that the function
\begin{equation}
u(\bx)\ =\ \int_{\R^2}\  \Gatom(\bx,\by)\ f(\by)\ d\by
\label{Latuf-sol}\end{equation}
solves
\begin{align}
\left(\ -\Delta + \lambda^2V_0(\bx) -E_0^\lambda\ \right)\ u(\bx)\ &=\ 
 f(\bx)\ -\ \left\langle p_0^\lambda,f\right\rangle_{L^2(\R)}\ p_0^\lambda(\bx),
\label{Latuf}\\
\left\langle p_0^\lambda\ ,\ u\ \right\rangle_{L^2(\R^2)}\ &=\ 0 \ .\label{orthog}
\end{align}

\subsubsection{$L^2$ bounds on $\bx\mapsto \Gatom(\bx,\by)$ and 
$\by\mapsto \Gatom(\bx,\by)$}

{\ }\medskip

By the spectral gap hypothesis on $\Hatom$, \eqref{EG}, we have that $u$ satisfies the bound:
\begin{equation}
\|u\|_{L^2(\R^2)}\ \le\ C\ \|f\|_{L^2(\R^2)}.
\label{ubound}
\end{equation}

We may next obtain pointwise bounds on $u(\bx)$ in terms of $\|f\|_{L^2(\R^2)}$. In particular, we claim that 
\begin{equation}
|u(\bx)|\ \le\ C\ \lambda^2\ \|f\|_{L^2(\R^2)}\ .
\label{u_ptbd}
\end{equation}
We prove this as follows:
\begin{align}
|u(\bx)| &\le\ C\left(\ \|\Delta u\|_{L^2(B_1(\bx))}\ +\ \|u\|_{L^2(B_1(\bx))} \right)\nn\\
&\le\ C\left(\ \Big\|(E_0^\lambda-\lambda^2V_0)u\ +\ f\ -\ \left\langle p_0^\lambda, f\right\rangle p_0^\lambda\ \Big\|_{L^2(B_1(\bx))}\ +\ \|u\|_{L^2(B_1(\bx))}\ \right)\nn\\
&\le\ C\ \lambda^2\ \|f\|_{L^2(\R^2)}\nn
\end{align}
which implies the bound \eqref{u_ptbd}. 

Therefore, by \eqref{Latuf-sol}, for all $f\in L^2(\R^2)$:
\begin{equation}
\Big|\ \int_{\R^2}\ \Gatom(\bx,\by)\ f(\by)\ d\by\ \Big|\ \le\ C\ \lambda^2\ \|f\|_{L^2(\R^2)}\ .
\label{u_ptbd1}\end{equation}

Consequently, 
 \begin{equation}
 \left(\ \int_{\R^2}\ |\Gatom(\bx,\by)|^2\ d\by\ \right)^{1\over2}\ \le\ C\lambda^2,\ \ \bx\in\R^2
 \label{Gat-L2y}
 \end{equation}
 and by symmetry of $\Gatom$
 \begin{equation}
 \left(\ \int_{\R^2}\ |\Gatom(\bx,\by)|^2\ d\bx\ \right)^{1\over2}\ \le\ C\lambda^2,\ \ \by\in\R^2\ .
 \label{Gat-L2x}
 \end{equation}
 
 We now use these $L^2$ bounds on $\Gatom(\bx,\by)$ to obtain pointwise bounds.
 
\subsubsection{Pointwise bounds on $\Gatom(\bx,\by)$}\ 
Recall that $\supp V_0\subset B_{r_0}(0)$.
\begin{theorem}[Pointwise bounds on $\Gatom(\bx,\by)$]\label{ptGat}
\begin{enumerate}
\item For all $R>0$,  there exist $\lambda_0=\lambda_0(R)$ and positive constants $c$, $C_R$ and $D_R$ such that for all $\lambda>\lambda_0$:
\begin{align}
\Big|\ \Gatom(\bx,\by)\ -\ c_0\ \log|\bx-\by|\ \Big|\ &\le\ C_R\ \lambda^4\ \ {\rm for}\ \ |\bx-\by|\ \le\ R,\quad \label{leR}
\end{align}
where $c_0=-(2\pi)^{-1}$.
\item There exist $R>10r_0$ and positive constants $\lambda^\prime$, $C$ and $c$, which depend on $R$ but not on $\lambda$, such that for all $\lambda>\lambda' (R)$: 
\begin{align}
|\Gatom(\bx,\by)|\ &\le\ C\ e^{-c\lambda}\ e^{-c\lambda|\bx-\by|},\ \ |\bx-\by|\ \ge\ R\ .\label{geR}
\end{align}
\item 
Choose $r_j$, $j=1,2,3$, such that $r_0<r_1<r_2<r_3<\frac{1}{10}R$. Assume $\by\in B_{r_1}(0)$ and $\bx\notin B_{r_3}(0)$. Then,
\begin{equation}\label{variant}
\Big|\ \Gatom(\bx,\by)\ \Big|\ \lesssim\ e^{-c\lambda}\ e^{-c\lambda|\bx-\by|}\ ,
\end{equation}
where the implied constants depend on $r_0, r_1$, $r_2$ and $r_3$. 
\end{enumerate}
\end{theorem}
\nit{\bf Proof of bound \eqref{leR}:}\\ Fix $\by\in\R^2$. By \eqref{Gat-def} we have
\begin{align}
-\Delta_\bx \Gatom(\bx,\by)\ &=\ \delta(\bx-\by)\ -\ p_0^\lambda(\bx)\ p_0^\lambda(\by)\ 
+\ \left(\ E_0^\lambda\ -\ \lambda^2 V_0(\bx)\ \right)\Gatom(\bx,\by)\nn\\
&\ =\ -\Delta_\bx\ c_0 \log|\bx-\by| \ -\ p_0^\lambda(\bx)\ p_0^\lambda(\by)\ +\ \left(\ E_0^\lambda\ -\ \lambda^2 V_0(\bx)\ \right)\Gatom(\bx,\by)\ .\end{align}
  Hence,
\begin{equation}
-\Delta_\bx\ \left[\  \Gatom(\bx,\by)\ -\ c_0 \log|\bx-\by|\ \right]
\ =\ -\ p_0^\lambda(\bx)\ p_0^\lambda(\by)\ +\ \left(\ E_0^\lambda\ -\ \lambda^2 V_0(\bx)\ \right)\Gatom(\bx,\by)\ .
\end{equation}
Therefore, using that $|f(\bx)|\lesssim \|\Delta f(\bz)\|_{L^2(B_1(\bx);d\bz)}+\|f(\bz)\|_{L^2(B_1(\bx);d\bz)}$  we have for arbitrary fixed $\by\in\R^2$ and all $\bx\in\R^2$ satisfying
 $|\bx-\by|\le R$:
\begin{align}
&\left|\ \Gatom(\bx,\by)\ -\ c_0 \log|\bx-\by|\ \right|\nn\\
&\qquad  \le\ 
 \left\| -\ p_0^\lambda(\bz)\ p_0^\lambda(\by)\ +\ \left(\ E_0^\lambda\ -\ \lambda^2 V_0(\bz)\ \right)\Gatom(\bz,\by) \right\|_{L^2(B_1(\bx);d\bz)} \nn\\
 &\qquad\qquad   +\ \left\| \Gatom(\bz,\by)\ -\ c_0 \log|\bz-\by|\  \right\|_{L^2(B_1(\bx);d\bz)}\ .
  \label{diff-est}
\end{align}
To continue this bound, we use  that
\begin{align}
&|p_0^\lambda(\by)|\ \lesssim\ \lambda\ (\textrm{see \eqref{p0-bd-cpam})},\ \ \|p_0^\lambda\|_{L^2}=1,\ \ |E_0^\lambda-\lambda^2 V_0(\bz)|\lesssim\lambda^2,\nn\\
& \|\Gatom(\bz,\by)\|_{L^2(B_1(\bx);d\bz)}
\ \ \lesssim\lambda^2\ \ \textrm{and}\ \  \left\|\ \log|\bz-\by|\  \right\|_{L^2(B_1(\bx);d\bz)}\le C^\prime_R\ .
\label{4bounds}\end{align}
The bounds \eqref{4bounds} follow since $|E_0^\lambda|\lesssim \lambda^2$ (since $\|V_0\|_\infty<\infty$) and by \eqref{p0-bd-cpam} and \eqref{Gat-L2x}.
 We obtain for any $R>0$ that there exists $C_R<\infty$ such that 
\begin{align}
&\left|\ \Gatom(\bx,\by)\ -\ c_0 \log|\bx-\by|\ \right|\ \le\ C_R\ \lambda^4,\ \ \textrm{for all}\ |\bx-\by|\le R,\ \textrm{with}\ \bx\ne\by\ .
\label{Gat-diff}\end{align}
\medskip

\nit{\bf Proof of bound \eqref{geR}:} Recall that the support of $V_0$ is contained in $B_{r_0}(0)$.  Assume $|\bx-\by|>R$, and choose constants:
\begin{equation}
 r_0<r_1<r_2<r_3<\frac{1}{10}R\ .\label{r_js}
 \end{equation}
 Thus, we require $R>10r_0$.
Without any loss of generality, we assume $|\by|\le|\bx|$. Therefore, $R<|\bx-\by|\le |\bx|+|\by|\le2|\bx|$ and therefore
\begin{equation}\label{wlog}
|\bx|\ge\frac12 |\bx-\by|>\frac12 R>r_3.
\end{equation}
  Let  $\Thout=\Thout(\bx)$ denote a smooth function of $r=|\bx|$, defined for all $\bx\in\R^2$, such that $0\le\Thout(\bx)\le1$ and
\begin{equation}
\Thout(\bx)\ \equiv\ 
\begin{cases} 
1\ , & |\bx|\ge r_2\\
0\ , & |\bx|\le r_1
\end{cases}
\end{equation}
We note that  $\Thout\cdot V_0\equiv0$. 

Using the defining equation for $\Gatom$, \eqref{Gat-def}, we obtain:
\begin{align}
&\left(\ -\Delta_\bz\ -\ E_0^\lambda\ \right)\ \left[\ \Thout(\bz)\ \Gatom(\bz,\by)\ \right]\nn\\
&\qquad =\
\Thout(\bz)\ \left\{\ -p_0^\lambda(\bz)\ p_0^\lambda(\by)\ \right\}\ +\ \Thout(\bz)\cdot\delta(\bz-\by)\nn\\
&\qquad\quad +\ 2\nabla_\bz\Thout(\bz)\cdot\nabla_\bz \Gatom(\bz,\by)\ +\ \left(\ \Delta_\bz\Thout(\bz)\ \right)\ \Gatom(\bz,\by)\ .\nn\\
 &\label{yo}\end{align}
 We next use the Green's function $\Gfree$ (see \eqref{Gfree-def}) to represent $\Thout(\bx)\ \Gatom(\bx,\by)$. Multiplication of \eqref{yo} by $\Gfree(\bx-\bz)$ and integration with respect to $\bz$ yields
 \begin{align}
 \Thout(\bx)\ \Gatom(\bx,\by)\ &=\ \int_{\R^2}\Gfree(\bx-\bz)\ \left(\ -\Delta_\bz-E_0^\lambda\ \right)
 \left[\ \Thout(\bz)\ \Gatom(\bz,\by)\ \right]\ d\bz\nn\\
  &=\  \Thout(\by)\Gfree(\bx-\by)\ -\int_{\R^2}\Gfree(\bx-\bz)\ 
  \Thout(\bz)\ p_0^\lambda(\bz)\ d\bz\ p_0^\lambda(\by) \nn\\
   &\quad +\  2\  \int_{\R^2}\Gfree(\bx-\bz)\ \nabla_\bz\Thout(\bz)\cdot\nabla_\bz \Gatom(\bz,\by)\ d\bz\nn\\
&\quad  +\ \int_{\R^2}\Gfree(\bx-\bz)\  \left(\ \Delta_\bz\Thout(\bz)\ \right)\ \Gatom(\bz,\by)\ d\bz\ ,\nn\end{align}
which, since $\Thout(\bx)=1$ for $|\bx|>r_2$, we write as
\begin{align} \Gatom(\bx,\by)\ = \Thout(\by)\Gfree(\bx-\by)\ +\ {\rm Term}_1(\bx,\by)\ +\ {\rm Term}_2(\bx,\by)\ +\ {\rm Term}_3(\bx,\by)\ . \label{Term123}
  \end{align}
 Since $|\bx-\by|>R$, by  \eqref{Gf-allx}  we have $\left| \Thout(\by)\Gfree(\bx-\by)\right|\lesssim e^{-c\lambda|\bx-\by|}$.
  We next estimate the latter three terms in \eqref{Term123} individually. \medskip
 
 \nit{\it Bound on ${\rm Term}_1(\bx,\by)$\ of \eqref{Term123}:}\ Consider the integral
 \begin{equation}
 {\rm Term}_1(\bx,\by)\ \equiv -\int_{\R^2}\Gfree(\bx-\bz)\ 
  \Thout(\bz)\ p_0^\lambda(\bz)\ d\bz\  p_0^\lambda(\by)\ .\label{T1}\end{equation}
   Due to the factor of $\Thout(\bz)$ in the integrand of  \eqref{T1}, only $\bz$ such that  $|\bz|\ge r_1$.
 are relevant. On this set  we have $p_0^\lambda(\bz)\lesssim e^{-c_1\lambda}\ e^{-c\lambda|\bz|}$  by \eqref{p0-bd-cpam}, for some constants $c_1, c>0$. Furthermore, by \eqref{Gf-allx}, there exists $c^\prime>0$ such that 
$\Gfree(\bx-\bz)\lesssim e^{-c^\prime\lambda|\bx-\bz|}\ \left(\ \Big|\log\lambda|\bx-\bz|\Big|\ {\bf 1}_{_{\{|\bx-\bz|\le1\}}}+1\ \right)$. 
Therefore, for some constant $\tc$ (smaller than the minimum of $c_1,c, c^\prime$) we have
 \begin{align}
 \Big| {\rm Term}_1(\bx,\by)\ \Big|\ &\lesssim\  e^{-\tc\lambda}\  \int_{|\bz|\ge r_1}
e^{-\tc\lambda|\bx-\bz|}\ \left(\ \Big|\log\lambda|\bx-\bz|\Big|\ {\bf 1}_{_{\{|\bx-\bz|\le1\}}}+1\ \right)\
 e^{-\tc\lambda|\bz|}\   d\bz\ p_0^\lambda(\by)\nn\\
 &=\ 
e^{-\tc\lambda}\ \int_{|\bz|\ge r_1}
e^{-\frac{\tc}{2}\lambda\left(|\bx-\bz|+|\bz|\right) }\ e^{-\frac{\tc}{2}\lambda\left(|\bx-\bz|+|\bz|\right) } \ \left(\ \Big|\log\lambda|\bx-\bz|\Big|\ {\bf 1}_{_{\{|\bx-\bz|\le1\}}}+1\ \right)\   d\bz\ p_0^\lambda(\by)
\nn\\
&\le\ e^{-\tc\lambda}\ \ e^{-\frac{\tc}{2}\lambda |\bx|}\ \int_{|\bz|\ge r_1}
 e^{-\frac{\tc}{2}\lambda\left(|\bx-\bz|+|\bz|\right) } \ \left(\ \Big|\log\lambda|\bx-\bz|\Big|\ {\bf 1}_{_{\{|\bx-\bz|\le1\}}}+1\ \right)\   d\bz\ p_0^\lambda(\by)\nn\\
& \lesssim\ e^{-c\lambda}\ e^{-c\lambda |\bx|}\ p_0^\lambda(\by)\ . \nn\end{align}
For $|\by|<r_0+\delta_0$, with small $\delta_0>0$, we have $p_0^\lambda(\by)\lesssim\lambda$.
For such $\by$, $|\bx|=|\bx-\by+\by|\ge |\bx-\by|-r_0-\delta_0\ge\frac12|\bx-\by|+\frac{R}{2}-r_0-\delta_0\ge\frac12|\bx-\by|$.
Therefore, for $|\bx-\by|>R$ and $|\by|<r_0+\delta_0$ we have 
$ \Big| {\rm Term}_1(\bx,\by)\ \Big|\lesssim e^{-c\lambda}\ e^{-c\lambda |\bx|}\ p_0^\lambda(\by)\lesssim e^{-c\lambda}\ e^{-c\lambda |\bx|}\ \lambda \lesssim e^{-c^\prime\lambda}\  e^{-c^\prime \lambda|\bx-\by|}$. 

Therefore, for $|\by|\ge r_0+\delta_0$ and  $|\bx-\by|>R$,  we have 
$ \Big| {\rm Term}_1(\bx,\by)\ \Big|\lesssim e^{-c\lambda}\ e^{-c\lambda |\bx|}\ p_0^\lambda(\by)\lesssim
e^{-c\lambda}\ e^{-c\lambda (|\bx|+|\by|)}\lesssim e^{-c^\prime\lambda}\  e^{-c^\prime \lambda|\bx-\by|}$\ .
\bigskip 

\medskip

\nit{\it Bound on ${\rm Term}_2(\bx,\by)$\ of \eqref{Term123}:}\ We first note that  $\nabla_\bz\Thout(\bz)=0$ for $|\bz|>r_2$. Since $|\bx|>\frac12 R>r_2$, the 
 integrand of ${\rm Term}_2(\bx,\by)$\ is supported away from $\bz=\bx$. 
Integration by parts yields
 \begin{equation}
{\rm Term}_2(\bx,\by)\ =\   -2\  \int_{\R^2}\ \nabla_\bz\cdot \left[\ \Gfree(\bx-\bz)\ \nabla_\bz\Thout(\bz)\ \right]\  \Gatom(\bz,\by)\ d\bz\ .
\label{int-b-pts}
\end{equation}
We note this integration by parts can be justified even though there is a weak singularity of the integrand at $\bz=\by$, and we remark on this at the conclusion of the proof.
Bounding ${\rm Term}_2(\bx,\by)$ using the Cauchy-Schwarz inequality we obtain:
\[
\Big|{\rm Term}_2(\bx,\by)\Big|\ \le\   2\  \left(\ \int_{\R^2}\ \Big|\nabla_\bz\cdot \left[\ \Gfree(\bx-\bz)\ \nabla_\bz\Thout(\bz)\ \right]\Big|^2\ d\bz\ \right)^{\frac12}\cdot \left( \int_{\R^2}  \Big|\Gatom(\bz,\by)\Big|^2\ d\bz\ \right)^{\frac12} \ .
\]
The second factor is bounded by a constant times $\lambda^2$ thanks to the  $L^2$ bound on $\Gatom$ given in \eqref{Gat-L2x}.
To bound the first factor note, due to the properties of $\Thout(\bz)$, that the support of the integrand is contained in: $r_1\le |\bz|\le r_2$ and $|\bx|\ge r_3$. Therefore, $|\bx-\bz|\ge|\ |\bx|-|\bz|\ |\ge r_3-r_2>0$.
 Therefore, by \eqref{Gf-allx} and \eqref{DGf-allx}, for all $|\bx|\ge r_3$:
 \[ \Big|\ \nabla_\bz\cdot \left[\ \Gfree(\bx-\bz)\ \nabla_\bz\Thout(\bz)\ \right]\ \Big|\ \lesssim e^{-c\lambda|\bx-\bz|}\ {\bf 1}_{_{\{r_1\le|\bz|\le r_2\}}}\ \lesssim\ e^{-c^\prime\lambda}\ 
  e^{-c^\prime\lambda|\bx|}.
 \]
 It follows from\eqref{wlog} that
 \begin{align}
\Big|\ {\rm Term}_2(\bx,\by)\ \Big|\ &\lesssim \ 
e^{-c^\prime\lambda}\ 
  e^{-c^\prime\lambda|\bx|}\ \left(\ \int_{_{|\bz|\le r_2}} \ |\Gatom(\bz,\by)|^2\ d\bz\ \right)^{1\over2}\ \nn\\
   &\lesssim \ e^{-c^\prime\lambda}\ e^{-c^\prime\lambda|\bx|}\ \lambda^2\ 
   \lesssim
     e^{-c\lambda}\ e^{-c\lambda|\bx-\by|}\ .
  \end{align}
The bound on ${\rm Term}_3(\bx,\by)$ is obtained in a manner similar to the bound on ${\rm Term}_2(\bx,\by)$,
but there is no need to integrate by parts.  

We conclude the proof of \eqref{geR} by remarking on the technical point raised above concerning the integration by parts leading to \eqref{int-b-pts}. Recall that
\[\left(-\Delta_\bz+\lambda^2V_0(\bz)\right)\Gatom(\bz,\by)=\delta(\bz-\by)+E_0^\lambda\Gatom(\bz,\by)-p_0^\lambda(\bz)p_0^\lambda(\by).\]
For fixed $\by$, $\bz\mapsto p_0^\lambda(\bz)p_0^\lambda(\by)$ is $C^\infty$ by elliptic regularity because
 $(-\Delta_\bz+\lambda^2V_0(\bz)-E_0^\lambda)p_0^\lambda(\bz)=0$ and $V_0\in C^\infty$. Furthermore, 
 \[
 (-\Delta_\bz+\lambda^2V_0(\bz)-E_0^\lambda)\left[\Gatom(\bz,\by)-\Gfree(\bz-\by)\right]=
 -\lambda^2V_0(\bz)\Gfree(\bz,\by)-p_0^\lambda(\bz)p_0^\lambda(\by).\]
 Since $V_0\in C^\infty$, $\bz\mapsto\Gfree(\bz,\by)\in H^{1-\eps}(\R^2)$ ($\eps>0$ arbitrary), we have by elliptic regularity
  that $\Gatom(\bz,\by)-\Gfree(\bz-\by)\in H_{\rm loc}^{3-\eps}(\R^2)$. Furthermore by \eqref{fg-Gfree}, for fixed $\by$
  \[ \Gatom(\bz,\by)=c_0\log|\bz-\by|\ +\ j(\bz,\by)\ \  \textrm{for $\bz$ near $\by$},\]
  where $\bz\mapsto j(\bz,\by)\in H^{2-\eps}_{\rm loc}(\R^2)$. 
   This makes it easy to justify the integration by parts. For example, replace $\Gatom(\bz,\by)$ by 
   $\frac12 c_0\log\left[|\bz-\by|^2+\tau^2\right]+j(\bz,\by)$, integrate by parts and pass to the limit $\tau\to0^+$.
This concludes the proof of \eqref{geR}. Since the proof of the bound \eqref{variant} follows from a very similar argument, we omit it. This completes the proof of  Theorem \ref{ptGat}.

  \subsection{Kernels}\label{kernels}

  Our goal will be to construct the Green's kernel for a Hamiltonian $H_\Gamma^{^\lambda}=-\Delta +\VGamma(\bx)-E_0^\lambda$, with potential $\VGamma$ defined via superposition involving translates of the atomic potential, $V_0$, centered at the sites of a discrete set $\Gamma$. 
  The construction of this Green's function, $G_\lambda^\Gamma(\bx,\by)$ makes use of some technical tools developed in this section. 
  \medskip
  
We work with integral operators of the form
\begin{equation}
f\ \mapsto\ A_\lambda[f](\bx)\ \equiv\ \int_{\R^2}\ A_\lambda(\bx,\by)\ f(\by)\ d\by\ .
\label{TA}\end{equation}
We shall use the notation  $A_\lambda f$ and $A_\lambda[f]$ to denote such operators and occasionally omit the $\lambda$ dependence. 

\begin{definition}[Main Kernel]\label{mainK}
The function $A_\lambda(\bx,\by):\R^2\times\R^2\to\R$ is called a main kernel if there exist positive constants 
 $R, c, C_1, C_2$ and $\lambda_0$ such that for all $\bx, \by\in\R^2$ with $\bx\ne\by$ we have
 \begin{equation}
 |A_\lambda(\bx,\by)|\ \le\ C_1\left[\ \lambda^4\ +\ \Big|\ \log|\bx-\by|\ \Big|\ \right]\ {\bf 1}_{|\bx-\by|\le R}
  \ +\ C_2\ e^{-c\lambda}\ e^{-c\lambda|\bx-\by|}
 \label{mainK-est} \end{equation}
 for all $\lambda\ge\lambda_0$.
\end{definition}
By Theorem \ref{ptGat},  the atomic Green's function $\Gatom(\bx,\by)$ is a main kernel.\medskip

\begin{definition}[Error Kernel]\label{errorK}
The function $\mathcal{E}_\lambda(\bx,\by):\R^2\times\R^2\to\R$ is called a {\it error kernel} if there exist positive constants 
 $c, C$ and $\lambda_0$ such that for all $\bx, \by\in\R^2$  
 \begin{equation}
 |\mathcal{E}_\lambda(\bx,\by)|\ \le\ C\ e^{-c\lambda}\ e^{-c\lambda|\bx-\by|}
 \label{errorK-est} \end{equation}
 for all $\lambda\ge\lambda_0$.
\end{definition} 
\medskip

If $A$ and $B$ are operators with kernels given by $A(\bx,\by)$ and $B(\bx,\by)$, respectively, 
 then $AB$ is defined to be the operator with kernel $(AB)(\bx,\by)$ given by
 \begin{equation}
 (AB)(\bx,\by)\ \equiv\ \int_{\R^2}\ A(\bx,\bz)\ B(\bz,\by)\ d\bz
 \label{AB-kernel}
 \end{equation}
 
 \begin{remark} If $\mathcal{E}(\bx,\by)$ is an error kernel, then  $\lambda^p\ \mathcal{E}(\bx,\by)$ is an error kernel for any $p\ge0$. To see this,  replace the constant $c$ in \eqref{errorK-est} by a slightly smaller positive constant, $c^\prime$.
 \end{remark}

   \begin{lemma}\label{error-main}
   {\ }\\
   Let $K_\lambda$ arise from a main kernel and $\mathcal{E}_\lambda$ arise from an error kernel.
   \begin{enumerate}
 \item  Then, 
  \begin{equation}
  \widetilde{\mathcal{E}}_\lambda=\ I\ -\ (I-\mathcal{E}_\lambda)^{-1}\ =\ \sum_{l\ge1}\mathcal{E}_\lambda^l
 \label{e-sum} \end{equation}
arises from an error kernel.
  \item The operators  $\mathcal{E}_\lambda\ K_\lambda$ and $K_\lambda\ \mathcal{E}_\lambda$ arise from error kernels.
  \item  The operator $e^{-c\lambda}\ K_\lambda^2$, where $c>0$,  arises from an error kernel.
  \end{enumerate}
     \end{lemma}
     
      \nit The proof of Lemma \ref{error-main} is presented in Appendix \ref{app:error-main}

 \subsection{Green's kernel for a set of atoms centered on points of a discrete set, $\Gamma$}
 
  Let $\Gamma$ denote a discrete subset of $\R^2$, which we refer to as a set of {\it nuclei}.
  The set $\Gamma$ may be finite or infinite. We assume that 
  \begin{equation}
   \inf\{ |\bv-\bw|\ :\ \bv,\bw\in\Gamma,\ \bv\ne\bw\ \}\ge r_{\rm min}>2r_0.\label{r_min-def}
   \end{equation}
   At sites $\omega\in\Gamma$ we center identical {\it atoms} described by the atomic potential $V_0$:
  \begin{equation}
  \VGamma(\bx)\ =\ \sum_{\omega\in\Gamma}\ \lambda^2\ V_\omega(\bx),\ \ {\rm where}\ \ V_\omega(\bx)\equiv V_0(\bx-\omega)\ .
  \label{VGamma}
  \end{equation}

  \begin{example}
  Some choices of $\Gamma$ which are of interest to us are:
  \begin{enumerate}
  \item $\Gamma=\mathbb{H}=\Lambda_A\cup \Lambda_B$, the bulk honeycomb structure.
  \item $\Gamma=\Lambda_I$, $\bI=A,B$, the $A-$ and $B-$ sublattices.
  \item $\Gamma=\mathbb{H}_\sharp\ =\ \{\bv_I+n_1\vtilde_1+n_2\vtilde_2\ :\ n_1\ge0,\ n_2\in\Z\ \}$,
   the set of lattice points in a zigzag- terminated honeycomb structure.
   \end{enumerate}
  \end{example}

  Our goal will be to construct the Green's kernel $G_\Gamma^{^\lambda}(\bx,\by)$ 
  associated with the operator 
    \begin{equation}
\HGamma\ =\ -\Delta\ +\ \VGamma(\bx)\ -\ E_0^\lambda\ ,
\label{HGamma}\end{equation}
where $E_0^\lambda$ is the ground state energy of $\Hatom=-\Delta+\lambda^2V_0$; see \eqref{GS}.

Recall $\Gatom$ which satisfies
\begin{align}
&\left(\ -\Delta\ +\ \lambda^2V_0(\bx)\ -\ E_0^\lambda\ \right)\Gatom(\bx,\by)\ =\ \delta(\bx-\by)-p_0^\lambda(\bx)\ p_0^\lambda(\by),\nn\\
&\int_{\R^2}\ \Gatom(\bx,\by) p_0^\lambda(\bx)\ d\bx\ =\ 0, \nn\\
&\Gatom(\bx,\by)\ =\ \Gatom(\by,\bx)\ .\nn
\end{align}

Recalling $r_j, j=1,2,3$ specified in \eqref{r_js}, we further introduce $r_4$ such that 
\begin{equation}
0<r_0<r_1<r_2<r_3<r_4<\ \frac12\ r_{\rm min}\ ,\qquad (r_{\rm min}>2r_0),
\label{r_jsA}
\end{equation}
where $r_{\rm min}$ is a lower bound for the minimum distance between points in $\Gamma$; see \eqref{r_min-def}. 
Introduce the  smooth cutoff function $\Theta_0(\bx)$ satisfying:
\\ 
$0\le\Theta_0\le1$ on $\R^2$, $\Theta_0(\bx)=1$ for $\bx\in B_{r_3}(0)$,
 and  $\Theta_0(\bx)=0$ for $\bx\notin B_{r_4}(0)$. 

 For $\omega\in\Gamma$, define $\Theta_\omega(\bx)=\Theta_0(\bx-\omega)$. Finally, let 
 \begin{equation}
 \Tfree(\bx)\ \equiv\ 1\ -\ \sum_{\omega\in\Gamma}\ \Theta_\omega(\bx).
 \label{theta-free}
 \end{equation}
 Then,  $0\le\Tfree\le1$ on $\R^2$; $\Tfree$ is smooth and supported away from $\Gamma$. In particular for all $\omega\in\Gamma$,  $\Tfree=0$ in $B_{r_3}(\omega)$. 
 
 We write $p_\omega^\lambda(\bx)\equiv p_0^\lambda(\bx-\omega)$, where $p_0^\lambda(\bx)$ is the ground state of $\Hatom=-\Delta+\lambda^2 V_0(\bx)$. Thus, $p_\omega^\lambda(\bx)$ is the ground state of $-\Delta+\lambda^2 V_\omega(\bx)$. We also express the translated atomic
  Green's kernel as 
  \begin{equation}
   G_{\lambda,\omega}^{\rm atom}(\bx,\by)\ =\ \Gatom(\bx-\omega,\by-\omega)\ .
  \label{Gomega}
  \end{equation}
 
For any $f\in L^2(\R^2)$ we may write:
 \begin{equation}
  f(\bx)\ =\ \sum_{\omega\in\Gamma}\ \left(\Theta_\omega f\right)(\bx)\ +\ (\Tfree f)(\bx)\ ,
 \label{f-split}\end{equation}
 and for each $\omega\in\Gamma$, we have by \eqref{Gat-def}
 \begin{align}
 \Theta_\omega(\bx) f(\bx)\ &=\
  \left(\ -\Delta_\bx+\lambda^2V_\omega(\bx)-E_0^\lambda\ \right)\ 
  \int_{\R^2} \ G_{\lambda,\omega}^{\rm atom}(\bx,\by)\  \left(\ \Theta_\omega(\by) f(\by)\ \right) d\by\nn\\
 & \qquad +\ \left\langle p_\omega^\lambda,\Theta_\omega f\right\rangle_{L^2(\R^2)}\ p_\omega^\lambda(\bx)\ ,
 \label{Tof} \end{align}
  and by \eqref{Gat-orth}
  \begin{equation}
  \int_{\R^2}\ p_\omega^\lambda(\bx)\ \left[\
   \int_{\R^2}\ G_{\lambda,\omega}^{\rm atom}(\bx,\by)\  
   \left(\ \Theta_\omega(\by) f(\by)\ \right)\ d\by\ \right]\ d\bx\ =\ 0\ .
\label{orth-Tof}  \end{equation}
   Next we express $\VGamma$ as:
   \[ \VGamma(\bx)\ =\ \lambda^2V_\omega(\bx)\ +\ 
   \sum_{\omega^\prime\in\Gamma\setminus\{\omega\}}\lambda^2 V_{\omega^\prime}(\bx)\ ,\]
   and therefore by \eqref{Tof}
   \begin{align}
 \Theta_\omega(\bx) f(\bx)\ &=\
  \left(\ -\Delta_\bx+\VGamma(\bx)-E_0^\lambda\ \right)\ 
  \int_{\R^2} \ G_{\lambda,\omega}^{\rm atom}(\bx,\by)\Theta_\omega(\by)\cdot f(\by)\ d\by \ d\by\nn\\
 & \qquad -\sum_{\omega^\prime\in\Gamma\setminus\{\omega\}}\ \lambda^2V_{\omega^\prime}(\bx)\
\int_{\R^2}G_{\lambda,\omega}^{\rm atom}(\bx,\by) \Theta_\omega(\by)\cdot f(\by)\ d\by
  +\ \left\langle p_\omega^\lambda,\Theta_\omega f\right\rangle_{L^2(\R^2)}\ p_\omega^\lambda(\bx)\ .
 \label{Tof1} \end{align}
 \medskip
 
 Similarly, 
  \begin{align}
 \Tfree(\bx) f(\bx)\ &=\ \left(\ -\Delta_\bx\ -\ E_0^\lambda\  \right)\
  \int_{\R^2}\ \Gfree(\bx-\by)\ \left( \Tfree(\by) f(\by) \right)\ d\by\nn\\
  &=  \left(\ -\Delta_\bx\ +\ \VGamma(\bx)\ -\ E_0^\lambda\  \right)\
  \int_{\R^2}\ \Gfree(\bx-\by)\Tfree(\by)\cdot f(\by)\ d\by\nn\\
  &\qquad -\ \VGamma(\bx)\  \int_{\R^2}\ \Gfree(\bx-\by)\Tfree(\by)\cdot f(\by)\ d\by\ .
 \label{Tfreef} \end{align}
  We note that $\VGamma(\bx)\equiv0$ on the support of $\Tfree$. 
  
Now summing  \eqref{Tof1} over $\omega\in\Gamma$  and adding the result to \eqref{Tfreef}, we have  by\eqref{f-split} the following:
 \begin{align}
 f(\bx)\ &=\ \left(\ -\Delta_\bx\ +\ \VGamma(\bx)\ -\ E_0^\lambda\  \right)\cdot\nn\\
 &\qquad  \int_{\R^2}\ \Big[\ \sum_{\omega\in\Gamma}\ G_{\lambda,\omega}^{\rm atom}(\bx,\by)\Theta_\omega(\by)\ +\ \Gfree(\bx-\by)\Tfree(\by)
  \Big]\cdot f(\by)\ d\by\nn\\
  &\qquad -  \int_{\R^2} \ \left[\  \sum_{\substack{\omega,\omega^\prime\in\Gamma \\
 \omega\ne\omega^\prime}}  \lambda^2\ V_{\omega^\prime}(\bx) G_{\lambda,\omega}^{\rm atom}(\bx,\by)\ 
\Theta_\omega(\by)\ +\ \VGamma(\bx)\Gfree(\bx-\by)\Tfree(\by)\ \right]\ \cdot f(\by)\ d\by\nn\\
 &\qquad +\ \sum_{\omega\in\Gamma}\ \left\langle \Theta_\omega \ p_\omega^\lambda,\ f\right\rangle_{L^2(\R^2)}\ p_\omega^\lambda(\bx)\ .\label{pre-inv}
 \end{align}
Introduce the kernels $ K_0^\lambda$ and $\mathcal{E}_0^\lambda$:
 \begin{align}
 K_0^\lambda(\bx,\by)\ &\equiv\ 
 \sum_{\omega\in\Gamma}\ G_{\lambda,\omega}^{\rm atom}(\bx,\by)\Theta_\omega(\by)\ +\ \Gfree(\bx-\by)\Tfree(\by)
\label{K0-def}\\
\mathcal{E}_0^\lambda(\bx,\by)\ &\equiv\ 
\sum_{\substack{\omega,\omega^\prime\in\Gamma \\
 \omega\ne\omega^\prime}}  \lambda^2\ V_{\omega^\prime}(\bx) G_{\lambda,\omega}^{\rm atom}(\bx,\by)\ 
\Theta_\omega(\by)\ +\ \VGamma(\bx)\Gfree(\bx-\by)\Tfree(\by)\ .
\label{E0-def}
 \end{align}
 
Equation \eqref{pre-inv} is equivalent to 
 \begin{align}
 f(\bx)\ &=\  \left(\ -\Delta_\bx\ +\ \VGamma(\bx)\ -\ E_0^\lambda\  \right)\ \int_{\R^2}\ K_0^\lambda(\bx,\by)\ f(\by)\ d\by\nn\\
 &\qquad +\ \sum_{\omega\in\Gamma}\ \left\langle \Theta_\omega \ p_\omega^\lambda,\ f\right\rangle_{L^2(\R^2)}\ p_\omega^\lambda(\bx)\ -\ \int_{\R^2}\ \mathcal{E}_0^\lambda(\bx,\by)\ f(\by)\ d\by\ .
 \label{appinv1}
 \end{align}
and in any even more compact form :
  \begin{align}
 f(\bx)\ &=\  \left(\ -\Delta_\bx\ +\ \VGamma(\bx)\ -\ E_0^\lambda\  \right)\ K^\lambda_0[f](\bx)\ -
 \  \mathcal{E}_0^\lambda[f](\bx)\nn\\
 &\qquad +\ \sum_{\omega\in\Gamma}\ \left\langle \Theta_\omega \ p_\omega^\lambda,\ f\right\rangle_{L^2(\R^2)}\ p_\omega^\lambda(\bx).
 \label{appinv1a}
 \end{align}

 \begin{proposition}\label{K0E0}
 $ K_0^\lambda(\bx,\by)$ is a main kernel in the sense of Definition \ref{mainK} and $\mathcal{E}_0^\lambda(\bx,\by)$ is an error kernel
  in the sense of \eqref{errorK}.
 \end{proposition}
 
\nit{\it Proof of Proposition \ref{K0E0}:}\ We first prove that $K_0^\lambda(\bx,\by)$, displayed in \eqref{K0-def}, is a main kernel. Note that for each $\by\in\R^2$ there is at most one $\omega=\omega_\by\in\Gamma$ with $\by\in {\rm supp}\ \Theta_\omega\ \subset\ \{\by:|\by-\omega|\le r_4\}$. Therefore, for the first term in \eqref{K0-def} we have by Theorem \ref{ptGat} the bound
\begin{align}
\Big|\ \sum_{\omega\in\Gamma}\ G_{\lambda,\omega}^{\rm atom}(\bx,\by)\Theta_\omega(\by)\
\Big|\ &\le\ \left|\ G^{\rm atom}_{\lambda,\omega_\by}(\bx,\by)\ \right|\nn\\
&  \lesssim C\left[\ \lambda^4\ +\ \left|\ \log|\bx-\by|\ \right|\ \right]\ {\bf 1}_{_{\{|\bx-\by|\le R\}}}\ +\ e^{-c\lambda}\ e^{-c\lambda|\bx-\by|}\ .\nn
\end{align}
Furthermore by \eqref{Gf-allx},  the second term in \eqref{K0-def} satisfies  the bound
\begin{align}
\Big|\ \Gfree(\bx-\by)\Tfree(\by)\ \Big|\ &\le  \Big|\ \Gfree(\bx-\by)\ \Big|\nn\\
&\lesssim\ 
 C\left[\ \lambda^4\ +\ \left|\ \log|\bx-\by|\ \right|\ \right]\ {\bf 1}_{_{\{|\bx-\by|\le R\}}}\ +\ e^{-c\lambda}\ e^{-c\lambda|\bx-\by|}\nn
\end{align} 
Adding the two previous bounds we conclude that $K_0^\lambda(\bx,\by)$ is a main kernel.\medskip

We now prove that $\mathcal{E}_0^\lambda(\bx,\by)$ given by \eqref{E0-def}
%
%
 is an error kernel. Consider the sum in \eqref{E0-def}. This sum is non-zero at 
 $(\bx, \by)\in\R^2\times\R^2$, if there are distinct points $\omega^\prime_\bx, \omega_\by\in\Gamma$  with $\bx\in{\rm supp}\ V_{\omega^\prime_\bx}$ and  $\by\in {\rm supp}\ \Theta_{\omega_\by}$ .
The choice of points $\omega^\prime_\bx, \omega_\by\in\Gamma$ is unique.
 We have $\by\in B_{r_4}(\omega_\by)$ and $\bx\notin  B_{r_4+\delta_1}(\omega_\by)$,
 where $\delta_1>0$. Therefore, part (3) of Theorem \ref{ptGat} implies
 %
%
\begin{align} 
\Big| \sum_{\substack{\omega,\omega^\prime\in\Gamma \\
 \omega\ne\omega^\prime}}  \lambda^2\ V_{\omega^\prime}(\bx) G_{\lambda,\omega}^{\rm atom}(\bx,\by)\ 
\Theta_\omega(\by)\ \Big|\ 
&\le\ \lambda^2\ |V_{\omega^\prime_\bx}(\bx)|\  |G_{\lambda,\omega_\by}^{\rm atom}(\bx,\by)|\ 
\Theta_{\omega_\by}(\by)\nn\\
&\le\ \lambda^2\ e^{-c\lambda}\ e^{-c\lambda|\bx-\by|}\ \lesssim\ e^{-\cp\lambda}\ e^{-\cp\lambda|\bx-\by|}\ . \nn
\end{align}
For the second term in \eqref{E0-def}, if $\bx\in {\rm supp}\ V_\Gamma$ and $\by\in {\rm supp}\
 \Tfree$, then $|\bx-\by|\ge r_3-r_0>0$. Therefore, 
 $\Gfree(\bx-\by)\lesssim e^{-c\lambda|\bx-\by|}\lesssim e^{-\cp\lambda}\ e^{-\cp\lambda|\bx-\by|}$. It follows that for some $\omega=\omega_\bx\in\Gamma$:
 \[
 \Big|\ \VGamma(\bx)\Gfree(\bx-\by)\Tfree(\by)\ \Big|\ \lesssim 
\lambda^2 \Big|\ V_{\omega_\bx}(\bx)\Gfree(\bx-\by)\ \Big|\ \lesssim\ e^{-c\lambda}\ 
e^{-c\lambda|\bx-\by|}.
 \]
 The latter two bounds imply that $\mathcal{E}_0^\lambda(\bx,\by)$, defined in \eqref{E0-def}, is an error kernel. 
  The proof of Proposition \ref{K0E0} is now complete.

 \begin{remark}\label{trans-invar}
 At this stage we wish to remark that if $\Gamma$ is translation invariant by some vector, then $K_0^\lambda$ and $\mathcal{E}_0^\lambda$ inherit this invariance. In particular, for $\Gamma=\mathbb{H}_\sharp$,
   the zigzag  truncation of the honeycomb $\mathbb{H}$, 
we have $ K_0^\lambda(\bx+\vtilde_2,\by+\vtilde_2)=K_0^\lambda(\bx,\by)$ and $ \mathcal{E}_0^\lambda(\bx+\vtilde_2,\by+\vtilde_2)=
 \mathcal{E}_0^\lambda(\bx,\by)$.
 \end{remark}

Introduce the orthogonal subspaces $\mathscr{X}_\Gamma$:
\begin{align}
\mathscr{X}_\Gamma\ &\equiv\ {\rm span}\Big\{p_\omega^\lambda:\omega\in\Gamma\Big\}^\perp
 =\ \Big\{f\in L^2(\R^2)\ :\ \left\langle p_\omega^\lambda, f\right\rangle_{_{L^2(\R^2)}}=0,\ \ \omega\in\Gamma\ \Big\}\ , 
\label{XGam1}
\end{align}
and  the orthogonal projections: 
\begin{equation}
 \Pi_\Gamma^\lambda:L^2(\R^2)\to \mathscr{X}_\Gamma,\qquad \widetilde\Pi_\Gamma^\lambda=I-\Pi_\Gamma^\lambda:
L^2(\R^2)\to{\rm span}\Big\{p_\omega^\lambda:\omega\in\Gamma\Big\}.
\label{PXGam1}
\end{equation} 
We seek the integral kernel for the inverse of the operator $\Pi^\lambda_\Gamma\ \left(\ H_\Gamma^\lambda-E_0^\lambda-\Omega\ \right)\ \Pi^\lambda_\Gamma$ on $\mathscr{X}_\Gamma$.  

 The operator $f\mapsto K_0^\lambda f$ (see \eqref{K0-def}, \eqref{appinv1})  defines an approximate inverse of $H^\lambda_\Gamma-E_0^\lambda-\Omega$ on the range of $\Pi_\Gamma^\lambda$ but we do not have
that $\Pi_\Gamma^\lambda K_0^\lambda[f]=K_0^\lambda[f]$.  Our next step is to correct  $K_0^\lambda$
in order achieve the desired projection. 
\medskip

Recall that the set $\{p_\omega^\lambda\ :\ \omega\in \Gamma\ \}$ is not orthonormal, but only nearly so; see Proposition \ref{prop:pkpar}.  The following lemma gives a representation for 
$\widetilde\Pi_\Gamma^\lambda$, defined in \eqref{PXGam1}.

\begin{lemma}\label{Pi-Gamma}
$\widetilde{\Pi}_\Gamma^\lambda=I-\Pi_\Gamma^\lambda$, the orthogonal projection of $L^2(\R^2)$ onto ${\rm span}\{p_\omega^\lambda\ :\ \omega\in \Gamma\ \}$, is given by 
\begin{align}
\widetilde{\Pi}_\Gamma^\lambda[g](\bx)\ &=\ 
\sum_{\omega,\hat\omega\in\Gamma}\  M^{\omega,\hat\omega}\left\langle p_{\hat\omega}^\lambda,g\right\rangle\ p^\lambda_{\omega}(\bx)\ ,
\label{PiTilde}
\end{align}
where $M^{\omega,\hat\omega}$  satisfies the estimate
\begin{align}
\left|\ M^{\omega,\hat\omega}\ -\ \delta_{\omega,\hat\omega}\ \right|\ &\lesssim\ e^{-c\lambda}\ e^{-c\lambda\ |\hat\omega-\omega|}\ .
\label{M-est}
\end{align}
\end{lemma}

\nit{\it Proof of Lemma \ref{Pi-Gamma}:}\ 
 If we define $\widetilde{\Pi}_\Gamma^\lambda[g]$ by \eqref{PiTilde}, then for all $g\in L^2(\R^2)$
\begin{align}
\left\langle p^\lambda_{\omega^\prime}, g\right\rangle\ =\left\langle p^\lambda_{\omega^\prime}, \widetilde{\Pi}_\Gamma^\lambda[g]\right\rangle\ &=\ 
\sum_{\omega,\hat\omega\in\Gamma}\  M^{\omega,\hat\omega}
\left\langle p^\lambda_{\hat\omega},g\right\rangle\ \left\langle p^\lambda_{\omega^\prime},p^\lambda_{\omega}\right\rangle\nn\\
&=\ \sum_{\omega,\hat\omega\in\Gamma}\ 
\Big(\ \sum_{\omega\in\Gamma}\
\left\langle p^\lambda_{\omega^\prime},p^\lambda_{\omega}\right\rangle\ M^{\omega,\hat\omega}\ \Big)\ \left\langle p^\lambda_{\hat\omega},g\right\rangle\
\end{align}
Therefore, $\widetilde{\Pi}_\Gamma^\lambda$ is as required provided:
\[
\sum_{\omega\in\Gamma}\
\left\langle p^\lambda_{\omega^\prime},p^\lambda_{\omega}\right\rangle\ M^{\omega,\hat\omega}\ =\
 \delta_{\omega^\prime\hat\omega}\ .
\]
We claim that if $\omega^\prime,\omega\in \Gamma$ are distinct, then 
\begin{equation} 
\left| \ \left\langle p^\lambda_{\omega^\prime},p^\lambda_{\omega}\right\rangle \ \right|\ \lesssim\ e^{-\cp\lambda|\omega-\omega^\prime|}\ e^{-\cp\lambda}. 
\label{nr-orth}
\end{equation}
Indeed, if $\omega\ne \omega^\prime$
\begin{align}
\left| \ \left\langle p^\lambda_{\omega^\prime},p^\lambda_{\omega}\right\rangle \ \right|
&\le\ \int_{B_{_{r_4}}(\omega)}\ p^\lambda_{\omega^\prime}(\bx)\ p^\lambda_{\omega}(\bx)\ d\bx
\ +\  \int_{B_{_{r_4}}(\omega^\prime)}\ p^\lambda_{\omega^\prime}(\bx)\ p^\lambda_{\omega}(\bx)\ d\bx\nn\\
&\qquad +\ \int_{_{\R^2\setminus B_{_{r_4}}(\omega)\cup B_{_{r_4}}(\omega^\prime)}}\
 p^\lambda_{\omega^\prime}(\bx)\ p^\lambda_{\omega}(\bx)\ d\bx\nn\\
 &\le\ \int_{B_{_{r_4}}(\omega)}\ \left[\ e^{-c\lambda|\bx-\omega^\prime|}\ \right]\cdot \left[\ \lambda^2\ \right] d\bx\ 
 +\ \int_{B_{_{r_4}}(\omega^\prime)}\ \left[\ \lambda^2\ \right]\cdot  \left[\ e^{-c\lambda|\bx-\omega|}\ \right] d\bx\nn\\
&\qquad +\ \int_{_{\R^2\setminus B_{_{r_4}}(\omega)\cup B_{_{r_4}}(\omega^\prime)}}\
 e^{-c\lambda|\bx-\omega|}
\cdot e^{-c\lambda|\bx-\omega^\prime|} \ d\bx\ \lesssim \ e^{-\cp\lambda|\omega-\omega^\prime|}\ e^{-\cp\lambda}. \nn
\end{align}
Since also $p_\omega^\lambda(\bx)=p_0^\lambda(\bx-\omega)$ is normalized in $L^2(\R^2)$, we have 
\begin{equation}
\Big|\ \left\langle p^\lambda_{\omega^\prime},p^\lambda_{\omega}\right\rangle\ -\ \delta_{\omega,\omega^\prime}\ \Big|
\ \lesssim\ e^{-c\lambda}\ e^{-c\lambda|\omega-\omega^\prime|}\ .
\label{ip-diff}\end{equation}
\medskip

Let $P=\left(\ \left\langle p^\lambda_{\omega^\prime},p^\lambda_{\omega}\right\rangle\ \right)_{_{\omega,\omega^\prime\in\Gamma}}$ and for any  $\nu\in\R^2$, $|\nu|=1$, let 
$D=\left(\ e^{\bar{c}\lambda \nu\cdot\omega}\ \delta_{\omega,\omega^\prime}\ \right)_{_{\omega,\omega^\prime\in\Gamma}}$, with $\bar{c}$ smaller than the constant $c$ appearing in \eqref{ip-diff}.
 Then, $ D\ P\ D^{-1}\ =\ 
 \left(\ e^{\bar{c}\lambda \nu\cdot(\omega-\omega^\prime)}\  \left\langle p^\lambda_{\omega^\prime},p^\lambda_{\omega}\right\rangle\  \right)_{_{\omega,\omega^\prime\in\Gamma}}\ =\ (\tilde{p}_{\omega,\omega^\prime})$ with 
 \[ 
  \Big| \ \tilde{p}_{\omega,\omega^\prime}\ -\  \delta_{\omega,\omega^\prime}\ \Big|\ \lesssim\ 
  e^{-\cp\lambda|\omega-\omega^\prime|}\ e^{-\cp\lambda}.
  \]
  by \eqref{ip-diff}. Hence, $D\ P^{-1}\ D^{-1}\ =\ \left(\ DPD^{-1}\ \right)^{-1}$ has an $(\omega,\omega^\prime)-$ entry
   that differs from $\delta_{\omega,\omega^\prime}$ by at most $e^{-\tc\lambda}$. That is,
   $
   \Big|\ \left[\ e^{\bar{c}\lambda \nu\cdot(\omega-\omega^\prime)}\ M^{\omega,\omega^\prime}\ \right]\ -\
    \delta_{\omega,\omega^\prime}\ \Big|\ \lesssim\ e^{-c\lambda}$ and hence
    \[
    \Big|\  e^{\bar{c}\lambda \nu\cdot(\omega-\omega^\prime)}\ \left[\ M^{\omega,\omega^\prime}\ -\
    \delta_{\omega,\omega^\prime}\ \right] \Big|\ \lesssim\ e^{-c\lambda}
    \]
    for all $\omega, ,\omega^\prime\in\Gamma$ and all unit vectors $\nu\in\R^2$. Optimizing over $\nu$
     gives
      \[
    \Big|\   M^{\omega,\omega^\prime}\ -\
    \delta_{\omega,\omega^\prime}\  \Big|\ \lesssim\ e^{-c\lambda}\ \ e^{-\bar{c}\lambda |\omega-\omega^\prime|} \ .   \]
   This completes the proof of   Lemma \ref{Pi-Gamma}.\medskip
   \medskip
 
By \eqref{appinv1a}, after subtracting and adding $\widetilde{\Pi}^\lambda_\Gamma\ K_0^\lambda$,  we have
 \begin{align}
 f(\bx)\ 
&=\  \left(\ -\Delta_\bx\ +\ \VGamma(\bx)\ -\ E_0^\lambda\  \right)\ \Big[\ K^\lambda_0[f](\bx)\ -\
  \left(\widetilde{\Pi}^\lambda_\Gamma\ K_0^\lambda\right)[f]\ \Big] \nn\\
 &\quad +\ \left(\ -\Delta_\bx\ +\ \VGamma(\bx)\ -\ E_0^\lambda\  \right)\
  \left(\widetilde{\Pi}^\lambda_\Gamma\ K_0^\lambda\right)[f] \nn\\  
 &\quad - \  \mathcal{E}_0^\lambda[f](\bx)\ +\ \sum_{\omega\in\Gamma}\ \left\langle \Theta_\omega \ p_\omega^\lambda,\ f\right\rangle_{L^2(\R^2)}\ p_\omega^\lambda(\bx)\ .
\label{appinv2} \end{align}

\nit Here, we have arranged for the expression within the square brackets in \eqref{appinv2}:
\begin{equation}
K^\lambda_1[f]\ \equiv\ K^\lambda_0[f]\ -\
  \left(\widetilde{\Pi}^\lambda_\Gamma\ K_0^\lambda\right)[f], 
  \label{K1-def}
  \end{equation}
 to be orthogonal to the translated atomic ground states $p_\omega^\lambda$, for all $\omega\in\Gamma$. 
 Our next task is to show that the remaining terms in \eqref{appinv2} comprise an error kernel.
 
 \begin{proposition}\label{HGP-ek} The operators $\widetilde{\Pi}^\lambda_\Gamma\ K_0^\lambda$ and 
$ \left(\ -\Delta_\bx\ +\ \VGamma(\bx)\ -\ E_0^\lambda\  \right)\
  \left(\widetilde{\Pi}^\lambda_\Gamma\ K_0^\lambda\right)$
  derive from error kernels   in the sense of Definition \ref{errorK}.
 \end{proposition}

\nit{\it Proof of Proposition \ref{HGP-ek}:}\ By \eqref{PiTilde}
\begin{align}
&\left(\ \widetilde{\Pi}^\lambda_\Gamma K_0^\lambda \right)[f](\bx) =\ 
\sum_{\omega,\hat\omega\in\Gamma}\  M^{\omega,\hat\omega}\left\langle p_\omega^\lambda\ ,\ K_0^\lambda[f]\ \right\rangle\ p^\lambda_{\hat\omega}(\bx)\ 
\label{tPiK0}\\
&=\sum_{\omega,\hat\omega\in\Gamma}\  M^{\omega,\hat\omega}\ 
\int_{\R^2}\ p_\omega^\lambda(\by)\  \int_{\R^2}\ K_0^\lambda(\by,\bz) f(\bz)\ d\bz\ d\by\ p^\lambda_{\hat\omega}(\bx)\nn\\
&=\sum_{\omega,\hat\omega\in\Gamma}\  M^{\omega,\hat\omega}\ 
\int_{\R^2}\ \left[\ \int_{\R^2}\ p_\omega^\lambda(\by)\  K_0^\lambda(\by,\bz)\ d\by\ p^\lambda_{\hat\omega}(\bx)\ \right] f(\bz)\ d\bz\ \nn\\
&=
\int_{\R^2}\ \left[\ \sum_{\omega,\hat\omega\in\Gamma}\  M^{\omega,\hat\omega}\ \int_{\R^2}\ p_\omega^\lambda(\by)\  K_0^\lambda(\by,\bz)\ d\by\ p^\lambda_{\hat\omega}(\bx)\ \right] f(\bz)\ d\bz\ \nn
\end{align}
Thus,
\begin{align}
\left(\ \widetilde{\Pi}^\lambda_\Gamma K_0^\lambda \right)(\bx,\bz)\ &=\
   \int_{\R^2}\ \left[\ \sum_{\omega,\hat\omega\in\Gamma}\ M^{\omega,\hat\omega}\ p^\lambda_{\hat\omega}(\bx)\ p_\omega^\lambda(\by)\ \right] \
 K_0^\lambda(\by,\bz)\ d\by\ , \nn
\end{align} 
where $K_0^\lambda$ is given by \eqref{K0-def}:
\begin{align} K_0^\lambda(\by,\bz)\ 
&\equiv\ 
 \sum_{\omega^\prime\in\Gamma}\ G_{\lambda,\omega^\prime}^{\rm atom}(\by,\bz)\Theta_{\omega^\prime}(\bz)\ +\ \Gfree(\by-\bz)\Tfree(\bz)\ .
\label{K0atfr} \end{align}

\nit Now decompose $\left(\ \widetilde{\Pi}^\lambda_\Gamma K_0^\lambda \right)(\bx,\bz)$ has follows:
\begin{align}
&\left(\ \widetilde{\Pi}^\lambda_\Gamma K_0^\lambda \right)(\bx,\bz)\ =\
   \int_{\R^2}\ \left[\ 
   \sum_{\substack{\omega,\hat\omega\in\Gamma\\ \omega\ne\hat\omega}}\ M^{\omega,\hat\omega}\ p^\lambda_{\hat\omega}(\bx)\ p_\omega^\lambda(\by)\ \right] \
 K_0^\lambda(\by,\bz)\ d\by\ \nn\\
 &\qquad\qquad \qquad\qquad +\  
 \int_{\R^2}\ \left[\ 
   \sum_{\omega\in\Gamma}\ M^{\omega,\omega}\ p^\lambda_{\omega}(\bx)\ p_\omega^\lambda(\by)\ \right] \
 K_0^\lambda(\by,\bz)\ d\by\nn\\
 &\equiv\ \ \left(\ \widetilde{\Pi}^\lambda_\Gamma K_0^\lambda \right)_{1}(\bx,\bz)\ +\
   \left(\ \widetilde{\Pi}^\lambda_\Gamma K_0^\lambda \right)_{2}(\bx,\bz)\ . \label{PiK012}
\end{align} 

We prove that each term in \eqref{PiK012} is an error kernel,  {\it i.e.}   $\Big| \left(\ \widetilde{\Pi}^\lambda_\Gamma K_0^\lambda \right)_{j}(\bx,\bz) \Big|\ \lesssim\ 
e^{-c\lambda}\ e^{-c\lambda|\bx-\bz|}$ for $j=1,2$. 
   For $\omega\ne\hat\omega$ we have by \eqref{M-est} that 
 \[| M^{\omega,\hat\omega} |\lesssim e^{-\cp\lambda|\omega-\hat\omega|}\ e^{-\cp\lambda}.\]
 We may therefore write:
\begin{align}
&| \ M^{\omega,\hat\omega}\ p^\lambda_{\hat\omega}(\bx)\ p_\omega^\lambda(\by)\  |\ 
\le\ e^{-\cp\lambda|\omega-\hat\omega|}\ e^{-\tc\lambda}p^\lambda_{\hat\omega}(\bx)\cdot 
e^{-\tc\lambda}p^\lambda_{\omega}(\by) .
\label{Mpp-bd}
\end{align}
Next, using \eqref{p0-bd-cpam} we bound $e^{-\tc\lambda}p^\lambda_{\hat\omega}(\bx)$ and 
$e^{-\tc\lambda}p^\lambda_{\omega}(\by)$ as follows:
\begin{align}
 e^{-\tc\lambda}p^\lambda_{\hat\omega}(\bx)\ &\lesssim\ \left(\  e^{-\cp\lambda}{\bf 1}_{\{|\bx-\hat\omega|\le r_1\}}+e^{-\cp\lambda} e^{-c\lambda|\bx-\hat\omega|}\ \right)\nn\\
 &\lesssim\ \left(\  e^{-\frac{\cp}{2}\lambda}\ e^{-\frac{\cp}{2r_1}\lambda|\bx-\hat\omega|}{\bf 1}_{\{|\bx-\hat\omega|\le r_1\}}+e^{-\cp\lambda} e^{-c\lambda|\bx-\hat\omega|}\ \right)
\end{align}
Therefore,  $e^{-\tc\lambda}p^\lambda_{\hat\omega}(\bx)\lesssim e^{-c\lambda} e^{-c\lambda|\bx-\hat\omega|}$
and similarly  $e^{-\tc\lambda}p^\lambda_{\omega}(\by)\lesssim e^{-c\lambda} e^{-c\lambda|\by-\omega|}$.
Substituting these bounds into \eqref{Mpp-bd}, we obtain for some $c>0$
\begin{align}
|\ M^{\omega,\hat\omega}\ p^\lambda_{\hat\omega}(\bx)\ p_\omega^\lambda(\by)\  |\ 
&\lesssim\ e^{-c\lambda}\ e^{-c\lambda|\omega-\hat\omega|}\ e^{-c\lambda|\bx-\hat\omega|}\ e^{-c\lambda|\by-\omega|}\nn\\
&\lesssim e^{-c\lambda}\  e^{-\frac{c}2\lambda|\bx-\by|}\ \times\ e^{-\frac{c}2\lambda|\omega-\hat\omega|}\ 
e^{-\frac{c}2\lambda|\bx-\hat\omega|}\ e^{-\frac{c}2\lambda|\by-\omega|}\  ,\nn
\end{align}
since $|\bx-\by|\le |\bx-\hat\omega|+|\omega-\hat\omega| + |\by-\omega|$. Therefore, for some $c^\prime$ which is independent of $\lambda$:
\[
\sum_{\substack{ \omega,\hat\omega\in\Gamma \\ \omega\ne\hat\omega}}\   |\ M^{\omega,\hat\omega}\ p^\lambda_{\hat\omega}(\bx)\ p_\omega^\lambda(\by)\  |\ \lesssim\ e^{-c^\prime\lambda}\  e^{-c^\prime\lambda|\bx-\by|}
\]
and therefore $\sum_{\substack{ \omega,\hat\omega\in\Gamma \\ \omega\ne\hat\omega}}\ M^{\omega,\hat\omega}\ p^\lambda_{\hat\omega}(\bx)\ p_\omega^\lambda(\by)$ is therefore an error kernel. And since $K_0^\lambda$ is a main kernel we have, by the expression for $\left(\ \widetilde{\Pi}^\lambda_\Gamma K_0^\lambda \right)_{1}(\bx,\bz)$ in \eqref{PiK012}, and by part 2 of Lemma 
\ref{error-main}, that $\left(\ \widetilde{\Pi}^\lambda_\Gamma K_0^\lambda \right)_{1}(\bx,\bz)$  is an error kernel.

We next prove that $\left(\ \widetilde{\Pi}^\lambda_\Gamma K_0^\lambda \right)_{2}(\bx,\bz)$, defined in \eqref{PiK012} is an error kernel. 
Using \eqref{K0atfr} we have
\begin{align}
\left(\ \widetilde{\Pi}^\lambda_\Gamma K_0^\lambda \right)_{2}(\bx,\bz)\ &\equiv\ \sum_{\omega\in\Gamma} \ M^{\omega,\omega}\ p^\lambda_{\omega}(\bx)\ \int_{\R^2}\   p_\omega^\lambda(\by) \
 K_0^\lambda(\by,\bz)\ d\by  \nn\\
&=\  \sum_{\omega\in\Gamma} \ M^{\omega,\omega}\ p^\lambda_{\omega}(\bx)\ \int_{\R^2}\   p_\omega^\lambda(\by) \
 \left[\   \sum_{\omega^\prime\in\Gamma\setminus\{\omega\}}\ G_{\lambda,\omega^\prime}^{\rm atom}(\by,\bz)\Theta_{\omega^\prime}(\bz) \ \right]\ d\by\nn\\
 &\qquad +\ \sum_{\omega\in\Gamma} \ M^{\omega,\omega}\ p^\lambda_{\omega}(\bx)\ \int_{\R^2}\   p_\omega^\lambda(\by) \
\Gfree(\by-\bz)\Tfree(\bz)\ d\by \nn\\
&\equiv\ \left(\ \widetilde{\Pi}^\lambda_\Gamma K_0^\lambda \right)_{2a}(\bx,\bz)\ +\ \left(\ \widetilde{\Pi}^\lambda_\Gamma K_0^\lambda \right)_{2b}(\bx,\bz). \label{I-IIdef}
\end{align}
Note the absence of the $\omega^\prime=\omega$ term in the inner sum just above since the atomic Green's function,
 $G_{\lambda,\omega^\prime}^{\rm atom}$, projects onto the orthogonal complement of the function $p_{\omega^\prime}^\lambda$. 

We prove  that the kernels 
$\left(\ \widetilde{\Pi}^\lambda_\Gamma K_0^\lambda \right)_{2a}(\bx,\bz)$ and $\left(\ \widetilde{\Pi}^\lambda_\Gamma K_0^\lambda \right)_{2b}(\bx,\bz)$, defined in \eqref{I-IIdef} are both bounded in absolute value by 
 $e^{-c\lambda}\ e^{-c\lambda|\bx-\bz|}$. We first recall the following relations and definitions:
\begin{align*}
&G^{\rm atom}_{\lambda,\omega}(\by,\bz) = \Gatom(\bx-\omega,\by-\omega),\quad \\
& \left(\ \Hatom-E_0^\lambda\ \right)\Gatom(\bx,\by) = \delta(\bx-\by)\ -\ p_0^\lambda(\bx)\ p_0^\lambda(\by)\\
& \Theta_0(\bx)\equiv
\begin{cases} 
1 , & |\bx|\le r_3\\
0 , & |\bx|\ge r_4
\end{cases},\qquad {\rm and}\\
&\Theta_\omega(\bx)\ =\ \Theta(\bx-\omega),\quad\textrm{for $\omega\in\Gamma$, and }\ 
\Theta_{\rm free}(\bx)\ =\ 1-\sum_{\omega\in\Gamma}\Theta_\omega(\bx)\ .
\end{align*}
\medskip

\nit{\it Estimation of $\left(\ \widetilde{\Pi}^\lambda_\Gamma K_0^\lambda \right)_{2a}(\bx,\bz)$; see \eqref{I-IIdef}}:\ Suppose first that $|\bz-\omega^\prime|\ge r_4$, for all $\omega^\prime\in\Gamma\setminus\{\omega\}$. Then, $\bz$ is outside the support of  $\Theta_{\omega^\prime}(\bz)$ for all $\omega^\prime\in\Gamma\setminus\{\omega\}$.
and we have: $\left(\ \widetilde{\Pi}^\lambda_\Gamma K_0^\lambda \right)_{2a}(\bx,\bz)\equiv0$. 

Suppose now that $\bz$ is such that $|\bz-\omega^\prime|\le r_4$ for some $\omega^\prime=\omega^\prime_\bz\in\Gamma\setminus\{\omega\}$. Therefore, the bracketed expression in the definition of $\left(\ \widetilde{\Pi}^\lambda_\Gamma K_0^\lambda \right)_{2a}$ (see  \eqref{I-IIdef}) is given by:
 $\left[\cdots\right](\by,\bz)=G^{\rm atom}_{\lambda,\omega^\prime_\bz}(\by,\bz)\Theta_{\omega^\prime_\bz}(\bz)$. Therefore, for $|\bz-\omega_\bz^\prime|\le r_4$, we have
 \begin{align}
\ \int\ p_\omega^\lambda(\by)\ \left[\cdots\right](\by,\bz)\ d\by\ 
 &=\ 
  \int\ p_\omega^\lambda(\by)\ G^{\rm atom}_{\lambda,\omega^\prime_\bz}(\by,\bz)\Theta_{\omega^\prime_\bz}(\bz)\ d\by\ \nn\\
  & \le\ \int_{|\by-\omega|\le r_1}\ p_\omega^\lambda(\by)\ G^{\rm atom}_{\lambda,\omega^\prime_\bz}(\by,\bz)\ d\by\ \nn\\
  &\qquad\qquad +\ \int_{|\by-\omega|\ge r_1}\ p_\omega^\lambda(\by)\ G^{\rm atom}_{\lambda,\omega^\prime_\bz}(\by,\bz)\ d\by\ .
\label{i-split}  \end{align}
We bound the latter two integrals individually by using the pointwise bounds on $p_\omega^\lambda(\by)=p_0^\lambda(\by-\omega)$  given in \eqref{p0-bd-cpam} and the pointwise bounds on $G^{\rm atom}_{\lambda,\omega^\prime_\bz}(\by,\bz)=G_\lambda^{\rm atom}(\by-\omega^\prime_\bz,\bz-\omega^\prime_\bz)$ of Theorem \ref{ptGat}.

With $|\bz-\omega_\bz^\prime|\le r_4$, we first consider the integral over the set $|\by-\omega|\le r_1$.  
For such $\by$, we have by \eqref{p0-bd-cpam}: $|p_\omega^\lambda(\by)|\lesssim\lambda^2$. 
Furthermore, note that $|\by-\omega_\bz^\prime|\ge|\omega-\omega_\bz^\prime|-|\by-\omega|\ge r_{\rm min}-r_1>r_4$;
 see\eqref{r_jsA}.  Because $|\by-\omega_\bz^\prime|>r_{\rm min}-r_1$, while $|\bz-\omega_\bz^\prime|<r_4$, it follows from \eqref{variant} (part 3 of Theorem \ref{ptGat}) that 
 $|G^{\rm atom}_{\lambda,\omega^\prime_\bz}(\by,\bz)| \lesssim e^{-c\lambda} e^{-c\lambda|\by-\bz|}$.
 The first integral in \eqref{i-split} therefore satisfies
  \begin{align*}
\int_{|\by-\omega|\le r_1}\ p_\omega^\lambda(\by)\ G^{\rm atom}_{\lambda,\omega^\prime_\bz}(\by,\bz)\ d\by
\lesssim\ \lambda^2\ \int_{|\by-\omega|\le r_1}\ e^{-c\lambda} e^{-c\lambda|\by-\bz|}\ d\by\ \lesssim\ e^{-c\lambda} e^{-c\lambda|\bz-\omega|}\ .
  \end{align*}
  
  Next, with $|\bz-\omega_\bz^\prime|\le r_4$, we consider the integral over the set $|\by-\omega|\ge r_1$. On this set, we have $|p_\omega^\lambda(\by)|\lesssim e^{-\cp\lambda}e^{-\cp\lambda |\by-\omega|}$ and, by the bounds of 
  Theorem \ref{ptGat}:
  \begin{align*}
&\int_{|\by-\omega|\ge r_1}\ p_\omega^\lambda(\by)\ G^{\rm atom}_{\lambda,\omega^\prime_\bz}(\by,\bz)\ d\by\\
&\quad \lesssim\ \int_{|\by-\omega|\ge r_1}\ e^{-\cp\lambda}e^{-\cp\lambda|\by-\omega|}\ \left[\ \left(c_0 \left| \log|\bz-\by| \right|\ +\ \lambda^4\ \right)\ {\bf 1}_{|\by-\bz|\le R}\ +\ e^{-c\lambda}\ e^{-c\lambda|\bz-\by|}\ \right]
\ d\by\\
&\quad \lesssim\ e^{-\tc\lambda} e^{-\tc\lambda|\bz-\omega|}\ .\end{align*}
Therefore, the integral expression in the definition of  $\left(\ \widetilde{\Pi}^\lambda_\Gamma K_0^\lambda \right)_{2a}(\bx,\bz)$ satisfies the bound:
{\footnotesize{
\[\ \int\ p_\omega^\lambda(\by)\ \left[\cdots\right](\by,\bz)\ d\by\ \ =\ \int\ p_\omega^\lambda(\by)\ G^{\rm atom}_{\lambda,\omega^\prime_\bz}(\by,\bz)\Theta_{\omega^\prime_\bz}(\bz)\ d\by\ \lesssim\ e^{-\tc\lambda} e^{-\tc\lambda|\bz-\omega|}\ =\ e^{-\tc\lambda}\ e^{-\frac12\tc\lambda|\bz-\omega|}\ e^{-\frac12\tc\lambda|\bz-\omega|}\ .
\]
}}
We next multiply this estimate by $p_\omega^\lambda(\bx)$ and once again use the pointwise
 bound \eqref{p0-bd-cpam}:
\begin{align*}
 p_\omega^\lambda(\bx)\ \int\ p_\omega^\lambda(\by)\ \left[\cdots\right](\by,\bz)\ d\by\ &\lesssim\ 
 \left(\ \lambda^2\ {\bf 1}_{|\bx-\omega|\le R}\ +\ e^{-c\lambda}e^{-c\lambda|\bx-\omega|}\ \right)\ e^{-\tc\lambda}\ e^{-\frac12\tc\lambda|\bz-\omega|}\ e^{-\frac12\tc\lambda|\bz-\omega|}\\
&\lesssim\ e^{-c\lambda} e^{-c\lambda|\bx-\bz|}\ e^{-\frac12\tc\lambda|\bz-\omega|}\ .
\end{align*}
Finally, we multiply the previous bound by $M^{\omega,\omega}=1 +\mathcal{O}(e^{-c\lambda})$ (see \eqref{M-est})
 and sum over all $\omega\in\Gamma$ to obtain:
 \begin{align*}
&\left(\ \widetilde{\Pi}^\lambda_\Gamma K_0^\lambda \right)_{2a}(\bx,\bz)
 =\ \sum_{\omega\in\Gamma}  p_\omega^\lambda(\bx)\ \int\ p_\omega^\lambda(\by)\ \left[\cdots\right](\by,\bz)\ d\by\\
 &\qquad \lesssim\ \left(\ 1 +\mathcal{O}(e^{-c\lambda})\ \right)\ e^{-c\lambda}\ e^{-c\lambda|\bx-\bz|}\
   \sum_{\omega\in\Gamma}\  e^{-\frac12\tc\lambda|\bz-\omega|}\lesssim\ e^{-c\lambda} e^{-c\lambda|\bx-\bz|}\ .
 \end{align*}
 Therefore, the contribution to $\left(\ \widetilde{\Pi}^\lambda_\Gamma K_0^\lambda \right)_{2}(\bx,\bz)$
 from $\left(\ \widetilde{\Pi}^\lambda_\Gamma K_0^\lambda \right)_{2a}(\bx,\bz)$ is an error kernel.
\medskip 

\nit{\it Estimation of $\left(\ \widetilde{\Pi}^\lambda_\Gamma K_0^\lambda \right)_{2b}(\bx,\bz)$; see \eqref{I-IIdef}}:\ 
 From the expression \eqref{I-IIdef} we need only consider $\bz\in\supp(\Theta_{\rm free})$, that is $\bz$ bounded away 
from the all sites $\omega\in \Gamma$; in particular, 
 $|\bz-\omega|\ge r_3$ for all $\omega\in\Gamma$.  By  \eqref{p0-bd-cpam} and \eqref{Gf-allx}: 
\begin{align*}
&p_\omega^\lambda(\by) G_\lambda^{\rm free}(\by-\bz)\ \Theta_{\rm free}(\bz)\\
&\quad \lesssim\ \left(\lambda^2{\bf 1}_{|\by-\omega|\le r_1}+e^{-c\lambda} e^{-c\lambda|\by-\omega|}\right)\cdot
  e^{-c\lambda|\by-\bz|}\cdot\left(1\  +\ \Big|\log\lambda|\by-\bz| \Big|\ \right)\ \Theta_{\rm free}(\bz)\\
  &\quad \lesssim\ \left(\ e^{-\cp\lambda|\by-\bz|}\ {\bf 1}_{|\by-\omega|\le r_1}\ +\ 
    e^{-c\lambda|\by-\bz|}\ e^{-c\lambda|\by-\omega|}\ {\bf 1}_{|\by-\omega|\ge r_1}\ \right)\cdot\left(\ 1\ +\ \Big|\log\lambda|\by-\bz| \Big|\ \right)
    \ \Theta_{\rm free}(\bz)\ .
  \end{align*}
  Integrating over $\R^2$ with respect to $\by$, we find that
  \begin{align*}
&\int_{_{\R^2}}\ p_\omega^\lambda(\by) G_\lambda^{\rm free}(\by-\bz)\  \Theta_{\rm free}(\bz)\ d\by\  \lesssim\ e^{-c\lambda|\bz-\omega|}\ \Theta_{\rm free}(\bz)\ .
\end{align*}
Now multiply this bound by $M^{\omega,\omega}\ p_\omega^\lambda(\bx)$ and apply the pointwise bound for $p_\omega^\lambda(\bx)$, implied by \eqref{p0-bd-cpam}, and the expansion $M^{\omega,\omega}=1+\mathcal{O}(e^{-c\lambda})$ of \eqref{M-est}, to obtain 
\begin{align*}
&M^{\omega,\omega}\ p_\omega^\lambda(\bx)\ \int_{_{\R^2}}\ p_\omega^\lambda(\by) G_\lambda^{\rm free}(\by-\bz)\  \Theta_{\rm free}(\bz)\ d\by\
\\
&\quad  \lesssim\ 
 \left(\lambda^2{\bf 1}_{|\bx-\omega|\le r_1} e^{-\frac14 c\lambda|\bz-\omega|}\ e^{-\frac14 c\lambda|\bz-\omega|}\ +\ e^{-c\lambda} e^{-c\lambda|\bx-\omega|}  e^{-\frac12 c\lambda|\bz-\omega|}\ \right)\
  \Theta_{\rm free}(\bz)\ e^{-\frac12 c\lambda|\bz-\omega|}\\
  &\quad\lesssim\  \left(\ {\bf 1}_{|\bx-\omega|\le r_1}\ e^{-\tc\lambda |\bx-\omega|}\ e^{-\tc\lambda|\bz-\omega|}\
   +\ e^{-c\lambda} e^{-c\lambda|\bx-\omega|}  e^{-\frac12 c\lambda|\bz-\omega|}\ \right)\
  \Theta_{\rm free}(\bz)\ e^{-\frac12 c\lambda|\bz-\omega|}\\
  &\quad\lesssim\  e^{-\cp\lambda |\bx-\bz|}\  \Theta_{\rm free}(\bz)\ e^{-\frac12 c\lambda|\bz-\omega|}\ .
\end{align*}
Summing over $\omega\in\Gamma$ and using that on the support of $\Theta_{\rm free}(\bz)$, $\bz$ is uniformly bounded away from $\Gamma$, we have that 
\begin{equation}
\sum_{\omega\in\Gamma}\ M^{\omega,\omega}\ p_\omega^\lambda(\bx)\ \int_{_{\R^2}}\ p_\omega^\lambda(\by) G_\lambda^{\rm free}(\by-\bz)\  \Theta_{\rm free}(\bz)\ d\by\ \lesssim e^{-c\lambda}\ e^{-c\lambda|\bx-\bz|}\ .
\nn\end{equation}
Hence, the contribution to $ \left(\ \widetilde{\Pi}^\lambda_\Gamma K_0^\lambda \right)_{2}(\bx,\bz)$ of 
$ \left(\ \widetilde{\Pi}^\lambda_\Gamma K_0^\lambda \right)_{2b}(\bx,\bz)$ is also an error kernel. Therefore, 
  $\left(\ \widetilde{\Pi}^\lambda_\Gamma K_0^\lambda \right)_{2}(\bx,\bz)$ is an error kernel, and since
  we have already verfied that  $\left(\ \widetilde{\Pi}^\lambda_\Gamma K_0^\lambda \right)_{1}(\bx,\bz)$ is an error kernel, we conclude that $ \left(\ \widetilde{\Pi}^\lambda_\Gamma K_0^\lambda \right)(\bx,\bz)$ is  an error kernel. Furthermore, it is straightforward to show by arguments similar to those above that
  $\HGamma \left(\ \widetilde{\Pi}^\lambda_\Gamma K_0^\lambda \right)(\bx,\bz)$ is an error kernel, where $\HGamma$ is defined in \eqref{HGamma}. Indeed, we just replace $p^\lambda_\omega(\bx)$ by 
  $\HGamma p^\lambda_\omega(\bx)$ in the previous discussion. Note that 
   $\HGamma p^\lambda_\omega(\bx)=\lambda^2\ \sum_{\omega^\prime\in\Gamma\setminus\{\omega\}}V_0(\bx-\omega^\prime)p^\lambda_0(\bx-\omega)$ and therefore 
   $|\HGamma p^\lambda_\omega(\bx)|\lesssim \lambda^2 \|V_0\|_{_{L^\infty}}\ p_\omega^\lambda(\bx)$. Hence, the estimates lose at worst one power of $\lambda^2$, which can be absorbed by our exponentials $e^{-c\lambda}$.
     This completes the proof of Proposition \ref{HGP-ek}.
     \medskip

 From \eqref{appinv2}, Proposition \ref{K0E0} and Proposition \ref{HGP-ek} we have
 \begin{align}
 f(\bx)\ &=\ \left(\ -\Delta_\bx\ +\ \VGamma(\bx)\ -\ E_0^\lambda\  \right)\ K_1^\lambda[f](\bx)\nn\\
  &\qquad +\ \sum_{\omega\in\Gamma}\ \left\langle \Theta_\omega \ p_\omega^\lambda,\ f\right\rangle_{L^2(\R^2)}\ p_\omega^\lambda(\bx)\  +\ \mathcal{E}^\lambda_1[f](\bx)\ ,
 \label{appinv2a} \end{align}
  where 
 \begin{equation}
 K_1^\lambda\ \equiv  K_0^\lambda\ -\ \widetilde{\Pi}^\lambda_\Gamma\ K_0^\lambda\ =\ \Pi^\lambda_\Gamma\ K_0^\lambda
  \quad \textrm{is a main kernel\ ,}
\label{K1-def-a}  \end{equation}
  \begin{equation}
\left\langle\ p_\omega^\lambda,K_1^\lambda[f]\ \right\rangle\ =\ 0,\quad \textrm{for all}\ \omega\in\Gamma, 
\label{K1-orth}\end{equation}
 and
\begin{equation}
 \mathcal{E}^\lambda_1\ =\ - \mathcal{E}_0^\lambda\ +\ 
 \left(\ -\Delta_\bx\ +\ \VGamma(\bx)\ -\ E_0^\lambda\  \right)\
  \left(\widetilde{\Pi}^\lambda_\Gamma\ K_0^\lambda\right)\ =\ - \mathcal{E}_0^\lambda\ +\ \HGamma\ \left(\widetilde{\Pi}^\lambda_\Gamma\ K_0^\lambda\right).
\label{E1-def}\end{equation}
 is derived from an error kernel.\medskip
  
  Now let $|\Omega|<e^{-\hat{c}\lambda}$, where $\hat{c}$ is a constant that was introduced in Remark \ref{hatc}, and thus 
 $(\rho_\lambda)^{-1}|\Omega|\le e^{-(\hat{c}-c_-)\lambda}\to0$, as  $\lambda\to\infty$.   Then, from \eqref{appinv2a} we have
   \begin{align}
 f(\bx)\ &=\ \left(\ -\Delta_\bx\ +\ \VGamma(\bx)\ -\ E_0^\lambda\ -\ \Omega\ \right)\ K_1^\lambda[f](\bx)\nn\\
  &\qquad +\ \sum_{\omega\in\Gamma}\ \left\langle \Theta_\omega \ p_\omega^\lambda,\ f\right\rangle_{L^2(\R^2)}\ p_\omega^\lambda(\bx)\  +\ \left(\ \mathcal{E}^\lambda_1\ +\ \Omega K_1^\lambda\ \right)[f](\bx)\ 
 \label{appinv2b} \end{align}
 and hence
  \begin{align}
\left(\ I\ -\ (\mathcal{E}_1^\lambda+\Omega K_1^\lambda)\ \right) f(\bx)\ &=\ \left(\ -\Delta_\bx\ +\ \VGamma(\bx)\ -\ E_0^\lambda\ -\ \Omega\ \right)\ K_1^\lambda[f](\bx)\nn\\
  &\qquad +\ \sum_{\omega\in\Gamma}\ \left\langle \Theta_\omega \ p_\omega^\lambda,\ f\right\rangle_{L^2(\R^2)}\ p_\omega^\lambda(\bx)\ .
 \label{appinv2c} \end{align}
 For $\lambda$ large, the operator $\mathcal{E}^\lambda_1\ +\ \Omega K_1^\lambda$ has small norm as a bounded operator on $L^2(\R^2)$. Hence, $I-\left(\mathcal{E}^\lambda_1\ +\ \Omega K_1^\lambda\right)$ is invertible. 
 Applying \eqref{appinv2b} to $\tilde{f}=\left(\ I-(\mathcal{E}_1^\lambda+\Omega K_1^\lambda)\ \right)^{-1}f$
  yields
   \begin{align}
 f(\bx)\ &=\ 
 \left(\ -\Delta_\bx\ +\ \VGamma(\bx)\ -\ E_0^\lambda\ -\ \Omega\ \right)\ 
 \Big(\ K_1^\lambda\left(\ I\ -\ (\mathcal{E}_1^\lambda+\Omega K_1^\lambda)\ \right)^{-1}\Big)[f](\bx)\nn\\
  &\qquad +\ \sum_{\omega\in\Gamma}\ \left\langle \Theta_\omega \ p_\omega^\lambda,\ \tilde{f}\right\rangle_{L^2(\R^2)}\ p_\omega^\lambda(\bx)\ .
 \label{appinv2d} \end{align}
  
 From \eqref{appinv2d} we see that for all $f\in L^2(\R^2)$ and $|\Omega|\lesssim e^{-\hat{c}\lambda}$
 \begin{align}
&\left(\ -\Delta_\bx\ +\ \VGamma(\bx)\ -\ E_0^\lambda\ -\ \Omega\ \right)\ \Big(\ K_1^\lambda\left(\ I\ -\ (\mathcal{E}_1^\lambda+\Omega K_1^\lambda)\ \right)^{-1}\Big)\ f\ =\ f\nn\\ 
&\qquad\qquad \textrm{modulo the  span of}\quad \{p_\omega^\lambda:\omega\in\Gamma\}\ .
\label{K_1-orth} \end{align}
Here, $K_1^\lambda$, defined in \eqref{K1-def},\ is derived from a  main kernel, $\mathcal{E}_1^\lambda$ is derived from an error kernel.

\begin{proposition}\label{E2} For $\lambda$ sufficiently large and $\Omega$ such that $|\Omega|\lesssim e^{-\hat{c}\lambda}$, 
\begin{equation}
K_2^\lambda\ \equiv\ K_1^\lambda\left(\ I\ -\ (\mathcal{E}_1^\lambda+\Omega K_1^\lambda)\ \right)^{-1}
\ \equiv\ K_1^\lambda\ +\ \mathcal{E}_2^\lambda .
\label{K2-def}\end{equation}
Here, $K_1^\lambda$ is derived from a main kernel, $\mathcal{E}_2^\lambda$  from an error kernel and therefore
  $K_2^\lambda$ is derived from a main kernel. 
Moreover, for all $f\in L^2(\R^2)$:
\begin{align}
& \left(\ -\Delta_\bx\ +\ \VGamma(\bx)\ -\ E_0^\lambda\ -\ \Omega\ \right)\ K_2^\lambda\ 
\  f\ =\ f, \nn\\
&\qquad\qquad \textrm{modulo the  span of}\quad \{p_\omega^\lambda:\omega\in\Gamma\}\ ,
\label{K1E2a}\\
&K_2^\lambda[f]\  \perp\ {\rm span}\Big\{p_\omega^\lambda:\omega\in\Gamma\Big\}.
\label{K1E2b}\end{align}
\end{proposition} 
\medskip

\begin{proof}
\nit{\it Proof of Proposition \ref{E2}:} Set $A=\Omega\ K_1^\lambda\ +\ \mathcal{E}_1^\lambda$, where $\lambda$ is taken sufficiently large.
 First note that by Lemma \ref{error-main} that the operator  $A^2$ is derived from an error kernel.
  As an operator on $L^2(\R^2)$ we have
  $\left(I-A\right)^{-1}= \left(I+A\right)\ \left(I-A^2\right)^{-1}=
  \left(I+A \right)\ \left(I+A_1 \right)$, where $A_1$ is an error kernel, again by Lemma \ref{error-main}. Therefore, 
$ (I-A)^{-1}=I\ +\ A\ +\  A_2\ =\ I\ +\ \Omega\ K_1^\lambda\ +\  A_3
$, where $A_j$ ($j=2,3$) arise from error kernels. Another application of Lemma \ref{error-main} completes the proof that 
$\mathcal{E}_2^\lambda$ is derived from an error kernel . That  \eqref{K1E2a}-\eqref{K1E2b} hold follows from  \eqref{K_1-orth}
 and \eqref{K1-orth}.
\end{proof}
\medskip

Recall the  subspace $\mathscr{X}_\Gamma$, the orthogonal complement of 
 ${\rm span}\Big\{p_\omega^\lambda:\omega\in\Gamma\Big\}$:
\begin{align}
\mathscr{X}_\Gamma\ &\equiv\ {\rm span}\Big\{p_\omega^\lambda:\omega\in\Gamma\Big\}^\perp
 =\ \Big\{f\in L^2(\R^2)\ :\ \left\langle p_\omega^\lambda, f\right\rangle_{_{L^2(\R^2)}}=0,\ \ \omega\in\Gamma\ \Big\}\ , 
\label{XGam}
\end{align}
and  the orthogonal projections: $\Pi_\Gamma^\lambda:L^2(\R^2)\to \mathscr{X}_\Gamma$
and $\widetilde{\Pi}_\Gamma^\lambda:L^2(\R^2)\to {\rm span}\Big\{p_\omega^\lambda:\omega\in\Gamma\Big\}$
; see \eqref{XGam1}.\ 
We now write
\[ K_2^\lambda\ =\ K_3^\lambda\ +\ \mathcal{E}_3^\lambda\]
where
\begin{equation}
 K_3^\lambda\ \equiv\ K_2^\lambda \Pi^\lambda_\Gamma, \quad 
\textrm{and}\quad \mathcal{E}_3^\lambda \equiv K_2^\lambda \widetilde{\Pi}^\lambda_\Gamma
\label{K3-def}\end{equation}
Note that 
\[ K_3^\lambda[f]=0\ \textrm{in $L^2(\R^2)$ if}\  f\in{\rm span}\{p_\omega^\lambda\ :\ \omega\in\Gamma\}\ ,\]
and by Proposition \ref{E2}:
 \begin{equation}
 \left(\ -\Delta_\bx\ +\ \VGamma(\bx)\ -\ E_0^\lambda\ -\ \Omega\ \right)\mathcal{E}_3^\lambda\ 
\in\ {\rm span}\Big\{p_\omega^\lambda:\omega\in\Gamma\Big\}\ .
\label{E3-prop}\end{equation} 
Hence,  for all $f\in L^2(\R^2)$:
\begin{align*}
&\left(\ -\Delta_\bx\ +\ \VGamma(\bx)\ -\ E_0^\lambda\ -\ \Omega\ \right)\  
 K_3^\lambda
  f=\ f\ \quad \textrm{modulo the  span of}\quad \{p_\omega^\lambda:\omega\in\Gamma\}.
\end{align*}
We therefore have
 \begin{proposition}\label{K-op} 
 Let $|\Omega|\le e^{-c\lambda}$ with $\lambda$ chosen sufficiently large. 
 Then,  the  operator $\Pi_\Gamma^\lambda(H^\lambda_\Gamma-E_0^\lambda-\Omega)=\Pi_\Gamma^\lambda\left(-\Delta+V_\Gamma^\lambda-E_0^\lambda-\Omega\right)$ is invertible on $\mathscr{X}_\Gamma$, the orthogonal complement of 
 $ {\rm span}\Big\{p_\omega^\lambda:\omega\in\Gamma\Big\}$. Its inverse is given by  $K_3^\lambda\Big|_{\mathscr{X}_\Gamma} $ and we write
 \[
\Res_\Gamma^\lambda(\Omega)\Big|_{\mathscr{X}_\Gamma}\ \equiv\ K_3^\lambda\Big|_{\mathscr{X}_\Gamma} :\ \mathscr{X}_\Gamma\ \to\ \mathscr{X}_\Gamma\ .
 \]
 \end{proposition}

 The following proposition characterizes the operator kernel we seek:
 
 \begin{proposition}\label{summary-K} 
  Let $|\Omega|\lesssim e^{-\hat{c}\lambda}$ with $\lambda$ chosen sufficiently large. Then, $\Res_\Gamma^\lambda(\Omega)$ defined in  Proposition \ref{K-op} satisfies the following properties:
  \begin{enumerate} 
  \item  \begin{equation}
 \mathcal{K}_\Gamma^\lambda(\Omega)[f]=0\ \textrm{in $L^2(\R^2)$ if}\  f\in{\rm span}\{p_\omega^\lambda\ :\ \omega\in\Gamma\}.\label{prop1K}
 \end{equation}
\item \begin{align}
&\mathcal{K}_\Gamma^\lambda(\Omega)[f]\ \perp\ {\rm span}\Big\{p_\omega^\lambda:\omega\in\Gamma\Big\}\  \textrm{in $L^2(\R^2)$}\ .\label{prop2K}
\end{align}
\item
\begin{align}
& \Pi_\Gamma^\lambda\ \left(-\Delta+V_\Gamma^\lambda-E_0^\lambda-\Omega\right)\mathcal{K}_\Gamma^\lambda(\Omega)[f]\ =\ f\ .
\label{smry1}
\end{align}
\item The operator $\Res_\Gamma^\lambda(\Omega)$ is derived from a kernel:
\begin{align}
&\mathcal{K}_\Gamma^\lambda(\Omega)[f](\bx)\ =\ \int_{\R^2}\ \mathcal{K}_\Gamma^\lambda(\bx,\by;\Omega)\ f(\by)\ d\by\quad
 \textrm{for all}\ f\in L^2(\R^2),\ \ \textrm{where}\label{smry2}\\
 & \left|\ \mathcal{K}_\Gamma^\lambda(\bx,\by;\Omega)\ \right|\ \le\ C\Big[ \left|\ \log|\bx-\by|\ \right|\ +\ \lambda^6\ \Big]\ {\bf 1}_{|\bx-\by|\le C}\ +\ e^{-c\lambda}\ e^{-c\lambda|\bx-\by|}\label{smry3}\\
 &\qquad \textrm{for all}\quad \bx, \by\ \in\R^2 .\nn
 \end{align}
%
 \end{enumerate}
\end{proposition}

The only assertion in Proposition \ref{summary-K} that requires proof is part $(4)$. Recall that 
$\mathcal{K}_\Gamma^\lambda(\Omega)=K_3^\lambda=K_2^\lambda \Pi^\lambda_\Gamma=K_2^\lambda\ -\  K_2^\lambda \widetilde{\Pi}^\lambda_\Gamma$. Since $K_2^\lambda$ is derived from a main kernel, it suffices to study
 the kernel of $K_2^\lambda \widetilde{\Pi}^\lambda_\Gamma$. 
 We begin with a bound on the kernel of $\widetilde{\Pi}^\lambda_\Gamma$, which we derive using Lemma \ref{Pi-Gamma}.
 The kernel of  $\widetilde{\Pi}^\lambda_\Gamma$, $K_{\widetilde{\Pi}}^\lambda(\bx,\by)$, is given by (see \eqref{PiTilde}):
 \begin{equation}
 K_{\widetilde{\Pi}}^\lambda(\bx,\by)\ =\ \sum_{\omega,\omega^\prime}\ M^{\omega,\omega^\prime} p^\lambda_\omega(\bx) p^\lambda_\omega(\by),
 \end{equation}
 and we have
 \begin{equation}
 \widetilde{\Pi}^\lambda_\Gamma[g](\bx)\ =\ \int_{\R^2}\ K_{\widetilde{\Pi}}^\lambda(\bx,\by) g(\by)\ d\by .
 \label{KtPi}
 \end{equation}
Our goal is to bound
 \begin{equation}
 \mathcal{K}_\Gamma^\lambda(\bx,\by;\Omega)\ =\ K_2^\lambda(\bx,\by)\ -\ \Big(K_2^\lambda\circ K_{\widetilde{\Pi}}^\lambda\Big)(\bx,\by)\ =\ 
K_2^\lambda(\bx,\by)\ -\ \int_{\R^2} K_2^\lambda(\bx,\bz)\ K_{\widetilde{\Pi}}^\lambda(\bz,\by) d\by
 \label{knl}\end{equation}
 
Note that 
 \begin{align}
 K_{\widetilde{\Pi}}^\lambda(\bx,\by) \ =\ \sum_{\omega} p^\lambda_\omega(\bx) p^\lambda_\omega(\by) 
 \ +\ \sum_{\omega,\omega^\prime} \left[M^{\omega,\omega^\prime}-\delta_{\omega,\omega^\prime}\right]\ p^\lambda_\omega(\bx) p^\lambda_\omega(\by) .
\label{KtPi-re} \end{align}
Recall from \eqref{M-est} that $\Big|M^{\omega,\omega^\prime}-\delta_{\omega,\omega^\prime}\Big|\lesssim e^{-c\lambda}e^{-c\lambda|\omega-\omega^\prime|}$. Also, from the pointwise bounds, \eqref{p0-bd-cpam}, on $p_0^\lambda$ we have:
 \[
 |p_\omega(\bx)|\ \lesssim\  \lambda {\bf 1}_{|\bx-\omega|\le R}+e^{-c\lambda|\bx-\omega|},\quad
\ |p_{\omega^\prime}(\by)|\ \lesssim\  \lambda {\bf 1}_{|\by-\omega|\le R}+e^{-c\lambda|\by-\omega^\prime|} ,
\]
which it follows that
  \begin{align}
  \Big|\ \sum_{\omega}p^\lambda_\omega(\bx) p^\lambda_\omega(\by)\ \Big|\ &\lesssim\ \lambda^2 {\bf 1}_{|\by-\omega|\le 2R}\
   +\ e^{-c^\prime|\bx-\by|},\nn\\
    \Big|\ \sum_{\omega,\omega^\prime}\left[M^{\omega,{\omega^\prime}}-\delta_{\omega,\omega^\prime}\right]\ p^\lambda_\omega(\bx) p^\lambda_{\omega^\prime}(\by)\ \Big|\ &\lesssim\ e^{-c^\prime\lambda}\ \left[\  {\bf 1}_{|\bx-\by|\le 2R}\ +\ 
    e^{-c^\prime\lambda|\bx-\by|}\ \right]
 \nn   \end{align}
 Substitution  into \eqref{KtPi-re}, we obtain
 \begin{equation}
  \Big| K_{\widetilde{\Pi}}^\lambda(\bx,\by) \Big|\ \lesssim\ {\bf 1}_{|\bx-\by|\le 2R}\ \lambda^2\ +\ e^{-c^\prime\lambda|\bx-\by|}\ .
  \label{Ktb-bd}
  \end{equation}
Now since $K_2(\bx,\by;\Omega)$ is a main kernel we have 
   \begin{equation}
  |K_2(\bx,\by;\Omega)|\ \lesssim\ \left[\ \lambda^4\ +\ \Big|\ \log|\bx-\by|\ \Big|\ \right]\ {\bf 1}_{|\bx-\by|\le R}
  \ +\  e^{-c\lambda}\ e^{-c\lambda|\bx-\by|}
  \label{K2bd}
  \end{equation}
   Inserting the bounds \eqref{Ktb-bd} and \eqref{K2bd} into \eqref{knl} we find that $\mathcal{K}_\Gamma^\lambda(\bx,\by;\Omega)$
   satisfies the bound:
   \begin{equation}
  |K_\Gamma(\bx,\by;\Omega)|\ \lesssim\ \left[\ \lambda^6\ +\ \Big|\ \log|\bx-\by|\ \Big|\ \right]\ {\bf 1}_{|\bx-\by|\le 3R}
  \ +\  e^{-c\lambda}\ e^{-c\lambda|\bx-\by|}
  \label{KGamma-bd}
  \end{equation}
  The proof is complete of Proposition \ref{summary-K} is complete.

\subsection{$\mathcal{K}_\Gamma^\lambda(\Omega)$ for the case where $\Gamma$, the set of nuclei, is translation invariant}\label{invGam}

We now suppose that our discrete set of nuclei, $\Gamma$, is translation invariant by a vector $\vtilde_2\in\R^2$. Of course, we have in mind, $\Gamma=\mathbb{H}_\sharp$, the zigzag truncation of $\mathbb{H}$; see \eqref{bbH+}. But our arguments would apply to other {\it rational truncations} of $\mathbb{H}$, for example along an armchair edge.
  For the particular choice 
$\Gamma=\mathbb{H}_\sharp$, we have  $V_\Gamma(\bx) = V_\sharp(\bx)$  and
\[ H_\Gamma^\lambda={H}_\sharp^{\lambda}\ \equiv\ -\Delta + \lambda^2 V_\sharp(\bx) - E_0^\lambda .\]

As commented upon in Remark \ref{trans-invar}, all our constructions of integral operators and kernels respect that translation invariance. Thus, at each stage our integral kernels
 $A(\bx,\by)$ satisfy: $A(\bx+\vtilde_2,\by+\vtilde_2)=A(\bx,\by)$.  It follows that 
 \begin{equation}
 \mathcal{K}_\Gamma^\lambda(\bx+\vtilde_2,\by+\vtilde_2)\ =\ \mathcal{K}_\Gamma^\lambda(\bx,\by)
 \quad \textrm{for all}\ \bx, \by\in\R^2. \label{KG-trin}
 \end{equation}
 
 \bigskip
 
\subsubsection{$ \mathcal{K}_\Gamma^\lambda$ as a bounded operator acting on $L_\kpar^2(\Sigma)$}\label{extendK}

{\ }\medskip

\nit  Let $\Gamma$ be invariant under translation by $\vtilde_2$.
We recall the setting discussed earlier. Associated with this translation invariance is a parallel quasi-momentum, $\kpar\in[0,2\pi)$. 
 We define the cylinder $\Sigma=\R^2/\R\vtilde_2$ and let $\Omega_\Sigma$ denote a fundamental domain for $\Sigma$. The space $L^2(\Sigma)$ consists of functions $f$ such that $f(\bx+\vtilde_2)=f(\bx)$
  for almost all $\bx\in\R^2$ and such that 
  $\|f\|_{L^2(\Sigma)}\equiv \left(\int_{_{\Omega_\Sigma}}|f(\bx)|^2d\bx\right)^{1\over2}<\infty$. 
  The space $L_\kpar^2(\Sigma)$ consists of functions $f$ such that
   $g(\bx)\equiv f(\bx)e^{-i\frac{\kpar}{2\pi}\ktilde_2\cdot\bx}$ satisfies $g(\bx+\vtilde_2)=g(\bx)$ almost everywhere in $\bx$ and 
   $g\in L^2(\Sigma)$. 
   
   We now show that  $ \mathcal{K}_\Gamma^\lambda$ also gives rise to a bounded operator $L_\kpar^2(\Sigma)$.  For any $f\in L_\kpar^2(\Sigma)$, we define 
\begin{equation} \mathcal{K}_\Gamma^\lambda[f](\bx)\ =\ \int_{\R^2}\ \mathcal{K}_\Gamma^\lambda(\bx,\by)\ f(\by)\ d\by.
\label{KonSigma}
\end{equation}
Similarly, $\Pi_\Gamma^\lambda$ may be defined on $L_\kpar^2(\Sigma)$ using Lemma \ref{Pi-Gamma}.

By our bounds on $\mathcal{K}_\Gamma^\lambda(\bx,\by)$, $\mathcal{K}_\Gamma^\lambda[f]$ is well-defined 
for all $f\in L^2_\kpar(\Sigma)$. Using
 \eqref{KG-trin} and our assumption that $f(\bx+\vtilde_2)=e^{i\kpar}f(\bx)$ almost everywhere, we obtain by change of variables:
\begin{align}
   \mathcal{K}_\Gamma^\lambda[f](\bx+\vtilde_2)\ &=\ \int_{\R^2}\ \mathcal{K}_\Gamma^\lambda(\bx+\vtilde_2,\by)\ f(\by)\ d\by\ =\ \int_{\R^2}\ \mathcal{K}_\Gamma^\lambda(\bx+\vtilde_2,\by+\vtilde_2)\ f(\by+\vtilde_2)\ d\by \nn\\
   &= \int_{\R^2}\ \mathcal{K}_\Gamma^\lambda(\bx,\by)\ f(\by+\vtilde_2)\ d\by\ =\ 
e^{i\kpar}\ \int_{\R^2}\ \mathcal{K}_\Gamma^\lambda(\bx,\by)\ f(\by)\ d\by\ =\ 
  e^{i\kpar}\ \mathcal{K}_\Gamma^\lambda[f](\bx).  \end{align}
Hence,  $e^{-i\frac{\ktilde_2\cdot\bx}{2\pi}\kpar}\ \mathcal{K}_\Gamma^\lambda[f](\bx)$ is a function defined on the cylinder $\Sigma$. Similarly, one shows easily that $\Pi_\Gamma^\lambda$ maps $L^2(\Sigma)$ into itself. Furthermore, we have
\begin{align}
 \Big(\ \Pi_\Gamma^\lambda\ \left(\ H^\lambda_\Gamma- E_D^0-\Omega\ \right)\Pi_\Gamma^\lambda\ \Big)\ \circ\ \mathcal{K}_\Gamma^\lambda\ f\ &=\ \Pi_\Gamma^\lambda\ f\ ,\qquad \mathcal{K}_\Gamma^\lambda\ f\ \in L^2_\kpar(\Sigma)
 \label{Klam-Sig}\end{align}
 thanks to Proposition \ref{summary-K}.
 That $ e^{-i\frac{\ktilde_2\cdot\bx}{2\pi}\kpar}\ \mathcal{K}_\Gamma^\lambda\ f\ \in L^2(\Sigma)$ is a consequence of the kernel bounds on 
 $\mathcal{K}_\Gamma^\lambda(\bx,\by)$ and Young's inequality. 
Therefore, we have\medskip

\begin{proposition}\label{KG-ext}
  Let $|\Omega|\le e^{-\hat{c}\lambda}$ with $\lambda$ chosen sufficiently large.  
Let the discrete set $\Gamma$ be invariant under translation by the vector $\vtilde_2$. Then, the
kernel  $\mathcal{K}_\Gamma^\lambda(\Omega)(\bx,\by)$, defined in Proposition \ref{summary-K} and  \eqref{KonSigma}, gives rise to a 
bounded  operator  on $L_\kpar^2(\Sigma)$. Furthermore, the operator 
\begin{equation}
 \mathcal{K}_\Gamma^\lambda(\Omega,\kpar)\ \equiv\ e^{-i\frac{\ktilde_2\cdot\bx}{2\pi}\kpar}\ \mathcal{K}_\Gamma^\lambda(\Omega)\ e^{i\frac{\ktilde_2\cdot\bx}{2\pi}\kpar}
 \label{KOmkp}
 \end{equation}
is  a bounded operator on 
 $L^2(\Sigma)$.
 \end{proposition}
 
\subsubsection{The operator $\mathcal{K}_\Gamma^\lambda(\Omega,\kpar)$ acting on periodized sums}\label{periodized}
{\ }\medskip

\nit Let $\Gamma$ be invariant under translates by integer multiples of $\vtilde_2$. 
We are interested in  $\mathcal{K}_\Gamma^\lambda(\Omega,\kpar): L^2(\Sigma)\to L^2(\Sigma)$ (see
\eqref{KOmkp})  applied to a sum over all $\vtilde_2-$ integer-translates of 
\begin{equation}
p^\lambda_{_{\kpar,\omega}}(\bx)\ =\ e^{i\frac{\kpar}{2\pi}\ktilde_2\cdot(\bx-\omega)}p_0(\bx-\omega)\ .
\label{p-kpar-om}
\end{equation}

For $\omega\in \Gamma$, let $[\omega]$ denote the equivalence class of all translates of $\omega$ by integer multiples of $\vtilde_2$. The set of such equivalence classes is 
\begin{equation}
 \Lambda_\Sigma
\equiv \{[\omega]:\omega\in\Gamma\}
\label{LamSig}\end{equation}
 For any $[\omega]\in\Lambda_\Sigma$ we set
  \begin{equation}
p_{_{\kpar,[\omega]}}^\lambda(\bx)\ \equiv\ \sum_{m\in\Z}\ p_{_{\kpar,\omega}}^\lambda(\bx+m\vtilde_2)\ .
\label{p[omega]}
\end{equation}
Our estimates on $p_\omega^\lambda\in L^2(\R^2)$ imply that $p_{_{\kpar,[\omega]}}^\lambda\in L^2(\Sigma)$, and by our discussion of the previous subsection $\mathcal{K}_\Gamma^\lambda[p_{_{\kpar,[\omega]}}^\lambda]\in L^2(\Sigma)$. Furthermore, we have

\begin{proposition}\label{Kp-orth}
 Let $|\Omega|\le e^{-\hat{c}\lambda}$ with $\lambda$ chosen sufficiently large. 
\begin{enumerate}
\item $\mathcal{K}_\Gamma^\lambda(\Omega,\kpar)[f]=0$ in $L^2(\Sigma)$ for all $f\in{\rm span}\{p_{_{\kpar,[\omega]}}^\lambda\ :\
\omega\in\Gamma\}$.
\item For all $\omega\in\Gamma$ and $f\in L^2(\Sigma)$, we have $\left\langle\ \mathcal{K}_\Gamma^\lambda(\Omega,\kpar)[f],p_{_{\kpar,[\omega]}}^\lambda\ \right\rangle_{L^2(\Sigma)}=0$.
\end{enumerate}
\end{proposition}
\medskip

\nit{\it Proof of claim $(1)$ of Proposition \ref{Kp-orth}:}\  We claim in fact for any $\omega\in\Gamma$, and for any $\bx\in\R^2$, we have $\mathcal{K}_\Gamma^\lambda(\Omega,\kpar)[p_{_{\kpar,[\omega]}}^\lambda](\bx)=0$. Indeed,
\begin{align}
\mathcal{K}_\Gamma^\lambda(\Omega,\kpar)[p_{_{[\kpar,\omega]}}^\lambda](\bx)\
&=\ \int_{\R^2} \mathcal{K}_\Gamma^\lambda(\Omega,\kpar)(\bx,\by)\sum_{m\in\Z} p^\lambda_{_{\kpar,\omega-m\vtilde_2}}(\by)\ d\by\nn\\
&=\ \lim_{N\to\infty}\ \int_{\R^2} \mathcal{K}_\Gamma^\lambda(\Omega,\kpar)(\bx,\by)\sum_{|m|\le N} p^\lambda_{_{\kpar,\omega-m\vtilde_2}}(\by)\ d\by\nn\\
&=\ \lim_{N\to\infty}\ \sum_{|m|\le N} \ \int_{\R^2} \mathcal{K}_\Gamma^\lambda(\Omega,\kpar)(\bx,\by)p^\lambda_{_{\kpar,\omega-m\vtilde_2}}(\by)\ d\by\nn\\
&=\ \lim_{N\to\infty}\ \sum_{|m|\le N} \ \mathcal{K}_\Gamma^\lambda(\Omega,\kpar)[p_{_{\kpar,\omega-m\vtilde_2}}^\lambda](\bx)\ =\ 0,\nn
\end{align}
by property \eqref{prop1K} of Proposition \ref{summary-K}. 
These formal manipulations are easily justified
thanks to our estimates on $\mathcal{K}_\Gamma^\lambda(\Omega,\kpar)(\bx,\by)$ and $p^\lambda_\omega(\bx)$.
 This completes the proof of the first claim of  Proposition \ref{Kp-orth}.
 \medskip
 
 \nit \nit{\it Proof of claim $(2)$ of Proposition \ref{Kp-orth}:}\  Let $\omega\in\Gamma$ and $f\in L^2(\Sigma)$. Then, 
 \begin{align}
& \left\langle\ \mathcal{K}_\Gamma^\lambda(\Omega,\kpar)[f],p_{_{\kpar,[\omega]}}^\lambda\ \right\rangle_{L^2(\Sigma)}
 \ =\ \int_{\bx\in \Omega_\Sigma}\  \sum_{m\in\Z}\ p^\lambda_{_{\kpar,\omega}}(\bx+m\vtilde_2)\cdot 
 \int_{\by\in\R^2}\ \mathcal{K}_\Gamma^\lambda(\Omega,\kpar)(\bx,\by)\ f(\by)\ d\by\  d\bx\nn\\
 \ &=\ \sum_{m\in\Z}\ \int_{\bx\in \Omega_\Sigma}\ p^\lambda_{_{\kpar,\omega}}(\bx+m\vtilde_2)\cdot 
 \int_{\by\in \R^2}\ \mathcal{K}_\Gamma^\lambda(\Omega,\kpar)(\bx,\by)\ f(\by)\ d\by\  d\bx\nn\\
 &=\ \sum_{m\in\Z}\ \int_{\bx\in \Omega_\Sigma}\ p^\lambda_{_{\kpar,\omega}}(\bx+m\vtilde_2)\cdot 
 \int_{\by\in \R^2}\ \mathcal{K}_\Gamma^\lambda(\Omega,\kpar)(\bx+m\vtilde_2,\by+m\vtilde_2)\ f(\by+m\vtilde_2)\ d\by\  d\bx
 \nn
 \end{align}
 The latter equality holds by properties of  $\mathcal{K}_\Gamma^\lambda(\Omega,\kpar)$ and $f$ under translation by $\vtilde_2$. Continuing, we have
 \begin{align}
 &\left\langle\ \mathcal{K}_\Gamma^\lambda(\Omega,\kpar)[f],p_{_{\kpar,[\omega]}}^\lambda\ \right\rangle_{L^2(\Sigma)}
 \ =\ \sum_{m\in\Z}\ \int_{\Omega_\Sigma+m\vtilde_2}\ p^\lambda_{_{\kpar,\omega}}(\bx^\prime)\ 
 \int_{\by^\prime\in\R^2}\ \mathcal{K}_\Gamma^\lambda(\Omega,\kpar)(\bx^\prime,\by^\prime)\ 
  f(\by^\prime)\ d\by^\prime\ d\bx^\prime \nn\\
 &=\  \int_{\bx\in\R^2}\ p^\lambda_{_{\kpar,\omega}}(\bx)\ \int_{\by\in\R^2}\ \mathcal{K}_\Gamma^\lambda(\Omega,\kpar)(\bx,\by)\ f(\by)\ d\by\ d\bx\nn\\
 &\ =\ \lim_{N\to\infty}\ 
 \int_{\bx\in\R^2}\ p^\lambda_{_{\kpar,\omega}}(\bx)\ \int_{|\by|\le N}\ \mathcal{K}_\Gamma^\lambda(\Omega,\kpar)(\bx,\by)\ f(\by)\ d\by\ d\bx\ =\ 0
\nn \end{align}
by property \eqref{prop2K} of Proposition \ref{summary-K}. Again, the formal manipulations are easily justified. This completes the proof of Proposition \ref{Kp-orth}.
\bigskip

\subsection{Green's kernel}\label{zz-ac-ker}
{\ }\medskip

We recall the cylinder $\Sigma=\R^2/\R\vtilde_2$ and the choice of fundamental domain $\Omega_\Sigma\subset\R^2$, given as the union of finite parallelograms, $\Omega_n,\ n\ge0$ together with one unbounded parallelogram, $\Omega_{-1}$,\ $\Omega_\Sigma=\cup_{n\ge0}\Omega_n\ \cup\ \Omega_{-1}$; see 
\eqref{OmSig}. In each finite parallelogram, $\Omega_n,\ n\ge0$, are two lattice points of $\mathbb{H}_\sharp$:\ $\bv_\bA^{(n)}$ and $\bv_\bB^{(n)}$ .
As our discrete set  we take $\Gamma=\mathbb{H}_\sharp$, our potential $V_\sharp(\bx)$ and 
our Hamiltonian $\cHeg$ acting on $L^2_\kpar(\Sigma)$.
  
  Next recall the subspace of $L^2(\Sigma)$
(see \eqref{XAB-def}):
\begin{align} 
\mathscr{X}^\lambda_{AB}(\kpar)\ &=\ \textrm{ orthogonal complement in $L^2(\Sigma)$ of span}\Big\{p_{_{\kpar,\bI}}^{\lambda}[n]:n\ge0,\ I=A, B\Big\}\nn
\end{align}
with orthogonal projection:
   \begin{equation} \Pi^\lambda_{_{AB}}(\kpar): L^2(\Sigma)\to \mathscr{X}^\lambda_{AB}(\kpar)\ . \nn
   \end{equation}  
By definition
\[ p_{_{\kpar,[\bv_\bI^{(n)}]}}^\lambda(\bx)\ =\  p^\lambda_{_{\kpar,\bI}}[n](\bx),\quad \bI=A, B , \]
where $p^\lambda_{_{\kpar,\bI}}[n]$ is defined in \eqref{pkpar}. 
\bigskip

%
\bigskip
Recall that $\Ressp^{\lambda}(\Omega,\kpar)$, the inverse of $\Pi_{_{AB}}(\kpar)\left(\ \cHeg(\kpar)-\Omega\ \right)\Pi_{_{AB}}(\kpar)$ (\ equivalently
 $ \Pi_{_{AB}}(\kpar)\circ$ $\left(\ \cHeg(\kpar)-\Omega\ \right)\ $) acting on $\mathscr{X}^\lambda_{_{AB}}(\kpar)$; see  Proposition \ref{resolvent}.
By  Propositions  \ref{summary-K} and \ref{KG-ext} this inverse is given by an integral operator 
\begin{equation} f\mapsto \mathcal{K}^\lambda_\sharp(\Omega,\kpar)[f]\ \equiv\
\int_{\R^2}\ \mathcal{K}^\lambda_\sharp(\bx,\by;\Omega,\kpar)\ f(\by)\ d\by\ ,
\label{k-sharp-int}
\end{equation}
with kernel 
\begin{equation}
\mathcal{K}^\lambda_\sharp(\bx,\by;\Omega,\kpar)= e^{-i\frac{\ktilde_2\cdot\bx}{2\pi}\kpar}\ \mathcal{K}_\sharp^\lambda(\bx,\by,\Omega)\ e^{i\frac{\ktilde_2\cdot\by}{2\pi}\kpar}\ 
\label{Kshp-ker}
\end{equation}
which satisfies the pointwise bounds: 
\begin{align}
&\left|\ \mathcal{K}_\sharp^\lambda(\bx,\by;\Omega,\kpar)\ \right|\ \le\ C\Big[ \left|\ \log|\bx-\by|\ \right|\ +\ \lambda^6\ \Big]\ {\bf 1}_{|\bx-\by|\le C}\ +\ e^{-c\lambda}\ e^{-c\lambda|\bx-\by|}
\label{K-ker-est} \\
 &\textrm{for all}\quad \bx, \by\ \in\R^2.
\nn\end{align}

 Now applying  Proposition \ref{Kp-orth} we obtain:
 
 \begin{proposition}\label{Kp-orthX}
 {\ }
  Let $|\Omega|\le e^{-\hat{c}\lambda}$ with $\lambda$ chosen sufficiently large.  
\begin{enumerate}
\item $\mathcal{K}_\sharp^\lambda(\Omega,\kpar)[f]=0$ in $L^2(\Sigma)$ for all 
$f\in{\rm span}\{p^\lambda_{_{\kpar,\bI}}[n]:\bI=A,B,\ n\ge0\}$.
\item Assume $f\in L^2(\Sigma)$. Then, for all $n\ge0$ and $\bI=A, B$, we have
 \[
 \left\langle\ \mathcal{K}_\sharp^\lambda(\Omega,\kpar)[f],p^{^{\lambda}}_{_{\kpar,I}}[n]\ \right\rangle_{L^2(\Sigma)}=0.
 \]
\item $[H^\lambda(\kpar)-\Omega]\mathcal{K}_\sharp^\lambda(\Omega,\kpar)[f]=f$ 
modulo ${\rm span}\{p^\lambda_{_{\kpar,\bI}}[n]:\bI=A,B,\ n\ge0\}$.
\end{enumerate}
\end{proposition}

A  consequence of the forgoing discussion is:
\begin{corollary}\label{corK}
 Let $|\Omega|\le e^{-\hat{c}\lambda}$ with $\lambda$ chosen sufficiently large.   The operator
$\mathcal{K}_\sharp^\lambda(\Omega,\kpar)$, the inverse of 
$\Pi_{AB}^\lambda(\kpar)\ \left(\ \cHeg(\kpar)-\Omega\ \right)\ 
\Pi_{AB}^\lambda(\kpar)$, arises from a kernel satisfying \eqref{k-sharp-int}-\eqref{K-ker-est}. $\mathcal{K}_\sharp^\lambda(\Omega,\kpar)$ is a bounded linear operator on $L^2(\Sigma)$.
 \end{corollary}

\section{Expansion and estimation of linear matrix elements:\\ proof of   Proposition \ref{zz-els}} \label{m-els-est}

Our first step in the proof of Proposition \ref{zz-els}  is to expand the inner products:
\[
\left\langle P_{_{\kpar,I}}^{^{\lambda}}[m],\cHeg\ P_{_{\kpar,J}}^{^{\lambda}}[n]\right\rangle_{_{L^2(\Sigma)}}
\ =\ \left\langle p_{_{\kpar,\bI}}^{^{\lambda}}[m],\cHeg(\kpar)\ p_{_{\kpar,\bJ}}^{^{\lambda}}[n]\right\rangle_{_{L^2(\Sigma)}},\] 
where $m,n\in \N_0$, in terms of overlap integrals of translates of the atomic potential, $V_0$, and the atomic ground state, $p_0^\lambda$. 
We have, by the definition of the $L^2(\Sigma)$ inner product:
\begin{align*}
\left\langle P_{_{\kpar,I}}^{^{\lambda}}[m],\cHeg\ P_{_{\kpar,J}}^{^{\lambda}}[n]\right\rangle_{_{L^2(\Sigma)}}
\ =\ \int_{_{\Omega_\Sigma}}\ \overline{P_{_{\kpar,I}}^{^{\lambda}}[m](\bx)}\ 
\cHeg\ P_{_{\kpar,J}}^{^{\lambda}}[n](\bx)\ d\bx\ .
\end{align*}

We first simplify the integrand: $\overline{ P_{_{\kpar,I}}^{^{\lambda}}[m]}\ \cHeg\ P_{_{\kpar,J}}^{^{\lambda}}[n]$.
 We recall the definition of $\cHeg$ (see \eqref{zz-evp}) and introduce the notation: 
\begin{equation} \textrm{
$J^\prime=A$ if $J=B$ and $J^\prime=B$ if $J=A$.}
\label{Jprime}\end{equation}
 For  $\bx\in\Omega_\Sigma$, the fundamental domain (see 
Figure \ref{zz-dimers}), we have for $J=A, B$:
{\footnotesize{
\begin{align}
\cHeg P_{_{\kpar,J}}^{^{\lambda}}[n](\bx)\ &=\ \sum_{\tm_2\in\Z}\ e^{i\tm_2\kpar}\nn\\
&\qquad\cdot \left[-\Delta+\lambda^2\sum_{n_1\ge0}V_0(\bx-\vtilde_J^{n_1})+\lambda^2\sum_{n_1\ge0}V_0(\bx-\vtilde_{J^\prime}^{n_1})\ -\ E_0^\lambda\right]\ p_0^\lambda(\bx-\vtilde_J^n-\tm_2\vtilde_2)\nn\\
&=\ \lambda^2\left[\ \sum_{\substack{n_1\ge0\\ n_1\ne n}}V_0(\bx-\vtilde_J^{n_1})+ \sum_{n_1\ge0}V_0(\bx-\vtilde_{J^\prime}^{n_1})\ \right]\ \cdot \Big[\ \sum_{\tm_2\in\Z}e^{i\tm_2\kpar}\ p_0^\lambda(\bx-\vtilde_J^n-\tm_2\vtilde_2)\ \Big]\nn\\
&\qquad\ +\ \lambda^2\ V_0(\bx-\vtilde_J^n)\ \sum_{\tm_2\in\Z\setminus\{0\}}\ e^{i\tm_2\kpar}\ p_0^\lambda(\bx-\vtilde_J^n-\tm_2\vtilde_2)\ .
\label{HsharpPJ}\end{align}
}}
To obtain \eqref{HsharpPJ} we 
use that $(-\Delta_\bx + \lambda^2 V_0(\bx)-E_0^\lambda)p_0^\lambda(\bx)=0$ and therefore
$
(-\Delta_\bx + \lambda^2 V_0(\bx-\bv)-E_0^\lambda)p_0^\lambda(\bx-\bv)=0\ 
$
for all $\bv\in\mathbb{H}$.  From \eqref{HsharpPJ} we obtain:
{\footnotesize{
\begin{align*}
& \overline{P_{_{\kpar,I}}^{^{\lambda}}[m](\bx)}\ \cHeg\ P_{_{\kpar,J}}^{^{\lambda}}[n](\bx)\nn\\
  \ &=\ \sum_{m_2\in\Z}\ \sum_{\tm_2\in\Z}\ e^{i(\tm_2-m_2)\kpar}\ \nn\\
  &\qquad\cdot p_0^\lambda(\bx-\vtilde_I^m-m_2\vtilde_2)\ \left[\ \sum_{\substack{n_1\ge0\\ n_1\ne n}}\ \lambda^2\ V_0(\bx-\vtilde_J^{n_1})+ \sum_{n_1\ge0}\ \lambda^2\ V_0(\bx-\vtilde_{J^\prime}^{n_1})\ \right]\ 
  p_0^\lambda(\bx-\vtilde_J^n-\tm_2\vtilde_2)\nn\\
  &\qquad\ +\  \sum_{m_2\in\Z}\ \sum_{\tm_2\in\Z\setminus\{0\}}\ e^{i(\tm_2-m_2)\kpar}\
   p_0^\lambda(\bx-\vtilde_I^m-m_2\vtilde_2)\ \lambda^2\ V_0(\bx-\vtilde_J^n)\ p_0^\lambda(\bx-\vtilde_J^n-\tm_2\vtilde_2)\ ,
\end{align*}
}}
for all $\bx\in\Omega_\Sigma$. 
Integrating the previous identity over $\Omega_\Sigma$, we obtain:
{\footnotesize{
\begin{align}
&\left\langle\ P_{_{\kpar,I}}^{^{\lambda}}[m](\bx)\ , \cHeg\ P_{_{\kpar,J}}^{^{\lambda}}[n]\ \right\rangle_{L^2(\Sigma)}\nn\\
 &\qquad =\ \sum_{m_2,\tm_2\in\Z}\  \sum_{\substack{n_1\ge0\\ n_1\ne n}}\ e^{i(\tm_2-m_2)\kpar}\
  \int_{_{\Omega_\Sigma}}\ p_0^\lambda(\bx-\vtilde^m_I-m_2\vtilde_2)\ \lambda^2\ V_0(\bx-\vtilde_J^{n_1})\ p_0^\lambda(\bx-\vtilde_J^n-\tm_2\vtilde_2)\ d\bx\nn\\
  &\qquad\ +\ \sum_{m_2,\tm_2\in\Z}\  \sum_{n_1\ge0}\ e^{i(\tm_2-m_2)\kpar}\
  \int_{_{\Omega_\Sigma}}\ p_0^\lambda(\bx-\vtilde^m_I-m_2\vtilde_2)\ \lambda^2\ V_0(\bx-\vtilde_{J^\prime}^{n_1})\ p_0^\lambda(\bx-\vtilde_J^n-\tm_2\vtilde_2)\ d\bx\nn\\
  &\qquad\ +\ \sum_{m_2\in\Z}\ \sum_{\tm_2\in\Z\setminus\{0\}}\ e^{i(\tm_2-m_2)\kpar}\
  \int_{_{\Omega_\Sigma}}\ p_0^\lambda(\bx-\vtilde^m_I-m_2\vtilde_2)\ \lambda^2\ V_0(\bx-\vtilde_J^n)\ p_0^\lambda(\bx-\vtilde_J^n-\tm_2\vtilde_2)\ d\bx\nn\\
&\equiv\  S_1^{IJ}(m,n)\ +\ S_2^{IJ}(m,n)\ +\ S_3^{IJ}(m,n)\ ,
\label{SipIJ}\end{align}
}}
where the three expressions $S_1^{IJ}(m,n)$, $ S_2^{IJ}(m,n)$, and $ S_3^{IJ}(m,n)$ denote the three sums in 
\eqref{SipIJ}. The dependence on $\lambda$ and $\kpar$ has been suppressed. We recall that in the expression for 
 $S_2^{IJ}(m,n)$, the index $J^\prime$ is defined in \eqref{Jprime}.\medskip
 
We now provide a general lemma, which will facilitate our determination of the leading terms and estimation of the error terms in the above sums. In preparation for the statement of this lemma we introduce some terminology.

\begin{definition}\label{lem-prep}
\begin{enumerate}
 \item For  $\Ione, \Jone\in\{A, B\}$, we write $\vtilde_{\Ione}-\vtilde_{\Jone}=\sigma (\vtilde_B-\vtilde_A)=\sigma \be$,
where $\sigma=1$ if $\Ione=B$ and $\Jone=A$, and $\sigma=-1$ if $\Ione=A$ and $\Jone=B$. We therefore write:
\begin{equation}
\sigma(B,A)=+1,\ \ \sigma(A,B)=-1,\ \ \textrm{and we define}\ \ \sigma(\Ione,\Ione)=0.
\label{sigmaIJ}\end{equation}
\item For $\sigma=+1,-1,0$ we define $N_b(\sigma)=\{\br=(r_1,r_2)\in\Z^2:|\sigma\be+\br\vec\vtilde|=|\be|\}$. Therefore
$N_b(+1)\ \equiv\ \{(0,0),(-1,0),(0,-1)\}$,\ 
$N_b(-1)\ \equiv\ \{(0,0),(1,0),(0,1)\}$, and $N_b(0)\ \equiv\ \emptyset$.
 \end{enumerate}
 \medskip
 
 Note that if $\bfm=(m_1,m_2)\in N_b(\sigma)$ with $\sigma=\pm1$, then there exists $l\in\{0,1,2\}$ such that 
 \begin{equation}
  \sigma\be\ +\ m_1\vtilde_1\ +\ m_2\vtilde_2\ =\ R^l\be
  \label{rotate}
  \end{equation}
 where $R$ denotes the $2\times2$ rotation in $\R^2$ by $2\pi/3$. 
 \end{definition}
 
 \begin{lemma}\label{OIsharp}
For $\Ione, \Jone, \tIone\in\{A,B\}$, $m, n, n_1\ge0$ and $m_2, \tm_2\in\Z$, consider the overlap integral
 \begin{equation}
\mathscr{I}_\sharp \equiv \int\ p_0^\lambda(\bx-\vtilde^m_{\Ione}-m_2\vtilde_2)\ \lambda^2\ |V_0(\bx-\vtilde_{\Jone}^{n_1})|\ p_0^\lambda(\bx-\vtilde_{\tIone}^n-\tm_2\vtilde_2)\ d\bx.
 \label{ol-sharp}
 \end{equation}
 Recall the hopping coefficient defined by: $\rho_\lambda=\int p^\lambda_0(\by)\ \lambda^2\ |V_0(\by)|\ p^\lambda_0(\by-\be)\ d\by$.
  Then we have the bound
 \begin{equation}
 \mathscr{I}_\sharp \ \lesssim\ e^{-c\lambda\left(\ |m-n_1|\ +\ |m_2|\ +\ |n-n_1|\ +\ |\tm_2|\ \right)}\ \rho_\lambda,
 \label{Isharp-est}
 \end{equation}
except in the following  cases of \underline{exceptional  indices} $(m,n,n_1,m_2,\tm_2)$:
 \begin{enumerate}
 \item[(a)]\ $\Ione=\tIone=\Jone$, $m=n=n_1$ and $m_2=\tm_2=0$. This case does not arise in the proof of Proposition \ref{zz-els} so we say nothing further about it\ .\\
 \item[(b)]\ $\tIone=\Jone$, $\Ione\ne\Jone$, $(m-n_1,m_2)\in N_b\left(\sigma(\Ione,\Jone)\right)$, $n=n_1$ and $\tm_2=0$,\\ in which case $ \mathscr{I}_\sharp=\rho_\lambda$\ .\\
 \item[(c)]\ $\Ione=\Jone$, $\tIone\ne\Jone$, $(n-n_1,\tm_2)\in N_b\left(\sigma(\tIone,\Jone)\right)$, $m=n_1$ and $m_2=0$,\\
in which case $ \mathscr{I}_\sharp=\rho_\lambda$\ .\\
 \end{enumerate}
 Furthermore, if $\Ione\ne\Jone$, $\tIone\ne\Jone$, then for all $m, n, n_1, m_2, \tm_2$:
 \begin{equation}
 \mathscr{I}_\sharp \ \lesssim\ e^{-c\lambda}\ e^{-c\lambda\left(\ |m-n_1|\ +\ |m_2|\ +\ |n-n_1|\ +\ |\tm_2|\ \right)}\ \rho_\lambda.
 \label{Isharp-est1}
 \end{equation}

 \end{lemma}
\nit  Lemma \ref{OIsharp} is proved in Appendix \ref{overlap}. It makes repeated use of the following pointwise decay estimates for the atomic ground state, $p_0^\lambda$:
\begin{lemma}[See Lemma 15.6 of \cite{FLW-CPAM:17}]\label{p0-bounds}
 There exists a constant $c$ such that for $\by\in {\rm supp}(V_0) \subset B_{r_0}({\bf 0})$, {\it i.e.} $|\by|\le r_0$, we have:
  \begin{align}
 p_0^\lambda(\by-\bn\vec\bv)\  &\lesssim\ e^{-c|\bn|\lambda}\ p_0^\lambda(\by)  ,\quad \bn\in\Z^2,   \label{p0-bound1} \\
p^\lambda_0\left(\by-(\sigma \be+\bn\vec\bv) \right)\ &\lesssim\ 
  e^{-c|\bn|\lambda}\ p^\lambda_0\left(\by-\sigma\be\right),\ \ \bn \notin N_{b}(\sigma),\ \sigma=\pm1, 
  \label{p0-bound2}\\
 p_0^\lambda(\by-\sigma\be)\ &\lesssim\ e^{-c\lambda}\ p_0^\lambda(\by) , \ \sigma=\pm1, \quad {\rm and}
 \label{p0-bound3}\\
 p_0^\lambda(\by-\bn\vec\bv)\ &\lesssim\ e^{-c\lambda|\bn|}\ p_0^\lambda(\by-\sigma\be) ,\ \bn\in\Z^2\setminus\{(0,0)\}.
 \label{p0-bound4}
 \end{align}
  \end{lemma}
\begin{remark}\label{ptwise}
  In  \cite{FLW-CPAM:17},  Lemma \ref{p0-bounds} was proved for all $r_0$ satisfying $0<r_0<r_{\rm critical}$,
  where $0.33|\be|\le r_{critical}<0.5|\be|$, and $|\be|=|\bv_B-\bv_A|=1/\sqrt3$.
  \end{remark}

To prove Proposition \ref{zz-els}, we now apply Lemma \ref{OIsharp} to the expansion of  the matrix elements: 
 $\left\langle\ P_{_{\kpar,I}}^{^{\lambda}}[m](\bx)\ , \cHeg\ P_{_{\kpar,J}}^{^{\lambda}}[n]\ \right\rangle_{L^2(\Sigma)}$, where $I, J=A, B$ and $m, n\in\N_0$,     for large $\lambda$. 
 
 \subsection{Expansion of the inner product
 $\left\langle\ P_{_{\kpar,B}}^{^{\lambda}}[m](\bx)\ , \cHeg\ P_{_{\kpar,A}}^{^{\lambda}}[n]\ \right\rangle_{L^2(\Sigma)}$}
 {\ }
 
 \nit We consider the summations $S_j^{IJ}(m,n), j=1,2,3$ in order (see \eqref{SipIJ}) with $I=B$ and $J=A$.\medskip
 
 \nit\underline{Estimation of $S_1^{BA}(m,n)$:} The expression to be summed over $m_2, \tm_2\in\Z$ and $n_1\ge0, n_1\ne n$ is: 
 \begin{equation}
 e^{i(\tm_2-m_2)\kpar}\ \int_{\R^2}\ p_0^\lambda(\bx-\vtilde^m_B-m_2\vtilde_2)\ \lambda^2\ V_0(\bx-\vtilde_A^{n_1})\ p_0^\lambda(\bx-\vtilde_A^n-\tm_2\vtilde_2)\ d\bx.
 \label{sum1-BA}
 \end{equation}
We apply Lemma \ref{OIsharp} with $\Ione=B$, $\Jone=A$ and $\tIone=A$. All summands \eqref{sum1-BA} of $S_1^{BA}(m,n)$, except for exceptional indices in case (b), defined by $\tIone=\Jone$, $\Ione\ne\Jone$, are bounded by $e^{-c\lambda(|m-n_1|+|m_2|+|n-n_1|+|\tm_2|)}\rho_\lambda$. The exceptional indices are characterized by the relations: $(m-n_1,m_2)\in N_b(\sigma(B,A))=N_b(+1)$, $n=n_1$ and $\tm_2=0$. Since the sum in the definition of $S_1^{BA}(m,n)$ is over $n_1\ge0$ with $n_1\ne n$, there are no relevant exceptional indices and we conclude for all $m, n\ge0$:
 \begin{equation}
 |S_1^{BA}(m,n)| \lesssim \rho_\lambda\  \sum_{m_2,\tm_2\in\Z}\  \sum_{\substack{n_1\ge0\\ n_1\ne n}}\ e^{-c\lambda\left(\ |m-n_1|\ +\ |m_2|\ +\ |n-n_1|\ +\ |\tm_2|\ \right)}\ \lesssim\ e^{-c^\prime\lambda}\ e^{-c^\prime\lambda|m-n|}\ \rho_\lambda\ , 
 \label{S1bound}\end{equation}
 for some strictly positive constant $c^\prime$. \medskip

\nit\underline{Expansion of $S_2^{BA}(m,n)$:} Since $I=B$, $J=A$ and $J^\prime=B$, the expression to be summed over $m_2, \tm_2\in\Z$ and $n_1\ge0$ is: 
 \begin{equation}
 e^{i(\tm_2-m_2)\kpar}\ \int_{\R^2}\ p_0^\lambda(\bx-\vtilde^m_B-m_2\vtilde_2)\ \lambda^2\ V_0(\bx-\vtilde_B^{n_1})\ p_0^\lambda(\bx-\vtilde_A^n-\tm_2\vtilde_2)\ d\bx.
 \label{sum2-BA}
 \end{equation}
 We apply Lemma \ref{OIsharp} with $\Ione=B$, $\Jone=B$ and $\tIone=A$.  
 All summands \eqref{sum2-BA} of $S_2^{BA}(m,n)$, except for exceptional indices in case (c), defined by  $\Ione=\Jone$ and  $\tIone\ne\Jone$,  are bounded by $e^{-c\lambda(|m-n_1|+|m_2|+|n-n_1|+|\tm_2|)}\rho_\lambda$. 
  The exceptional indices are characterized by the relations: $(n-n_1,\tm_2)\in N_b(\sigma(\tIone,\Jone))=N_b(\sigma(A,B))=N_b(-1)=\{(0,0),(1,0),(0,1)\}$, $m=n_1$ and $m_2=0$. We next simplify the expression \eqref{sum2-BA} in each of these three exceptional cases.
\medskip

\nit \underline{$(n-n_1,\tm_2)=(0,0)$,  $m=n_1$, $m_2=0$}: We have $n_1=m=n$ and $m_2=\tm_2=0$. 
For this case, the expression in \eqref{sum2-BA} is equal to $-\rho_\lambda$ and contributes to $S_2^{BA}(m,m)$. \medskip
 
 \nit \underline{$(n-n_1,\tm_2)=(0,1)$,  $m=n_1$, $m_2=0$}: We have $n_1=n=m$, $m_2=0$ and $\tm_2=1$.
 For this case, the expression \eqref{sum2-BA} is equal to
 $-e^{i\kpar}\ \rho_\lambda$
 and contributes to $S_2^{BA}(m,m)$.\medskip
 
  \nit \underline{$(n-n_1,\tm_2)=(1,0)$,  $m=n_1$, $m_2=0$}: We have $n_1=m$, $n=m+1$, $m_2=\tm_2=0$. 
  For this case, the expression in \eqref{sum2-BA}  is equal to  $-\rho_\lambda$
 and contributes to $S_2^{BA}(m,m+1)$.\medskip
 
 We conclude from the above  discussion of $S_2^{BA}(m,n)$  that:
 \begin{align}
 S_2^{BA}(m,m)\ &=\ -\left(1+e^{i\kpar}\right)\rho_\lambda\ +\ \mathcal{O}\left(\ e^{-c\lambda}\  \rho_\lambda \right),\ (n=m) \label{S2mm}\\
 S_2^{BA}(m,m+1)\ &=\ -\rho_\lambda\ +\ \mathcal{O}\left(\ e^{-c\lambda}\  \rho_\lambda\ \right),\ (n=m+1)\label{S2mm+1}\\
 S_2^{BA}(m,n)\ &=\ \mathcal{O}\left(e^{-c\lambda}\ e^{-c\lambda|m-n|}\ \rho_\lambda\right)\ ,\qquad \textrm{if}\quad n\ne m, m+1\ .
\label{S2mn} \end{align}
 The $\mathcal{O}(\cdot)$ error terms are bounds on contributions to $S_2^{BA}(m,n)$ arising from the summation over $m_2, \tm_2\in\Z$ and $n_1\ge0$ of the bound $e^{-c\lambda(|m-n_1|+|m_2|+|n-n_1|+|\tm_2|)}\rho_\lambda$ for 
  non-exceptional indices (as in  \eqref{S1bound}).
 \medskip
 
\nit\underline{Expansion of $S_3^{BA}(m,n)$:} Since $I=B$ and $J=A$, the expression to be summed over $m_2\in\Z$ and $\tm_2\in\Z\setminus\{0\}$ is: 
 \begin{equation}
 e^{i(\tm_2-m_2)\kpar}\ \int_{\R^2}\ p_0^\lambda(\bx-\vtilde^m_B-m_2\vtilde_2)\ \lambda^2\ V_0(\bx-\vtilde_A^n)\ p_0^\lambda(\bx-\vtilde_A^n-\tm_2\vtilde_2)\ d\bx\ .
 \label{sum3-BA}
 \end{equation}
We apply Lemma \ref{OIsharp} with $\Ione=B$, $\Jone=A$, $\tIone=A$ and $n_1=n$. All summands \eqref{sum3-BA} of $S_3^{BA}(m,n)$, except for exceptional indices in case (b), defined by  $\Ione\ne\Jone$ and  $\tIone=\Jone$,  are bounded by $e^{-c\lambda(|m-n|+|m_2|+|\tm_2|)}\rho_\lambda$ ($n_1=n$). 
 Now exceptional indices in case (b) of Lemma \ref{OIsharp} are such that $\tm_2=0$. However, in  $S_3^{BA}(m,n)$ we sum over $\tm_2\ne0$. Hence, there are no relevant exceptional indices and therefore all expressions  \eqref{sum3-BA} are bounded by $e^{-c\lambda(|m-n|+|m_2|+|\tm_2|)}\rho_\lambda$ . Summing over $m_2\in\Z$ and $\tm_2\in\Z\setminus\{0\}$ we obtain:
 \begin{equation}\label{S3bound}
 |S_3^{BA}(m,n)|\ \lesssim\ e^{-c\lambda}\ e^{-c\lambda|m-n|},\ \ m,n\ge0.
 \end{equation}
 \medskip
 
 Putting together the expression \eqref{SipIJ} for the inner product $\left\langle\ P_{_{\kpar,B}}^{^{\lambda}}[m](\bx)\ , \cHeg\ P_{_{\kpar,A}}^{^{\lambda}}[n]\ \right\rangle_{L^2(\Sigma)}$ with the expansions and bounds in \eqref{S1bound}, \eqref{S2mm}, \eqref{S2mm+1}, \eqref{S2mn} and \eqref{S3bound} we obtain:
 {\footnotesize{
 \begin{align}
 \left\langle\ P_{_{\kpar,B}}^{^{\lambda}}[m](\bx)\ , \cHeg\ P_{_{\kpar,A}}^{^{\lambda}}[n]\ \right\rangle_{L^2(\Sigma)}\ 
 &=\begin{cases}  
 -\left(1+e^{i\kpar}\right)\ \rho_\lambda\ +\ \mathcal{O}\left(e^{-c\lambda}\ e^{-c\lambda|m-n|}\ \rho_\lambda\right),\ & n=m\\
 &\\
 -\rho_\lambda\ +\ \mathcal{O}\left(e^{-c\lambda}\ e^{-c\lambda|m-n|}\ \rho_\lambda\right),\ & n=m+1 \\
 &\\
 \mathcal{O}\left(e^{-c\lambda}\ e^{-c\lambda|m-n|}\ \rho_\lambda\right),\ & n\ne m, m+1.
\end{cases}
\label{mm+1} 
 \end{align}
 }}
 
 By self-adjointness, 
 
  {\footnotesize{
 \begin{align}
 \left\langle\ P_{_{\kpar,A}}^{^{\lambda}}[m](\bx)\ , \cHeg\ P_{_{\kpar,B}}^{^{\lambda}}[n]\ \right\rangle_{L^2(\Sigma)}\ 
 &=\begin{cases}  
 -\left(1+e^{-i\kpar}\right)\ \rho_\lambda\ +\ \mathcal{O}\left(e^{-c\lambda}\ e^{-c\lambda|m-n|}\ \rho_\lambda\right),\ & n=m\\
 &\\
 -\rho_\lambda\ +\ \mathcal{O}\left(e^{-c\lambda}\ e^{-c\lambda|m-n|}\ \rho_\lambda\right),\ & n=m-1 \\
 &\\
 \mathcal{O}\left(e^{-c\lambda}\ e^{-c\lambda|m-n|}\ \rho_\lambda\right),\ & n\ne m, m-1.
\end{cases} 
 \label{mm-1}\end{align}
 }}
 \nit Equations \eqref{mm-1}  and \eqref{mm+1} imply assertions (1), (2) and (3) of Proposition \ref{zz-els}. \medskip
 
 Finally, we turn to the proof of  part (4) of Proposition \ref{zz-els}. By \eqref{SipIJ}, we have
 for $I=A, B$:
\[
\left\langle\ P_{_{\kpar,I}}^{^{\lambda}}[m](\bx)\ , \cHeg\ P_{_{\kpar,I}}^{^{\lambda}}[n]\ \right\rangle_{L^2(\Sigma)}
\ =\ S_1^{II}(m,n)\ +\ S_2^{II}(m,n)\ +\ S_3^{II}(m,n)\ .
\]
 We claim that $|S_j^{II}(m,n)|\lesssim e^{-c\lambda}\ e^{-c\lambda|m-n|}$ for $j=1,2,3$ and $I=A, B$.
  We consider the case $I=A$. The case $I=B$ is essentially the same.
 
 \nit\underline{Estimation of $S_1^{AA}(m,n)$:} 
 The expression to be summed over $m_2, \tm_2\in\Z$ for $n_1\ge0, n_1\ne n$ is: 
 \begin{equation}
 e^{i(\tm_2-m_2)\kpar}\ \int_{\R^2}\ p_0^\lambda(\bx-\vtilde^m_A-m_2\vtilde_2)\ \lambda^2\ V_0(\bx-\vtilde_A^{n_1})\ p_0^\lambda(\bx-\vtilde_A^n-\tm_2\vtilde_2)\ d\bx\ .
 \label{sum1-AA}
 \end{equation}
We apply Lemma \ref{OIsharp} with $\Ione=A$, $\Jone=A$ and $\tIone=A$. All summands in the expression for $S_1^{AA}(m,n)$, except for exceptional indices are bounded by $e^{-c\lambda(|m-n_1|+|m_2|+|n-n_1|+|\tm_2|)}\rho_\lambda$. The only possible exceptional indices are of case (a) in Lemma \ref{OIsharp}. This case requires $n_1= n$ and since the summation in  $S_1^{AA}(m,n)$ is over $n_1\ge0$ with $n_1\ne n$, there are no relevant exceptional indices. We conclude for all $m, n\ge0$:
 \begin{equation}
 |S_1^{AA}(m,n)| \lesssim \rho_\lambda\  \sum_{m_2,\tm_2\in\Z}\  \sum_{\substack{n_1\ge0\\ n_1\ne n}}\ e^{-c\lambda\left(\ |m-n_1|\ +\ |m_2|\ +\ |n-n_1|\ +\ |\tm_2|\ \right)}\ \lesssim\ e^{-c^\prime\lambda}\ e^{-c^\prime\lambda|m-n|}\ \rho_\lambda\ , 
 \label{S1-AA}\end{equation}
 for some strictly positive constant $c^\prime$.

 \medskip
 
 \nit\underline{Estimation of $S_2^{AA}(m,n)$:} The expression to be summed over $m_2, \tm_2\in\Z$ for $n_1\ge0$ is
 \[
e^{i(\tm_2-m_2)\kpar}\
  \int\ p_0^\lambda(\bx-\vtilde^m_A-m_2\vtilde_2)\ \lambda^2\ V_0(\bx-\vtilde_B^{n_1})\ p_0^\lambda(\bx-\vtilde_A^n-\tm_2\vtilde_2)\
   d\bx\ .
  \]
Since  $\Ione=A$, $\Jone=B$ and $\tIone=A$, we have that $\Ione\ne\Jone$ and $\tIone\ne\Jone$. Hence, the bound \eqref{Isharp-est1} applies.  Thus,  all summands in the expression for $S_2^{AA}(m,n)$ are bounded by $e^{-c\lambda}\ e^{-c\lambda(|m-n_1|+|m_2|+|n-n_1|+|\tm_2|)}\rho_\lambda$. Summing over all relevant indices we have: 
\begin{equation}
 |S_2^{AA}(m,n)| \lesssim\
  e^{-c\lambda}\  \sum_{m_2,\tm_2\in\Z}\ \sum_{n_1\ge0}\ e^{-c\lambda(|m-n_1|+|m_2|+|n-n_1|+|\tm_2|)}\ \rho_\lambda\  \lesssim\ e^{-c^\prime\lambda}\ e^{-c^\prime\lambda|m-n|}\ \rho_\lambda\ , 
 \label{S2-AA}\end{equation}
 for some strictly positive constant $c^\prime$.

  \medskip
 
 \nit\underline{Estimation of $S_3^{AA}(m,n)$:} The expression to be summed over $m_2\in\Z$ and $ \tm_2\in\Z\setminus\{0\}$ for $n\ge0$ is
\[
 e^{i(\tm_2-m_2)\kpar}\
  \int\ p_0^\lambda(\bx-\vtilde^m_A-m_2\vtilde_2)\ \lambda^2\ V_0(\bx-\vtilde_A^n)\ p_0^\lambda(\bx-\vtilde_A^n-\tm_2\vtilde_2)\ d\bx\ .\]
  Since  $\Ione=\Jone=\tIone=A$, the only possible exceptional case is case (a). However, note that $\tm_2=0$ is omitted in the summation and hence there are no relevant exceptional cases. Thus,  summands in the expression for $S_3^{AA}(m,n)$ are bounded by $e^{-c\lambda(|m-n|+|m_2||+|\tm_2|)}\rho_\lambda$, and we have: 
  \begin{equation}
 |S_3^{AA}(m,n)| \lesssim\  \sum_{m_2\in\Z}\sum_{\tm_2\in\Z\setminus\{0\}}\ 
 e^{-c\lambda(|m-n|+|m_2|+|\tm_2|)}\rho_\lambda
\ \lesssim\  e^{-c^\prime\lambda}\ e^{-c^\prime\lambda|m-n|}\ \rho_\lambda\ , 
 \label{S3-AA}\end{equation}
 for some strictly positive constant $c^\prime$.\medskip
 
 Finally, summing the bounds \eqref{S1-AA}, \eqref{S2-AA} and \eqref{S3-AA} implies the bound \eqref{ip-II}.
  This completes the proof of Proposition \ref{zz-els}.

\section{Estimation of the nonlinear matrix elements; Proof of Proposition \ref{zz-nl-els}}\label{nl-els} 

Recall our decomposition of $\mathcal{M}^{\lambda}[m,n](\Omega,\kpar)$ into its {\it linear} and {\it nonlinear} contributions:
\begin{align}
&\mathcal{M}^{\lambda}[m,n](\Omega,\kpar)\ =\ \mathcal{M}^{^{\lambda,{\rm l}}}[m,n](\Omega;\kpar)\ -\ \mathcal{M}^{^{\lambda,{\rm nl}}}[m,n](\Omega;\kpar)\ ,
\end{align}
%
where the latter nonlinear matrix elements are given by (see \eqref{cM-def}):
\begin{align}
&\mathcal{M}_{JI}^{^{\lambda,{\rm nl}}}[m,n](\Omega;\kpar)\nn\\
&\quad \equiv\    \left\langle\  H^\lambda_\sharp(\kpar)\ p_{_{\kpar,J}}^\lambda[m]\ ,\ \Pi_{_{AB}}(\kpar)\   \Ressp^\lambda(\Omega,\kpar)\ \Pi_{_{AB}}(\kpar)\  H^\lambda_\sharp(\kpar)\ p_{_{\kpar,I}}^\lambda[n]
\right\rangle_{L^2(\Sigma)}\  .
\label{cM-nl}
\end{align}
Here,  we recall (from Section \ref{XAB}) $\Pi_{_{AB}}(\kpar)$ denotes the projection onto 
\begin{align}
 \mathscr{X}_{AB}(\kpar)\ &=\ \textrm{the orthogonal complement in $L^2(\Sigma)$ of}\ {\rm span}\Big\{p_{_{\kpar,I}}^{^{\lambda}}[n]\ :\ \bI=A,B,\ n\ge0\Big\}\ ,
 \nn\end{align}
 and $\Pi_{_{AB}}(\kpar)\   \Ressp(\Omega,\kpar)\ \Pi_{_{AB}}(\kpar): \mathscr{X}_{_{AB}}(\kpar)\to \mathscr{X}_{_{AB}}(\kpar) $ is the inverse 
of\\ $\Pi_{_{AB}}(\kpar)\ \left(-\left(\nabla_\bx+i\frac{\kpar}{2\pi}\ktilde_2\right)^2+V_\sharp-E_0^\lambda-\Omega\right)\ \Pi_{_{AB}}(\kpar)$ .\\
Furthermore,  the operator $\Pi_{_{AB}}(\kpar)\Ressp(\Omega,\kpar)\Pi_{_{AB}}(\kpar)$ arises from a kernel $\Ressp(\bx,\by,\Omega,\kpar)$;
 see Corollary \ref{corK}. And finally we recall the projection operators $\Pi_\Gamma^\lambda$ (see \eqref{XGam}) which projects onto the orthogonal complement of the set of atomic ground states, centered at nuclei of the discrete set $\Gamma$, 
\[\mathscr{X}_\Gamma\ \equiv\ {\rm span}\Big\{p_\omega^\lambda:\omega\in\Gamma\Big\}^\perp\]
and $\tilde{\Pi}_\Gamma^\lambda=I-\Pi_\Gamma^\lambda$; see \eqref{XGam} and Proposition \ref{summary-K}. In the following discussion we shall be interested in the choice $\Gamma=\mathbb{H}_\sharp$, the zigzag truncation of $\mathbb{H}$. Finally, we recall the notation: $F_\omega(\bx)=F(\bx-\omega)$.

Given  $F(\bx)$, a rapidly decaying function on $\R^2$, define
 \begin{equation}
 F_{_{[\omega]}}(\bx)\ \equiv\ \sum_{n\in\Z}\ F(\bx-\omega+n\vtilde_2)\ 
 =\ \sum_{n\in\Z}\ F_\omega(\bx+n\vtilde_2)\label{F[om]} .
 \end{equation}
The functions $p_{_{\kpar,J}}^\lambda[m]$ in \eqref{cM-nl} are of this type and we now seek to bound
 inner products in  $L^2(\Sigma)$ of the form \eqref{cM-nl} .

 For a constant $\gamma>0$ to be fixed, we   introduce the weighted $L^2(\R^2)-$ spaces:
\begin{equation}
\mathscr{H}^{(\omega)}\ \equiv\ L^2\left(\ \R^2;e^{\gamma |\bx-\omega|}\ d\bx\ \right) .
\label{Hom}\end{equation}
\begin{proposition}\label{sig2R2}
Fix $\Gamma=\mathbb{H}_\sharp$, which is translation-invariant by the vector $\vtilde_2\in\mathbb{H}$.   Let $[\omega], [\omega^\prime]$
 denote equivalence classes (see \eqref{LamSig} with $\Gamma=\mathbb{H}_\sharp$),  and   $\omega_0\in[\omega]\cap\Omega_\Sigma$  and $\omega^\prime_0\in[\omega^\prime]\cap\Omega_\Sigma$. 
 
 \begin{enumerate}
\item For any rapidly decaying functions $F$ and $G$  on $\R^2$ we have
 \begin{align}
& \left\langle\ F_{_{[\omega]}}\ ,\ \Pi_{_{AB}}(\kpar)\ \mathcal{K}^\lambda_\sharp(\Omega,\kpar)\ \Pi_{_{AB}}(\kpar)
 \ G_{_{[\omega^\prime]}}\ \right\rangle_{L^2(\Sigma)}\nn\\
 & =\ 
\sum_{l\in\Z}\ \int_{\bx\in\R^2}\ F_{_{\omega_0}}(\bx)\ \int_{\by\in\R^2}\  \mathcal{K}^\lambda_\sharp(\bx,\by+l\vtilde_2;\Omega,\kpar)\ G_{_{\omega^\prime_0}}(\by)\ d\by\ d\bx\ .
 \label{ipSigR2} \end{align}

\item The expression in \eqref{ipSigR2} may be bounded in exponentially weighted norms as follows:
\begin{align}
& \Big|\ \left\langle\ F_{_{[\omega]}}\ ,\ \Pi_{_{AB}}(\kpar)\ \mathcal{K}^\lambda_\sharp(\Omega,\kpar)\ \Pi_{_{AB}}(\kpar)
 \ G_{_{[\omega^\prime]}}\ \right\rangle_{L^2(\Sigma)}\ \Big|\nn\\
 & \le\  \Big[\sum_{l\in\Z}\ \|\mathcal{K}_\sharp^{\lambda,\omega_0,\omega_0^\prime,l}(\Omega,\kpar)\|_{_{L^2(\R^2)\to L^2(\R^2)}}\Big]\ \|F_{\omega_0}\|_{\mathcal{H}^{(\omega_0)}}\
\|G_{\omega^\prime_0}\|_{\mathcal{H}^{(\omega^\prime_0)}}\ .
  \label{bd-sker}\end{align}
where 
\begin{equation}
\left( \mathcal{K}_\sharp^{\lambda,\omega_0,\omega_0^\prime,l}f \right)(\bx)
 =\ \int_{\R^2}\ e^{-\frac{\gamma}{2}|\bx-\omega_0|}\  
 \mathcal{K}^\lambda_\sharp(\bx,\by+l\vtilde_2)\ 
 e^{-\frac{\gamma}{2}|\by-\omega^\prime_0|}\ f(\by)\ d\by. 
   \label{Kshift1}\end{equation}
\end{enumerate}
Note: The above may be formulated for an arbitrary discrete set $\Gamma$ satisfying $\inf\{|\omega-\omega^\prime|:\omega,\omega^\prime\in\Gamma\ {\rm distinct}\ \}>r_4$, which is translation invariant by the vector  $\vtilde_2$.
  \end{proposition}
  
  \nit{\it Proof of Proposition \ref{sig2R2}:} By Corollary \ref{corK} we have that 
  the operator $ \Pi_{_{AB}}(\kpar)\ \Ressp(\Omega,\kpar)\ \Pi_{_{AB}}$ arises 
  from a kernel $\mathcal{K}^\lambda_\Gamma(\bx,\by;\Omega,\kpar)$. We have 
  \begin{align}
  & \left\langle\ F_{_{[\omega]}}\ ,\ \ \Pi_{_{AB}}(\kpar)\ \Ressp(\Omega,\kpar)\ 
 \Pi_{_{AB}}(\kpar)\ G_{_{[\omega^\prime]}}\ \right\rangle_{L^2(\Sigma)}\nn\\
 & =\ \int_{_{\Omega_\Sigma}}\ F_{_{[\omega]}}(\bx)\ \int_{\by\in\R^2}\ \mathcal{K}^\lambda_\sharp(\bx,\by;\Omega,\kpar)\ G_{_{[\omega^\prime]}}(\by)\ d\by\ d\bx\nn\\
 & =\ \int_{_{\Omega_\Sigma}}\ \sum_{n\in\Z}\ F(\bx-\omega_0+n\vtilde_2)\ \int_{\by\in\R^2}\ 
\mathcal{K}_\sharp^\lambda(\bx,\by;\Omega)\ \sum_{n^\prime\in\Z}\ G(\by-\omega_0^\prime+n^\prime\vtilde_2)\ d\by\ d\bx\nn\\
&=\ \sum_{n,n^\prime\in\Z}\ \int_{_{\Omega_\Sigma}}\  F_{\omega_0}(\bx+n\vtilde_2)\ \int_{\by\in\R^2}\ \mathcal{K}_\sharp^\lambda(\bx,\by;\Omega,\kpar)\ G_{\omega_0^\prime}(\by+n\vtilde_2)\ d\by\ d\bx\nn\\
&=_{\Big[\substack{\tx=\bx+n\vtilde_2\\ \ty=\by+n^\prime\vtilde_2}\Big]}\ \sum_{n,n^\prime\in\Z}
\ \int_{_{\tx\in\Omega_\Sigma+n\vtilde_2}}\ F_{\omega_0}(\tx)\  \int_{\ty\in\R^2}
\mathcal{K}_\sharp^\lambda(\tx-n\vtilde_2,\ty-n^\prime\vtilde_2;\Omega,\kpar)\ G_{\omega_0^\prime}(\ty)\ d\ty\ d\tx\nn\\
 &=\qquad\  {\textrm{by equation \eqref{KG-trin}}}\nn\\
 &\qquad\  \sum_{n,n^\prime\in\Z}
\ \int_{_{\tx\in\Omega_\Sigma+n\vtilde_2}}\ F_{\omega_0}(\tx)\ \int_{\ty\in\R^2}
\mathcal{K}_\sharp^\lambda(\tx,\ty+(n-n^\prime)\vtilde_2;\Omega,\kpar)\ G_{\omega_0^\prime}(\ty)\ d\ty\ d\tx\nn\\
 &= \sum_{n\in\Z}
\ \int_{_{\tx\in\Omega_\Sigma+n\vtilde_2}}\ F_{\omega_0}(\tx)\
 \sum_{n^\prime\in\Z}\  \int_{\ty\in\R^2}\ \mathcal{K}^\lambda_\sharp(\tx,\ty+(n-n^\prime)\vtilde_2;\Omega,\kpar)\ G_{\omega_0^\prime}(\ty)\ d\ty\ d\tx\nn\\
 &= 
\ \int_{_{\tx\in\R^2}}\ F_{\omega_0}(\tx)\
  \int_{\ty\in\R^2}\ \sum_{l\in\Z}\ \mathcal{K}^\lambda_\sharp(\tx,\ty+l\vtilde_2;\Omega,\kpar)\ G_{\omega_0^\prime}(\ty)\ d\ty\ d\tx \nn\\
  &=\ 
\sum_{l\in\Z}\  \int_{_{\tx\in\R^2}}\ F_{\omega_0}(\tx)\
  \int_{\ty\in\R^2}\  \mathcal{K}^\lambda_\sharp(\tx,\ty+l\vtilde_2;\Omega,\kpar)\ G_{\omega_0^\prime}(\ty)\ d\ty\ d\tx \nn
   \end{align}
  This completes the proof of part (1) of Proposition \ref{sig2R2}. 
 To prove part (2) of Proposition \ref{sig2R2}, we bound the expression in \eqref{ipSigR2}. 
 Write 
  $\mathcal{K}_\sharp^{\lambda,\omega_0,\omega_0^\prime,l}$ for the operator:  \begin{equation}
\left( \mathcal{K}_\sharp^{\lambda,\omega_0,\omega_0^\prime,l}f \right)(\bx)
 =\ \int_{\R^2}\ e^{-\frac{\gamma}{2}|\bx-\omega_0|}\  
 \mathcal{K}^\lambda_\sharp(\bx,\by+l\vtilde_2)\ 
 e^{-\frac{\gamma}{2}|\by-\omega^\prime_0|}\ f(\by)\ d\by. 
\label{Kshift1a}\end{equation}
Then, by part (1) of Proposition \ref{sig2R2}, we have 
\begin{align*}
& \Big|\ \left\langle\ F_{_{[\omega]}}\ ,\ \Pi_{_{AB}}(\kpar)\ \mathcal{K}^\lambda_\sharp(\Omega,\kpar)\ \Pi_{_{AB}}(\kpar)
 \ G_{_{[\omega^\prime]}}\ \right\rangle_{L^2(\Sigma)}\ \Big|\nn\\
 & \le\ 
\sum_{l\in\Z}\ \Big|\ \int_{\bx\in\R^2}\ \left[\ e^{\frac{\gamma}{2}|\bx-\omega_0|}\ F(\bx-\omega_0)\ \right]\nn\\ &\qquad\qquad \int_{\by\in\R^2}\ \left[ e^{-\frac{\gamma}{2}|\bx-\omega_0|}\mathcal{K}^\lambda_\sharp(\bx,\by+l\vtilde_2;\Omega,\kpar)\ e^{-\frac{\gamma}{2}|\by-\omega_0|}\ \right]\  \left[\ e^{-\frac{\gamma}{2}|\by-\omega^\prime_0|}\ G_{_{\omega^\prime_0}}(\by)\ \right] d\by\ d\bx\  \Big|\\
&\le\ \sum_{l\in\Z}\ \|F_{\omega_0}\|_{\mathcal{H}^{(\omega_0)}}
\ \|\mathcal{K}_\sharp^{\lambda,\omega_0,\omega_0^\prime,l}\|_{_{L^2(\R^2)\to L^2(\R^2)}} 
\ \|G_{\omega^\prime_0}\|_{\mathcal{H}^{(\omega^\prime_0)}}\\
& =\ \Big[\sum_{l\in\Z}\ \|\mathcal{K}_\sharp^{\lambda,\omega_0,\omega_0^\prime,l}\|_{_{L^2(\R^2)\to L^2(\R^2)}}\Big]\ \|F_{\omega_0}\|_{\mathcal{H}^{(\omega_0)}}\
\|G_{\omega^\prime_0}\|_{\mathcal{H}^{(\omega^\prime_0)}}\ .
  \end{align*}
This completes the proof of part (2) of Proposition \ref{sig2R2}.\medskip

  We shall apply conclusion (2) of Proposition  \ref{sig2R2} with 
  $F_{_{[\omega]}}\ =\ \cHeg\ p^\lambda_{\kpar,J}[n]$
  and $G_{_{[\omega^\prime]}}=\Pi^\lambda_{_{AB}}(\kpar)\ \cHeg\ p^\lambda_{\kpar,I}[m]$,\ $J,I\in\{A,B\}$. 
  Two more tasks remain in this section:
  \begin{enumerate}
 \item  Bound the sum of norms on the right hand side of  \eqref{bd-sker} using 
  our pointwise kernel bounds, \eqref{K-ker-est},  on $\mathcal{K}_\sharp^\lambda(\bx,\by;\Omega,\kpar)$,
   and 
   \item Bound 
   $ \|F_{\omega_0}\|_{_{\mathscr{H}^{(\omega_0)}}}$ and $ \|G_{\omega^\prime_0}\|_{_{\mathscr{H}^{(\omega^\prime_0)}}}$, where $F_{\omega_0}\ =\ \cHeg\ p_{\omega_0}^\lambda$ and 
   $G_{\omega_0^\prime}=\cHeg\  
   p_{\omega_0^\prime}^\lambda$.
   \end{enumerate}
   
   This will enable us to bound the nonlinear contributions to matrix $\mathcal{M}[m,n](\Omega,\kpar)$, displayed in \eqref{cM-nl}, thereby proving Proposition \ref{zz-nl-els}.
 \medskip
 
 The following two propositions will do the trick:
 
 \begin{proposition}\label{Kshift-bd}
 Let $\omega_0$ and $\omega^\prime_0$ be as in the statement of Proposition \ref{sig2R2}. 
 There exist constants $\lambda_1>0$ and $c>0$ such that for all $\lambda\ge\lambda_1$
and  $|\Omega|\le e^{-c\lambda}$:
 \begin{equation}
 \sum_{l\in\Z}
 \|\mathcal{K}_\sharp^{\lambda;{\omega_0},{\omega_0^\prime},l}(\Omega,\kpar)\|_{_{L^2(\R^2)\to L^2(\R^2)}}\
 \lesssim\ \lambda^{10}\ e^{-c|\omega_0-\omega_0^\prime|}.
 \label{Kshift-bd1}
 \end{equation}
 \end{proposition}
 
  \begin{proposition}\label{FG-bd}
  Then, 
 \[ \|\ \cHeg\ p_{\omega_0}^\lambda\|_{_{\mathscr{H}^{(\omega_0)}}}\ \le\ e^{-c\lambda}\ \sqrt{\rho_\lambda}\quad\ {\rm and}\quad  \|\ 
\cHeg\ p_{\omega_0^\prime}^\lambda\|_{_{\mathscr{H}^{(\omega^\prime_0)}}} \lesssim\ e^{-c\lambda}\ \sqrt{\rho_\lambda}\ .\]
  \end{proposition}
  \medskip
  
  The proofs of Propositions \ref{Kshift-bd} and \ref{FG-bd} are presented in the following two subsections. 
  We first apply them to conclude the proof of Proposition \ref{zz-nl-els}, which gives our bound on nonlinear matrix elements. 
Estimate \eqref{bd-sker} with $F_{\omega_0}= \cHeg p_{\omega_0}^\lambda$ and 
   $G_{\omega_0^\prime}=
  \cHeg  p_{\omega_0^\prime}^\lambda$ implies
    \begin{align}
  &\Big|\ \left\langle\  \cHeg p^\lambda_{\kpar,J}[n]\ ,\  \Pi^\lambda_{_{AB}}(\kpar)\ \Ressp(\Omega,\kpar)\ 
 \Pi^\lambda_{_{AB}}(\kpar)\ \cHeg\ p^\lambda_{\kpar,I}[m]\ \right\rangle_{L^2(\Sigma)}\ \Big|\nn\\
 &\ \le\  \Big[\ \sum_{l\in\Z}
 \|\mathcal{K}_\sharp^{\lambda;{\omega_0},{\omega_0^\prime},l}\|_{_{L^2(\R^2)\to L^2(\R^2)}}\Big]\cdot 
 \|\cHeg\ p_{\omega_0}^\lambda\|_{_{\mathscr{H}^{(\omega_0)}}}\cdot \|\cHeg\ p_{\omega_0^\prime}^\lambda\|_{_{\mathscr{H}^{(\omega^\prime_0)}}}\ .
 \label{sum-ker1}
\end{align}
Now apply Propositions \ref{Kshift-bd} and \ref{FG-bd} to obtain
 \begin{align}
  &\Big|\ \left\langle\  \cHeg\ p^\lambda_{\kpar,J}[n]\ ,\  \Pi^\lambda_{_{AB}}(\kpar)\ \Ressp(\Omega,\kpar)\ 
 \Pi^\lambda_{_{AB}}(\kpar)\ \cHeg\ p^\lambda_{\kpar,I}[m]\ \right\rangle_{L^2(\Sigma)}\ \Big|\nn\\
 &\ \lesssim\  \lambda^{10}\ e^{-c|\omega_0-\omega_0^\prime|}\cdot
 e^{-c\lambda}\ \sqrt{\rho_\lambda}\cdot e^{-c\lambda}\ \sqrt{\rho_\lambda}
 \nn\\
 & \lesssim\  \rho_\lambda\ e^{-c\lambda}\ e^{-c|\omega_0-\omega_0^\prime|}
 \ .  \label{nl-est} \end{align}
 
 We have proved  Proposition \ref{zz-nl-els} for the case $j=0$. From this, the case $j=1$ follows by analytic dependence of the  inner product on $\Omega$; see the remark just prior to the statement of Proposition \ref{zz-nl-els}. 
This completes the proof of Proposition \ref{zz-nl-els}.

 \subsubsection{Proof of Proposition \ref{Kshift-bd}:}
 From the expression for the  integral kernel, displayed in \eqref{Kshift1}, we have
 \begin{align}
   \|\mathcal{K}_\sharp^{\lambda,{\omega_0},{\omega_0^\prime},l}\|_{_{L^2(\R^2)\to L^2(\R^2)}}
\ &\le\  \sup_{\bx\in\R^2}\ \int_{\ty\in\R^2}\ e^{-\frac{\gamma}{2}|\bx-\omega_0|}\ |\mathcal{K}_\sharp^\lambda(\bx,\ty+l\vtilde_2)|\ e^{-\frac{\gamma}{2}|\ty-\omega^\prime_0|}\ d\ty\nn\\
&\qquad\qquad +\ 
   \sup_{\ty\in\R^2}\ \int_{\bx\in\R^2}\ e^{-\frac{\gamma}{2}|\bx-\omega_0|}\ 
   |\mathcal{K}_\sharp^\lambda(\bx,\ty+l\vtilde_2)|\ 
   e^{-\frac{\gamma}{2}|\ty-\omega^\prime_0|}\ d\bx\nn\\
  &=\  \sup_{\bx\in\R^2}\ \int_{\ty\in\R^2}\ e^{-\frac{\gamma}{2}|\bx-\omega_0|}\ |\mathcal{K}_\sharp^\lambda(\bx,\by)|\ e^{-\frac{\gamma}{2}|\by-l\vtilde_2-\omega^\prime_0|}\ d\by\nn\\
&\qquad\qquad +\ 
   \sup_{\ty\in\R^2}\ \int_{\bx\in\R^2}\ e^{-\frac{\gamma}{2}|\bx-\omega_0|}\ 
   |\mathcal{K}_\sharp^\lambda(\bx,\by)|\ 
   e^{-\frac{\gamma}{2}|\by-l\vtilde_2-\omega^\prime_0|}\ d\bx
   \nn\\
   &=\ \sup_{\bx\in\R^2}\ \mathscr{I}^\lambda(\bx;l)\ +\ \sup_{\by\in\R^2}\ \mathscr{J}^\lambda(\by;l)
\label{Knorm} \end{align}
Recall that the kernel $\mathcal{K}_\sharp^\lambda(\bx,\by;\Omega,\kpar)$ satisfies the pointwise bound
 \eqref{K-ker-est}:
\begin{align}
&\left|\ \mathcal{K}_\sharp^\lambda(\bx,\by;\Omega,\kpar)\ \right|\ \le\ C\Big[ \left|\ \log|\bx-\by|\ \right|\ +\ \lambda^{10}\ \Big]\ {\bf 1}_{|\bx-\by|\le R}\ +\ e^{-c\lambda}\ e^{-c\lambda|\bx-\by|}
\label{K-ker-est1} \\
 &\textrm{for all}\quad \bx, \by\ \in\R^2\ .
\nn\end{align}
\nit The bounds on $\sup_{\bx\in\R^2}\ \mathscr{I}^\lambda(\bx;l)$ and $ \sup_{\by\in\R^2}\ \mathscr{J}^\lambda(\by;l)$ are obtained very similarly. We present the argument for $\sup_{\bx\in\R^2}\ \mathscr{I}^\lambda(\bx;l)$.
To bound $\mathscr{I}^\lambda(\bx;l)$, we bound the $d\by$ integral over $\R^2$ separately over the sets $|\bx-\by|\le R$ and $|\bx-\by|\ge R$. Call these parts:  $\mathscr{I}_{_{\le R}}^\lambda(\bx;l)$
and $\mathscr{I}_{_{\ge R}}^\lambda(\bx;l)$.

 First assume $|\bx-\by|\le R$. By \eqref{K-ker-est1}
\begin{align*}
\mathscr{I}_{_{\le R}}^\lambda(\bx;l)\ &\le\ e^{-\frac{\gamma}{2}|\bx-\omega_0|}\ \int_{|\bx-\by|\le R}\  |\mathcal{K}_\sharp^\lambda(\bx,\by;\Omega,\kpar)|\ e^{-\frac{\gamma}{2}|\by-l\vtilde_2-\omega^\prime_0|}\ d\by\nn\\
&\quad\lesssim\ e^{-\frac{\gamma}{2}|\bx-\omega_0|}\ \int_{|\bx-\by|\le R}\  \Big[ \left|\ \log|\bx-\by|\ \right|\ +\ \lambda^{10}\ \Big]\ e^{-\frac{\gamma}{2}|\by-l\vtilde_2-\omega^\prime_0|} d\by\nn\\
&\quad \lesssim\ e^{-\frac{\gamma}{2}|\bx-\omega_0|}\  \int_{|\bz|\le R}\  \Big[ \left|\ \log|\bz|\ \right|\ +\ \lambda^{10}\ \Big]\ e^{-\frac{\gamma}{2}|\bx-\bz-l\vtilde_2-\omega^\prime_0|} d\bz\nn\\
&\quad \lesssim\ e^{-\frac{\gamma}{2}|\bx-\omega_0|}\  \int_{0\le|\bz|<\rho}\  \Big[ \left|\ \log|\bz|\ \right|\ +\ \lambda^{10}\ \Big]\ e^{-\frac{\gamma}{2}|\bx-\bz-l\vtilde_2-\omega^\prime_0|} d\bz\nn\\
&\quad\quad +\ e^{-\frac{\gamma}{2}|\bx-\omega_0|}\ \int_{\rho\le|\bz|\le R}\ \Big[ \left|\ \log|\bz|\ \right|\ +\ \lambda^{10}\ \Big]\ e^{-\frac{\gamma}{2}|\bx-\bz-l\vtilde_2-\omega^\prime_0|} d\bz\nn\\
&\quad \lesssim\ e^{-\frac{\gamma}{2}|\bx-\omega_0|}\  e^{-c_1|\bx-l\vtilde_2-\omega^\prime_0|}\ \int_{0\le|\bz|\le\rho}\ \Big[ \left|\ \log|\bz|\ \right|\ +\ \lambda^{10}\ \Big]\ d\bz\nn\\
&\quad\quad +\ e^{-\frac{\gamma}{2}|\bx-\omega_0|}\ \Big[\ C_{\rho,R}\ +\ \lambda^{10}\ \Big]\ \ 
\int_{\rho\le|\bz|\le R}\ e^{-\frac{\gamma}{2}|\bx-\bz-l\vtilde_2-\omega^\prime_0|} d\bz\ .\nn\\
\end{align*}
The latter two terms are each $\lesssim\ \lambda^{10}\ e^{-c_2|\omega_0-\omega_0^\prime|}\ e^{-c_3|l|}$.
Therefore, 
 \begin{align}
& \sup_{\bx\in\R^2}\ \mathscr{I}_{_{\le R}}^\lambda(\bx;l) \lesssim\ \lambda^{10}\ e^{-c_2|\omega_0-\omega_0^\prime|}\  e^{-c_3|l|}\ .
\label{eyeleR}\end{align}
A similar argument yields  a  bound of this type for  $\sup_{\bx\in\R^2}\ \mathscr{J}_{_{\le R}}^\lambda(\bx;l)$.\medskip

Next assume $|\bx-\by|\ge R$. By \eqref{K-ker-est1},
\begin{align*}
\mathscr{I}_{_{\ge R}}^\lambda(\bx;l)\ &\lesssim\ e^{-\frac{\gamma}{2}\lambda}\ e^{-c|\bx-\omega_0|}\ \int_{|\bx-\by|\ge R}\  e^{-c\lambda|\bx-\by|}\ e^{-\frac{\gamma}{2}|\by-l\vtilde_2-\omega^\prime_0|}\ d\by\ .\nn
\end{align*}
Note that $|\bx-\omega_0|+|\by-l\vtilde_2-\omega^\prime_0|\ge |(\bx-\omega_0)-(\by-l\vtilde_2-\omega^\prime_0)|=|\bx-\by- (\omega_0-\omega^\prime_0)+l\vtilde_2|\ge c_3\left(|\omega_0-\omega^\prime_0| + |l|\right)-|\bx-\by|$. Thus, 
\begin{align}
\mathscr{I}_{_{\ge R}}^\lambda(\bx;l)\ &\lesssim\   e^{-c\lambda} \int_{|\bx-\by|\ge R}\  e^{-c_4\lambda|\bx-\by|}\  d\by\ e^{-c_3|\omega_0-\omega^\prime_0|}\ e^{-c_3|l|}\lesssim e^{-c\lambda} \ e^{-c_3|\omega_0-\omega^\prime_0|}\ e^{-c_3|l|}\ .
\label{eyegeR}
\end{align}
The bounds \eqref{eyeleR} and \eqref{eyegeR} imply that 
\[\sup_{\bx\in\R^2}\ \mathscr{I}^\lambda(\bx;l)\ \le\ 
e^{-c_3|\omega_0-\omega^\prime_0|}\ e^{-c_3|l|}\ \lambda^{10}.
\]
and similarly
\[\sup_{\by\in\R^2}\ \mathscr{J}^\lambda(\by;l)\ \le\ 
e^{-c_3|\omega_0-\omega^\prime_0|}\ e^{-c_3|l|}\ \lambda^{10}.
\]
Therefore, by \eqref{Knorm} it follows that 
$ \|\mathcal{K}_\sharp^{\lambda;{\omega_0},{\omega_0^\prime},l}\|_{_{L^2(\R^2)\to L^2(\R^2)}}
\lesssim e^{-c_3|\omega_0-\omega^\prime_0|}\ e^{-c_3|l|} \ \lambda^{10}$. Finally, summing over $l\in\Z$ we deduce
\eqref{Kshift-bd1}. The proof of Proposition \ref{Kshift-bd} is now complete.
\bigskip

\subsubsection{Proof of Proposition \ref{FG-bd}} We need to verify that there are constants $\gamma, \lambda_1>0$, such that  for all $\omega\in\mathbb{H}_\sharp$ and all $\lambda\ge\lambda_1$:
\begin{equation}
\|\ {\mathcal H}^\lambda_\sharp p_\omega^\lambda\ \|_{_{\mathscr{H}^{(\omega)}}}\ =\
\|\ e^{\frac{\gamma}{2}|\bx-\omega|}\ (-\Delta+V_\sharp^\lambda(\bx)-E_0^\lambda)p_\omega^\lambda(\bx)\ 
\|_{_{L^2(\R^2_\bx)}}\ \lesssim\ e^{-c\lambda}\ \sqrt{\rho_\lambda}\ .
\label{Fdef}
\end{equation}

\nit Since $(-\Delta+\lambda^2V_\omega(\bx))p_\omega^\lambda(\bx)=E_0^\lambda p_\omega^\lambda(\bx)$, it follows that 
\[ {\mathcal H}^\lambda_\sharp p_\omega^\lambda(\bx) \equiv\ (-\Delta+V_\sharp^\lambda(\bx)-E_0^\lambda)p_\omega^\lambda(\bx)
\ =\ \sum_{\omega^\prime\in\mathbb{H}_\sharp\setminus\{\omega\}}\ \lambda^2V_0(\bx-\omega^\prime)p_\omega^\lambda(\bx)\ .
\]
By invariance  of  $\cHeg$ under translation by $\vtilde_2$, we may assume $\omega\in\Omega_\Sigma$. Thus, $\omega=\bv_I+n\vtilde_1$ for $I=A$ or $B$ and $n\ge0$ . Fix $I=A$; the argument for $I=B$ is similar. Then, $p_\omega^\lambda(\bx)=p_0^\lambda(\bx-\bv_A-n\vtilde_1)$.
 Recall, for $I=A, B$ and $n_1, n_2\in\Z$: $\bv_I^{n_1,n_2}=\bv_I+n_1\bv_1+n_2\bv_2$.  
 Thus, 
\begin{align*} \cHeg\ p_\omega^\lambda(\bx) 
\ &=\ \sum_{n_1\ge0,\ n_2\in\Z }\ \lambda^2V_0(\bx-\bv_B^{n_1,n_2})p_0^\lambda(\bx-\bv_A-n\vtilde_1)\nn\\
&\qquad +\ 
\sum_{\substack{n_1\ge0,n_2\in\Z\\ (n_1,n_2)\ne(n,0)} }\ \lambda^2V_0(\bx-\bv^{n_1,n_2}_A)p_0^\lambda(\bx-\bv_A-n\vtilde_1)\ .
\end{align*}
\begin{align*}
\end{align*}
For the $\mathscr{H}^{(\omega)}$ norm ($\omega=\bv_I+n\vtilde_1$) we have
\begin{align}
&\|\ \cHeg\ p_\omega^\lambda\ \|_{_{\mathscr{H}^{(\omega)}}}\ =\ 
\Big\|\ e^{\frac{\gamma}{2}|\bx-\omega|}\ \cHeg\ p_\omega^\lambda(\bx)\ \Big\|_{_{L^2(\R^2_\bx)}}\nn\\
&\qquad\ \le\ \lambda^2\
 \sum_{n_1\ge0, n_2\in\Z}\ \left(\ \int e^{\gamma|\bx-(\bv_A+n\vtilde_1)|}\ |V_0(\bx-\bv_B^{n_1,n_2})|^2\ |p_0^\lambda(\bx-\bv_A-n\vtilde_1)|^2\ d\bx\ \right)^{1\over2}\nn\\
 &\qquad\qquad +\ \lambda^2\ \sum_{\substack{n_1\ge0,n_2\in\Z\\ (n_1,n_2)\ne(n,0)} }\
  \left(\ \int e^{\gamma|\bx-(\bv_A+n\vtilde_1)|}\ |V_0(\bx-\bv_A^{n_1,n_2})|^2\ |p_0^\lambda(\bx-\bv_A-n\vtilde_1)|^2\ d\bx\ \right)^{1\over2}\nn\\
  &\qquad \equiv\  \sum_{n_1\ge0, n_2\in\Z}\ A_{n_1,n_2}^\lambda\ +\ \sum_{\substack{n_1\ge0,n_2\in\Z\\ (n_1,n_2)\ne(n,0)} }\ B_{n_1,n_2}^\lambda\ .
  \label{ABdef}
\end{align}
\nit Consider $A_{n_1,n_2}^\lambda$,  for any fixed $n_1\ge0$ and $n_2\in\Z$.
\begin{align*}
&|\ A_{n_1,n_2}^\lambda\ |^2\ =\ \lambda^4\ \int_{|\bx-\bv_B^{n_1,n_2}|\le r_0}\ 
e^{\gamma|\bx-(\bv_A+n\vtilde_1)|}\ |V_0(\bx-\bv_B^{n_1,n_2})|^2\ |p_0^\lambda(\bx-\bv_A-n\vtilde_1)|^2\ d\bx\\
&\quad =\ \lambda^4\ \int_{|\by|\le r_0}\ 
e^{\gamma|\by+\bv_B^{n_1,n_2}-\bv_A^{n,0}|}\ |V_0(\by)|^2\ |p_0^\lambda(\by+\bv_B^{n_1,n_2}-\bv^{n,0}_A)|^2\ d\by\\
&\quad =\ \lambda^4\ \int_{|\by|\le r_0}\ 
e^{\gamma|\by+\bv_B-\bv_A+(n_1-n)\vtilde_1+n_2\vtilde_2|}\ |V_0(\by)|^2\ 
|p_0^\lambda(\by+\bv_B-\bv_A+(n_1-n)\vtilde_1+n_2\vtilde_2)|^2\ d\by\nn\\
&\quad =\ \lambda^4\ \int_{|\by|\le r_0}\ 
e^{\gamma|\by+\be+(n_1-n)\vtilde_1+n_2\vtilde_2|}\ |V_0(\by)|^2\ 
|p_0^\lambda\left(\by-[-\be+(n-n_1)\vtilde_1-n_2\vtilde_2]\right)|^2\ d\by.
\end{align*}
As in  Section \ref{m-els-est} we divide index pairs $(n-n_1,-n_2)$ into those in the set $N_b(-1)=\{(0,0),(1,0),(0,1)\}$ and 
those not in $N_b(-1)$. 
Those in $N_b(-1)$, ``bad index pairs'' ,  correspond to the cases:
(i) $(n_1,n_2)=(n-1,0)$ with $n\ge1$, (ii) $(n_1,n_2)=(n,0)$ with $n\ge0$ or (iii) $(n_1,n_2)=(n,-1)$
with $n\ge0$.  By the remark immediately following Definition \ref{lem-prep}, we then have for some $l=0,1$ or $2$
 \[
 p_0^\lambda\left(\by-[-\be+(n-n_1)\vtilde_1-n_2\vtilde_2]\right)=p_0^\lambda(\by-[-R^l\be]),
 \]  
 where $R$ is a $2\pi/3$ rotation matrix.
 Therefore, by orthogonality of the matrix $R$ and symmetry assumptions on $V_0$, we have:
 \begin{align*}
&|\ A_{n_1,n_2}^\lambda\ |^2\ = \lambda^4\ \int_{|\by|\le r_0}\ 
e^{2c|\by-[-R^l\be]|}\ |V_0(\by)|^2\ 
|p_0^\lambda(\by-[-R^l\be])|^2\ d\by\\
&\ =\ \lambda^4\ \int_{|\by|\le r_0}\ 
e^{2c|R^{-l}\by+\be|}\ |V_0(R^{-l}\by)|^2\ 
|p_0^\lambda(R^{-l}\by+\be)|^2\ d\by\\
&\ =\ \lambda^4\ \int_{|\bz|\le r_0}\ 
e^{2c|\bz+\be|}\ |V_0(\bz)|^2\ 
|p_0^\lambda(\bz+\be)|^2\ d\bz\ .
\end{align*}
Next, applying the bound \eqref{p0-bound3} to one factor of $p_0^\lambda(\bz+\be)$ yields
\begin{align*}
|\ A_{n_1,n_2}^\lambda\ |^2\ &\ \lesssim \lambda^4\ \|V_0\|_\infty\ \int_{|\bz|\le r_0}\ 
e^{2c|\bz+\be|-c\lambda}\ |V_0(\bz)|\ 
 p_0^\lambda(\bz)\ p_0^\lambda(\bz+\be)\ d\bz\\
&\lesssim e^{-c^\prime\lambda}\rho_\lambda\ .
\end{align*}

Next consider $n_1\ge0$ and $n_2\in\Z$, for which $(n-n_1,-n_2)\notin N_{\rm bad}(-1)$.
  By Proposition \ref{p0-bounds}, in particular \eqref{p0-bound2}, we have 
  \begin{equation}
  p_0^\lambda\left(\by-[-\be+(n-n_1)\vtilde_1-n_2\vtilde_2]\right)
  \lesssim
   e^{-c\lambda(|n-n_1|+|n_2|)}\ p_0^\lambda(\by + \be) .
    \label{pAB-bd}  \end{equation}
    Therefore, for $|\by|\le r_0$ and $\lambda$ sufficiently large:
   \begin{align}
 &e^{\gamma|\by-[-\be+(n-n_1)\vtilde_1-n_2\vtilde_2]|}\ \ p_0^\lambda(\by-[-\be+(n-n_1)\vtilde_1-n_2\vtilde_2])\nn\\
 &\qquad \lesssim
 e^{\gamma|\by-[-\be+(n-n_1)\vtilde_1-n_2\vtilde_2]|}\ \ \   e^{-c\lambda(|n_1-n|+|n_2|)}\ p_0^\lambda(\by + \be)\nn\\
 &\qquad  \lesssim\  e^{-c^\prime\lambda(|n_1-n|+|n_2|)} p_0^\lambda(\by + \be)\lesssim\  
e^{-c^\prime\lambda(|n_1-n|+|n_2|)}\ p_0^\lambda(\by)\ ,
   \label{pAB-bd1}
   \end{align} 
   where the last inequality uses \eqref{p0-bound3}.
 Therefore, for good index pairs $(n-n_1,-n_2)$ we have
 \begin{align*}
|\ A_{n_1,n_2}^\lambda\ |^2\ &\lesssim\ \lambda^4\  \|V_0\|_{_\infty}\  \ e^{-c\lambda(|n-n_1|+|n_2|)}\ \int\ |V_0(\by)| p_0^\lambda(\by)\ p_0^\lambda(\by + \be)\ d\by\\
&\lesssim\ e^{-c\lambda(|n-n_1|+|n_2|)}\ \rho_\lambda.
\end{align*}
Taking the square root and summing over good index pairs $(n_1,n_2)$ we have:
\begin{equation} \sum_{\substack{n_1, n_2\\ (n-n_1,-n_2)\ {\rm good}}}\  A_{n_1,n_2}^\lambda\ \ \lesssim\ e^{-c\lambda}\ \sqrt{\rho_\lambda}.
\label{g-sum}\end{equation}

Taken together with our bound on $|\ A_{n_1,n_2}^\lambda\ |$ for the three cases of bad indices,
 this tells us that  
\begin{equation} \sum_{n_1\ge0, n_2\in\Z}\ A_{n_1,n_2}^\lambda\ \lesssim\ e^{-c\lambda}\ \sqrt{\rho_\lambda}.
\label{A-sum}\end{equation}
The proof that 
\begin{equation} \sum_{n_1\ge0, n_2\in\Z}\ B_{n_1,n_2}^\lambda\ \lesssim\ e^{-c\lambda}\ \sqrt{\rho_\lambda}.
\label{B-sum}\end{equation}
is similar, so this completes the proof of \eqref{FG-bd}.

\appendix

  \section{Error and Main Kernels; Proof of  Lemma \ref{error-main}} \label{app:error-main}
  
We prove that if $\mathcal{E}$ is an operator derived from an error kernel $\mathcal{E}(\bx,\by)$
in the sense of Definition \ref{errorK}, then $\widetilde{\mathcal{E}}=I-(I-\mathcal{E})^{-1}$ is an operator
derived from an error kernel $\widetilde{\mathcal{E}}(\bx,\by)$. 

\subsection{Elementary integrals in 1d}\label{1d-int}
Let $f\in L^1(\R)$.  We define  
$f^{*0}=\delta$, the Dirac delta function and $f^{*1}=f$. Let $f^{*n}$ denote the $n-$ fold convolution of $f$ with itself:

For $f$ and $g$ in $L^1(\R)$,
\begin{equation}
\left(\ f\ +\ g\ \right)^{*n}\ =\ \sum_{k=0}^n\ {n \choose k}\  f^{*k}\ g^{*(n-k)}\ .
\label{fgn}
\end{equation}
Let $f(t)=a e^{-\gamma|t|}$, where $a$ and $\gamma$ are positive constants with $\gamma>a$. We may write
 \[ f(t)=f_+(t)+f_-(t),\quad  f_+(t)=a e^{-\gamma t}\ {\bf 1}_{\{t>0\}},\ \quad f_-(t)=a e^{-\gamma |t|}\ {\bf 1}_{\{t<0\}} .\]
 Induction on $k$ gives:
 \[ f_+^{*k}(t)\ =\ a^k e^{-\gamma t}\ \frac{t^{k-1}}{(k-1)!}\ {\bf 1}_{\{t>0\}}\ 
 \le\ a e^{-\gamma t}\ \sum_{l=0}^\infty 
 \frac{(at)^l}{l!} {\bf 1}_{\{t>0\}}\ =\ a e^{-(\gamma-a)t} {\bf 1}_{\{t>0\}},\ \ k\ge1\ .\]
A similar bound holds for $f_-$. Therefore, for all $0<a<\gamma$:
  \[f_+^{*k}(t)\ \le\ a\ e^{-(\gamma-a)t} {\bf 1}_{\{t>0\}}\qquad {\rm and}\qquad   f_-^{*k}(t)\ \le  a\ e^{-(\gamma-a)|t|} {\bf 1}_{\{t<0\}},\ \ k\ge1\ .  \]
  Therefore, for $m\ge1$, we have from \eqref{fgn} that 
  \begin{align}
 f^{*m}(t)\ =\  \left(\ f_+\ +\ f_-\ \right)^{*m}(t)\ &\le\ 
  a^2\ \sum_{k=1}^{m-1}\ {m \choose k}\  \left[\ e^{-(\gamma-a)t} {\bf 1}_{\{t>0\}}\ \star\ 
  e^{-(\gamma-a)|t|} {\bf 1}_{\{t<0\}} \ \right](t)\nn\\
  &\qquad +\ a e^{-(\gamma-a)t}\ {\bf 1}_{\{t>0\}}\ + \ a e^{-(\gamma-a)|t|}\ {\bf 1}_{\{t<0\}}\ .
  \label{n-conv}\end{align}
  The last two terms, which sum to $a e^{-(\gamma-a)|t|}$,  correspond to $k=0$ and $k=m$ in the binomial formula. We calculate the convolution in \eqref{n-conv}. For $t>0$, 
  \begin{align*}
 & \left[\ e^{-(\gamma-a)|t|} {\bf 1}_{\{t>0\}}\ \star\ 
  e^{-(\gamma-a)|t|} {\bf 1}_{\{t<0\}} \ \right](t)\\
  \ &\qquad =\ \int_0^\infty e^{-(\gamma-a)s}\ e^{-(\gamma-a)|t-s|} {\bf 1}_{\{t-s<0\}}\ ds\ =\ \int_t^\infty e^{-2(\gamma-a)s}\ e^{(\gamma-a)t} \ ds\ =\ \frac{e^{-(\gamma-a)t}}{2(\gamma-a)}.
  \end{align*}
  Similarly, if $t<0$ then this convolution is $\frac{e^{-(\gamma-a)|t|}}{2(\gamma-a)}$. Therefore, 
   \begin{align*}
 & \left[\ e^{-(\gamma-a)|t|} {\bf 1}_{\{t>0\}}\ \star\ 
  e^{-(\gamma-a)|t|} {\bf 1}_{\{t<0\}} \ \right](t) \ =\ \frac{e^{-(\gamma-a)|t|}}{2(\gamma-a)},
  \ \ \textrm{for all}\ t\in\R\ .
  \end{align*}
  Substituting into \eqref{n-conv} we have
  \begin{align}
   f^{*m}(t) &= \left(\ a e^{-\gamma|t|}\ \right)^{*m}(t)\ \le\ a^2\ 
    \frac{e^{-(\gamma-a)|t|}}{2(\gamma-a)}\ \sum_{k=1}^{m-1}\ {m \choose k} + a e^{-(\gamma-a)|t|}
\le \left[ a +\ \frac{ 2^m\ a^2 }{2(\gamma-a)} \right] e^{-(\gamma-a)|t|}  .
  \nn \end{align}
  Therefore,
  \begin{equation} \left(\ \frac{a}4\ e^{-\gamma|t|}\ \right)^{*m}(t)\ \le\ 
  \left[ 4^{-m}\ a +\ \frac{ 2^{-m}\ a^2 }{2(\gamma-a)} \right] e^{-(\gamma-a)|t|}\quad  \textrm{for}\ m\ge1. 
  \label{m-fold}\end{equation}
  
  \subsection{Elementary integrals in  $n$ dimensions}\label{nd-int}
  
 For $(x_1,\dots,x_n)\in\R^n$,  let 
  \[K(x_1,\dots,x_n)=\frac{a^n}{4^n}\ e^{-\gamma\left(|x_1|+\dots+|x_n|\right)}\ \textrm{with $0<a<\gamma$}.\] 
  We now apply \eqref{m-fold} to the $l-$ fold convolution of $K(x_1,\dots,x_n)$:
  \begin{align}
  K^{*l}(x_1,\dots,x_n)\ \ &\equiv\ \underbrace{ K\star K\star\cdots\star K }_{\textrm{$l$- times} }(x_1,\dots, x_n)
   \ \le\ \Pi_{j=1}^n\ \Big\{\ \left[\ 4^{-l}\ a +\ \frac{ 2^{-l}\ a^2 }{2(\gamma-a)} \right] e^{-(\gamma-a)|x_j|}\ \Big\}
   \nn\\
 &\ =\  \left[\ 4^{-l}\ a +\ \frac{ 2^{-l}\ a^2 }{2(\gamma-a)} \right]^n e^{-(\gamma-a)(|x_1|+\dots+|x_n|)}
  \ .
 \label{K-conv-l} \end{align}
    
\subsection{Proof of  part (1) of Lemma \ref{error-main}} 
 For $\bx=(x_1,\dots,x_n)\in\R^n$ we write $|\bx|_{l^1}$ to denote $|x_1|+\dots+|x_n|$.
Suppose that $E(\bx,\by)$ satisfies the bound:
\begin{equation}
 |E(\bx,\by)|\ \le\ (a/4)^ne^{-\gamma|\bx-\by|_{l^1}},\ \ \textrm{for all}\ \ \bx, \by\in\R^n
 \label{Ebound}\end{equation}
and gives rise to the integral operator:
\begin{equation}
\left(\ Ef\ \right)(\bx)\ =\ \int_{\R^n}\ E(\bx,\by)\ f(\by)\ d\by,
\label{cEdef}\end{equation}
 then for all $l\ge1$ the $l^{th}$ power of the operator $E$: $f\mapsto E^l[f]$, is given by 
\[
E^l[f](\bx)\ =\ \int_{\R^n}\ E_l(\bx,\by)\ f(\by)\ d\by\ ,
\]
where by \eqref{K-conv-l}, $E_l$ satisfies the bound
\[ |\ E_l(\bx,\by)\ |\ \le\ \left[\ 4^{-l}\ a +\ \frac{ 2^{-l}\ a^2 }{2(\gamma-a)} \right]^n e^{-(\gamma-a)\ |\bx-\by|_{l^1}}.
\]
If $\gamma>2a$, then $\frac{a^2}{2(\gamma-a)}\le \frac{a}2$. Therefore,  for $l\ge1$:
\[\left[\ 4^{-l}\ a +\ \frac{ 2^{-l}\ a^2 }{2(\gamma-a)} \right]\le 
\left[\ 4^{-l}\ a +\ 2^{-l}\ \frac{a}{2} \right]\ =\ 2^{-l} a\ \left[ 2^{-l}+2^{-1}\right]\le 2^{-l} a. 
\]
Hence, 
\[ |\ E_l(\bx,\by)\ |\ \le\ 2^{-ln} a^n e^{-(\gamma-a)\ |\bx-\by|_{l^1}},\ \ l\ge1.
\]

 \medskip

Let's now apply these observations to  $E(\bx,\by)=\mathcal{E}(\bx,\by)$, where $\mathcal{E}(\bx,\by)$ is an error kernel which by Definition \ref{errorK} satisfies
 $|\mathcal{E}(\bx,\by)|\lesssim e^{-c\lambda}\ e^{-c\lambda|\bx-\by|}$ for $\bx, \by\in\R^2$; here $n=2$.   Note that $e^{-c^{\prime\prime}\lambda|\bx-\by|_{l^1}}\le e^{-c\lambda|\bx-\by|}\ \le\ e^{-c^\prime\lambda|\bx-\by|_{l^1}}$. Therefore, 
$|\mathcal{E}(\bx,\by)|\lesssim e^{-c\lambda}\ e^{-c^\prime\lambda|\bx-\by|_{l^1}}$. It follows
that $\mathcal{E}(\bx,\by)$ satisfies the bound \eqref{Ebound} with $n=2$, $(a/4)^2=e^{-c\lambda}$
 and $\gamma=c^\prime\lambda$. Therefore, the operator $\mathcal{E}^l$ is given by a kernel
  $\mathcal{E}_l(\bx,\by)$:
 \[
 \mathcal{E}^l[f](\bx)\ =\ \int_{\R^n}\ \mathcal{E}_l(\bx,\by)\ f(\by)\ d\by\ ,
\]
where $\mathcal{E}_l$ satisfies the bound
\begin{equation}
|\mathcal{E}_l(\bx,\by)\ |\ \le\ 2^{-2l} e^{-c\lambda} e^{-c\lambda |\bx-\by|},\ \ l\ge1
\label{cal-El}
\end{equation}
for some $c>0$, which is independent of $l$.   Consequently, $f\mapsto \widetilde{\mathcal{E}}f=\left(\ I\ - \left(I-\mathcal{E}\right)^{-1}\ \right) f=\sum_{l\ge1}\ \mathcal{E}^l f$
is given by the kernel $\widetilde{\mathcal{E}}(\bx,\by)=\sum_{l\ge1}\ \mathcal{E}_l(\bx,\by)$, which by
 \eqref{cal-El} satisfies the bound  $|\widetilde{\mathcal{E}}(\bx,\by)|\lesssim e^{-c\lambda} e^{-c\lambda |\bx-\by|}$. Thus,  $\widetilde{\mathcal{E}}$ is an error kernel and
 \begin{align}
\widetilde{\mathcal{E}}f(\bx) =\ \int_{\R^2}\ \widetilde{\mathcal{E}}(\bx,\by)\ f(\by)\ d\by \ .
\end{align} 
The  proof of part (1) of  Lemma \ref{error-main} is now complete.
 
\subsubsection{Proof of part (2) of Lemma \ref{error-main}} We need to prove that if $\mathcal{E}^\lambda$ derives from an error kernel and $K^\lambda$ from a main kernel, then $K^\lambda\mathcal{E}^\lambda$ and 
$\mathcal{E}^\lambda K^\lambda$  derive from error kernels $(K\mathcal{E}^\lambda)(\bx,\by)$ and 
$(\mathcal{E}^\lambda K^\lambda)(\bx,\by)$. 
 We begin with the following bounds on  $\mathcal{E}^\lambda(\bx,\bz)$ and $K^\lambda(\bz,\by)$: 
\begin{align}
|\ \mathcal{E}^\lambda(\bx,\bz)\ |\ &\lesssim\ e^{-c\lambda}\ e^{-c\lambda|\bx-\bz|}\nn\\
|\ K^\lambda(\bz,\by)\ |\ &\lesssim\ \left[\ \lambda^4\ +\ \left|\ \log|\bz-\by|\ \right|\ \right]\ {\bf 1}_{\{|\bz-\by|\le R\}}\ +\ 
e^{-c\lambda}\ e^{-c\lambda|\bz-\by|}\ .\nn
\end{align}
Thus, 
\begin{align}
&|\ (\mathcal{E}^\lambda K^\lambda)(\bx,\by)\ |\nn\\
 &\qquad \lesssim\  \int_{|\bz-\by|\le R}\ \left(\lambda^4+\Big|\log|\bz-\by|\Big|\right)\ e^{-c\lambda}\ e^{-c\lambda|\bx-\bz|}\ d\bz\ +\ 
 \int_{|\bz-\by|\ge R}\  e^{-c\lambda}\ e^{-c\lambda|\bx-\bz|}\  e^{-c\lambda|\bz-\by|} d\bz\nn\\
 &\qquad \lesssim\ e^{-c^\prime\lambda}\ e^{-c^\prime\lambda|\bx-\by|}\ .\ 
 \nn
\end{align}
Thus, $\left(\mathcal{E}^\lambda K^\lambda\right)(\bx,\by)$ is an error kernel. A similar bound shows that 
$\left(K^\lambda\mathcal{E}^\lambda\right)(\bx,\by) $ is an error kernel. 

\subsubsection{Proof of part (3) of Lemma \ref{error-main}} 
We show that if $K_\lambda$ arises from a main kernel, 
then $e^{-c\lambda}K_\lambda^2$ arises from an error kernel. Since $K_\lambda(\bx,\by)$ is bounded by the sum of a first term:
  $\sim\ (\lambda^4\ +\ \Big|\log|\bx-\by|\Big|)\ 1_{|\bx-\by|<R}$  and 
a second term $\lesssim\ e^{-c\lambda}\ e^{-c\lambda|\bx-\by|}$ (an error kernel), by part (2) we need only consider the contribution to 
$e^{-c\lambda}\ (K^2_\lambda)(\bx,\bz)\ =\ e^{-c\lambda}\ \int K_\lambda(\bx,\by)K_\lambda(\by,\bz)d\by$ arising from the first term.
The size of this contribution is $\lesssim\ \lambda^8e^{-c\lambda} 1_{|\bx-\bz|<2R}\lesssim e^{-c^\prime\lambda}\ 1_{|\bx-\bz|<2R}$. Hence,
 \[
 e^{-c\lambda}\ (K^2_\lambda)(\bx,\bz)\lesssim e^{-c^\prime\lambda}\ 1_{|\bx-\bz|<2R}\ +\ e^{-c\lambda}\ e^{-c\lambda|\bx-\bz|}\
  \lesssim\ e^{-c^{\prime\prime}\lambda}\ e^{-c^{\prime\prime}\lambda|\bx-\bz|}\ .
  \]
  Hence, $e^{-c\lambda} K^2_\lambda$ derives from an error kernel.

\section{Overlap integrals; proof of Lemma \ref{OIsharp}}\label{overlap}

In this section we prove Lemma \ref{OIsharp}, which we restate here for convenience:
 
\nit {\it For $\Ione, \Jone, \tIone\in\{A,B\}$, $m, n, n_1\ge0$ and $\tm_2\in\Z$, consider the overlap integral
 \begin{equation}
\mathscr{I}_\sharp \equiv \int\ p_0^\lambda(\bx-\vtilde^m_{\Ione}-m_2\vtilde_2)\ \lambda^2\ |V_0(\bx-\vtilde_{\Jone}^{n_1})|\ p_0^\lambda(\bx-\vtilde_{\tIone}^n-\tm_2\vtilde_2)\ d\bx\ .
 \label{ol-sharp1}
 \end{equation}
 Note that the overlap integral in \eqref{ol-sharp1}, although taken over $\R^2$, has an integrand  supported 
 on the disc $B_{r_0}(\vtilde_{\Jone}^{n_1})$. 
 Recall the hopping coefficient defined by: 
 \[\rho_\lambda=\int p^\lambda_0(\by)\ \lambda^2\ |V_0(\by)|\ p^\lambda_0(\by-\be)\ d\by\ .\]
 \nit We also recall from Lemma \ref{lem-prep} that for  $I,J\in\{A,B\}$, we define  $\sigma(I,J)$ so that: $\vtilde_I-\vtilde_J=\sigma\ (\vtilde_B-\vtilde_A)\equiv \sigma\be$. Thus, 
 $\sigma(A,B)=-1$, $\sigma(B,A)=1$, and $\sigma(A,A)=\sigma(B,B)=0$.\\
 
 \nit  Further, 
 for $\sigma=+1,-1,0$ we define $N_b(\sigma)=\{(r_1,r_1): |\sigma\be+r_1\vtilde_1+r_2\vtilde_2|=|\be|\}$. Hence,
$N_b(+1)\ \equiv\ \{(0,0),(-1,0),(0,-1)\}$,\ 
$N_b(-1)\ \equiv\ \{(0,0),(1,0),(0,1)\}$, and $N_b(0)\ \equiv\ \emptyset$}\ .
 \medskip

\nit Lemma \ref{OIsharp} asserts  the bound:
 \begin{equation}
 \mathscr{I}_\sharp \ \lesssim\ e^{-c\lambda\left(\ |m-n_1|\ +\ |m_2|\ +\ |n-n_1|\ +\ |\tm_2|\ \right)}\ \rho_\lambda,
 \label{Isharp-est}
 \end{equation}
except in the following  cases of \underline{exceptional  indices} $(m,n,n_1,m_2,\tm_2)$:
 \begin{enumerate}
 \item[(a)]\ $\Ione=\tIone=\Jone$, $m=n=n_1$ and $m_2=\tm_2=0$. This case does not arise in the proof of Proposition \ref{zz-els},  so we say nothing further about it.\\
 \item[(b)]\ $\tIone=\Jone$, $\Ione\ne\Jone$, $(m-n_1,m_2)\in N_b\left(\sigma(\Ione,\Jone)\right)$, $n=n_1$ and $\tm_2=0$,\\ in which case $ \mathscr{I}_\sharp=\rho_\lambda$.\\
 \item[(c)]\ $\Ione=\Jone$, $\tIone\ne\Jone$, $(n-n_1,\tm_2)\in N_b\left(\sigma(\tIone,\Jone)\right)$, $m=n_1$ and $m_2=0$,\\
in which case $ \mathscr{I}_\sharp=\rho_\lambda$.\\
 \end{enumerate}
 \medskip
 
\nit Lemma \ref{OIsharp} further asserts that if $\Ione\ne\Jone$, $\tIone\ne\Jone$, then 
 \begin{equation}
 \mathscr{I}_\sharp \ \lesssim\ e^{-c\lambda}\ e^{-c\lambda\left(\ |m-n_1|\ +\ |m_2|\ +\ |n-n_1|\ +\ |\tm_2|\ \right)}\ \rho_\lambda,
 \label{Isharp-est1A}
 \end{equation}

\nit We shall occasionally use the notation: $\bfm\vec\vtilde=m_1\vtilde_1+m_2\vtilde_2$, where $\bfm=(m_1,m_2)\in\Z^2$.
\medskip
 
  To prove Lemma \ref{OIsharp} we begin with a change of variables: $\by=\bx-\vtilde_{\Jone}^{n_1}$. Therefore,
 {\footnotesize{
  \begin{align}
\mathscr{I}_\sharp &\equiv \int\ p_0^\lambda(\by-[\sigma(\Ione,\Jone)\be+(m-n_1)\vtilde_1+m_2\vtilde_2])\ 
\lambda^2\ |V_0(\by)|\ p_0^\lambda(\by-[\sigma(\tIone,\Jone)\be+(n-n_1)\vtilde_1+\tm_2\vtilde_2])\ d\by . \nn\\
 &\label{ol-sharp2}
 \end{align}
 }}
 Thus, our task is to consider integrals of the form
 {\footnotesize{
 \begin{equation}
 \mathscr{I}=\ \int\ p_0^\lambda(\by-[\sigma\be+r_1\vtilde_1+r_2\vtilde_2])\ 
\lambda^2\ |V_0(\by)|\ p_0^\lambda(\by-[\tsigma\be+\tr_1\vtilde_1+\tr_2\vtilde_2])\ d\by\ .
\label{OI}
\end{equation}
}}

\begin{lemma}\label{I-bound}
Consider the overlap integral \eqref{OI}, which depends on $\sigma, \tsigma\in \{0,+1,-1\}$ and 
 $\br=(r_1, r_2)$, $\btr=(\tr_1,\tr_2)\in\Z^2$. 
  The expression $\mathscr{I}$ satisfies the bound:
\begin{equation}
 \mathscr{I}(\sigma,\br,\tsigma,\btr)\ \lesssim\ e^{-c\lambda(|r_1|+|r_2|+|\tr_1|+|\tr_2|)}\ \rho_\lambda\ 
 \label{I-bound1}
 \end{equation}
 except in the following cases:
 \begin{enumerate}
 \item[($\alpha$)]\ $\sigma=\tsigma=0$, $\br=\bzero$, $\btr=\bzero$.\\ This case does not arise in the proof of Proposition \ref{zz-els} so we say nothing further about it.\\
 \item[($\beta$)]\ $\tsigma=0$, $\sigma\ne0$, $\br\in N_b(\sigma)$, $\btr=\bzero$, in which case 
 $\mathscr{I}=\rho_\lambda$.\\
 \item[($\gamma$)]\ $\tsigma\ne0$, $\sigma=0$, $\btr\in N_b(\tsigma)$, $\br=\bzero$, in which case $\mathscr{I}=\rho_\lambda$. 
 \end{enumerate}
 \end{lemma}
 We shall also make use of 
 \begin{lemma}\label{sigtsig}
Suppose   $\tsigma\ne0$ and $\sigma\ne0$. Then, 
\begin{enumerate}
\item If $\br\in N_b(\sigma)$ and $\btr\in N_b(\tsigma)$, then 
 \begin{equation}
 \mathscr{I}(\sigma,\br,\tsigma,\btr)\ \lesssim\ e^{-c\lambda}\  \rho_\lambda\  .
 \label{I-bound2a}
 \end{equation}
\item If $\br\in N_b(\sigma)$ and $\btr\notin N_b(\tsigma)$, then 
 \begin{equation}
 \mathscr{I}(\sigma,\br,\tsigma,\btr)\ \lesssim\ e^{-c\lambda}\  e^{-c\lambda(|\tr_1|+|\tr_2|)}\rho_\lambda\ .
 \label{I-bound2b}
 \end{equation}
 The analogous bound  holds with $\br$ and $\btr$ interchanged. \\
\item If $\br\notin N_b(\sigma)$ and $\btr\notin N_b(\tsigma)$ (and therefore $\br,\btr\ne (0,0)$),  then 
 \begin{equation}
 \mathscr{I}(\sigma,\br,\tsigma,\btr)\ \lesssim\  e^{-c^\prime\lambda}\ e^{-c^\prime \lambda(|r_1|+|r_2|+|\tr_1|+|\tr_2|)}\rho_\lambda\ .
 \label{I-bound2c}
 \end{equation}
 \end{enumerate}
\end{lemma}
Note that Lemma \ref{OIsharp} is an immediate consequence of Lemma \ref{I-bound} and Lemma \ref{sigtsig} since $\mathscr{I}_\sharp=\mathscr{I}(\sigma,\br,\tsigma,\btr)$ (see \eqref{OI}), for the choices:
$\sigma=\sigma(\Ione,\Jone)$, $\tsigma=\sigma(\tIone,\Jone)$, $(r_1,r_2)=(m-n_1,m_2)$ and $(\tr_1,\tr_2)=(n-n_1,\tm_2)$.
 Hence it suffices to prove Lemma \ref{I-bound} and Lemma \ref{sigtsig}.
 \subsection{Proof of Lemma \ref{I-bound} and Lemma \ref{sigtsig}:} We estimate the overlap integral \eqref{OI} by considering the two cases:
 {\bf Case 1: $\tsigma=0$} and {\bf Case 2: $\tsigma\ne0$}, and a number of subcases within each.
 \medskip
 
 \nit{\bf Case 1:\ $\tsigma=0$.}\quad In this case, for all $\by\in B_{r_0}(0)$, we have by \eqref{p0-bound1}:
  \begin{equation}
 p_0^\lambda(\by-\tsigma\be-\tr_1\vtilde_1-\tr_2\vtilde_2)\ =\ p_0^\lambda(\by-\tr_1\vtilde_1-\tr_2\vtilde_2)\ 
 \lesssim\ e^{-c\lambda(|\tr_1|+|\tr_2|)}\ p_0^\lambda(\by)\ .
 \label{b1}
 \end{equation}
 Thus, 
 \begin{align}
 \mathscr{I}(\sigma,\br,\tsigma,\btr)\ &=\ \int\ p_0^\lambda(\by-[\sigma\be+r_1\vtilde_1+r_2\vtilde_2])\ 
\lambda^2\ |V_0(\by)|\ p_0^\lambda(\by-[\tr_1\vtilde_1+\tr_2\vtilde_2])\ d\by\nn\\
&\lesssim\ e^{-c\lambda(|\tr_1|+|\tr_2|)}\ \int\ p_0^\lambda(\by-[\sigma\be+r_1\vtilde_1+r_2\vtilde_2])\ 
\lambda^2\ |V_0(\by)|\ p_0^\lambda(\by)\ d\by\ . \label{p1}
\end{align}
We next consider two subcases: 
\[\textrm{Subcase 1A:\ $\tsigma=0$ and $\sigma=0$ and Subcase 1B:\ $\tsigma=0$ and $\sigma\ne0$}\]
 \nit {\bf Subcase 1A: \ $\tsigma=0$ and  $\sigma=0$}\ \ For any $(r_1,r_2)\ne(0,0)$,  we have by \eqref{p0-bound4}
\begin{equation}
 p_0^\lambda(\by-\sigma\be-r_1\vtilde_1-r_2\vtilde_2)\ =\ p_0^\lambda(\by-[r_1\vtilde_1+r_2\vtilde_2])\lesssim\ e^{-c\lambda(|r_1|+|r_2|)}\ p_0^\lambda(\by-\be)\ .
 \label{b2}
 \end{equation}
 Therefore, in subcase 1A  we have after substitution of \eqref{b2}
 into \eqref{p1}, that
 \begin{align}
 \mathscr{I}(\sigma,\br,\tsigma,\btr)\ \lesssim\ e^{-c\lambda(|r_1|+|r_2|+|\tr_1|+|\tr_2|)}\ \rho_\lambda\ .
 \label{b3}
 \end{align}
 Interchanging the roles of $\br$ and $\btr$ in the case where $\tsigma=\sigma=0$, we also have that \eqref{b3}
 holds unless $\btr=0$. Hence when $\sigma=\tsigma=0$, we have \eqref{b3} unless $r_1=r_2=\tr_1=\tr_2=0$.
 \medskip
 
 \nit{\bf Subcase 1B, $\tsigma=0$ and  $\sigma\ne0$:}\ Then, by \eqref{p0-bound2} we have
\begin{equation}
 p_0^\lambda(\by-\sigma\be-r_1\vtilde_1-r_2\vtilde_2)\ \lesssim\ e^{-c\lambda(|r_1|+|r_2|)} p_0^\lambda(\by-\sigma\be)
 \label{b4}
 \end{equation}
 unless $(r_1,r_2)\in N_b(\sigma)$. Substituting \eqref{b4} into \eqref{p1}, we obtain the bound \eqref{b3} unless $(r_1,r_2)\in N_b(\sigma)$.
 
 Now consider the case where $(r_1,r_2)\in N_b(\sigma)$. Then, for some $l\in\{0,1,2\}$ which depends on $\sigma$, $r_1$ and $r_2$ we have:\ $p_0^\lambda(\by-(\sigma\be+r_1\vtilde_1+r_2\vtilde_2))\ =\ p_0^\lambda(\by-\sigma R^{-l}\be)$, where $l=0,1$ or $2$ and  $R$ is a $2\pi/3$ rotation matrix. Substituting into \eqref{p1}, we conclude that $\mathscr{I}(\sigma,\br,\tsigma,\btr)\ \lesssim\ e^{-c\lambda(|\tr_1|+|\tr_2|)}\rho_\lambda$. Indeed,  using symmetry we obtain for $(r_1,r_2)\in N_b(\sigma)$:
 \begin{align}
 \mathscr{I}(\sigma,\br,\tsigma,\btr)\ &\lesssim\ e^{-c\lambda(|\tr_1|+|\tr_2|)}\ \int p_0^\lambda(\by) \lambda^2|V_0(\by)| p_0^\lambda(R^l\by-\sigma\be)\ d\by\nn\\
 &=\ e^{-c\lambda(|\tr_1|+|\tr_2|)}\ \int p_0^\lambda(R^l\by) \lambda^2|V_0(R^l\by)| p_0^\lambda(R^l\by-\sigma\be)\ d\by\nn\\
 &=\ e^{-c\lambda(|\tr_1|+|\tr_2|)}\ \int p_0^\lambda(\bz) \lambda^2|V_0(\bz)| p_0^\lambda(\bz-\sigma\be)\ d\bz\ =\  e^{-c\lambda(|\tr_1|+|\tr_2|)}\ \rho_\lambda. \label{b5}
 \end{align}
 Since $|r_1|+|r_2|=0$ or $1$ for $(r_1,r_2)\in N_b(\sigma)$, it follows that \eqref{b3} holds (with a smaller constant, also denoted $c$, than appearing on the right hand side of \eqref{b5}), unless $\tr_1=\tr_2=0$. Therefore, if $\tsigma=0$ and $\sigma\ne0$, the bound \eqref{b3} holds provided $(\tr_1,\tr_2)\ne(0,0)$. 
 \medskip
 
 Now consider the case where $\tsigma=0, \sigma\ne0$,  $(r_1,r_2)\in N_b(\sigma)$ and $(\tr_1,\tr_2)=(0,0)$.  Then, 
 \begin{align*} \mathscr{I}(\sigma,\br,\tsigma,\btr)\ &=\ \int\ p_0^\lambda(\by-[\sigma\be+r_1\vtilde_1+r_2\vtilde_2])\ 
\lambda^2\ |V_0(\by)|\ p_0^\lambda(\by)\ d\by\\
&=\ \int\ p_0^\lambda(\by-[\sigma R^{-l}\be])\ 
\lambda^2\ |V_0(\by)|\ p_0^\lambda(\by)\ d\by\ =\ \rho_\lambda,
\end{align*}
where $R$ is a $2\pi/3$ rotation matrix and we have used the symmetry assumptions on $V_0$.

Summarizing, for  Case 1 we have  proved:\medskip

\underline{Claim 1:} \nit{\it  $\tsigma=0$, then  \eqref{b3} holds
unless
\begin{enumerate}
\item $\sigma=0$ and $r_1=r_2=\tr_1=\tr_2=0$, a case we address no further since it does not arise in the proof of Proposition \ref{zz-els}\\

or\\

\item $\sigma\ne0$ and $\tr_1=\tr_2=0$, $(r_1,r_2)\in N_b(\sigma)$, in which case $\mathscr{I}(\sigma,\br,\tsigma,\btr)=\rho_\lambda$.
\end{enumerate} 
}
\medskip

Furthermore, because $\tsigma$ and $\sigma$ play symmetric roles as do $\br$ and $\tilde{\bf r}$, we have\bigskip

\underline{Claim 2:}
\nit{\it if $\sigma=0$, then  the bound \eqref{b3} on $\mathscr{I}(\sigma,\br,\tsigma,\btr)$ holds 
unless
\begin{enumerate}
\item $\tsigma=0$ and $r_1=r_2=\tr_1=\tr_2=0$, a case we address no further since it does not arise in the proof of Proposition \ref{zz-els}\\
or\\
\item $\tsigma\ne0$ and $r_1=r_2=0$, $(\tr_1,\tr_2)\in N_b(\tsigma)$, in which case $\mathscr{I}(\sigma,\br,\tsigma,\btr)=\rho_\lambda$.
\end{enumerate} 
}
\bigskip

We now turn to bound on $\mathscr{I}(\sigma,\br,\tsigma,\btr)$ in\\
\nit{\bf Case 2: $\sigma\ne0$ and $\tsigma\ne0$}\\
 \nit{\it Case 2a:}\ \underline{ $\br\in N_b(\sigma)$ and $\btr\in N_b(\tsigma)$}: We claim that 
 \begin{align}
 \mathscr{I}(\sigma,\br,\tsigma,\btr)\  &\lesssim\ e^{-c\lambda}\ \rho_\lambda\quad {\rm for}\ \ \br\in N_b(\sigma),\ \ \btr\in N_b(\tsigma).
  \label{2a}
\end{align}
 By \eqref{rotate}, there exist $l,\tilde{l}\in\{0,1,2\}$ such that 
  $p_0^\lambda(\by-[\sigma\be+\br\vec\vtilde])=p_0^\lambda(\by-\sigma R^l\be)$ and 
  $p_0^\lambda(\by-[\tsigma\be+\btr\vec\vtilde])=p_0^\lambda(\by-\tsigma R^{\tilde l}\be)$. Therefore,
  \begin{align}
 \mathscr{I}(\sigma,\br,\tsigma,\btr)\  &=\ \int\ p_0^\lambda(\by-\sigma R^l\be)\ 
\lambda^2\ |V_0(\by)|\ p_0^\lambda(\by-\tsigma R^{\tilde l}\be)\ d\by\nn\\
&\ \lesssim\ e^{-c\lambda}\ \int\ p_0^\lambda(\by-\sigma R^l\be)\ 
\lambda^2\ |V_0(\by)|\ p_0^\lambda(\by)\ d\by\qquad \textrm{\ (\ by \eqref{p0-bound3}\ ) }\nn\\
&\ \lesssim\ e^{-c\lambda}\ \int\ p_0^\lambda(R^{-l}\by-\sigma \be)\ 
\lambda^2\ |V_0(R^{-l}\by)|\ p_0^\lambda(R^{-l}\by)\ d\by\ =\ e^{-c\lambda}\ \rho_\lambda\ .\nn
\end{align}
   \nit{\it Case 2b:}\ \underline{ $\br\in N_b(\sigma)$ and $\btr\notin N_b(\tsigma)$}: We claim that 
 \begin{align}
 \mathscr{I}(\sigma,\br,\tsigma,\btr)\  &\lesssim\ e^{-c\lambda}\ e^{-c\lambda|\btr|}\ \rho_\lambda\quad
  {\rm for}\ \ \br\in N_b(\sigma),\ \ \btr\notin N_b(\tsigma).
  \label{2b}
\end{align}
By \eqref{rotate} $p_0^\lambda(\by-[\sigma\be+\br\vec\vtilde])=p_0^\lambda(\by-\sigma R^l\be)$, and  by 
\eqref{p0-bound2} and \eqref{p0-bound3}
$p_0^\lambda(\by-[\tsigma\be+\btr\vec\vtilde])\lesssim e^{-c|\btr|\lambda}\ p_0^\lambda(\by-\sigma\be)\lesssim
e^{-c\lambda}\ e^{-c|\btr|\lambda}\ p_0^\lambda(\by) $. These observations together with symmetry imply:
\begin{align}
\mathscr{I}(\sigma,\br,\tsigma,\btr)\  &\lesssim\ e^{-c\lambda}\ e^{-c\lambda|\btr|}\ \int\ p_0^\lambda(\by-\sigma R^l\be)\ 
\lambda^2\ |V_0(\by)|\ p_0^\lambda(\by)\ d\by\ = e^{-c\lambda}\ e^{-c\lambda|\btr|}\ \rho_\lambda\ .
\nn\end{align}
This proves \eqref{2b}. Similarly, if $\br\notin N_b(\sigma)$ and $\btr\in N_b(\tsigma)$ we have
$\mathscr{I}(\sigma,\br,\tsigma,\btr)\lesssim e^{-c\lambda}\ e^{-c\lambda|\br|}\ \rho_\lambda$.
\medskip

\nit{\it Case 2c:}\ \underline{ $\br\notin N_b(\sigma)$ and $\btr\notin N_b(\tsigma)$}:  We claim that 
\begin{equation}
\mathscr{I}(\sigma,\br,\tsigma,\btr)\ \lesssim\  e^{-c\lambda}\ e^{-c\lambda(|\br|+|\btr|)}\ \rho_\lambda\quad 
\ {\rm for}\ \ \br\notin N_b(\sigma),\ \  \btr\notin N_b(\tsigma).
\label{2c}\end{equation}
 
By
\eqref{p0-bound2} and \eqref{p0-bound3}, 
$p_0^\lambda(\by-[\tsigma\be+\btr\vec\vtilde])\lesssim e^{-c|\btr|\lambda}\ p_0^\lambda(\by-\tsigma\be)$
 and $p_0^\lambda(\by-[\sigma\be+\br\vec\vtilde]) \lesssim
e^{-c\lambda}\ e^{-c\lambda|\br|}\ p_0^\lambda(\by)$ . Therefore,
\[\mathscr{I}(\sigma,\br,\tsigma,\btr)\lesssim 
e^{-c|\btr|\lambda}\ \ e^{-c\lambda}\ e^{-c\lambda|\br|}\ \int\ p_0^\lambda(\by-\tsigma\be)\lambda^2|V_0(\by)| p_0^\lambda(\by)\ d\by
\ =\ e^{-c\lambda}\ e^{-c\lambda(|\br|+|\btr|)}\ \rho_\lambda\ .\]

The bounds \eqref{2a}, \eqref{2b} and \eqref{2c} imply Lemma \ref{sigtsig}, and together with 
 Claim 1 and Claim 2 above Lemma \ref{I-bound} follows.
This also completes the proof of Lemma \ref{OIsharp}.

\bibliographystyle{amsplain}
\bibliography{hard_edge-jams}

\end{document}